\documentclass[a4paper,11pt]{amsart}

\usepackage{amsmath}
\usepackage{amssymb,amsbsy,amsmath,amsfonts,amssymb,amscd}
\usepackage{latexsym}
\usepackage{graphics}
\usepackage{color}
\input xy
\xyoption{all}

\newcommand\sC{{\mathcal C}}
\newcommand\sT{{\mathcal T}}
\newcommand\sD{{\mathcal D}}
\newcommand\sE{{\mathcal E}}
\newcommand\sA{{\mathcal A}}
\newcommand\sF{{\mathcal F}}
\newcommand\sG{{\mathcal G}}
\newcommand\sI{{\mathcal I}}
\newcommand\sJ{{\mathcal J}}
\newcommand\sL{{\mathcal L}}

\newcommand\sN{{\mathcal N}}
\newcommand\sK{{\mathcal K}}
\newcommand\sX{{\mathcal X}}
\newcommand\sY{{\mathcal Y}}
\newcommand\sH{{\mathcal H}}

                 \newcommand\sM{{\mathcal M}}
                 \newcommand\XX{{\mathfrak X}}

\newcommand\om{\omega}
\newcommand\la{\lambda}
\newcommand\Lam{\Lambda}
\newcommand\al{\alpha}
\newcommand\be{\beta}
\newcommand\e{\epsilon}
\newcommand\s{\sigma}

\newcommand\Ga{\Gamma}
\newcommand\De{\Delta}
\newcommand\ga{\gamma}
\newcommand\de{\delta}
\newcommand\fie{\varphi}
\DeclareMathOperator{\Pic}{Pic}
\DeclareMathOperator{\Mat}{Mat}

\DeclareMathOperator{\Hom}{Hom}
\DeclareMathOperator{\Alb}{Alb}

\DeclareMathOperator{\Tors}{Tors}

\DeclareMathOperator{\Def}{Def}

\def\Bbb{\bf}

\newcommand{\CC}{\ensuremath{\mathbb{C}}}
\newcommand{\RR}{\ensuremath{\mathbb{R}}}
\newcommand{\ZZ}{\ensuremath{\mathbb{Z}}}
\newcommand{\QQ}{\ensuremath{\mathbb{Q}}}

\newcommand{\sS}{\ensuremath{\mathcal{S}}}

\newcommand{\NN}{\ensuremath{\mathbb{N}}}
\newcommand{\hol}{\ensuremath{\mathcal{O}}}

\newcommand{\HH}{\ensuremath{\mathbb{H}}}
\newcommand{\BB}{\ensuremath{\mathbb{B}}}
\newcommand{\PP}{\ensuremath{\mathbb{P}}}

\newcommand{\FF}{\ensuremath{\mathbb{F}}}
\newcommand{\HHH}{\ensuremath{\mathcal{H}}}

\newcommand{\ra}{\ensuremath{\rightarrow}}

\newcommand{\F}{\ensuremath{\mathbb{F}}}

\def\eea{\end{eqnarray*}}
\def\bea{\begin{eqnarray*}}

\def\C{{\Bbb C}}

\def\La{\langle\langle}
\def\Ra{\rangle\rangle}

\DeclareMathOperator{\Id}{Id}
\DeclareMathOperator{\Aut}{Aut}
\DeclareMathOperator{\Out}{Out}
\DeclareMathOperator{\Inn}{Inn}
\DeclareMathOperator{\Sing}{Sing}
\DeclareMathOperator{\diag}{diag}
\DeclareMathOperator{\Van}{Van}
\DeclareMathOperator{\rank}{rank}
\DeclareMathOperator{\Inv}{Inv}

\newcommand{\Proof}{{\it Proof. }}

\newcommand\dual{\mathrel{\raise3pt\hbox{$\underline{\mathrm{\thinspace d
\thinspace}}$}}}
\newcommand\qe{\ifhmode\unskip\nobreak\fi\quad $\Box$}       

\def\BOX{\hfill\lower.5\baselineskip\hbox{$\Box$}}

\newtheorem{theorem}[equation]{Theorem}
\newtheorem{theo}[equation]{Theorem}
\newtheorem{remark}[equation]{Remark}
\newenvironment{rem}{\begin{remark}\rm}{\end{remark}}

\newtheorem{defin}[equation]{Definition}

\newtheorem{question}[equation]{Question}
\newtheorem{prop}[equation]{Proposition}
\newtheorem{cor}[equation]{Corollary}
\newtheorem{lemma}[equation]{Lemma}
\newtheorem{example}[equation]{Example}
\newenvironment{ex}{\begin{example}\rm}{\end{example}}
\newtheorem{conj}[equation]{Conjecture}

\newcommand{\SSS}{\ensuremath{\mathcal{S}}}
\newcommand{\sR}{\ensuremath{\mathcal{R}}}

\newcommand{\X}{\ensuremath{\mathcal{X}}}
\newcommand{\Y}{\ensuremath{\mathcal{Y}}}

\newcommand{\T}{\ensuremath{\mathbb{T}}}

\DeclareMathOperator{\Aff}{Aff}
\DeclareMathOperator{\im}{Im}
\DeclareMathOperator{\Irr}{Irr}
\DeclareMathOperator{\Sp}{Sp}
\DeclareMathOperator{\GL}{GL}

\def\C{{\Bbb C}}

\def\La{\langle\langle}
\def\Ra{\rangle\rangle}
\DeclareMathOperator{\Jac}{Jac}
\DeclareMathOperator{\id}{id}

\begin{document}

\title[Topological methods for  moduli]{ Topological methods in moduli theory}
\author{ F. Catanese}
\address {Lehrstuhl Mathematik VIII\\
Mathematisches Institut der Universit\"at Bayreuth\\
NW II,  Universit\"atsstr. 30\\
95447 Bayreuth}
\email{fabrizio.catanese@uni-bayreuth.de}

\thanks{The present work took place in the realm of the DFG
Forschergruppe 790 ``Classification of algebraic
surfaces and compact complex manifolds''.
Part of the article was written when the author was  visiting KIAS, Seoul,  as KIAS research scholar. }

\date{\today}
\maketitle

\tableofcontents

\newpage

\section*{Introduction}

The interaction of algebraic geometry and topology has been  such, in the last three centuries,  that it is often difficult to say when
does a result belong to one  discipline or to the other, the archetypical example being the B\'ezout theorem, first conceived through 
a process of geometrical degeneration (algebraic hypersurfaces degenerating  to union of hyperplanes), and later clarified through  topology  and through algebra.

This `caveat' is meant  to warn the reader that a more appropriate title for the present survey article could be: `Some topological methods
in moduli theory, and from the personal viewpoint, taste and understanding of the author'. In fact,  many topics are treated, some classical and some very recent, but with a choice converging towards 
some well defined research interests. 

 I considered for some time  the tempting and  appealing   title  `How can the angel of topology
live happily with the devil of abstract algebra', paraphrasing the motto by Hermann Weyl\footnote {In these days the angel of topology and the devil of abstract algebra fight for the
soul of each individual mathematical domain, \cite{Weyl}, p.500. 

My motto is instead: `Any good mathematical theory requires several good theorems. Conversely,
a really  good theorem requires several good theories. '}.

The latter title would have matched with my personal philosophical point of view:
while it is reasonable that researchers in mathematics develop with enthusiasm and
 dedication new promising mathematical tools and theories, it is important then that  the accumulated knowledge and cultural wealth
  (the instance of topology in the twentieth century being a major one)
 be  not lost  afterwards. This wealth must   indeed not only be  invested and   exploited, but also further developed by  
addressing  problems in other fields, problems which  often raise new and 
fascinating questions. In more down to earth words, the main body of the article is meant to be an invitation for algebraic geometers
to use more classical topology. This invitation  is not new, see for istance the work of Atiyah and Bott on the moduli spaces
of vector bundles on curves (\cite{a-b}); but explains the structure of the article which  is, in a sense, that of a protracted colloquium talk, and where we hope that also topologists, for which many of these notions are well known, will get new kicks coming from algebraic geometry, and especially moduli theory.

In this article we  mostly consider  moduli theory as the fine part 
of classification theory of complex varieties: and we want to show how in some lucky cases topology helps also  for the fine classification,
allowing the  study of  the structure of moduli spaces: as we have done quite concretely  in several papers (\cite{burniat1}, \cite{keumnaie},
\cite{burniat2}, \cite{burniat3}, \cite{bc-inoue}, \cite{bc-CMP}, \cite{bcf}) .

We have already  warned the reader about the inhomogeneity of the level assumed in the text: usually many sections start with very elementary arguments
but, at a certain point, when we deal with current problems, the required knowledge may raise considerably.

Let us try to summarize the logical thread of the article.

Algebraic topology  flourished from some of its applications (such as Brouwer's fixed point theorem, or the theorem of Borsuk-Ulam)   inferring the  non existence of certain continuous maps from the observation that  their existence would imply the existence of homomorphisms satisfying   algebraic properties which are manifestly impossible to be verified.

Conversely, the theory of fibre bundles and homotopy theory give a topological incarnation of a group $G$ through its classifying space $BG$. The theory of classifying spaces translates then group homomorphisms into continuous maps to classifying spaces.
For instance, in algebraic geometry, the theory of Albanese varieties can be understood as dealing with the case where G is free abelian and the classifying maps are holomorphic.

For more general $G$, an important question is the one of the regularity of these classifying maps, such as harmonicity, addressed by Eells and Sampson, and their complex analyticity addressed by Siu and others. These questions, which were at the forefront of mathematical research  in the last 40 years, have powerful applications to moduli theory.

After a general introduction directed towards a broader public, starting  with classical theorems by Zeuthen-Segre and Lefschetz,
proceeding to classifying spaces and their properties,  I shall concentrate on some classes of projective varieties which are classifying spaces for some group, 
providing several explicit examples. I  discuss then locally symmetric varieties, and at a certain length the quotients of Abelian varieties by a cyclic group acting freely, which are here called Bagnera-De Franchis varieties.

At this point the article becomes instructional, and oriented towards graduate students, and several important topics, like orbifold fundamental groups,
Teichm\"uller spaces, moduli spaces of curves, group cohomology and homology are treated in detail (and a new proof of a classical theorem of Hopf is sketched).

Then some applications are given to concrete problems in moduli theory, in particular a new construction of surfaces with $p_g=q=1$ is given.

The next section is devoted to a preparation for the  rigidity and quasi-rigidity properties of projective varieties which are classifying spaces
 (meaning that their moduli spaces are completely determined by their  topology); in the section are recalled the by now classical results of Eells and Sampson,
and Siu's results about complex analyticity of harmonic maps, with particular emphasis on bounded domains and locally symmetric varieties.

Other more elementary results, based on Hodge theory, the theorem of Castelnuovo-De Franchis, and on the explicit constructions of classifying spaces
are explained in detail because of their importance for K\"ahler manifolds. We then briefly discuss Kodaira's problem and Voisin's counterexamples,
then we dwell on fundamental groups of projective varieties, and on the   Shafarevich conjecture.

Afterwards we deal with several concrete  investigations of  moduli spaces, which in fact  lead to some group theoretical questions, and to the investigation of moduli spaces of varieties with symmetries.

Some key examples are: varieties isogenous to a product, and the Inoue-type varieties introduced in recent work with Ingrid Bauer: for these the moduli space is determined by the topological type. I shall present new results and open questions concerning this class of varieties.

In   the final part, after recalling basic results on complex moduli theory,  we shall also illustrate 
 the concept of symmetry marked varieties and their moduli,
discussing the  several reasons why  it is interesting to consider moduli spaces of triples $(X,G,\al)$ where $X$ is a projective variety, $G$ is a finite group, and $\al$ is an effective action of $G$  on $X$. If $X$ is the canonical model of a variety of general type, then $G$ is acting linearly on some pluricanonical model, and we have a moduli space which is a finite covering of a closed subspace $\mathfrak M^G$  of the moduli space.

In the case of curves we show how  this investigation is related to the description  of the singular locus of the moduli space $\mathfrak M_g$ (for instance of its irreducible components, see \cite{cornalba}),  and of  its compactification $\overline{\mathfrak M_g}$, (see \cite{singMg}).

In the case of surfaces there is another occurrence of Murphy's law, as shown in my joint work with Ingrid Bauer (\cite{burniat3}): the deformation equivalence for minimal models $S$ and for canonical models differs drastically (nodal Burniat surfaces being the easiest example).
This shows how appropriate it is to work with Gieseker's moduli space of canonical models of surfaces.

In the case of curves, there are interesting relations with topology. Moduli spaces $\mathfrak M_g (G)$ of curves with a group $G$  of automorphisms of a fixed topological type have a description by Teichm\"uller theory, which  naturally leads to conjecture genus stabilization for rational homology groups. I will then describe  two equivalent descriptions of the  irreducible components of $\mathfrak M_g (G)$, surveying known
irreducibility results  for some special groups. A  new fine homological invariant was introduced in our joint work  with L\"onne and Perroni: it allows to prove  genus stabilization in the ramified case, extending a beautiful  theorem due to Livingston (\cite{Liv}) and  Dunfield-Thurston (\cite{DT}), who dealt with  the easier unramified case.

Another  important application is  the following one, in the direction of arithmetic: in the 60's J. P. Serre (\cite{serre}) showed   that there exists a field automorphism  $\s$  in the absolute Galois group $Gal(\bar{\QQ} / \QQ)$,
and a variety $X$ defined over a number field, such that $X$ and the Galois conjugate variety $X^{\s}$  have non isomorphic fundamental groups, in particular they are not homeomorphic.

In a joint paper with I. Bauer  and F. Grunewald  we proved  a  strong sharpening of this phenomenon discovered by Serre, namely, 
that if $\s$   is not in the conjugacy class of the complex conjugation then there exists a surface (isogenous to a product)  $X$  such that $X$ and the Galois conjugate variety $X^{\s}$ have non isomorphic fundamental groups.

In the end we finish with an extremely quick mention of several interesting topics which we do not have the time to
describe properly, among these, the stabilization results for the cohomology of moduli spaces and of arithmetic varieties.

\newpage

\section{Prehistory and beyond}

The following discovery belongs to the 19-th century: consider the complex projective plane $\PP^2$ and two general
homogeneous polynomials $ F, G \in \CC[x_0,x_1,x_2]$ of the same degree $d$. Then  $F,G$ determine a {\bf linear pencil}
of curves $ C_{\la} , \ \forall \la = (\la_0, \la_1) \in \PP^1$,
$$  C_{\la} := \{ x = (x_0,x_1,x_2) \in \PP^2 | \la_0 F(x) +  \la_1 G(x) = 0 \}.$$
One sees that the curve $ C_{\la}$ is singular for exactly $ \mu = 3 (d-1)^2 $ values of $\la$, as can be verified 
by an elementary argument which we now sketch.

In fact, $x$ is a singular point of some $ C_{\la}$  iff the following system of three homogeneous linear equations
in $ \la = ( \la_0, \la_1)$ has a nontrivial solution:
$$ \la_0 \frac{\partial F}{\partial x_i} (x) +  \la_1 \frac{\partial G}{\partial x_i} (x) = 0, \ \forall i = 0,1,2.$$

 By generality of $F,G$ we may assume that the curves $ C_0 : = \{x | F(x) = 0\} $ and $ C_1 : = \{x | G(x) = 0\} $ are smooth and  intersect
 transversally (i.e., with distinct tangents) in $d^2$ distinct points; hence  if a curve of the pencil $ C_{\la}$ has a singular point $x$,
 then we may assume that for this point we have 
 $F(x) \neq 0 \neq G(x)$, and then  $\la$ is uniquely determined.

If now $\frac{\partial F}{\partial x_0} $ and $\frac{\partial G}{\partial x_0} $ do not vanish simultaneously in $x$, then the above system
has a nontrivial solution if and only if

$$[ \frac{\partial G}{\partial x_0} \cdot \frac{\partial F}{\partial x_i } -  \frac{\partial F}{\partial x_0} \cdot \frac{\partial G}{\partial x_i }](x)= 0, \ i=1,2. $$

By the theorem of B\'ezout (see \cite{walker}) the above two equations have  $ ( 2 (d-1))^2 = 4 (d-1)^2$ solutions, including among these
the $(d-1)^2$ solutions of the system of two equations $\frac{\partial F}{\partial x_0}(x) = \frac{\partial G}{\partial x_0} (x) = 0$.

One sees that, for $F, G$ general, there are no common solutions of the system
 $$\frac{\partial F}{\partial x_0}(x) = \frac{\partial G}{\partial x_0} (x) = [ \frac{\partial G}{\partial x_1} \cdot \frac{\partial F}{\partial x_2 } -  \frac{\partial F}{\partial x_1} \cdot \frac{\partial G}{\partial x_2 }](x)= 0,$$
hence  the solutions of the  above system are indeed $$ 3 (d-1)^2 = 4 (d-1)^2 -  (d-1)^2 .$$
 
 It was found indeed that, rewriting   $ \mu =  3 (d-1)^2 = d^2 + 2d (d-3)+ 3  $,  the above  formula generalizes  to
 a beautiful formula, valid for any smooth algebraic surface $S$, and which is  the content of  the so-called theorem
 of  Zeuthen-Segre; this goes as follows: observe in fact  that $d^2$ is the number of points where the curves of the pencil meet, 
 while  the genus $g$ of a plane curve of degree $d$ equals $\frac{ (d-1)(d-2)}{2}$.
 
 \begin{theo}{\bf (Zeuthen-Segre, classical)}
Let $S$  be a smooth projective surface, and let $ C_{\la} , \  \la  \in \PP^1$ be a linear pencil of curves of genus $g$ which meet transversally 
in $\de $ distinct points. If $\mu$ is the number of singular curves in the pencil (counted with multiplicity), then
$$ \mu -  \de - 2 (2g-2) = I + 4, $$ where the integer $I$ is an invariant of the algebraic surface, called Zeuthen-Segre invariant.
 
 \end{theo}
 
 Here, the integer $\de$ equals the self-intersection number $C^2$ of the curve $C$, while in modern terms the number $2 g -2 =
 C^2 + K_S \cdot C$, $K_S$ being the divisor (zeros minus poles) of a rational differential 2-form.
 
 In particular, our previous calculation shows that for  $\PP^2$ the invariant $I = -1$.
 
 The interesting part of the discovery is that the integer $ I + 4$ is not only an algebraic invariant, but is indeed a topological invariant.
  
  Indeed, for a compact topological space $X$ which can be written as the disjoint union of locally closed sets 
 $X_i, \  i = 1, \dots r,$
 homeomorphic to an Euclidean space $\RR^{n_i}$, one can define
 $$  e(X) := \sum_i^r (-1)^{n_i},$$
 and indeed this definition is compatible with the more abstract definition
 $$ e(X) =  \sum_{j=0}^{\dim (X)} (-1)^j \rank H_j(X, \ZZ). $$ 
 
 For example,  the plane $\PP^2 = \PP^2_{\CC}$ is obtained from a point attaching $\CC = \RR^2$  and then $\CC^2 = \RR^4$,
 hence $ e(\PP^2) = 3$ and we verify that $  e(\PP^2) = I + 4$. 
 
 While for an algebraic curve $C$ of genus $g$ its topological Euler-Poincar\'e characteristic, for short Euler number, equals  $ e(C) = 2 - 2g$,
 since $C$  is obtained as the disjoint
 union of one point, $2g$ arcs, and a 2-disk (think of the topological realization as the quotient of a polygon with $4 g$ sides).
 
 The Euler Poincar\'e characteristic is multiplicative for products: $$ e (X \times Y ) = e (X) \cdot e (Y),$$ and more generally for fibre bundles
 (a concept we shall introduce in the next section), and accordingly there is a generalization of the theorem of Zeuthen-Segre:
 
  \begin{theo}{\bf (Zeuthen-Segre, modern)}\label{ZS}
Let $S$  be a smooth compact complex surface, and let $ f : S \ra B$ be a fibration onto a projective curve $B$ of genus $b$
(i.e., the fibres $f^{-1}(P)$, $ P \in B$, are connected), and denote by $g$ the genus of 
the smooth fibres of $f$. Then 
$$  e(S) =   (2b- 2 )(2g-2) + \mu, $$ where  $\mu \geq 0$, and $\mu = 0$ if and only if all the fibres of $f$ are either 
smooth or, in the case where $g=1$,
a multiple of a smooth curve of genus $1$.
 
 \end{theo}
 
 The technique of studying linear pencils turned out to be an invaluable tool for the study of the topology of projective varieties.
 In fact, Solomon Lefschetz in the beginning of the 20-th century was able to describe the relation holding between a smooth projective variety $X \subset \PP^N$ of dimension $n$
 and its hyperplane section $ W = X \cap H$, where $H $ is a general linear subspace of codimension 1, a hyperplane.
 
 The work of Lefschetz deeply impressed the Italian algebraic geometer Guido Castelnuovo, who came to the conclusion that algebraic geometry
 could no longer be carried over without the new emerging techniques, and convinced Oscar Zariski  to go on setting the building of algebraic geometry on a more solid basis. The report of Zariski (\cite{Zar}) had a big influence and the results of Lefschetz were reproven and vastly extended
 by several authors: they say essentially that homology and homotopy groups of real dimension smaller than the complex dimension  $n$
 of $X$ are the same for $X$ ands its hyperplane section $W$. 
 
 In my opinion the nicest proofs of the theorems of Lefschetz are those given much later by Andreotti and Frankel (\cite{AF1} , \cite{AF2}).
 
 \begin{theo}\label{hyperplanesection}
 Let $X$ be a smooth projective variety of complex dimension $n$, let $W = X \cap H$ be a smooth hyperplane section of $X$,
 and let  further $Y = W \cap H'$ be  a smooth hyperplane section of $W$.
 
 {\bf First Lefschetz' theorem:} the natural homomorphism $ H_i (W, \ZZ) \ra H_i (X, \ZZ)$ is bijective for $ i < n-1$, and surjective for $i = n-1$;
 the same is true for the natural homomorphisms of homotopy groups $\pi_i (W) \ra \pi_i (X)$
 (the results hold more generally, see \cite{milnorMT},page 41,  even if $X$ is singular and $W$ contains the singular locus of $X$).
 
  {\bf Second Lefschetz' theorem:} The kernel of  $ H_{n-1} (W, \ZZ) \ra H_{n-1} (X, \ZZ)$ is the subgroup $\Van H_{n-1} (W, \ZZ)$ generated by the
  vanishing cycles, i.e., those cycles which are mapped to $0$ when $W$ tends to a singular hyperplane section $W_{\la}$ in a pencil
  of hyperplane sections of $X$.
  
   {\bf Generalized Zeuthen-Segre theorem:} if $\mu$ is the number of singular hyperplane sections
   in a general linear pencil
  of hyperplane sections of $X$, then
$$  e(X) = 2 e(W) - e(Y) + (-1)^n \mu.$$
   
    {\bf Third Lefschetz' theorem or Hard Lefschetz' theorem:} 
    
    The first theorem and the universal coefficients theorem  imply for the cohomology groups 
     that  $ H^i (X, \ZZ) \ra H^i (W, \ZZ)$ is bijective for $ i < n-1$, while  $ H^{n-1} (X, \ZZ) \ra H^{n-1} (W, \ZZ)$ is injective.
  Defining $ \Inv H_{n-1} (W, \ZZ)$ as the Poincar\'e dual of the image of  $ H^{n-1} (X, \ZZ)$, then we have a  direct sum decomposition
  (orthogonal for the cup product)
  after tensoring with $\QQ$:
  $$ H_{n-1} (W, \QQ) =  \Inv H_{n-1} (W, \QQ) \oplus  \Van H_{n-1} (W, \QQ).$$
  Equivalently, the operator $L$ given by cup product with the cohomology class $h \in H^2 (X, \ZZ)$ of a hyperplane,
  $  L \colon H^i (X, \ZZ) \ra H^{i+2} (X, \ZZ)$, induces an isomorphism
  $$ L^j \colon H^{n-j} (X, \QQ) \ra H^{n+j} (X, \QQ), \ \forall j \leq n. $$
 \end{theo}
 
 Not only the theorems of Lefschetz  play an important role for our particular purposes, but we  feel  that we should also spend  a few words
 sketching how they lead to   some very  interesting and still widely open conjectures, the Hartshorne conjectures (see \cite{hartshorne}).
 
 Assume now that the smooth projective variety $X \subset \PP^N$  is the  complete intersection of $N-n$ hypersurfaces (this means the the sheaf $\sI_X$ of
 ideals of functions vanishing on $X$ is generated by polynomials $F_1, \dots, F_c$, $c : = N - n$ being the codimension of $X$). 
 Then the theorems of Lefschetz \footnote{one uses here the following trick: the intersection of a projective variety $X$ with a hypersurface is equal to the
 intersection of $X$ with a hyperplane, but  for a different embedding of  $X$; this trick is used several times, starting with $X = \PP^N$}
 imply that the homology groups of $X$ equal those of $\PP^N$ for $i \leq n-1$ (recall that  $H_i(\PP^N, \ZZ) = 0$ for $i $ odd,
 while $H_i(\PP^N, \ZZ) = \ZZ$ for $i \leq 2N$, $i$ even). Similarly holds true for the homotopy groups, and we recall that,
 since $ \PP^N = S^{2N+1} / S^1$, then $ \pi_i (\PP^N) = 0$ for $ i \leq 2N, i \neq 0,2$.
 
 It was an interesting discovery by Barth  ( see \cite{Barth1, bl, bvdv, fh} for the following and related results)  that a similar (but weaker)  result holds true for each smooth subvariety of $\PP^N$, provided the codimension $ c = N - n$ of $X$  is smaller than the dimension.  
 
 \begin{theo}{\bf (Barth-Larsen)}\label{BL}
 Let $X$ be a smooth subvariety of dimension $n$ in $\PP^N$: then the homomorphisms
 $$    H_i (X, \ZZ)  \ra H_i (\PP^N, \ZZ) , \     \pi_i (X) \ra \pi_i (\PP^N) $$
 are bijective for $ i \leq  n - c \Leftrightarrow i < 2n - N + 1$, and surjective for $ i = n - c + 1 = 2n - N + 1$.
 
 \end{theo}
 Observe that, if $N= n+1$, then the above result  yields exactly the one of Lefschetz, hence the theorem is  sharp in this trivial case
 (but much weaker for complete intersections of higher codimension). The case of the Segre embedding $ X :=  \PP^1 \times \PP^2 \ra \PP^5$
 is a case which shows how the theorem is sharp since $  (\ZZ)^2   \cong H_2 (X, \ZZ) \ra H_2 (\PP^8, \ZZ)   \cong  \ZZ $ is  surjective
 but not bijective.
 
 The reader might wonder why the theorem of Barth and Larsen is a generalization of the theorem of Lefschetz. 
 First of all, while Barth used originally methods of holomorphic convexity in complex analysis
 (somehow reminiscent of Morse theory in the real case)  Hartshorne showed (\cite{hartshorne}) how the third Lefschetz Theorem implies the result of Barth for cohomology
 with coefficients in $\QQ$. Moreover a strong similarity with the Lefschetz situation  
 follows from the fact that one may view it (as shown by Badescu, see \cite{badescu}) as an application of the classical Lefschetz theorem
 to the intersection $$  ( X \times \PP^N) \cap \De \cong X,$$
 where $\De \subset   \PP^N \times \PP^N$ is the diagonal. In turn, an idea of Deligne (\cite{del}, \cite{f-l})  shows that the diagonal
 $\De \subset   \PP^N \times \PP^N$ behaves `like' a complete intersection, essentially because, under the standard birational
 map $\PP^N \times \PP^N \dashrightarrow \PP^{2N}$, it maps to a linear subspace of $\PP^{2N}$.

The philosophy is then that smooth subvarieties of small codimension behave like complete intersections.
This could be no accident if the well known Hartshorne conjecture (\cite{hartshorne}) were true.

\begin{conj}{\bf (On subvarieties of small codimension, Hartshorne)}\label{HC1}
Let  $X$ be a smooth subvariety of dimension $n$ in $\PP^N$, and assume that the dimension is bigger than twice the codimension,
$ n > 2 (N-n)$: then $X$ is a complete intersection.
\end{conj}

While the conjecture says nothing in the case of curves and surfaces and is trivial in the case where $ n \leq 4$, since a codimension
1 subvariety is defined by a single equation, it starts to
have meaning for $n \geq 5$ and $ c = N - n \geq  2$. 

In the case where $c = N - n = 2$, then by a result of Serre one knows (see \cite{ellia}, page 143, also  for a general survey of the Hartshorne  conjecture for codimension 2 subvarieties) that $X$ is the zero set of a section $ s$ of a rank 2 holomorphic
vector bundle  $V$ on $\PP^N$ (observe that in this  codimension one has that $\Pic  (X) \cong \Pic  (\PP^N) $, so that $X$ satisfies the condition of
being subcanonical: this means that $\omega_X = \hol_X(d)$ for some $d$): in the case $c=2$ the conjecture by Hartshorne is then equivalent to the conjecture

\begin{conj}{\bf (On vector bundles of rank 2 on projective space, Hartshorne)}\label{HC2}
Let  $V$ be a rank 2 vector bundle  on  $\PP^N$, and assume that  $ N \geq 7$: then $V$ is a direct sum of line bundles.
\end{conj}

The major evidence for the conjecture on subvarieties of small codimension comes from the concept of positivity of vector bundles
(\cite{f-l}, also \cite{laz}). In fact, many construction methods of subvarieties $X$ of small codimension involve a realization of  $X$ as the locus where a vector bundle map drops
rank to an integer $r$, and $X$ becomes singular if there are points of $X$ in the locus $\Sigma$ where the rank drops further down to $(r-1)$.

The expected dimension of $\Sigma$ is positive in the range of Hartshorne's conjecture, but nevertheless this is not sufficient
to show that $\Sigma \cap X $ is non empty.

Hartshorne's  conjecture \ref{HC1} is related to projections: in fact, for each projective variety $X \subset \PP^N$ there exists a linear projection
$ \PP^N \dashrightarrow \PP^{2n+1}$ whose restriction to $X$ yields an embedding  $X \subset \PP^{2n+1}$: in other words,
an embedding where the codimension is equal to the dimenion $n$ plus 1. The condition of being embedded as a subvariety where
the codimension is smaller or equal  than the dimension is already a restriction (for instance, not all curves are plane curves, and smooth surfaces
in $\PP^3$ are simply connected by Lefschetz's theorem),
and the smaller the codimension gets, the stronger the restrictions are (as shown by theorem \ref{BL}).

Speaking now in more technical terms, a necessary condition for $X$ to be a complete intersection is that the sheaf of ideals $\sI_X$ be arithmetically Cohen-Macaulay
(ACM, for short),
which means that all higher cohomology groups $H^i(\PP^N, \sI_X(d)) = 0$ vanish for $n \geq  i >0$ and $\forall d \in \ZZ$. In view of the exact sequence
$$ 0 \ra \sI_X \ra \hol_{\PP^N} \ra \hol_X \ra 0 $$ the CM condition amounts to two conditions:

1) $X$ is projectively normal, i.e., the linear system cut on $X$ by polynomials of degree $d$ is complete ( $H^0(\PP^N, \hol_{\PP^N}(d) ) \ra 
H^0 (\hol_X(d))$ is surjective for all $d \geq 0$)

2)  $H^i(X , \hol_X(d)) = 0$ for all $d \in \ZZ$, and for all  $ 0 < i < n = \dim (X)$.

In the case where $X$ has codimension 2, and $ N \geq 6$ (see \cite{ellia}, cor. 4.2, page 165), 
 the condition of being a complete intersection is equivalent to projective normality.
 
Linear normality is the case $d = 1$ and means that $X$ is not obtained as the projection of a non-degenerate variety  from a higher dimensional  projective space
$\PP^{N+1}$. This part of Hartshorne's conjecture is the only one which has been verified: a smooth subvariety with $ n \geq 2 c -1 = 2 (N-n) -1$ is linearly normal
 (theorem of Zak,  \cite{zak}).
 
 For  higher $d$, one considers the so-called formal neighbourhood of $X$: denoting by $N_X^{\vee} : =  \sI_X /  \sI_X^2$ the conormal bundle
 of $X$, one sees that, in order to show projective normality, in view of the exact sequence
 $$ 0 \ra  \sI_X^{m+1}(d) \ra   \sI_X^m (d)\ra Sym^m ( N_X^{\vee})(d) \cong    (\sI_X^m  /  \sI_X^{m+1})(d)  \ra    0,$$
 a crucial role is played by  the cohomology groups
 $$ H^q( Sym^m ( N_X^{\vee})(d))).$$ We refer the reader to \cite{PPS} for a discussion of more general  Nakano type vanishing
 statements of the form $ H^q( Sym^m ( N_X^{\vee})(d) \otimes \Omega^p_X)$ (these could be implied by very strong curvature properties
 on the normal bundle, see also \cite{laz}).

Observe finally that the theorems of Lefschetz have been also extended to the case of singular varieties (see \cite{g-mp}, \cite{fultonAMS}),
but we shall not need to refer to these extensions in the present paper.

\section{Algebraic topology: non existence and existence of continuous maps}

The first famous achievements of algebraic topology were based on functoriality, which was used to infer the nonexistence of certain
 continuous maps.
 
 The Brouwer's fixed point theorem says that every continuous self map $ f : D^n \ra D^n$, where $D^n = \{ x \in \RR^n | |x| \leq 1\}$
 is the unit disk, has a fixed point. The argument is by contradiction: otherwise, letting $\phi(x)$ be the intersection of the boundary $S^{n-1}$ of $D^n$
 with the half line stemming from $f (x)$ in the direction of $x$, $\phi$ would be a continuous map 
 $$ \phi : D^n \ra  S^{n-1}, \  s.t. \ \phi|_{S^{n-1}} = \Id_{S^{n-1}}.$$ 
 The key point is to show that the reduced \footnote{the reduced homology group differs from the ordinary one only for $i=0$,
 and for $i=0$ is defined as the kernel of the degree surjection onto $\ZZ $: this distinction is only needed  in order to treat the case $n=1$
 on an equal footing.  }
 homology group $H_{n-1} ({S^{n-1}} , \ZZ) \cong \ZZ$, while $H_{n-1} (D^n , \ZZ) = 0$, the disc being 
 contractible; after that, denoting by $\iota : S^{n-1} \ra D^n$ the inclusion,  functoriality of homology groups,
 since $ \phi \circ \iota = \Id_{S^{n-1}}$,  would imply  $0 = H_{n-1} ( \phi ) \circ H_{n-1} ( \iota )= H_{n-1} ( \Id_{S^{n-1}}) = \Id_{\ZZ} $,
 the desired contradiction.
 
 Also well known is the Borsuk-Ulam theorem, asserting that there is no odd continuous  function $ F : S^n \ra S^m$ for $ n > m$ (odd means that $ F ( -x) = - F(x)$).
 
 Here there are two ingredients, the main one being the cohomology algebra, and its contravariant functoriality: to any continuous map
 $ f : X \ra Y$ there corresponds an algebra homomorphism $$ f^* : H^* (Y, R) = \oplus_{i=0}^{\dim(Y)}  H^i (Y, R)  \ra H^* (X, R) ,$$
 for any ring $R$ of coefficients.

 In our case one  takes as $X : = \PP^n_{\RR} = S^n /\{\pm 1\}$, similarly $Y : = \PP^m_{\RR} = S^m /\{\pm 1\}$ and lets $f$ be the continuus map induced by $F$.
 One needs to show that, choosing $ R = \ZZ / 2 \ZZ$, then the cohomology algebra of real projective space is a truncated polynomial algebra,
 namely:
 $$ H^* (  \PP^n_{\RR},   \ZZ / 2 \ZZ) \cong ( \ZZ / 2 \ZZ ) [\xi_n] / (\xi_n^{n+1}). $$
 
 The other ingredient consists in showing that $$ f^* ([\xi_m])  = [\xi_n],$$ $[\xi_m]$ denoting the residue class in the quotient algebra.
 
 One gets then the desired contradiction since,  if $ n > m$, $$ 0 =  f^*(0) = f^* ([\xi_m]^{m+1}) =  f^* ([\xi_m])^{m+1} = [\xi_n]^{m+1} \neq 0.$$
 Notice that up to now we have mainly used that $f$ is a continuous map  $f  : =\PP^n_{\RR} \ra  \PP^m_{\RR},$ 
 while precisely in order to obtain that 
 $ f^* ([\xi_m])  = [\xi_n]$ we must make use of the hypothesis that $f$ is induced by an odd function $F$.
 
 This property can be interpreted as the property that one has a commutative diagram
 
 \[ 
  \begin{matrix}
    S^n &    \ra   &   S^m \\
     \downarrow &  &  \downarrow  \\
     \PP^n_{\RR}&\ra  & \PP^m_{\RR}
  \end{matrix} 
\] 

which exhibits the two sheeted covering of  $\PP^n_{\RR} $ by $S^n$ as the pull-back of the analogous two sheeted
cover for $\PP^m_{\RR}$. Now, as we shall digress soon, any such two sheeted covering is given by a homomorphism
of $H_1(X,  \ZZ / 2 \ZZ) \ra  \ZZ / 2 \ZZ$, i.e., by an element in $H^1(X,  \ZZ / 2 \ZZ) $, and this element is trivial if and only if the
covering is trivial (that is, homeomorphic to $ X \times ( \ZZ / 2 \ZZ) $, in other words a disconnected cover).

This shows that the pull back of the cover, which is nontrivial, corresponds to  $ f^* ([\xi_m])$ and is nontrivial, hence  $ f^* ([\xi_m]) = [\xi_n]$. 

As we saw already in the first section, algebraic topology attaches to a good topological space homology groups $H_i (X, R)$, which are
covariantly functorial, a cohomology algebra   $H^* (X, R)$ which is contravariantly functorial, and these groups can be calculated,
by virtue of the Mayer Vietoris exact sequence and of excision (see any textbook), by chopping the space in smaller pieces. In particular,
these groups vanish when $i > dim (X)$. But to $X$ are also attached the homotopy groups $\pi_i(X)$. The common feature is that homotopic maps
induce the same homomorphisms on homology, cohomology, and homotopy.

We are, for our purposes, more interested in the more mysterious homotopy groups, which, while not necessarily vanishing 
for  $ i > dim (X)$, enjoy however  a fundamental property.

Recall the definition due to Whitney and Steenrod (\cite{steenrod}) of a fibre bundle.  In the words of Steenrod, the notion of a fibre bundle is a weakening
of the notion
of a product, since a product $ X \times Y$ has two continuous projections $p_X : X \times Y \ra X$, and  $p_Y : X \times Y \ra Y$,
while a fibre bundle $E$ over $B$ with fibre $F$  has only one projection, $p = p_B : E \ra B$ and its similarity to a product lies in the fact that for each point $ x \in B$ there
is an open set $ U$ containing $x$, and a homeomorphism  $ p_B^{-1} (U) \cong U \times F$ compatible with both projections onto $U$.

The fundamental property of fibre bundles is that there is a long  exact sequence of homotopy groups 
$$ \dots \ra  \pi_i(F)  \ra  \pi_i(E)  \ra  \pi_i(B )  \ra  \pi_{i-1}(F )  \ra  \pi_{i-1}(E )  \ra  \pi_{i-1}(B ) \ra \dots $$
where one should observe that $\pi_i (X)$ is  a group for $ i \geq 1$, an abelian group for $ i \geq 2$, and for
$i=0$ is just the set of arc-connected components of $X$ (we assume the spaces to be good, that is, locally arcwise connected,
 semilocally simply connected, see \cite{greenberg}, and, most of the times, connected). 

The special case where the fibre $F$ has the discrete topology is the case of a {\bf covering space}, which is called the universal covering
if moreover $\pi_1(E)$ is trivial.  

Special mention deserves the following more special case.

\begin{defin}\label{kp1}
Assume that $E$ is arcwise connected, contractible (hence all homotopy groups $\pi_i(E)$ are trivial), and that the fibre $F$ is discrete, so that
all  the higher homotopy groups  $ \pi_i(B ) = 0$ for $i \geq 2$, while  $ \pi_1(B ) \cong \pi_0(F) = F$.
Then one says that $B$ is a classifying space $K(\pi, 1)$ for the group $\pi =  \pi_1(B )$.

In general, given a group $\pi$,  a CW complex $B$ is said to be a $K(\pi, 1)$ if  $ \pi_i(B ) = 0$ for $i \geq 2$, while  $ \pi_1(B ) \cong \pi$.

\end{defin}

\begin{ex}
The easiest examples are the following ones, where (2) is a case where we have a complex projective variety (see section 3 for more such examples):
\begin{enumerate}
\item
the real torus $T^n : = \RR^n / \ZZ^n$ is a classifying space $ K (\ZZ^n, 1)$ for the group $\pi = \ZZ^n$;
\item
a complex projective curve $C$ of genus $ g \geq 2$ is a classifying space $ K (\pi_g, 1)$, since by the uniformization theorem
its universal covering is the Poincar\'e upper half plane $ \HHH : = \{ z  \in \CC |  Im (z) > 0 \}$ and
its  fundamental group $\pi_1(C)$ is isomorphic to the group
$$\pi_g  : = \langle \al_1, \be_1, \dots \al_g,\be_g | \Pi_1^g [ \al_i, \be_i]  = 1\rangle ,$$
quotient of a free group with $2g$ generators by the normal subgroup generated by the relation $ \Pi_1^g [ \al_i, \be_i] $; 
\item
a classifying space $ K (\ZZ/ 2 \ZZ, 1)$ is given by the inductive limit $ \PP^{\infty}_{\RR} : = \lim_{n \ra \infty} \PP^n_{\RR}$. To show this, it suffices
to show  that  $ S^{ \infty}  : =  \lim_{n \ra \infty} S^n$
is contractible or, equivalently, that the identity map is homotopic to a constant map.

We do this as follows \footnote{following a suggestion of Marco Manetti}: let $\s : \RR^{ \infty}   \ra \RR^{ \infty}  $ be the shift operator,
and first define a homotopy of the identity map of $\RR^{ \infty}  \setminus \{ 0 \}$ to the constant map with value $e_1$.
The needed homotopy is the composition of two homotopies:
$$ F (t,v) : =  (1-t) v +  t  \s (v),  0 \leq t \leq 1 , \ \ $$
$$  F (t,v) : =  (2-t) \s (v) + (t-1) e_1, \   1 \leq t \leq 2, \ \forall v \in \RR^{ \infty} .$$ 

Then we simply project the homotopy from $\RR^{ \infty}  \setminus \{ 0 \}$ to $ S^{ \infty} $ considering $ \frac{F(t,v)}{| F(t,v)|}$.

\end{enumerate}

\end{ex}

These classifying spaces, although not unique, are unique up to {\bf homotopy-equivalence} (we use the notation $X \sim_{h.e.} Y $ 
to denote homotopy equivalence: this means  that there exist continuous maps $ f : X \ra Y, \ g : Y \ra X$ such that both compositions $ f\circ g$ and $ g \circ f$ are homotopic to the identity).

Therefore, given two classifying spaces for the same group, they not only do have the same homotopy groups, but also
the same homology and cohomology groups. Thus the following definition is well posed.

\begin{defin}
Let $\Ga$ be a finitely presented group, and let $ B\Ga$ be a classifying space for $\Ga$: then the homology and cohomology groups and algebra of $\Ga$ are defined as
$$ H_i (\Ga, \ZZ) : =  H_i ( B \Ga, \ZZ), \  H^i (\Ga, \ZZ) : =  H^i ( B \Ga, \ZZ), \  H^* (\Ga, \ZZ) : =  H^*( B \Ga, \ZZ), $$
and similarly for other rings of coefficients instead of $\ZZ$. 
\end{defin}

\begin{rem}
The concept of a classifying space $BG$ is indeed more general: the group $G$ could also be a Lie group, and then, if $EG$
is a contractible space over which $G$ has a free (continuous) action, then one defines $ BG : = EG / G$.

The typical example is the simplest compact Lie group $G = S^1$: then, keeping in mind that $ S^1 = \{ z \in \CC | |z| = 1 \}$,
we take as $EG $ the space $$ E S^1 : = S (\CC^{\infty}) = lim_{n \ra \infty}  S (\CC^{n}) = lim_{n \ra \infty}  \{ w \in  \CC^{n} | |w| = 1 \}. $$
\end{rem}

Classifying  spaces, even if often quite difficult to construct explicitly, are very important because they guarantee the existence of
continuous maps! We have more precisely the following  (cf. \cite{Spanier}, Theorem 9, page 427, and Theorem 11, page 428)

\begin{theo}\label{classifying map}
Let $Y$ be a `nice'  topological space, i.e., $Y$ is homotopy-equivalent to a CW-complex, and let $X$ be a nice space which is a  $ K (\pi, 1)$
space: then, choosing base points $y_0 \in Y, x_0 \in X$,  one has a bijective  correspondence 
$$[ (Y, y_0 ), (X, x_0 )] \cong \Hom (\pi_1(Y, y_0), \pi_1(X, x_0)) , [f]  \mapsto   \pi_1 (f) ,$$
where $ [ (Y, y_0 ), (X, x_0 )] $ denotes the set of homotopy classes $[f]$ of continuous maps $ f : Y \ra X$ such that $ f(y_0 ) = x_0$ (and where the homotopies 
$ F(y,t)$ are also required to satisfy $ F(y_0,t) = x_0, \ \forall t \in [0,1]$).

In particular, the free homotopy classes $[Y, X]$ of continuous maps are in bijective correspondence with the conjugacy classes of homomorphisms
$  \Hom (\pi_1(Y, y_0), \pi )$ (conjugation is here inner conjugation by $\Inn (\pi )$ on the target). 
\end{theo}

Observe that, quite generally, the universal covering $E _{\pi}$ of a classifying space $ B \pi : = K (\pi, 1)$ 
associates (by the lifting property) to a continuous map $f : Y \ra B \pi$ a $\pi_1(Y) $-equivariant map $ \tilde{f} $ 
$$ \tilde{f} :  \tilde{Y} \ra E _{\pi},$$
where the action of $\pi_1(Y) $ on $E _{\pi}$ is determined by the homomorphism $\varphi : = \pi_1(f) : \pi_1(Y) \ra \pi = \pi_1 (B \pi )$.

Moreover, any $\varphi :  \pi_1(Y) \ra \pi = \pi_1 (B \pi )$ determines a fibre bundle  $E_{\varphi}$ over $Y$ with fibre $E _{\pi}$:
$$  E_{\varphi} : = ( \tilde{Y} \times E _{\pi} ) / \pi_1(Y)  \ra Y =  ( \tilde{Y} ) / \pi_1(Y) = Y,$$
where the action of $\ga \in \pi_1(Y)$ is as follows:  $ \ga (y' , v) =  (\ga (y'), \varphi(\ga) (v))$.

While topology deals with continuous maps, when dealing with manifolds more regularity is wished for.
For instance, when we choose for $Y$   a differentiable manifold $M$, and the group $\pi$ is
abelian and torsion free, say $\pi = \ZZ^r$, then a more precise incarnation of the above theorem is given by the De Rham theory. 

In fact,  a homomorphism $\varphi :  \pi_1(Y) \ra  \ZZ^r$ factors through the Abelianization $H_1 (Y, \ZZ)$ of the fundamental group.
Since $H^1 (Y, \ZZ)= \Hom ( H_1 (Y, \ZZ), \ZZ)$, 
 $\varphi$ is equivalent to giving an element in $$ \varphi \in H^1 (Y, \ZZ)^r \subset H^1 (Y, \RR)^r \cong H^1_{DR} (Y, \RR)^r,$$
 where  $H^1_{DR} (Y, \RR)$ is the quotient space of the space of closed differentiable 1-forms modulo exact 1-forms.

In this case the classifying space is a real torus $$T_r : = \RR^r /  \ZZ^r.$$
Observe however that to give  $\varphi :  \pi_1(Y) \ra  \ZZ^r$ it is equivalent to give its $r$ components $\varphi_i , i = 1, \dots, r$, which are homomorphisms into $\ZZ$,
 and giving a map to $T_r : = \RR^r /  \ZZ^r$ is equivalent to giving $r$ maps to $T_1 : = \RR /  \ZZ$: hence  we may restrict ourselves to consider the case $r=1$.

Let us sketch the basic idea of the previous theorem \ref{classifying map} in this special case.
Let us assume that $Y$ is a cell complex, and define as usual $Y^j$ to be its j-th skeleton,
the union of all the cells of dimension $i \leq j$.

Since the fundamental group of $Y$ is generated by the free group $ F : = \pi_1 (Y^1)$, we get a homomorphism $\Phi : F \ra \ZZ$
inducing $\varphi$.
For each 1-cell $\ga \cong S^1$ we send $\ga \ra S^1$ according to the map $ z \in S^1 \mapsto z^m  \in S^1$, where $ m = \Phi (\ga)$.

In this way we get a continuous map $f^1 : Y^1 \ra S^1$, and we want to extend it inductively to $Y^j$ for each $j$.
Now,  assume that $f$ is already defined on $Z$, and that you are attaching an n-cell to $Z$, according to a continuous map $ \psi : \partial (D^n)= S^{n-1}  \ra Z$.
In order to extend $f$ to $ Z \cup_{\psi} D^n$ it suffices to extend the map $ f \circ {\psi}  $ to the interior of the disk $D^n$. This is possible
once the map $ f \circ {\psi}  : S^{n-1} \ra S^1$ is homotopic to a constant map. 
Now, for $n=2$, this condition holds by assumption: since $ {\psi}  (S^1)$ yields a relation for $\pi_1(Y)$, therefore its image under $\phi$ must be equal to zero.

For higher $n$, $n \geq 3$,  it suffices to observe that a continuous map $ h : S^{n-1} \ra S^{1}$ extends to the interior always: since $S^{n-1}$
is simply connected, and $S^1= \RR / \ZZ $,  $h$ lifts to a continuous map $ h' : S^{n-1} \ra \RR$, and  we can extend $h'$ to $D^n $ by setting 
$$ h' (x) : = |x|  h' ( \frac{x}{|x|}) , \forall x \neq 0, \  h'(0) : = 0.$$
($h$ is then the composition of $h'$ with the projection $ \RR \ra \RR / \ZZ = S^1$). 

In general,  when both $Y$ and  the classifying space $X$ ( as $T_1$ here) are differentiable  manifolds,
then each continuous map $ f : Y \ra X$ is homotopic to a differentiable map $f'$.
Take in fact  $ X \subset \RR^N$ and observe that the implicit function theorem implies that there is a tubular neighbourhood
$ X \subset  \sT_X  \subset \RR^N$ diffeomorphic to a tubular neighbourhood  of $X$ embedded as the $0$-section
of the normal bundle $N_X$ of the embedding $ X \subset \RR^N$.

Therefore, approximating the function $ f : Y \ra  \RR^N$ by a differentiable function $f''$ with values in $\sT_X$, 
we can use the bundle projection $ N_X \ra X$ to project $ f''$ to a differentiable function $ f' : Y \ra X$,
and similarly we can project the natural homotopy between $f$ and $ f''$, $ f(y) + t f'' (y)$ to obtain a homotopy between $f$
and $f'$.

Once we have a differentiable map $ f' : Y \ra T_1 = \RR / \ZZ$, we simply take the lift $\tilde{f'} : \tilde{Y} \ra \RR$,
and
the differential $d \tilde{f'}$ descends to a closed differential form $\eta$ on $Y$ such that its integral over a closed loop $\ga$ is just $\phi (\ga)$.

We obtain  the

\begin{prop}
Let $Y$ be a differentiable manifold, and let  $X$ be a differentiable manifold  that is a  $ K (\pi, 1)$
space: then, choosing base points $y_0 \in Y, x_0 \in X$,  one has a bijective  correspondence 
$$[ (Y, y_0 ), (X, x_0 )]^{diff} \cong \Hom (\pi_1(Y), \pi ) , [f]  \mapsto   \pi_1 (f) ,$$
where $ [ (Y, y_0 ), (X, x_0 )] ^{diff} $ denotes the set of differential homotopy classes $[f]$ of differentiable maps $ f : Y \ra X$ such that $ f(y_0 ) = x_0$.

In the case where $ X$ is a torus $T^r = \RR^r / \ZZ^r$, then $f$ is obtained as the projection  onto $T^r$ of 
$$ \tilde{\phi} (y) : = \int_{y_0}^y (\eta_1, \dots, \eta_r), \eta_j \in H^1 (Y, \ZZ) \subset H^1_{DR} (Y, \RR). $$

\end{prop}

\begin{rem}
In the previous proposition, $\eta_j$ is indeed a closed $1$-form, representing a certain De Rham cohomology class with integral periods
( i.e., $ \int_{\ga} \eta_j = \varphi (\ga) \in \ZZ$, $\forall \ga \in \pi_1 (Y)$).  Therefore $f $ is defined by $ \int_{y_0}^y (\eta_1, \dots, \eta_r) \mod (\ZZ^r)$.
Moreover, changing  $\eta_j$ with another form 
$\eta_j + d F_j$ in the same cohomology class, one finds a homotopic map, since   $ \int_{y_0}^y  (\eta_j  + t dF_j)= \int_{y_0}^y  (\eta_j ) + t( F_j(y) - F_j (y_0))$.
\end{rem}

 Before we dwell into a review of results concerning higher regularity of the classifying maps,
 we consider in the next section the basic examples of projective varieties that are classifying spaces.
 
\section{Projective varieties which are $K(\pi, 1)$}

The following are the easiest examples of projective varieties which are  $K(\pi,1)$'s. 
\begin{enumerate}
\item
Projective curves $C$ of genus $g \geq 2$.

By the {\bf Uniformization theorem}, these have  the Poincar\'e upper half plane
$ \HHH : =\{ z \in \CC  | Im (z) > 0 \}$ as universal covering, hence they are  compact quotients $C = \HHH /\Ga$, where
$ \Ga \subset \PP SL (2, \RR)$ is a discrete subgroup isomorphic to the fundamental group of $C$, $\pi_1(C) \cong \pi_g $.
Here $$\pi_g  : = \langle \al_1, \be_1, \dots \al_g,\be_g | \Pi_1^g [ \al_i, \be_i]  = 1\rangle $$
contains no elements of finite order.  Hence, given a faithful action of $\pi_g $ on $\HHH$,  it follows that necessarily $\Ga$ acts freely on $ \HHH $. Moreover, the quotient
must be compact, otherwise $C$ would be homotopically equivalent to a bouquet of circles, hence $H_2(C, \ZZ) = 0$,
a contradiction, since  $H_2(C, \ZZ) \cong  H_2(\pi_g, \ZZ) \cong \ZZ$, as one sees taking the standard realization of
a classifying space for $\pi_g$ by glueing the $2g$ sides of a polygon in the usual pattern.

Moreover, the complex  orientation of $C$ induces a standard generator $[C] $ of  $H_2(C, \ZZ) \cong  \ZZ$,
the so-called fundamental class.
\item
AV : = Abelian varieties.

More generally, a complex torus $ X = \CC^g / \Lam$, where $\Lam $ is a discrete subgroup of maximal rank (isomorphic then to $\ZZ^{2g}$),
is a K\"ahler classifying space $ K (\ZZ^{2g},1)$, the K\"ahler metric being induced by the translation invariant Euclidean metric
$ \frac{i}{2} \sum_1^g dz_j \otimes d \overline{ z_j}$.

For $g=1$ one gets in this way all projective curves of genus $g=1$;
but, for $ g > 1$, $X$ is in general not projective: it is projective, and called then an Abelian variety, if it satisfies
the Riemann bilinear relations.
These amount to the existence of a   positive definite  Hermitian form $H$ on  $\CC^g$ whose imaginary part $A$  ( i.e., $H = S + i A$), 
takes integer values on $\Lam \times \Lam$. In modern terms, there exists a positive line bundle $L$ on $X$, with Chern class $A \in H^2 (X, \ZZ) =
H^2 (\Lam , \ZZ) = \wedge^2 (\Hom (\Lam, \ZZ))$, whose curvature form, equal to  $H$, is positive (the existence of a positive line bundle 
on a compact complex manifold $X$ implies that $X$ is projective algebraic, by
 Kodaira's
theorem, \cite{kodemb}). 

We shall indeed see (\ref{abelianKp1}) that Abelian varieties are exactly the projective $K(\pi, 1)$ varieties, for which $\pi$ is an abelian group.
\item
LSM : = Locally symmetric manifolds.

These are the quotients of a {\bf bounded symmetric domain} $\sD$ by a cocompact discrete subgroup $\Ga \subset Aut (\sD)$
acting freely.  Recall that a bounded symmetric domain $\sD$ is a bounded domain $\sD \subset \subset \CC^n$ such that its
group $Aut (\sD)$ of biholomorphisms contains for each point $p \in \sD$, a holomorphic automorphism 
$ \s_p$ such that $ \s_p (p)= p$, and such that the derivative of $ \s_p $ at $ p$ is equal to $ - Id$. This property implies that $\sigma$ is an involution (i.e., it has order 2), and that
$Aut (\sD)^0$ (the connected component of the identity) is transitive on $\sD $ , and one can write $\sD = G/K $, where $G$ is a 
connected Lie group, and $K$ is a maximal compact subgroup.

The two important properties are:

(3.1) $\sD$ splits uniquely as the product of irreducible bounded symmetric domains.

(3.2) each such $\sD$ is contractible, since there is a Lie subalgebra $\sL$ of the Lie algebra $\mathfrak G$ of $G$ such that the exponential
map is a homeomorphism $\sL \cong \sD$. Hence $X$ is a classifying space for the group $\Ga \cong \pi_1(X)$.

Bounded symmetric  domains were classified by Elie Cartan in
\cite{Cartan}, and  there is only   a finite number of them (up to isomorphism) for each dimension $n$.

Recall the notation for the irreducible domains:

\begin{enumerate}
\item[(i)] 
   $I_{n,p} $ is the domain $ \sD = \{ Z \in Mat (n,p, \mathbb{C}) :
\mathrm{I}_p - ^tZ \cdot \overline{Z} > 0 \}$.\\
\item[(ii)] 
   $II_{n} $ is the intersection of the domain  $I_{n,n} $ with the
subspace of skew symmetric matrices.
 \item[(iii)] 
   $III_{n} $ is instead the intersection of the domain  $I_{n,n} $
with the subspace of  symmetric matrices.

\item[(iv)]   The Cartan - Harish Chandra realization of a domain of type $IV_{n}$ in $\mathbb{C}^{n}$ is the subset $\sD$ defined by the inequalities
(compare \cite{Helgason2}, page 527)

\vspace{.3cm}

\begin{tabular}{c}
 $|z_1^2 + z_2^2 + \cdots + z_n^2 | < 1$ , \\
  $1 + | z_1^2 + z_2^2 + \cdots + z_n^2 |^2 - 2\left( |z_1|^2 + |z_2|^2 + \cdots + |z_n|^2 \right) > 0 \, .$ \\
\end{tabular}\\

\item[(v)] $\sD_{16}$ is the exceptional domain of dimension $d=16$.

\item[(vi)] $\sD_{27}$ is the exceptional domain of dimension $d=27$.

\end{enumerate}
We refer the reader to
\cite{Helgason2}, Theorem 7.1, page 383 and exercise D, pages 526-527, and \cite{Roos} page 525 for a list of these irreducible bounded symmetric domains,
and a description of all of them as homogeneous spaces $G/ K$. In this context the domains are also called {\bf Hermitian symmetric spaces of non compact type}. Each of these 
is contained in the so-called {\bf compact dual}, which is a Hermitian symmetric spaces of  compact type. 

The easiest example is, for type I, the Grassmann manifold.
For type IV, the compact dual of $\sD$ is the hyperquadric $Q^n \subset \mathbb{P}^{n+1}$ defined by the polynomial $\sum_{j=0}^{n-1} X_j^2 - X_n^2 - X_{n+1}^2$. Notice that $\mathrm{SO}_0(n,2) \subset \mathrm{Aut}(Q^n)$.
The Borel embedding $j : \sD \rightarrow Q^n$ is given by \[ j(z_1,\cdots,z_n) = [2z_1 : 2z_2 : \cdots : 2z_n : \mathrm{i}(\Lambda - 1) : \Lambda + 1 ] \, ,\]
where $\Lambda := z_1^2 + \cdots + z_n^2$.
The map $j$ identifies the domain $\sD$ with the $\mathrm{SO}_0(n,2)$-orbit of the point $[0:0:\cdots:1:\mathrm{i}] \in Q^n$, i.e. $\sD \cong \mathrm{SO}_0(n,2)/\mathrm{SO}(n)\times \mathrm{SO}(2)$.

Among the   bounded symmetric domains are the so called {\bf bounded symmetric domains of tube
type}, those which are biholomorphic to a {\bf tube domain}, a generalized
Siegel
upper half-space
$$ T_{\sC} =   \mathbb{V} \oplus  \sqrt -1  \sC$$ where $\mathbb{V}$ is
a real vector
space and $\sC \subset \mathbb{V}$ is a { \em symmetric cone}, i.e., a self dual homogeneous convex cone
containing
no full lines.

In the case of type III domains, the tube domain is Siegel's  upper half space:
$$  \sH_g : = \{ \tau \in \Mat(g,g,\CC)| \tau = \ ^t\tau, \im (\tau) > 0 \},$$ 
a generalisation of the upper half-plane of Poincar\'e.

Borel proved in  \cite{Bo63} that for each bounded symmetric domain
$\sD$ there exists  a compact free
quotient $ X
= \sD / \Ga$,  called a compact Clifford-Klein form  of the
symmetric domain $\sD$.

 A classical result of J. Hano (see \cite{Hano} Theorem IV, page 886, and Lemma 6.2, page 317 of  \cite{milnorcurv}) asserts that a bounded homogeneous domain
that is the universal cover of  a compact complex manifold is symmetric.

\item
A particular, but very explicit case of  locally symmetric manifolds  is given by the
VIP : = Varieties isogenous to a product.

These were  studied in \cite{isogenous}, and they are defined as quotients 
$$X = (C_1 \times C_2 \times \dots \times C_n) / G$$ 
of the product of projective curves $C_j$ of respective genera $g_j \geq 2$ by the action of a finite 
group $G$  acting freely on the product.

In this case the fundamental group of  $X$ is not so mysterious and fits into an exact sequence
$$1 \ra  \pi_1 (C_1 \times C_2 \times \dots \times C_n) \cong \pi_{g_1} \times \dots \times \pi_{g_n}    \ra \pi_1 (X) \ra G  \ra 1.   $$ 

Such varieties are said to be of the {\bf unmixed type} if the group $G$ does not permute the factors, i.e., 
there are actions of $G$ on each curve such that 
$$ \ga (x_1, \dots, x_n ) =   ( \ga x_1, \dots, \ga x_n  ) , \forall \ga \in G.$$ 

Equivalently, each individual subgroup $\pi_{g_j}$ is normal in $\pi_1 (X)$.
\item
Kodaira fibrations $ f : S \ra B$.

 Here $S$ is a smooth projective surface and all the fibres of $f$ are smooth curves of genus $ g \geq 2$,
 in particular $f$ is a differentiable fibre bundle. Unlike  the examples above, where 
 $ C= (C_1 \times C_2) / G \ra C_2 / G$ is a holomorphic fibre bundle with fibre $C_1$ if $G$ acts freely on $C_2$,
 the second defining property for Kodaira fibrations is that the fibres are not all biholomorphic to each other.
 
 These Kodaira fibred  surfaces $S$ are very interesting topological objects: they were constructed by Kodaira (\cite{Kodsurf}) as a counterexample
to the conjecture that the index (of the cup product in middle cohomology) would  be multiplicative for fibre bundles.
In fact, for curves the pairing $H^1(C, \ZZ) \times H^1(C, \ZZ) \ra H^2(C, \ZZ) \cong \ZZ$ is skew symmetric, hence it has index zero; while Kodaira showed that
the index of the cup product $H^2(S, \ZZ) \times H^2(S, \ZZ) \ra H^4(S, \ZZ) \cong \ZZ$ is strictly positive.

The fundamental group of $S$ fits obviously into an exact sequence:

$$1 \ra   \pi_{g}    \ra \pi_1 (S) \ra  \pi_{b}  \ra 1,  $$
where $g$ is the fibre genus and $b$ is the genus of the base curve $B$, and it is known that $b \geq 2$, $g \geq 3$ (see \cite{CR}, also for more constructions
and a thorough discussion of their moduli spaces).

By simultaneous uniformization (\cite{Bers}) the universal covering $\tilde{S}$ of a Kodaira fibred surface $S$ is biholomorphic
to a bounded domain in $\CC^2$ (fibred over the unit disk $\De : = \{ z \in \CC | |z| < 1\}$ with fibres isomorphic to $\De$), which is not homogeneous.

\item
Hyperelliptic surfaces: these are the quotients of a complex torus of dimension 2  by a finite 
group $G$ acting freely, and in such away that the quotient is not again a complex  torus.

These surfaces were classified by Bagnera and de Franchis (\cite{BdF}, see also \cite{Enr-Sev} and \cite{bpv}) and they are obtained as quotients $(E_1 \times E_2)/G$
where  $E_1, E_2$ are  two elliptic curves, and $G$ is an abelian group acting on $E_1$ by translations, and on $E_2$ effectively and in such a way
that $ E_2/G \cong \PP^1$.  

\item
In higher dimension we define the 
Generalized Hyperelliptic Varieties (GHV) as quotients $A/G$ of an Abelian Variety $A$ by a finite group $G$ acting freely,
and with the property  that $G$ is not a subgroup of the group of  translations. Without loss of generality one can then assume
that $G$ contains no translations, since the subgroup $G_T$ of translations in $G$ would be a normal subgroup, and
if we denote $G' = G/G_T$, then $A/G = A' / G'$, where $A'$ is the Abelian variety $ A' : = A/G_T$.

We propose instead  the name {\bf Bagnera-de Franchis (BdF) Varieties} for those quotients $X = A/G$ were $G$ contains no translations,
and $G$ is a cyclic group of order $m$, with generator $g$ (observe that, when $A$ has  dimension $n=2$, the two notions coincide, thanks
to the classification result of Bagnera-De Franchis in \cite{BdF}).

A concrete description of such Bagnera-De Franchis varieties  shall be given in the following section.

We end this section giving an example of a projective variety which is not a $K(\pi,1)$, thus showing that the property of being a projective classifying space is lost after
taking hyperplane sections.

\begin{prop}
Let $n \geq 3$, and consider $X = C_1 \times \dots \times C_n$, the product of $n$ projective curves $C_i$ of respective genera $g_i \geq 1$.
Let $X \subset \PP^N$ be a projective embedding, and let $S$ be a smooth surface, obtained taking  the complete intersection of $X$ with
$n-2$ hypersurfaces. Then $\pi_2 (S) \neq 0$, in particular, $S$ is not a projective classifying space.
\end{prop}

\Proof
By the theorem of Lefschetz $\pi_1(S) \cong \pi_1(X)$, hence the universal covering $\tilde{S}$ of $S$ is a closed complex submanifold of 
$\tilde{X} = \CC^r \times \HH^{n-r}$. Hence $\tilde{S}$ is a Stein manifold, therefore (see \cite{AF1}) it has the homotopy type of a CW complex of real dimension $\leq 2$.

Since $\pi_1 (\tilde{S}) = 0$, we claim that $H_2 (\tilde{S}, \ZZ) = \pi_2 (\tilde{S}) = \pi_2 (S) \neq 0$ (the first isomorphism follows from the theorem of Hurewicz).

Otherwise, all homology and homotopy groups of $\tilde{S}$ would be trivial and $\tilde{S}$ would be contractible, $S$ would be a $K(\pi,1)$,
hence $H^* (S, \ZZ) = H^* (\pi_1(S), \ZZ) = H^* (X, \ZZ) $. This is a contradiction, since $H^{2n} (X, \ZZ) \neq 0, \ H^{2n}  (S, \ZZ) = 0$,
by our hypothesis that $ n \geq 3$.

\qed

\section{A trip around   Bagnera-de Franchis Varieties and group actions on Abelian Varieties}

\subsection{Bagnera-de Franchis Varieties}
Let $A/G$  be a Generalized Hyperelliptic Variety.
An easy observation is that any $g \in G$  is induced by an affine transformation $ x \mapsto \al x + b $ on the universal cover, hence it 
does not have a fixed point on $ A =  V/ \Lambda$ if and only if there is no solution of the equation
$$ \exists x \in V , \  g (x) \equiv  x  \ (\mod  \Lambda)$$ $$ \Leftrightarrow  \exists x \in V, \la \in \Lambda, \  (\al - \Id) x = \la - b.$$
This remark implies that $1$ must be an eigenvalue of $\al$ for all non trivial transformations $ g \in G$.

Since $ \al \in GL (\Lambda) \cong GL (2n, \ZZ)$, the eigenspace $V_1 = Ker ( \al - \Id) $ is a complex subspace defined over $\QQ$,
hence we have an Abelian subvariety $$A_1 \subset A, A_1 := V_1 / \Lambda_1,  \   \Lambda_1 : = (V_1 \cap \Lambda).$$

While what we said up to now was valid for any complex torus, we replace this assumption by the stronger assumption that $A$ is an Abelian variety.
This assumption allows us to invoke Poincar\'e' s complete reducibility theorem, so that we can split 
$$  V = V_1 \oplus V_2, \ V_2 : = (V_1)^{\perp}, A_2 : = V_2 /  \Lambda_2, \ \Lambda_2 : = (V_2 \cap \Lambda).$$

We get then an isogeny $  A_1 \times A_2 \ra A,$ with kernel $ T: = \Lambda / ( \Lambda_1 \oplus  \Lambda_2)$.
Observe that $ T \cap A_j = \{ 0 \}, j=1,2$.

In view of the splitting $  V = V_1 \oplus V_2$, we can write, after a change of the origin in the affine space $V_2$, 
$$ \al (v_1 + v_2) = v_1 + \al_2 v_2, g(v) = v_1 + b_1 +  \al_2 v_2, b_1 \in V_1, b_1  \notin \Lambda.$$

If $\al_2$ has order equal to $m$, then necessarily the image of $b_1$ has order exactly $m$ in $A$,
by virtue of our  assumption that $G$ contains no translations. In other words, $b_1$ induces a translation on $A_1$
of order exactly $m$. 

Now, $g $ lifts naturally to $A_1 \times A_2$, by 
$$ g (a_1, a_2) = (a_1 + [b_1], \al_2 a_2),$$
where $[b_1] $ is the class of $b_1$ in $A_1$.

We reach the conclusion that $ X = A /G = ((A_1 \times A_2)/ T) /G$, where $T$ is the finite group of translations
$ T =  \Lambda / ( \Lambda_1 \oplus  \Lambda_2)$. 

Conversely, given such an automorphism $g$  of $A_1 \times A_2$, it descends to $ A : = (A_1 \times A_2)/ T$
if and only if the linear part of $g$ sends $T$ to $T$.

Denote now by $T_j$ the (isomorphic) image of $T \ra A_j$: then $T \subset T_1 \times T_2$ is the graph of an isomorphism 
$\phi : T_1 \ra T_2$, hence the condition which allows $g$ to descend to $A$ is that:
$$ (**)  (\Id \times \al_2)  (T) = T \Leftrightarrow \al_2 \circ \phi = \phi  \Leftrightarrow$$
$$ \Leftrightarrow  ( \al_2 - \Id) \circ \phi = 0  \Leftrightarrow  ( \al_2 - \Id) T_2 = 0.$$

We are in the position to  illustrate  the standard  example, before we give the more general  more complete description of BdF varieties.

The basic example is the one where $m=2$, hence $\al_2$ is scalar multiplication by $-1$. Then $\phi = - \phi$ implies that
$T_2, T_1$ are 2-torsion subgroups. Then also $ 2 [b_1] = 0$ implies that $[b_1]$ is a 2-torsion element.
However $ [b_1]$ cannot belong to $T_1$, else $g : A \ra A$ would be induced by 
$$  (a_1, a_2 ) \mapsto ( a_1 + [b_1], \al_2 a_2) \equiv ( a_1,  \al_2 a_2 - \phi ([b_1]) ) \ (mod \ T) ,$$
which has a fixed point on $A$.

We conclude that the standard example of Bagnera-de Franchis Varieties of order $m=2$  is the following:
$$ X = A / g, \ {\rm where} \  A : = (A_1 \times A_2 )/ T , \ {\rm and} \ $$ 
$$  \ T = \{ (t_1, \phi (t_1))\}, T_j \subset A_j[2] , \phi : T_1 \cong T_2,  \be_1 \in A_1[2] \setminus T_1,$$
$$  g : A \ra A , \ \  g(a_1, a_2) = (a_1 + \be_1, - a_2). $$  
  
\end{enumerate}

In order to conclude appropriately the above discussion, we give some useful definition.

\begin{defin}
We define first a  Bagnera-de Franchis manifold (resp.: variety)  {\bf  of product type } as a quotient $ X= A/G $ where $A = A_1 \times A_2$,  $A_1, A_2$ are complex tori
(resp.: Abelian Varieties), 
and $G \cong \ZZ/m$ is a cyclic group operating freely on $A$, generated by an automorphism of the form
$$ g (a_1, a_2 ) = ( a_1 + \be_1, \al_2 (a_2)) ,$$
where $\be_1 \in A_1[m]$ is an element of order exactly $m$, and similarly $\al_2 : A_2 \ra A_2$ is a linear automorphism
of order exactly $m$ without $1$ as eigenvalue (these conditions guarantee that the action is free). If moreover all eigenvalues  of $\al_2 $
are primitive $m$-th roots of $1$, we shall say that  $ X= A/G $ is a {\bf primary}  Bagnera-de Franchis manifold.

\end{defin}

We have the following proposition, giving a characterization of Bagnera- De Franchis varieties.

\begin{prop}\label{BdF}
Every Bagnera-de Franchis variety $ X= A/G $, where  $G \cong \ZZ/m$ contains no translations, is the quotient of a  Bagnera-de Franchis variety
of product type,  $(A_1 \times A_2)/ G$ by any finite subgroup $T$ of $A_1 \times A_2$ which satisfies the following properties:

1) $T$  is the graph of an isomorphism between two respective 
subgroups $ T_1 \subset A_1,  T_2 \subset A_2,$

2) $(\al_2 - \Id) T_2 = 0$

3) if $ g (a_1, a_2 ) = ( a_1 + \be_1, \al_2 (a_2)) ,$ then the subgroup of order $m$ generated by $\be_1$ intersects $T_1$ only in $\{0\}$.

In particular, we may write $X$ as the quotient $ X =  (A_1 \times A_2)/ (G \times T)$ by the abelian group $G \times T$. 

\end{prop}

\subsection{Actions of a finite group on an Abelian variety}

Assume that we have the  action of a finite group $G$ on  a complex torus 
 $ A = V / \Lam$. Since every holomorphic map between complex tori lifts to a complex  affine map
 of the respective universal covers, we can attach  to the group $G$  the group of affine transformations
 $\Ga$ which fits into the exact sequence:
 
 $$0 \ra \Lam \ra \Ga \ra G \ra 1.  $$
 $\Ga$ consists of all affine transformations of $V$ which lift transformations of the group $G$.
 
 Define now $G^0$ to be the subgroup of $G$ consisting of all the translations in $G$.
 
 \begin{prop}\label{notranslations}
 The abstract group $\Ga$ determines an exact sequence 
 $$0 \ra \Lam' \ra \Ga \ra G' \ra 1  $$
 such that $\Lam' \subset  \Lam$, $\Lam' $ is a lattice in $V$, $\Lam' / \Lam = G^0 $, $ G' \subset Aut (V/ \Lam')$
 contains no translations.
  \end{prop}
 
 \Proof
 
 It is clear that $ V = \Lam \otimes_{\ZZ} \RR$ as a real vector space, and we denote by $V_{\QQ}: = 
 \Lam \otimes \QQ$.  Consider the homomorphism $\al_L$ associating to an affine transformation its linear part, and 
 let $$\Lam' : = \ker (  \alpha_L : \Ga \ra \GL  (V_{\QQ}) \subset \GL  (V)) ,$$
$$\overline{G}_1 : = \im (  \alpha_L : \Ga \ra \GL  (V_{\QQ})).$$
$\Lam'$ is a subgroup of the group of translations in $V$, hence it is obviously Abelian,  and maps isomorphically onto a lattice of $V$ 
which contains $\Lambda$.  We shall identify this lattice with $\Lambda' $, writing with shorthand notation $\Lambda' \subset V$.

In turn $ V = \Lam'  \otimes_{\ZZ} \RR$, and, if $ G' : = \Ga / \Lam'$, then $ G' \cong \overline{G}_1 $ and we have  the exact sequence 
$$0 \ra \Lam' \ra \Ga \ra G' \ra 1  ,$$
yielding  an embedding $ G' \subset GL(\Lam')$.
 
There remains  to show that $\Lam'$ is determined by $\Ga$ as an abstract group, independently
of the exact sequence we started with. In fact, one property of $\Lam'$ is that it is a maximal abelian subgroup,
 normal and of finite index.

Assume that $\Lam''$ has the same property: then their intersection $\Lam^0 : = \Lam' \cap \Lam''$ is a normal subgroup of
 finite index, in particular $\Lam^0  \otimes_{\ZZ} \RR = \Lam'  \otimes_{\ZZ} \RR = V $;
hence $\Lam'' \subset \ker (  \alpha_L : \Ga \ra \GL  (V))= \Lam'$,
where $\al_L$ is induced by conjugation on $\Lam^0$
By maximality $\Lam' = \Lam''$.

\qed

In the case where $A$ is an Abelian variety, we can say a little bit more about the affine representation of the group $\Ga$.

\begin{lemma}\label{torsiontranslations}
Let $A$ an Abelian variety, and let $G$ be a finite group of automorphisms. Then there is a positive integer $r$ such that
the group $\Ga$ of affine transformations which are lifts of transformations in $G$ satisfies
  $$ \Ga \subset Aff (\frac{1}{r} \Lam).$$
  In other words, we can write any transformation $g \in G$ as induced by the affine transformation of $V$:
  $ v \mapsto \al (v) + \frac{\la}{r}, \la \in \Lam$.
\end{lemma}

\Proof
Let $L$ be an ample divisor class on $A$. Since $G$ is finite, there is a $G$-invariant such class $L$
(simply replace $L$ by $\sum_{g\in G} g^* (L') = : L'$).

We can further assume that the class of $L$ is not only $G$-invariant , but also indivisible.
There is a homomorphism $\phi_L : A \ra A^{\vee}$, such that $\phi_L(a) : = T_a^* L - L \in Pic^0(A) = A^{\vee}$,
and where $ T_a$ is the translation by $a$, $ v \mapsto v + a$.

Let $K_L : = \ker (\Phi_L)$, which is a finite group, and denote by $r$ its exponent (hence,  $ r K_L = 0$
and $K_L \subset  (\frac{1}{r} \Lam)/ \Lam$). 

Represent now $L$ by a line bundle whose cocycle is in Appell-Humbert normal form
$ f_{\la} (z) =  \rho (\la) exp ( \pi H(z,\la) + \frac{\pi}{2} H (\la, \la))$, and write a transformation $g\in G$
as induced by  $ v \mapsto \al (v) + b$. The condition that the Chern class of $L$ is $G$-invariant implies that 
the Hermitian form $H$ on $V$ is left invariant by $\al$ (i.e., we have a group of isometries of $V$ for the metric given by the Hermitian form $H$).

Hence also the translation part leaves $L$ invariant, therefore the class of $b$ lies in  $K_L \subset  (\frac{1}{r} \Lam)/ \Lam$.

\qed

Observe now that, given an affine group $\Ga \subset Aff ( \Lam \otimes_{\ZZ} \RR)$ as above,  in order to obtain the structure of a complex torus on $V / \Lam'$, we must give a complex structure
on $V$ which makes the action of $G' \cong \overline{G}_1$ complex  linear. 

 In order to study the moduli spaces of the associated complex manifolds, we introduce  therefore a further invariant, called Hodge type, according to
the following
definition. 

\begin{defin}\label{HodgeT}
(i) Given a  faithful representation $ G \ra Aut (\Lam)$, where $\Lam$ is a free abelian group  of even rank $2n$,
 a {\bf $G$- Hodge decomposition} is a $G$-invariant decomposition 
 $$\Lam \otimes \CC = H^{1,0} \oplus H^{0,1}, \ H^{0,1} = \overline{H^{1,0}}.$$

(ii)
Write   $\Lam \otimes \CC$
as the sum of
isotypical components $$\Lam \otimes \CC = \oplus_{\chi \in \Irr (G)}
U_{ \chi}.$$
Write also  $U_{ \chi} = W_{ \chi} \otimes M_{ \chi}$, where $W_{ \chi} $ is the  irreducible representation corresponding to the character $\chi$,
and $ M_{ \chi} $ is a trivial representation whose  dimension is denoted $n_{ \chi} $.

Write accordingly $V : = H^{1,0} =  \oplus_{\chi \in \Irr (G)}
V_{ \chi},$ where $V_{ \chi} = W_{ \chi} \otimes M^{1,0}_{ \chi}$.

 Then the {\bf Hodge type} of the decomposition
is the datum 
of the dimensions   $$\nu ( \chi): = dim_{\CC} M^{1,0}_{ \chi} $$
corresponding to  the Hodge summands for non real representations (observe in fact  that one must have: $\nu ( \chi) + \nu ( \bar{\chi}) = \dim ( M_{ \chi})$).
\end{defin}

\begin{rem}\label{Hodgetype}
Given a  faithful representation $ G \ra Aut (\Lam)$, where $\Lam$ is a free abelian group  of even rank $2n$,
 all the  $G$- Hodge decompositions of a fixed Hodge type are parametrized by an open set in a product of Grassmannians.
 Since, for a non real irreducible representation $\chi$ one may simply choose $M^{1,0}_{ \chi}$ to be a complex subspace of dimension  $\nu ( \chi)$ of  $ M_{ \chi}$, and for $M_{ \chi} = \overline{M_{ \chi}}$,
one simply chooses a complex subspace $M^{1,0}_{ \chi}$ of middle dimension. Then the open condition is just that (since $  M^{0,1}_{ \chi} : = \overline {M^{1,0}_{ \chi} }$) we want 
$$ M_{ \chi} = M^{1,0}_{ \chi}  \oplus  M^{0,1}_{ \chi} ,$$ or, equivalently,  $ M_{ \chi} = M^{1,0}_{ \chi} \oplus  \overline{M^{1,0}_{\bar{ \chi}} } $.
\end{rem}

\subsection{The  general case where $G$ is Abelian.}

Assume now that $G$ is Abelian, and consider the linear  representation
$$ \rho : G \ra \GL (V).$$

We get a splitting 
$$ V = \oplus_{\chi \in G^{\vee}} V_{\chi}  = \oplus_{\chi \in \sX} V_{\chi}  ,$$ 
where $G^{\vee} : = Hom (G, \CC^*) = \{ \chi : G \ra \CC^* \}$ is the dual group of characters of $G$,
 $V_{\chi} $ is the $\chi$-eigenspace of $G$ ( for each $v \in V_{\chi}$, $ \rho(g) (v) = \chi (g) v$),
 and $\sX$ is the set of characters $\chi$ such that $ V_{\chi}  
 \neq 0$.
 
 \begin{prop}\label{affinetype}
 Assume that an abelian  group $G$ acts on a complex torus $ X : = V / \Lam$.
 
 Then the isomorphism class of the group $\Ga$ of affine transformations which are lifts of transformations in $G$  determines the real affine type of the action 
 of $\Ga$ on $V = ( \Lam \otimes \RR )$, in particular the above
 exact sequence determines the action of $G$ up to a real affine isomorphism of $A= 
( \Lam \otimes \RR )/ \Lam$.
 \end{prop}
 \Proof
 By proposition \ref{notranslations} we may reduce ourselves to the case where $G$ contains no translations.
 
 Step I): let us first treat the case where $G$ is cyclic, generated by an element $g'$ of order $m$.
 
 Then we claim that the power $g^m$ of a lift $g$ of $g'$ determines the affine type of the transformation $g$.
 
 In fact, if $ g(v) = \al v + b$, $V$ splits according to the eigenspaces of the linear map $\al$, which has order exactly $m$,
 and we can write $ V = V_1 \oplus V_2$, where $V_1 = \ker ( \al - Id)$.
 
 Then $ g (v_1 + v_2) = (v_1 + b_1 ) + (\al_2 (v_2) + b_2)$, and choosing as the origin in the affine space $V_2$ the unique fixed point
 of $g_2 : =  v_2 \mapsto \al_2 (v_2) + b_2$, we reach the normal form $ g (v_1 + v_2) = (v_1 + b_1 ) + \al_2 (v_2) $.
 Now,  $ g^m(v_1 + v_2) = (v_1 + m b_1 ) + v_2 $, hence $g^m$ determines $b_1$, and therefore also the normal form of $g$.
 
 Step II): split  $ V = \oplus _{\chi} V_{\chi} $ according the the characters of $G$, and choose, for each $\chi$,
 an element $g_{\chi}$ whose image under $\chi$ generates $ \chi (G) \subset \CC^*$.
 
 For each $\chi$, as in step I, choose as origin in the affine space $V_{\chi} $ a fixed point for $g_{\chi}$.
 
 Let now $g$ be an arbitrary element in $G$: then $g$ acts on $V_{\chi} $ by $ v \mapsto \chi(g) v + b_{\chi} $.
 
 We only need to show that $b_{\chi} $ is uniquely determined.
 
 This is clear, by step I, if $ \chi(g) =1$. Else, there is an integer $r$ such that $ \chi(g g_{\chi}^r) =1$,
 hence the affine action of $g g_{\chi}^r$ on $V_{\chi} $ is uniquely determined, therefore also the one of $g$.
 
 \qed
 
 \begin{rem}
 If $G$ is a general group, we can write  $ V = \oplus _{\rho} V_{\rho} $ where $V_{\rho} $ is irreducible (but the same irreducible representation may occur more than once)
 and let $H_{\rho} $ be a subgroup such that $\rho (H_{\rho} ) = \rho (G)$. Argueing as in step II above, the theorem holds for $G$
 if one can prove that the affine action of $H_{\rho} $ on $V_{\rho} $ is uniquely determined up to affine conjugation.
 
 I.e., it suffices to prove the result when $V$ is irreducible and  $ G \subset Aut (V)$.
 \end{rem}
 
  \begin{rem}\label{erratum}
  In the papers \cite{bc-inoue} and \cite{bcf} we claimed the validity of proposition \ref{affinetype} for any finite group $G$, but the proof was incorrect.
  Does the result hold for any finite group $G$, at least in the case of an Abelian variety $A$?
  We leave this question open, hoping to return on it.

 \end{rem}

In the case of generalized  hyperelliptic varieties, we require the condition that the action of $G$ is free: this  implies the following condition

$$ (1) \   \forall g \in G, \exists \chi  \in \sX, \ \ \chi (g) = 1, \Leftrightarrow \cup _{\chi \in \sX} \ker (\chi) = G, $$
while we may assume that $G$ contains no translations, i.e., that $\rho$ is injective (equivalently , $\sX$ spans  $G^{\vee}$).

With the above in mind, we pass to investigate in more detail  the case where $G$ is cyclic, the case where the quotient variety shall
be called a Bagnera- de Franchis variety.

\subsection{Bagnera-de Franchis Varieties of small dimension}

In view of the characterization given in proposition \ref{BdF}, we can see a Bagnera-de Franchis variety as the quotient of
 one of product type,  since the actions of $T$ and $g$ commute (by the property that $ \al_2 \circ \phi = \phi$).
 
 Dealing with appropriate choices of $T$ is the easy part, since,
as we saw,  the points $t_2$ of $T_2$ satisfy the property $ \al_2 (t_2)  = t_2$. It suffices to choose  $T_2 \subset  A_2[*]:= \ker (\al_2 - \Id_{A_2})$,
which is a finite  subgroup of $A_2$, and then  to pick an isomorphism $\psi : T_2 \ra T_1 \subset A_1$,
such that $ T_1 : = Im (\psi) \cap \langle \be_1 \rangle = \{ 0\}$.

Let us then restrict ourselves to consider Bagnera-de Franchis varieties of product type.

We show now how to further reduce to the investigation of primary Bagnera-de Franchis varieties.

In fact, in the case of a BdF variety of product type, $\Lambda_2$ is a $G $-module, hence a module over the group ring
$$ R : = R(m)  : = \ZZ[G]  \cong  \ZZ[x] / (x^m - 1).$$

This ring is in general far from being an integral domain, since indeed it can be written as a direct sum of cyclotomic rings, which 
are the integral domains
defined as  
$$R_m: = \ZZ[x] / P_m(x).$$
 Here $P_m(x)$ is the m-th cyclotomic polynomial
$$ P_m(x) = \Pi_{0 < j < m, (m,j)=1} (x - \e^j),$$ where $\e = \exp ( 2 \pi i /m)   $.

Then $$  R (m) = \oplus_ { k | m} R_k .$$ 

The following  elementary lemma, together with the splitting of the vector space $V$ as a direct sum of eigenspaces for $g$,
yields a decomposition of $A_2$ as a direct product $ A_2 =  \oplus_{ k | m} A_{2,k} $
of $G$ -invariant Abelian subvarieties $A_{2,k} $ on which $g$ acts with eigenvalues of order precisely $k$.

\begin{lemma}
Assume that $M$ is a module over a ring $R =  \oplus_k R_k .$

Then $M$ splits uniquely as a direct sum $M =  \oplus_k M_k $ such that $M_k$ is an $R_k$-module, and the
$R$-module structure is obtained through the projection $ R \ra R_k$.

\end{lemma}

\Proof
We can write the identity in $R$ as a sum of idempotents $ 1 = \Sigma_k e_k$, where $e_k$ is the identity of $R_k$,
and $ e_k e_j = 0$ for $ j \neq k$.

Then each element $w \in M$ can be written as
$$ w = 1 w = ( \Sigma_k e_k) w = \Sigma_k e_k w =  :  \Sigma_k w_k.$$ 

Hence $ M_k$ is defined as $ e_k M$.

\qed

In the situation where we have a primary Bagnera-de Franchis variety $\Lambda_2$ is a module over the integral domain $R : = R_m: = \ZZ[x] / P_m(x)$, 

Since $\Lambda_2$ is a projective module, a classical result (see \cite{Milnor}, lemmas 1.5 and 1.6) is that $\Lambda_2$  splits as the direct sum 
$\Lambda_2 = R^r \oplus I $ of a free module with an ideal
$ I \subset R$, and it is indeed free if the class number $ h(R)=1$ (to see for which  integers $m$ this occurs, see the table in \cite{washington}, page 353).
To give a complex structure to $A_2 : = (\Lambda_2 \otimes_{\ZZ} \RR )/ \Lambda_2$ 
it suffices to give a decomposition  $\Lambda_2 \otimes_{\ZZ} \CC = V \oplus \bar{V},$
such that the action of $x$ is holomorphic, which is equivalent to asking that $V$ is a direct sum of eigenspaces
$ V_{\la}$, for $\la =  \e^j$ a primitive m-th root of unity. 

Writing  $U: = \Lambda_2 \otimes_{\ZZ} \CC = \oplus U_{\la}$, the desired decomposition is obtained by choosing,
for each eigenvalue $\la$, a decomposition $  U_{\la} = U_{\la}^{1,0}  \oplus U_{\la}^{0,1} $ such that $ \overline{U_{\la}^{1,0}} = U^{0,1}_{\bar{\la}}$.

The simplest case (see \cite{cacicetraro} for more details) is the one where $I =0, r=1$, hence $  \dim (U_{\la} ) = 1$, so that we have only a finite number of complex structures,
depending on the choice of the $\varphi(m)/2$ indices $j$ such that $  U_{\e^j} = U^{1,0}_{\e^j}$.

The above discussion does however leave open the following  question (see  \cite{LR} for some partial results), which should have an affirmative
answer (in the sense that each family of complex structures should contain some Abelian variety).

{\bf Question: when is then the complex torus $A_2$ an Abelian variety?}

Observe that the classification in small dimension is possible thanks to the observation that the $\ZZ$- rank of $R $ (or of any ideal $I \subset R$) cannot exceed 
the real dimension of $A_2$: in other words we have 
$$\varphi(m) \leq 2 (n-1),$$
where $\varphi(m)$ is the Euler function, which is multiplicative for relatively prime numbers,
and satisfies $\varphi(p^r) = (p-1)  p^{r-1}$ when $p$ is a prime number.

For instance, when $n \leq 3$, then $\varphi(m) \leq 4$, and  $\varphi(p^r) \leq 4$ iff $p=3, 5 , r=1$, or $p=2, r \leq 3$.

Observe that the case $m=2$ was  completely described, in view of proposition \ref{BdF}, hence we may assume that $m \geq  3$
and ask for the classification for the Bagnera-de Franchis varieties (or manifolds) of order $m$.

Assuming then $m \geq  3$, if $ n \leq 3$, it is only possible to have $ m = 3,4,6, \varphi(m)= 2$, or  $ m = 5,8,10,12, \varphi(m)= 4$.
The classification is then made easier by the fact that, in the above range for $m$, $R$ is a P.I.D., hence every  torsion free module is free.
In particular $\Lam_2$ is a free $R$-module.

The classification for $ n =  4$, since we must have $\varphi(m) \leq 6$, is going to include also the case $m=7, 9$.

We are not aware of literature dedicated to a precise classification of Bagnera-de Franchis varieties, or generalised hyperelliptic varieties,
at least in dimension $>$ 3 (see however \cite{UchidaYoshihara} and \cite{Lange}  for  results in dimension 3).
Observe that the hypothesis that $G$ is a finite group allows to find a $G$-invariant Hermitian metric on $V$, hence the affine group
$\Ga$ extension of $\Lambda$ by $G$ is a Bieberbach group, and in each dimension we have only a finite number of those.

We end this section mentioning some elementary results which are useful to locate  the BdF varieties in the classification theory of algebraic varieties.

\begin{prop}
The Albanese variety of a Bagnera-de Franchis variety $X = A/G$ is the quotient $A_1 / (T_1 +  \langle \be_1 \rangle)$.
\end{prop}
\Proof
Observe that the Albanese variety  $H^0(\Omega^1_X)^{\vee}/ Im (H_1(X, \ZZ))$ of $X = A/G$ is a quotient of the vector space $V_1$ by the image of the fundamental group of $X$
(actually of its abelianization, the  first homology group $H_1(X, \ZZ)$):
since the dual of $V_1$ is the space of $G$-invariant forms on $A$, $H^0(\Omega^1_A)^G \cong H^0(\Omega^1_X)$.

We also observe that there is a well defined  map $X \ra A_1 / (T_1 + \La \be_1 \Ra)$, since $T_1$ is the first projection of $T$.
The image of the fundamental group of $X$ contains the image of $\Lambda$, which is precisely the extension of $\Lambda_1$
by the image of $T$, namely $T_1$. Since we have the exact sequence 
$$ 1 \ra \Lambda = \pi_1 (A) \ra   \pi_1 (X) \ra   G \ra 1$$
the image of the fundamental group of $X$ is generated by the image of $\Lambda$ and by the image of the transformation
$g$, which however acts on $A_1$ by translation by $\be_1 = [b_1]$.
\qed

\begin{rem}
Unlike the case of complex dimension $n=2$, there  are   Bagnera-de Franchis varieties $X = A/G$  with trivial canonical divisor,
for instance some examples are given by:

1) any  BdF variety which is standard (i.e., $m=2$) and is such that  $A_2$ has even dimension has trivial canonical divisor, as well as

2) some BdF varieties $ X = A/G = (A_1 \times A_2) / (T \times G)$ whose associated product type  BdF variety
$ Y = (A_1 \times A_2) /  G$ is primary and elementary (meaning that $\Lam_2 \cong R_m$). For instance,
when $m=p$ is an odd prime, there is a  primary and elementary  BdF variety with trivial canonical divisor  if and only if, given the set of  integers $\{ 1,2 ,  \dots ,  [\frac{p-1}{2}] \}$,
there is a partition of it into two sets whose respective sums are the same:
for $ m=7$, it suffices to choose the partition $\{ 1,2\} \cup \{3\}$.  This partition yields the sets $\{1,2,4\}$ and $\{3,5,6\}$, 
which determine $V= H^{1,0} $ as  $ V = H^{1,0}   : = U_1 \oplus U_2 \oplus U_4$.   Whereas, for $ m = 17$,
there are many choices, consisting  of the union of two doubletons of the form $ \{ j, 9 - j\}$.
\end{rem}

We do not further consider the question of determining exactly the case where the divisor $K_X$ is trivial, also because this question
can be asked also for the more general case where the action of $G$ is not free (see \cite{Oguiso} for a classification
in dimension $n=3$).

\section{ Orbifold fundamental groups and rational  $K(\pi, 1)$'s }

\subsection{Orbifold fundamental group of an action}

In the previous section we have considered quotients $ X = A/G$ of a complex torus  $A = V / \Lam $ by the free action of 
a finite group $G$. In this case the affine group $\Ga$ fitting into the exact sequence 
$$ 1 \ra \Lambda = \pi_1 (A) \ra   \Ga  \ra   G \ra 1$$
equals $\pi_1 (X)$.

We have also seen that we have the same exact sequence in the more general case where  the action of $G$ is no longer free,
and that  the group $\Ga$ determines in general the  affine action of $\Ga$ on the real affine space $V$ (resp. on the rational affine space $\Lam \otimes \QQ$)
and that, once the Hodge type of $V$ is fixed, these varieties are parametrized by a connected complex manifold.

In the case where the action of $G$ is no longer free, we would like to remember the group $\Ga$, in view of its importance,
even if it is no longer a `bona fide' fundamental group. This can be done through a more general
 correspondence which associates to the pair of $X$ and the group action of $G$ a group $\Ga$ which is called 
 the {\bf orbifold fundamental group} and  can be defined in many ways (see \cite{DM}, \cite{isogenous}). Here is one. 

\begin{defin}
Let $Z$ be a `good' topological space, i.e., arcwise connected and semi-locally 1-connected, so that there exists the
universal cover $\sD$ of $Z$. Then we have $ Z = \sD / \pi$ where $\pi : = \pi_1 (Z)$; denote $p :  \sD \ra Z$ the quotient projection.

Assume now that a group $G$ acts  properly discontinuously on $Z$, and set $ X : = Z/ G$.
Then we define the orbifold fundamental group of the quotient of $Z$ by $G$ as the group of all the possible liftings
of the action of $G$ on $\sD$, more precisely:

$$   \pi_1^{orb} ( Z, G): = \{ \ga : \sD \ra \sD | \exists g \in G, \ s.t. \ p \circ \ga = g \circ p  \}$$
\end{defin}

\begin{rem}
(1) We obtain an exact sequence, called orbifold fundamental group exact sequence, 
$$  (***)  \ 1 \ra \pi_1 (Z) \ra  \pi_1^{orb} ( Z, G) \ra G \ra 1 $$
since for each $ g \in G$ we have a lifting $\ga$ of $g$, because $\sD$ is the universal cover, so its associated
subgroup is the identity, hence the lifting property is ensured. Moreover the lifting is uniquely determined,
given base points $z_0 \in Z, y_0 \in \sD$, with $ p(y_0) = z_0$, by the choice of $y \in p^{-1} (\{ z_0\})$
with $\ga (y_0) = y$.

(2) In the case where the action of $G$ is free, then $\pi_1^{orb} ( Z, G) = \pi_1 ( Z/G) = \pi_1( X) $.

(3) If $G$ acts properly discontinuously on $Z$, then $\Ga : = \pi_1^{orb} ( Z, G) $ acts also 
properly discontinuously on $\sD$: in fact the defining property that, for each compact $K \subset \sD$,
$ \Ga(K,K) : = \{ \ga| \ga(K) \cap K \neq \emptyset \}$ is finite, follows since $p(K) $ is compact,
and since $\pi$ acts properly discontinuously.

(4) If the space $Z$ is a $ K ( \pi, 1)$, i.e., a classifying space $Z =  B \pi $, then (***) determines the topological type
of the $G$-action. The group $\Ga$ acts by conjugation on $\pi$, and this yields an homomorphism
$ G \ra \Out(\pi) = \Aut (\pi) / \Inn (\pi) $, and we just have to recall that the topological action of $\ga \in \Ga$
is determined by its action on $\pi$, taken up to $\Inn (\pi) $ if we do not keep track of the base point.

In fact  there is a bijection between homotopy classes of self maps of $Z$ and homomorphisms of $\pi$,
taken of course up to inner conjugation (inner conjugation is the effect of  changing the base point, and if we do not insist
on taking pointed spaces, i.e., pairs $(Z,z_0)$ the action of a continuous map on the fundamental group is only determined
up to inner conjugation). Clearly a homeomorphism $\fie : Z \ra Z$ yields then an associated element $\pi_1 (\fie) \in \Out (\pi)$.

(5) Even if not useful for computations, one can still interpret the exact sequence (***) as an exact sequence for 
the fundamental groups of a tower
of covering spaces, using the standard construction for equivariant cohomology that we shall discuss later.

Just let $ BG = EG / G$ be a classifying space, and consider the free product action  of $G$ on $Z \times EG$.

Then we have a tower of covering spaces
$$ (\sD \times EG) \ra ( Z \times EG) \ra   ( Z \times EG)/ G, $$ 
where $\sD \times EG$ is simply connected, and is the universal cover of $ Y : =  ( Z \times EG)/ G$,
such that $\pi_1 (Y) \cong \pi_1^{orb} ( Z, G)$.
\end{rem} 

Of course, the above definitions seem apparently unrelated to the fundamental group of $X$, 
however we can take the maximal open set $ U' \subset \sD$ where $\Ga : = \pi_1^{orb} ( Z, G) $ 
acts freely. Indeed $ U' = p^{-1} (U) $, where $U$ is the open set $ U : = Z \setminus \{ z | G_z : = Stab (z) \neq \{ \Id_G\} \}$.

The conclusion is that $ U' / \Ga = U / G =: U'' \subset X$, hence 
we get  exact sequences 

$$ 1 \ra \pi_1 (U) \ra  \pi_1 (U'')  \ra G \ra 1 , $$
$$1 \ra \pi_1 (U') \ra  \pi_1 (U'')  \ra \Ga \ra 1.$$

We make now the hypothesis that $\sD$, hence also $Z,X$ are normal complex spaces, in particular they are locally contractible,
and that the actions are holomorphic: then we have surjections  $$\pi_1 (U'') \ra \pi_1 (X), \ \pi_1 (U) \ra \pi_1 (Z),$$
moreover $\pi_1 (U') \ra \pi_1 (\sD) = \{1\}$, and finally we get 
$$\pi_1(Z) = \pi_1 (U) / \pi_1 (U'),   \ \pi_1^{orb} ( Z, G) \ra \pi_1(X) \ra 0.$$

\begin{ex}\label{orbicurves}
Assume that Z, X are smooth complex curves $Z = C$, $X = C'$. Then for each point $p_i \in X \setminus U''$ we get the conjugate $\ga_i \in \pi_1(U'')$ 
of a circle around the  point $p_i$, and  we have $ \ker (  \pi_1(U'') \ra \pi_1(X)) = \La \{ \ga_i \} \Ra$,
i.e., the kernel is normally generated by these loops.

Each of these loops $\ga_i$ maps to an element $c_i \in G$, and we denote by $m_i $ the order of the element $ (c_i)$ in the group $G$.
Then $\ga_i ^{m_i} \in \pi_1 (U)$, and it is the conjugate of  small circle around a point of $Z \setminus U$, hence it maps to
zero in $\pi_1(Z)$. One can see that in this case 

$$  \pi_1^{orb} ( Z, G) = \pi_1(U'') /  \La \{ \ga_i ^{m_i}\} \Ra$$ 
and the kernel of the natural surjection onto $\pi_1(X)$ is normally generated by the $\ga_i$'s.

This means that, if the genus of $X = C'$ is equal to $g'$, then the orbifold fundamental group is isomorphic to the abstract group
$$\pi_{g',m_1, \dots , m_d} : =\langle \al_1, \be_1, \dots ,  \al_{g'},\be_{g'}, \ga_1, \dots , \ga_d | \Pi_1^d \ga_j  \Pi_1^{g'} [ \al_i, \be_i]  = 1, $$
$$\ga_1^{m_1} = \dots =   \ga_d^{m_d} =1 \rangle . $$

And the orbifold exact sequence is an extension, called Nielsen extension, of  type 

$$ (NE)  \   1 \ra \pi_g \ra \pi_{g',m_1, \dots , m_d}  \ra G \ra 1.$$

\end{ex}

A similar description holds in general (see \cite{isogenous}, definition 4.4 and  proposition 4.5, pages 25-26 ),
at least when $Z$ is a complex manifold.

\begin{rem}
The notion of orbifold fundamental group has turned out to be quite useful also in real algebraic geometry,
where one considers the action of complex conjugation on a real algebraic variety (see \cite{c-f} for example).
\end{rem}

\subsection{Rational $K(\pi,1)$'s: basic examples}

An important role is also played by {\bf complex Rational $K(\pi,1)$'s}, i.e., quasi projective varieties (or complex spaces) $Z$ such that
$$ Z = \sD / \pi ,$$
where $\sD$ is a contractible manifold (or complex space) and the action of $\pi$ on $\sD$ is properly discontinuous but not necessarily free.

While for a $K(\pi,1)$ we have $ H^* (\pi , \ZZ) \cong H^* (Z, \ZZ)$,   $ H_* (\pi , \ZZ) \cong H_* (Z, \ZZ)$,
for a rational $K(\pi,1)$, as a consequence of proposition \ref{finite} and
of the spectral sequence of theorem  \ref{tohoku}, we  have $ H^* (\pi , \QQ) \cong H^* (Z, \QQ)$ and therefore also  $ H_* (\pi , \QQ) \cong H_* (Z, \QQ)$.

Typical examples of such rational $K(\pi,1)$'s are:

\begin{enumerate}
\item
quotients of a bounded symmetric domain $\sD$ by a subgroup $\Ga \subset \Aut (\sD)$ which is acting properly discontinuously
(equivalently, $\Ga$ is discrete); especially noteworthy are the case where $\Ga$ is {\bf cocompact}, meaning that $X = \sD/\Ga$
is compact, and the {\bf finite volume} case where the volume of $X$ via the invariant volume form for $\sD$ is finite.

 \item 
 the moduli space of principally polarized Abelian Varieties, which is an important special case of the above examples;
 \item

 the moduli space of curves $\sM_g$.

\end{enumerate}
For the first example the property is clear (we observed already that such a domain $\sD$ is contractible).

For the second example  it suffices to observe that the moduli space of Abelian Varieties of dimension $g$, and with
a polarization of type $ (d_1, d_2, \dots , d_g)$ is the quotient of Siegel's upper half space
$$  \sH_g : = \{ \tau \in \Mat(g,g,\CC)| \tau = ^t \tau, \im (\tau) > 0 \},$$ 
which is biholomorphic to a bounded symmetric domain of type III in E. Cartan's classification ({\cite{Cartan}),
by the properly discontinuous action of the group 

$$  \Sp (D, \ZZ) : = \{ M \in Mat(2g, \ZZ)| ^t M  D M = D \} , $$
where 
 \[ D' : = \diag (d_1, d_2, \dots , d_g), 
 D: =   \begin{pmatrix}
    0  &    D' \\
  - D' &    0
  \end{pmatrix} 
\] 
and, 
 \[  
 M: =  \left( \begin{matrix}
    \al  &    \be \\
   \ga &    \de
  \end{matrix}\right), \ \
  \tau \mapsto - D' (D' \al - \tau \ga)^{-1} ( D' \be - \tau \de).
\] 

In fact, these Abelian varieties are quotients $\CC^g / \Lam$, where $\Lam$ is generated by the columns of the matrices $D, - \tau,$
and the Hermitian form associated to the real matrix $ \Im m (\tau)^{-1}$ has imaginary part which takes integral values on
$ \Lam \times \Lam$, and its associated antisymmetric matrix in the chosen basis is  equal to $D$.

The quotient $\sA_{g,D'} : =  \sH_g /  \Sp (D, \ZZ) $ is not compact, but of finite volume, and there exist several compactifications 
which are  projective algebraic varieties (see \cite{Namikawa}, \cite{FaltingsChai}, \cite{gerardAV}).

\subsection{The moduli space of curves}

The most useful (and first fully successful) approach to the moduli space of curves of genus $g $ is to view it as a quotient

$$ (**) \   \mathfrak M_g  =  \sT _g  /  \sM ap _g $$
of a connected complex manifold $\sT _g $ of dimension $ 3g-3 +a(g)$, called {\bf Teichm\"uller space}, by the properly discontinuous action of 
the {\bf Mapping class group} $\sM ap _g $ (here $ a (0) = 3, a(1) = 1, a(g)= 0,  \forall g \geq 2$ is the complex dimension
of the group of automorphisms of a curve of genus $g$).
A key result (see \cite{kerckhoff}, \cite{tromba}, \cite{hubbard}) is the

\begin{theo}
Teichm\"uller space $\sT _g $ is diffeomorphic to a ball, and the action of $\sM ap _g $ is properly discontinuous.
\end{theo}

Denoting  as usual  by $\pi_g$ the fundamental group of a compact complex curve $C$ of genus $g$,  we have in fact a more concrete
description of the mapping class group (see \cite{Clebsch, Hurwitz, art1, art, Bir69, birman, h-t, fm} for general results on braid and mapping class groups):

$$ (M) \ \sM ap _g \cong \Out^+ (\pi_g). $$ 

The above superscript $^+$ refers to the orientation preserving property. 

There is a simple algebraic way to describe the orientation preserving property: any automorphism of a group induces an
automorphism of its abelianization, and any inner automorphism acts trivially. 

Hence $  \Out (G) $ acts on $G^{ab}$, in our particular case  $\Out (\pi_g)$ acts on $\pi_g^{ab} \cong \ZZ^{2g}$,
and  $\Out^+ (\pi_g) \subset  \Out (\pi_g)$ is the inverse image of $ SL (2g, \ZZ)$.

The above isomorphism (M)  is of course related to the fact that $C$ is a $ K (\pi_g,1)$, as soon as $ g \geq 1$.

As we already discussed,  there is a bijection between homotopy classes of self maps of $C$ and endomorphisms of $\pi_g$,
taken up to inner conjugation. Clearly a homeomorphism $\fie : C \ra C$ yields then an associated element $\pi_1 (\fie) \in \Out (\pi_g)$.

Such a homeomorphism acts then on the second homology group $H_2 (C, \ZZ) \cong \ZZ [C]$, where the generator $[C]$ corresponds
to the orientation associated to the complex structure; the condition $H_2 (\fie) = +1$ that $\fie$ is orientation-preserving
translates into the above  algebraic condition on $\psi : = \pi_1 (\fie) $.  That the induced action $\psi^{ab}$ on the Abelianization $\pi_1(C)^{ab} \cong \pi_g^{ab}
\cong \ZZ^{2g}$ satisfies:  $\wedge^{2g} (\psi^{ab})$ acts as the identity  on $ \wedge^{2g} (  \ZZ^{2g} ) \cong \ZZ$. 

In other words, the image of the product of commutators $\e : = \Pi_j [\al_j, \be_j] $ is sent to a conjugate of $\e$.

Turning now to the  definition of the Teichm\"uller space $\sT _g $, we observe that it is somehow conceptually easier 
to give the  definition of Teichm\"uller space for more general manifolds.

\subsection{Teichm\"uller space}
Let $M$ be an  oriented real differentiable manifold of real dimension $2n$, for simplicity let's assume 
  that $M$ is compact.
 
 Ehresmann (\cite{ACS}) defined an {\bf almost complex structure} on $M$ as the structure
 of a complex vector bundle on the real tangent bundle $TM_{\RR}$: namely, the action of 
 $\sqrt {-1}$ on $TM_{\RR}$ is provided by an endomorphism
 $$ J : TM_{\RR} \ra TM_{\RR}, {\rm{\  with  }}\  J^2 = - Id.$$

Equivalently,  as done for the complex tori, one  gives the decomposition of the complexified tangent bundle $TM_{\CC} : =  TM_{\RR} \otimes_{\RR}\CC$
as the direct sum of the $i$, respectively $-i$ eigenbundles:
$$  TM_{\CC}  =  TM^{1,0} \oplus  TM^{0,1} {\rm{ \  where  }} \  TM^{0,1} = \overline {TM^{1,0} }.$$

The space $\sA \sC (M)$ of almost complex structures, once  $TM_{\RR}$ (hence all associated bundles) is endowed with a
Riemannian metric, is  a subset of the  Fr\'echet space $H^0 (M, \sC^{\infty} (End (TM_{\RR}) ))$  .

A closed subspace of $\sA \sC (M)$ consists of the set $ \sC (M)$ of complex structures: these are the almost complex structures for
which there are at each point $x$ local holomorphic coordinates, i.e., $\CC$-valued functions $z_1, \dots , z_n$ whose differentials
span the dual $(TM^{1,0}_y)^{\vee}$ of $TM^{1,0}_y$ for each point $y$ in a neighbourhood of $x$.

In general, the  splitting 
$$TM_{\CC}^{\vee} = (TM^{1,0})^{\vee} \oplus (TM^{0,1})^{\vee}$$ yields a decomposition of exterior differentiation
 of functions as
$ df = \partial f + \bar{\partial} f$, and a function is said to be holomorphic if its
differential is complex linear, i.e., $ \bar{\partial} f = 0$. 

This decomposition $ d= \partial  + \bar{\partial} $
extends to higher degree differential forms.

The theorem of Newlander-Nirenberg (\cite{NN}), first proven by Eckmann and Fr\"olicher in the real analytic case (\cite{E-F})
characterizes the complex structures through an explicit equation:

\begin{theo}\label{NN}
{\bf (Newlander-Nirenberg)} An almost complex structure $J$ yields the structure of a complex manifold if and only if
it is integrable, which means $ \bar{\partial} ^2 = 0. $
\end{theo}

Obviously the group of oriented diffeomorphisms of $M$ acts on the space of complex structures, hence one can
define in few words some basic concepts.

\begin{defin}
Let $\sD iff ^+ (M)$ be the group of orientation preserving diffeomorphisms of $M$ , and let $\sC (M)$  be the space of complex structures
on $M$. Let $\sD iff ^0 (M) \subset \sD iff ^+ (M)$ be the connected component of the identity, 
the so called subgroup of diffeomorphisms which are isotopic
to the identity.

Then Dehn  (\cite{dehn}) defined the mapping class group of $M$ as 
$$\sM ap (M) : =  \sD iff ^+ (M) /  \sD iff ^0 (M),$$
while the Teichm\"uller space of $M$, respectively the moduli space of complex structures on $M$  are defined as
$$ \sT (M) : = \sC (M) / \sD iff ^0 (M) , \ \mathfrak M (M)  : = \sC (M) / \sD iff ^+ (M).$$
\end{defin}

From these definitions follows that 
$$  \   \mathfrak M (M) =  \sT (M)  /  \sM ap (M) .$$

The simplest examples here are two: complex tori and compact complex curves.
In the case of tori a connected component of Teichm\"uller space (see  \cite{cat02} and  also
\cite{cat04})
is  an open set $\mathcal T_n$ of the complex
Grassmann Manifold $Gr(n,2n)$, image of the open set of matrices

$\{ \Omega \in \Mat(2n,n; \CC) \ | \ i^n det  (\Omega \overline
{\Omega}) > 0 \}.$

This parametrization is very explicit: if we consider
a fixed lattice $ \Lam \cong \ZZ^{2n}$,
to each matrix $ \Omega $ as above we associate the subspace of $\CC^{2n} \cong \Lam \otimes \CC$ given as
$$ V =  \Omega \CC^{n},$$ so that
$ V \in Gr(n,2n)$ and $\Lam \otimes \CC \cong V \oplus \bar{V}.$

Finally, to $ \Omega $ we associate the torus 
$$Y_V : = V / p_V (\Lam) = (\Lam \otimes \CC)/ (\Lam \oplus  \bar{V} ),$$
$p_V : V \oplus
\bar{V} \ra V$ being the projection onto the first addendum.

It was observed however by Kodaira and Spencer already in their first article on deformation theory
(\cite {k-s58}, and volume II of Kodaira's
collected works)
that the mapping class group $ SL ( 2n, \ZZ)  $ does not act properly
discontinuously
on $\mathcal T_n$.

\bigskip

The case of compact complex curves $C$ is instead the one  which was originally 
considered by Teichm\"uller.

In this case,  if the genus $g$ is at least $2$, the
Teichm\"uller space $\sT_g$  is a bounded domain,   diffeomorphic to a ball,
contained in the vector space of quadratic differentials 
$H^0 ( C, \hol_C ( 2 K_C)) \cong \CC^{3g-3}$ on a fixed such curve $C$.

In fact, for each other complex structure on the oriented 2-manifold $M$ underlying 
$C$ we obtain a complex curve $C'$, and there is a unique extremal quasi-conformal
map $ f : C \ra C'$, i.e., a map such that the Beltrami distortion  $\mu_f : = \bar{\partial} f / \partial f$ 
has minimal norm (see for instance  \cite{hubbard} or \cite{ar-cor}). 

The fact that  the Teichm\"uller space $\sT_g$  is diffeomorphic to a ball (see \cite{tromba} for a simple proof)
is responsible for the fact that the moduli space of curves $\mathfrak M_g$ is a rational $K(\pi, 1)$. 

\subsection{Singularities of $\frak M_g$, I}

Teichm\"uller's theorem says that, for $g \geq 2$,  $\mathcal{T}_g \subset \CC^{3g-3}$ is an open subset diffeomorphic to a ball. 
Moreover ${\sM ap}_g \cong {\sM ap}(C)$ acts properly discontinuously on $\mathcal{T}_g$, but not freely.

The lack of freeness of this action is responsible of the fact that $\frak M_g$ is  a singular complex space (it is quasi-projective
by the results of Mumford and Gieseker, see  \cite{GIT}, \cite{EnsMath}).

Therefore one has to understand when a class of complex structure $C$  on the manifold $M$ is fixed by an element $\ga \in \sM ap_g$.
Since $\ga$ is represented by the class of an orientation preserving diffeomorphism $ \fie : M \ra M$, it is not difficult to see that this situation means
that the diffeomorphism $\fie$ carries the given complex structure to itself, or, in other words, the differential of $\fie$
preserves the splitting of the complexified tangent bundle of $M$, $  TM_{\CC}  =  TM^{1,0} \oplus  TM^{0,1} $.

This is then equivalent to $\bar{\partial} \fie = 0$, i.e., $\fie$ is a biholomorphic map, i.e., $ \fie \in \Aut (C)$.

Conversely,  let $\ga \in \Aut(C)$ be a nontrivial automorphism.  Consider then the graph of $\ga$, $ \Ga_{\ga} \subset C \times C$,
and intersect it with the diagonal $\De \subset C \times C$: their intersection number must be a nonnegative integer, because intersection
points of complex subvarieties carry  always a positive local intersection multiplicity.

This argument was used by Lefschetz to prove the following  

\begin{lemma}{\bf (Lefschetz' lemma)}\label{LefLemma} 
If we have an automorphism $\ga \in \Aut(C)$ of a projective curve of genus $ g \geq 2$, then 
$\ga$ is homotopic to the identity iff it is the identity:
$$
\ga \underset{h}{\sim} Id_C \Leftrightarrow \ga =Id_C \, .
$$
\end{lemma}
\Proof
The self intersection of the diagonal $\De^2$ equals the degree of the tangent bundle to $C$, which is $ 2 - 2g$.
If $\ga$ is homotopic to the identity, then $\Ga_{\ga} \cdot \De = 2- 2g < 0$,
and this is only possible if $ \Ga_{\ga}= \De$, i.e., iff $\ga=id_C.$

\qed

In particular, since $C$ is a classifying space $ B \pi_g$ for the group $\pi_g$, 
$
\ga \colon C \to C \, 
$
is such that the homotopy class of $\ga$ is determined by the conjugacy class of
$
\pi_1(\ga) \colon \pi_1(C,y_0) \to \pi_1(C, \ga \cdot y_0) \, .
$
Hence we get an injective  group homomorphism
$$
\rho_C  \colon \Aut(C)  \to \frac{{\rm Aut}(\pi_g)}{{\rm Inn}(\pi_g)}={\rm Out}(\pi_g) \, .
$$
Actually, since a holomorphic map is orientation preserving, we have that  $\rho_C  \colon \Aut(C) \to {\rm Out}^+(\pi_g)={\rm Map}_g \cong \frac{{\rm \sD iff}^+(C)}{{\rm \sD iff}^0(C)}$.

The conclusion is then that curves with a nontrivial group of automorphisms $ \Aut(C) \neq \{ Id_C\}$ correspond to points of $\sT _g$ with a non trivial stabilizer for the
action of  $\sM ap_g$.

\begin{rem}
There is of course the possibility that the action of $\sM ap_g$ is not effective, i.e., there is a normal subgroup acting trivially on $\sT _g$.
This happens exactly in genus $g=2$, where every curve is hyperelliptic and the normal subgroup of order $2$ generated by the hyperelliptic involution
is exactly the kernel of the action.
\end{rem}

The famous theorem of Hurwitz gives a precise upper bound $ | \Aut(C) | \leq 84 (g-1)$ for  the cardinality of the stabilizer $\Aut(C)$ of $C$ (that this is finite  follows 
also by the proper discontinuity 
of the action). 

Now, coming to $\Sing ( \frak M_g)$, we get that $\frak M_g$ is locally analytically isomorphic to a  quotient `singularity'
$\CC^{3g-3} / \Aut (C)$ (indeed, as we shall see, to the quotient singularity of $ H^1 (\Theta_C) / Aut (C)$ at the origin).

Now,  given a quotient singularity $\CC^{n} / G$, where $G$ is a finite group, then by  a well known lemma by H. Cartan (\cite{cartan}) we may assume
that $G$ acts linearly, and then a famous theorem by Chevalley (\cite{chevalley} and Shephard-Todd (\cite{s-t}) says that the quotient $\CC^{n} / G$ is smooth if and only if the action of $G \subset GL (n, \CC)$
is generated by {\bf pseudo-reflections}, i.e., matrices which are diagonalizable with $(n-1)$  eigenvalues equal to $1$, and the last  one  a root of unity.

\begin{theo}
Let $G \subset GL (n, \CC)$,  let $G^{pr}$ be the normal subgroup of $G$ generated by pseudo reflections, and let $G^q$
be the factor group. Then

(1)  the quotient $\CC^{n} / G^{pr}$ is smooth,

(2) We have a factorization  of the quotient map $\CC^{n} \ra \CC^{n} / G$ as $\CC^{n} \ra  (\CC^{n} / G^{pr} ) \cong \CC^{n} \ra \CC^{n} / G^q$  
and  $\CC^{n} / G$ is singular if $ G \neq G^{pr}$.

\end{theo}

This result is particularly interesting in this situation, since the only pseudo-reflection occurs for the case of the hyperelliptic involution in genus $g=3$.

Teichm\"uller theory can be further applied in order to analyse the fixed loci of finite subgroups $G$ of the mapping class group (see \cite{kerckhoff}, \cite{tromba}, 
\cite{isogenous})

\begin{theo} {\bf (Refined Nielsen realization)}\label{refinedNR}
Let $ G \subset \sM ap_g$ be a finite subgroup. Then $ Fix(G) \subset \sT _g$ is a non empty complex manifold, diffeomorphic to a ball.
\end{theo}

\begin{rem}\label{G-extension}

1) The name Nielsen realization comes from the fact that Nielsen (\cite{nielsen}, \cite{nielsen2},  \cite{nielsen3}) conjectured that any topological action of a finite group $G$
on a Riemann surface could be realized as  a group of biholomorphic automorphisms if and only if the action would yield
an embedding $\rho :  G \ra \sM ap_g$.

2) The question can be reduced to  a question of algebra, by the following observation.

The group $\pi_g$, for $ g \geq 2$, has trivial centre: in fact, two commuting hyperbolic transformations in $\PP SL(2, \RR)$ can be 
simultaneously written as M\"obius transformations  of the form $ z \mapsto \la z, z \mapsto \mu z$, hence if the first transformation were in the centre
of $\pi_g \subset \PP SL(2, \RR)$, the group $\pi_g$ would be abelian, a contradiction.

Then $\pi_g \cong Inn (\pi_g)$ and from  the exact sequence 
$$ 1 \ra   \pi_g \cong Inn (\pi_g) \ra \Aut(\pi_g ) \ra \Out (\pi_g) \ra 1 $$
and an injective homomorphism $\rho :  G \ra \sM ap_g =  \Out^+ (\pi_g)$ one can pull-back an extension \footnote{the same result holds
without the assumption that $\rho$ be injective, since it suffices to consider the image  $G' = Im ( \rho) $, apply the pull-back construction to
get  $1 \ra   \pi_g \ra \hat{G}' \ra G' \ra 1$, and then construct $\hat{G}$ as the inverse image of the diagonal under the epimorphism
$  \hat{G}'  \times G \ra G' \times G' $.}

$$ 1 \ra   \pi_g \ra  \hat{G} \ra G \ra 1.$$

4)  The main question is to show that the group $\hat{G}$ is isomorphic to an orbifold fundamental  group of the form $ \pi_{g', m_1, \dots , m_d} $.
Because then we have a Nielsen extension, and the epimorphism $\mu :  \pi_{g', m_1, \dots , m_d}  \ra G \ra 1$ combined with 
 Riemann's existence theorem finds for us a curve $C$ over which $G$
acts with the  given topological type: such a curve $C$ is given by the ramified covering of a curve $C'$ of genus $g'$, branched on $d$ points $y_1, \dots, y_d$
and with monodromy $\mu$.

5) Since $\pi_g$ is torsion free (hyperbolic elements in $\PP SL(2, \RR)$ have infinite order), we see that the branch points correspond to conjugacy classes
of cyclic subgroups  of finite order in $\hat{G}$, and that their order $m_i$  in $\hat{G}$ equals the order of their image in $G$.
\end{rem}

In turn, the question can be reduced to showing that there is a differentiable action of the group $G$ on the curve $C$ inducing
the topological action $\rho$. There is the following 

\begin{lemma}
Given $ \rho : G \ra Map_g = \sD iff^+(C) /  \sD iff^0(C)$, if there is a homomorphism $ \psi: G \ra \sD iff^+(C)$  whose projection to
$Map_g$ is $\rho$, then there is a complex curve with a $G$ action of topological type $\rho$.
\end{lemma}

\Proof
In fact, by Cartan's lemma \cite{cartan}, at each fixed point $x \in C$, there are local coordinates
such that in these coordinates the action of the stabilizer $G_x$ of $x$ is linear (in particular, if we assume that the action of $G$
is orientation preserving, then $G_x \subset SO(2, \RR)$, and $G_x$ is a cyclic group of rotations).

Therefore, for each $\ga \in G$, $\ga \neq 1_G$, the set of fixed points $ Fix (\ga)$, which  is closed, is either discrete (hence finite),
or it consists of a discrete set plus a set (the  points $ x \in Fix (\ga)$ where the derivative is the identity) which is open and closed.
Since $C$ is connected, the only possible alternative is that, for $\ga \neq 1_G$, the set of fixed points $ Fix (\ga)$ is finite.

Now, it is easy to see that the quotient $C/G = : C'$ is a smooth Riemann surface and we can take a complex structure on $C'$
and lift it to $C$, so that $C' = C/G$ as complex curves.

\qed

\begin{rem}
In general Markovi\'c has proven that the mapping class group cannot be realized, as hoped  by Thurston,
as a subgroup of the group of diffeomorphisms of a Riemann surface of genus $g$. There are other 
 recent results concerning non finite groups (see \cite{bcs13}).
\end{rem}

\begin{rem}
1) It is interesting to observe that one can view the orbifold fundamental group, in the case of curves, as the quotient of a bona fide
fundamental group.

Let $TC$ be the tangent bundle of the complex curve, and assume that $G$ has a differentiable action on $C$, so that $G$ acts on $TC$,
preserving the zero section. We can further assume, by averaging a Riemannian metric on $TC$, that $G$ preserves the sphere bundle
of $TC$, $$ SC : = \{ v \in TC | |v| = 1\}.$$ 

Since each stabilizer  $G_x$ is a cyclic group of rotations, the conclusion is that $G$ acts on $SC$ by a free action, and we get an exact sequence
$$ (Seifert) \   1 \ra \pi_1(SC) \ra \pi_1(SC/G ) \ra G \ra 1.$$

It is well known (see e.g. \cite{localpi1}) that the fundamental group of the fibre bundle $SC$ is a central extension
$$  1 \ra \ZZ \ra \pi_1(SC) \ra \pi_g \ra 1,$$
and is the quotient of the direct product $\ZZ \times \FF_{2g}$ by the subgroup normally generated by the relation
$$ c^{2g-2}   \Pi_1^{g} [ \al_i, \be_i]  = 1, $$
where $c$ is a generator of the central subgroup $\langle c \rangle \cong \ZZ$. 

Taking the quotient by the subgroup $\langle c \rangle $ generated
by $c$ one obtains the orbifold group exact sequence. The fibration $SC/G \ra C' = C/G$ is a so-called Seifert fibration,
with fibres all homemorphic to $S^1$, but the fibres over the branch points $y_j$ are multiple of multiplicity $m_j$. 

Now, an injection $ G \ra \Out^+(\pi_g)$ determines an action of $G$ on $\pi_1(SC)$ (sending $c$ to itself), an extension of the type
we denoted (Seifert), and a topological action on the classifying space $SC$.

2) The main point of the proof of the Nielsen realization  is however based on analysis, and, 
more precisely, on  the construction of  a Morse function $f$  on $\sT _g$ of class $C^1$ which is $G$-invariant,  strictly convex and proper,
and   therefore has only one minimum.

The unique minimum must be $G$-invariant, hence we get a point $C$ in $ Fix(G) $. Since the function is strictly convex and proper,
it does not have other critical points, and it has only one minimum: therefore Morse theory tells that $\sT _g$ is diffeomorphic
to a ball.

Now, in turn,  the locus $ Fix(G) $ is a non empty submanifold, and the group $G$ acts nontrivially on the normal bundle 
of $ Fix(G) $  in $\sT_g$. 

The same analysis applies to  the restriction of the Morse function to any connected component $Y$ of $ Fix(G) $.
The function has exactly one critical point $D$, which is an absolute minimum on $Y$, hence $Y$ is also  diffeomorphic to a ball.

To show that $ Fix(G) $ is connected it suffices to show  that $C=D$, i.e., that $D$ is also a critical point on the whole space $\sT _g$.
But the derivative of $f$ at $D$ is a $G$-invariant linear functional, vanishing on the subspace of $G$-invariants: therefore 
the derivative of $f$ at $D$ is zero.

Hence $ Fix(G) =Y $, as well as $\sT _g$,  is diffeomorphic to a ball. 

3) It may happen that $ Fix(G) $ may be also the fixed locus for a bigger subgroup $ H \subset \sM ap_g$. We shall see examples later in examples \ref{notfull1}, \ref{notfull2}.
To determine when this happens is an interesting question, fully answered by Cornalba in the case where $G$ is a cyclic group of prime order
(\cite{cornalba}) (this result  is the key to understanding the structure of the singular locus $\Sing ( \frak M_g)$).

A partial answer to the question whether $G$ is not a {\bf full} subgroup (i.e., $\exists H \supset G$,  $ Fix(H) = Fix(G) $ was given by several
authors, who clarified the orbifold fundamental groups associated to the quotients $ C \ra C/G \ra C/H$ (\cite{singer}, \cite{ries}, \cite{mssv}).
\end{rem}

\subsection{Group cohomology and equivariant cohomology}

The roots of group cohomology go back to Jacobi and to the study of periodic meromorphic functions as quotients of quasi-periodic holomorphic  functions;
these functions, which can more generally be taken as vector-valued functions $ s : \CC^n \ra \CC^m$, are solutions to the functional equation
$$ s(x + \ga ) =  f_{\ga}(x) s(x),$$ 
where $ \Lam \subset \CC^n$ is a discrete subgroup, $\forall \ga \in \Lam \subset \CC^n, x \in \CC^n, f_{\ga} : \CC^n \ra \GL (m, \CC),$
and clearly $f_{\ga}$ satisfies then the cocycle condition

$$   f_{\ga + \ga' }(x)  =  f_{\ga}(x + \ga' )  f_{\ga'}(x) , \ \forall  \ga, \ga' \in \Lam,   \forall x \in \CC^n .$$

Later on, the ideas of Jacobi were generalised through the concept of vector bundles, since the cocycle $f_{\ga} (x)$ can be used to
constructing a vector bundle on the quotient manifold $\CC^n / \Lam$ taking a quotient of the trivial bundle $  \CC^n \times \CC^m \ra \CC^n$ by the  relation
$$ (x,v) \sim ( x + \ga,  f_{\ga}(x) v) ,$$ 
and the cocycle relation just says that the above is an equivalence relation. The main difficulty classically consisted
in simplifying as much as possible the form of these cocycles, replacing $ f_{\ga}(x) $ by an equivalent one
$ \phi(x + \ga)  f_{\ga}(x) \phi (x)^{-1}$, for $\phi : \CC^n \ra \GL (m, \CC).$

The construction can be vastly generalized by viewing any bundle on a classifying space $ Y = B \Ga  $,
quotient of a contractible space $ X = E \Ga$, as the quotient of the trivial bundle $ X \times F$
by an action which covers the given action of $G$ on $X$, hence such that there exists 
$  f_{\ga }: X \ra Homeo (F)$ with 
$$ \ga (x, v) = ( \ga (x) , f_{\ga}(x)(v) ) \Rightarrow \forall \ga, \ga' ,  \forall x \in X, \ f_{\ga'  \ga }(x) =  f_{\ga'}(\ga (x)) f_{\ga}(x) .$$

As remarked by David Mumford \cite{mum-notices}, classical mathematics led to concrete understanding of 
1-cocycles and 2-cocycles, hence of first and second  cohomology groups. But the general machinery
of higher cohomology groups  is harder to understand concretely.

As always in mathematics, good general definitions help to understand what one is doing, but  explicit 
calculations remain often a hard task. Let us concentrate first on  the special case where we look at the singular cohomology groups of a classifying space.

In the case of a classifying space $ B \Ga$, we know that for each ring of coefficients $R$,
$$ H^* (\Ga, R)  =  H^*( B \Ga, R) ,$$ and  for this it suffices indeed that 
$$ H^* (\Ga, \ZZ)  =  H^*( B \Ga, \ZZ) .$$

Which side of the equation is the easier one to get hold of?  The answer depends of course on our knowledge
of the group $\Ga$ and  of the classifying space $ B \Ga$, which is  determined only up to homotopy equivalence.

For instance, we saw that  an incarnation of $ B (\ZZ / 2) $ is given by the inductive limit $\PP^{\infty}_{\RR} : =  \lim_{n \ra \infty} \PP^n_{\RR}$.

A topologist would just observe that $\PP^{\infty}_{\RR}$ is a CW-complex obtained adding a cell $\sigma_n$ in each dimension $n$,
and with boundary map $$ \partial ( \sigma_{2n}) = 2 \sigma_{2n-1}, \  \partial ( \sigma_{2n-1}) = 0 \cdot \sigma_{2n-2}, \forall n \geq 1 .$$
Therefore $ H_{2i}( \ZZ / 2, \ZZ) = 0, H_{2i-1}( \ZZ / 2, \ZZ) =  \ZZ / 2, \forall i \geq 1,$
and taking the dual complex of cellular cochains one gets:
$$ H^{2i}( \ZZ / 2, \ZZ) = \ZZ / 2, H^{2i-1}( \ZZ / 2, \ZZ) =  0,  \forall i \geq 1.$$

There are two important ways how this special calculation generalizes: 

i) a completely general construction of a CW complex  which is a classifying space $ B \Ga$ for a finitely presented group $\Ga$,

ii) a completely algebraic definition of group cohomology.

i):  In fact, assume we  are given a finite presentation of  $\Ga$ as
  $$ \Ga = \langle  x_1, \dots, x_n | R_1(x), \dots R_m (x) \rangle, $$
 which means that $\Ga$ is isomorphic to  the quotient of the free group $\F_n$ with generators $x_1, \dots, x_n$ by the minimal normal subgroup $R$ 
 containing the words $R_j \in \F_n$ ( $R$ is called the subgroup of relations, and it is said to be normally generated by the relations    $R_j$).
 
 Then the standard construction of the 2-skeleton of $B\Ga$ is the CW-complex $B\Ga^2 $, of dimension two, which is obtained by
 attaching, to a bouquet of $n$ circles (which correspond to the generators  $x_1, \dots, x_n$), $m$ 2-cells whose respective boundary is the closed path
 corresponding to the word $R_j, \forall j=1, \dots, m$.
 
 For instance , $$ \ZZ / 2 = \ZZ / ( 2 \ZZ )= \La  x_1 | x_1^2 \Ra, $$ and the above procedure produces $\PP^2_{\RR}$ as 2-skeleton. 
 
 However, the universal cover of $\PP^2_{\RR}$ is the sphere $S^2$, which has a homotopy group $\pi_2 (S^2) \cong \ZZ$, and we know that
 the higher ($i \geq 2$) homotopy groups $\pi_i (X)$ of a space equal the ones of its universal cover. Then a 3-cell is attached to  $\PP^2_{\RR}$,
 obtaining $\PP^3_{\RR}$, and one continues to attach cells in order to kill all the homotopy groups. The same is done more generally
 to obtain $B\Ga$ from the CW-complex $B\Ga^2 $: 3-cells are attached in order to kill the second homotopy groups,
 and one then obtains $B\Ga^3 $; one obtains then $B\Ga^4 $ adding 4-cells in order to kill $\pi_3(B\Ga^3)$, and so on.
 
 The disadvantage of this construction is that from  the third step onwards  it is no longer so explicit, since calculating homotopy 
 groups is  difficult. 
 
An easy but important remark used in  the proof of a theorem of Gromov (see \cite{Gromov}, \cite{Levico})   is that when the group $\pi$ has {\bf few relations}, which more precisely means
here $ n \geq m + 2$, then necessarily there is an integer $b$ with $0 \leq b \leq m$ such that the rank of $H^1 (B \pi, \ZZ)$ equals $  n-m +b$, while  $ \rank (H^2 (B \pi, \ZZ)) \leq b:$
since by the above construction $ H_2 (B \pi^2, \ZZ)$ is the kernel of a linear map $ \partial : \ZZ^m \ra  H_1 (B \pi^1, \ZZ) = \ZZ^n$, whose rank is denoted  by $ m-b$,  
and we have a surjection $ H_2 (B \pi^2, \ZZ) \ra H_2 (B \pi, \ZZ)$ .

Hence, if we consider the cup product  bilinear map $$ \cup :  H^1 (B \pi, \ZZ) \times H^1 (B \pi, \ZZ) \ra H^2 (B \pi, \ZZ),$$
for each $\psi \in H^1 (B \pi, \ZZ)$ the linear map $ \phi \mapsto \psi \cup \phi$ has a kernel of rank $ \geq n-m \geq 2$.
This means that each element $\psi \in H^1 (B \pi, \ZZ)$ is contained in an isotropic subspace (for the cup product map) of 
rank at least $2$.

It follows  that if the fundamental group of a space $X$, $\Ga : = \pi_1(X)$ admits a surjection onto $\pi$, then the induced classifying
continuous map $\phi : X \ra   B \pi$ (see \ref{classifying map}) has the property 

\begin{itemize}
\item
that its induced homomorphism on first cohomology
$$ H^1 (\phi) : H^1 (B \pi, \ZZ) \ra H^1 (X, \ZZ)$$ is injective (since $ H_1 (\phi) :  H_1 (X, \ZZ) \ra H_1 (B \pi, \ZZ) $ is surjective)
\item
the image $W$ of $ H^1 (\phi) $ is such that each element $ w \in W$ is contained in an  isotropic subspace of rank $\geq 2$.
\end{itemize}

We shall explain Gromov's  theorem, for the case where $X$ is a K\"ahler manifold, in the next section.

But let us return to the algebraic point of view (see for instance \cite{BAII}, pages 355-362  , also \cite{Brown}).

ii):  Assume now that $G$ is a  group, so that we have the group algebra $ \sA : = \ZZ[G]$.
Recall the functorial definition of group cohomology and of group homology, which gives a high brow explanation for
a rather concrete definition we shall give later.

\begin{defin}
Consider the category of modules $M$ over the group ring $ \sA : = \ZZ[G]$ (i.e., of abelian groups over which there is an action of $G$,
also called $G$-modules).

Then there are two functors, the functor of {\bf invariants}
$$ M \mapsto M^G : = \{  w \in M | g(w) = w \ \forall g \in G \} = \cap_{g \in G } \ker (g-Id), $$
associating to each module $M$ the submodule of elements which are left fixed by the action of $G$, 
and the functor of {\bf co-invariants}
$$ M \mapsto M_G : =M / ( \Sigma_{g \in G } Im (g-Id)), $$
associating to each module $M$ the minimal quotient module on which  the action of $G$ becomes trivial.

Both $M^G$ and $M_G$ are trivial $G$-modules, and one defines the cohomology groups $H^i(G,M)$ as the derived functors
of the functor of invariants, while   the homology groups $H_i(G,M)$ are the derived functors
of the functor of co-invariants.

\end{defin}

The relation of these concepts to homological algebra is furnished by the following elementary lemma.

\begin{lemma}
Let $\ZZ$ be the trivial $G$-module, i.e., we consider the abelian group $ \ZZ$ with  module structure such that
every $g \in G$ acts as the identity. 

 Then we have canonical isomorphisms
$$  \Hom_{\ZZ[G]}  ( \ZZ, M) \cong  M^G , \  M \otimes_{\ZZ[G]} \ZZ \cong M_G,$$
given by the homomorphisms $ f \mapsto f(1)$, respectively $ w \otimes n \mapsto [ n w]$.

In particular, we have 
$$ H^i(G,M) = Ext^i_{\ZZ[G]}  ( \ZZ, M), \   H_i(G,M) = Tor_i^{\ZZ[G]}  ( \ZZ, M). $$
\end{lemma}

We illustrate now the meaning of the above abstract definition for the case of group cohomology.

One can  construct an explicit free resolution (a projective resolution would indeed suffice) for 
the trivial $\sA$-module $M = \ZZ$
$$  \dots L_n \ra L_{n-1} \ra \dots \ra L_1 \ra L_0       \ra  \ZZ \ra 0  .$$

Then for each $\sA$-module $M$  the cohomology groups $ H^n (G, M) $ are computed as the cohomology groups of the complex
 $$ \Hom (L_0, M)  \ra \Hom (L_1, M)  \ra \dots \ra \Hom (L_{n-1}, M)  \ra \Hom (L_n, M)       \ra   \dots  .$$
 
 While the homology groups are calculated as the homology groups of the complex
 $$ \dots \ra M \otimes_{\ZZ[G]}   L_{i+1}   \ra M \otimes_{\ZZ[G]}  L_{i}  \ra  M \otimes_{\ZZ[G]}   L_{i-1}     \ra   \dots  .$$

 For instance, if $G$ is a cyclic group of order $m$, then $\sA = \ZZ[x]/ (x^m - 1)$ and the free resolution of $\ZZ$ is given by free modules of rank one:
 the homomorphisms  are the augmentation $ \e : L_0 \ra \ZZ$ (defined by setting $\e (\Sigma_g a_g \cdot g ) : =  \Sigma_g a_g$) and  the scalar multiplications:
 $$ ( x-1) : L_{2n+1} \ra  L_{2n} , \  ( 1 + x \dots + x^{m-1}) : L_{2n+2} \ra  L_{2n+ 1}, $$
 in view of the fact that $ (x-1)( 1 + x \dots + x^{m-1}) = x^m -1 \equiv 0 \in \sA $. 
 
 The  case $m=2$ we saw before is a very special case, since then $\sA = \ZZ \oplus \ZZ x $, with $ x^2 = 1$, hence on the trivial module
 $\ZZ = \Hom_{\ZZ[G]}  ( \ZZ[G], \ZZ )$ $(1-x)$ acts as multiplication by $0$, $1 + x$ acts as multiplication by $2$.
 But the same pattern happens for all $m$: $(1-x)$ acts as multiplication by $0$, $1 + x \dots + x^{m-1}$ acts as multiplication by $m$.
 
 In particular, the homology groups can be calculated 
 as $$H_{2i}(\ZZ/m, \ZZ) = 0, i \geq 1, \ H_{2i+1}(\ZZ/m, \ZZ) \cong \ZZ/m, \  i \geq 0.$$

 For more general groups there is a general complex, called the bar-complex, which  yields a resolution of the $\sA$-module $\ZZ$.
 
 It is easier to first give the concrete definition of the cohomology groups (see also \cite{godement}).
  
 \begin{defin}
 Given a group $G$ and a $\sA : = \ZZ [G]$-module $M$, one defines the group of $i$-cochains with values in $M$ as:
 $$ C^i (G, M) : = \{ f : G^{i+1}  \ra M| f (\ga g_0, \dots, \ga g_i) = \ga   f ( g_0, \dots,  g_i) \ \forall \ga \in G \},$$
 and the differential $ d_i :  C^i (G, M) \ra  C^{i+1} (G, M) $ through the familiar formula
 $$ df ( g_0, \dots,  g_{i+1} )  : = \Sigma_0^{i+1} (-1)^j f ( g_0, \dots, \hat{g_j}  ,\dots,   g_{i+1} ) .$$
 
 Then the groups $H^i (G, M)$ are defined as the cohomology groups of the complex of cochains, that is, $$H^i (G, M) = \ker (d_i) / Im (d_{i-1}).$$ 
 \end{defin}
 
 \begin{rem}
 1) In group theory one gives a different, but equivalent formula, obtained by considering a cochain as a function of $i$ instead of $i+1$ variables,
 as follows:
 $$ \varphi ( g_1, \dots, g_i) : = f (1, g_1, g_1 g_2, \dots, g_1 g_2 \dots g_i) .$$
 One observes in fact that, in view of equivariance of $f$ with respect to left translation on $G$, giving $f$ is equivalent to giving $\varphi$.
 
 Then the formula for the differential becomes asymmetrical
 $$ d \varphi ( g_1, \dots, g_{i+1}) = g_1  \varphi ( g_2, \dots, g_{i+1}) - \varphi ( g_1 g_2, g_3 , \dots, g_{i+1}) + $$
 $$ + \varphi ( g_1, g_2 g_3 , \dots, g_{i+1}) \dots (-1)^i  \varphi ( g_1, g_2, g_3 , \dots, g_{i}).$$
 
 2) The second formula is more reminiscent of the classical formulae, and indeed, for $i=0$, yields
 $$ H^0(G, M) = \{ x \in M | g x - x = 0 \ \forall g \in G \} = M^G, $$
  $$ Z^1(G, M) = \{ \varphi = \varphi (g)  | g _1  \varphi (g_2)  -  \varphi (g_1 g_2 ) +  \varphi (g_1)   = 0 \ \forall g_1, g_2  \in G \} , $$
  $$  B^1(G, M) = \{\varphi =  \varphi (g)  | \exists x \in M ,   \varphi (g)   = g x - x \}, \  H^1(G, M) = Z^1(G, M) / B^1(G, M).$$ 
  
 3) Hence, if $M$ is a trivial $G$-module, then $  B^1(G, M) =0$ and $  H^1(G, M) = \Hom (G, M)$.
 \end{rem}
 
 We have been mainly interested in the case of $\ZZ$-coefficients, however if we look at coefficients in a field $k$,
 we have:
 
 \begin{prop}\label{finite}
 Let $G$ be a finite group, and let $k$ be a field such that the characteristic of $k$ does not divide $|G|$. Then 
 $$ H^i(G, k) = 0 \ \forall i \geq 1. $$ 
 
 \end{prop}
 \Proof
 If $M$ is a $k[G]$-module,  we observe first  that $$M^G  \cong Hom_{\ZZ[G]} (\ZZ, M) \cong Hom_{k[G]} (k , M),$$
 and then that  the functor $M \mapsto M^G$ is exact, because $M^G$ is a direct summand 
 of $M$ (averaging $ v \mapsto \frac{1}{|G|} \sum_{g \in G} gv $ yields a projector with image $M^G$).
 Hence $ H^i(G, M) = 0 \ \forall i \geq 1. $

 \qed
 
 \subsection{Group homology, Hopf's theorem, Schur multipliers}
To give some explicit formula in order to calculate homology groups, we need to describe the bar-complex, which gives
a resolution of the trivial module $\ZZ$.

In order to relate it to the previous formulae for group cohomology, we should preliminary observe that
the group algebra $\ZZ[G]$ can be thought of as the subalgebra of the space of functions $ a: G \ra \ZZ$ generated
by the characteristic functions of elements of $G$ (so that $g \in  \ZZ[G]$ is the function such that $g(x) = 1 $ for $x=g$, else
$g(x) = 0, x\neq g $. As well known, multiplication on the group algebra corresponds to convolution of the corresponding functions,
$ f_1 * f_2 (x) = \int _G f_1 (x y^{-1} ) f_2 (y) dy $.

It is then clear that the tensor product $\ZZ[G] \otimes_{\ZZ} M$ yields a space of $M$-valued functions on $G$
(all of them, if the group $G$ is finite), and we get a space of $\ZZ$-valued functions on the $i$-th Cartesian product $G^{i+1}$
by considering the $(i+1)$-fold tensor product
$$ C_i : = \ZZ[G] \otimes_{\ZZ} \ZZ[G] \otimes \dots \otimes_{\ZZ} \ZZ[G]  ,$$
with $\ZZ[G] $-module structure given by
$$  g (x_0 \otimes  \dots \otimes x_i ): =  g (x_0) \otimes  \dots \otimes x_i .$$

\begin{defin}
The bar-complex of a group $G$ is the homology complex given by the  free  $\ZZ[G] $-modules
$C_i$ (a basis is given by
$$ \{ (g_1, \dots, g_i) : =  1 \otimes g_1 \otimes  \dots \otimes g_i \})$$
and with differential $ d_i : C_i \ra C_{i-1}$ 
(obtained by dualizing  the previously considered  differential for functions), defined by
  
  $$ d_i   (g_1, g_2, \dots, g_i) : = $$
  $$ g_1  (g_2, \dots, g_i) + \Sigma_1^{i-1} (-1)^j   (g_1, \dots g_{j-1},  g_j g_{j+1} , \dots, g_i) + (-1)^i (g_1, g_2, \dots, g_{i-1}).$$
  
  The augmentation map $ \e : C_0 = \ZZ[G] \ra \ZZ$ shows then the the bar-complex is a free resolution of the
  trivial $\ZZ[G] $-module $\ZZ$.
  
  \end{defin}
  
  As a consequence, the homology groups $H_n(G, \ZZ)$ are then the homology groups of the complex
  $( C_n \otimes_{\ZZ[G]} \ZZ, d_n)$, and one notices that, since $ \ZZ[G] \otimes_{\ZZ[G]} \ZZ \cong \ZZ$,
  (observe in fact that $ g\otimes 1 = 1 \otimes g 1 = 1 \otimes 1 $),
  then  $ C_n \otimes_{\ZZ[G]} \ZZ $ is a trivial $\ZZ[G] $-module and a free $\ZZ$-module with basis $\{ (g_1, \dots, g_n) \}$.
  
  The algebraic definition yields the expected result for $i=1$.
  
  \begin{cor}
  $H_1(G, \ZZ) = G^{ab} = G / [G,G] $, $[G,G] $ being  as usual the subgroup generated by commutators.
  \end{cor}
    \Proof
  First argument, from algebra: we have
  $$ \oplus_{g_1, g_2  \in G}  \ZZ (g_1, g_2)  \ra \oplus_{g \in G} \ZZ g \ra \ZZ  $$ 
  but the second homomorphism is zero, hence
  $$ H_1(G, \ZZ) = ( \oplus_{g \in G} \ZZ g ) / \langle g_2  - g_1 g_2 + g_1 \rangle = G^{ab}, $$
  since we have then the equivalence relation
 $  g_1 g_2 \equiv    g_1 + g_2 $.
 
 Second argument, from topology: $H_1(X, \ZZ)$ is the abelianization of the fundamental group,
 so we just apply this to the case of $ X = BG$, since $\pi_1 (BG) \cong G$.
  
  \qed
  
  The machinery of algebraic topology is also useful in order to show the following `duality' statement.
  
  \begin{theo}\label{Schur}
  If $G$ is a finite group there is an isomorphism between the group $H_2(G, \ZZ)$ and the group of `Schur multipliers' $H^2 (G, \CC^*)$.
  \end{theo}
  
  \Proof
Let us consider the exact sequence of groups
$$ 0 \ra \ZZ \ra \CC \ra \CC^* \ra 1 ,$$  
and use that $H^i(G, \CC ) = 0$ for $ i \geq 1$ (see \cite{BAII}).

Hence we have an isomorphism   $H^2 (G, \CC^*) \cong  H^3 (G, \ZZ)$.

Now, by the universal coefficient formula, the torsion subgroup of $H_n(X, \ZZ)$ equals 
 the torsion subgroup of $H^{n+1}(X, \ZZ)$. Apply this observation  to the case $ X = BG$, $n=2$, remarking 
that the groups $H^i(G, \ZZ)$, hence also the groups $H_i(G, \ZZ)$, are torsion abelian groups for $i \geq 1$.
 Indeed,   $H^i(G, \RR) =  H^i(BG, \RR) =  H^i(EG, \RR)^G = 0$).  
  
  \qed
  
  The calculation of the above  group $H_2(G, \ZZ)$ was achieved by Heinz Hopf (\cite{Hopf}), and we shall sketch the underlying topological idea.
  
  \begin{theo} {(\bf Hopf)}
  Assume that we have a presentation of a group $G$ as the quotient $ G = F/R$ of a free group.
  Then $$H_2(G, \ZZ) = ( [F,F] \cap R) / [F,R].$$  
  \end{theo}
  
  \Proof {(\em Sketch)}
  
  {\bf Proof I (currently fashionable argument)}
  
  To the exact sequence of groups
  $$   1 \ra R \ra F \ra G \ra 1$$
  there correspond continuous maps of classifying spaces. These maps are only well defined up to homotopy,
  and using a small trick we can arrange the corresponding maps to be the standard inclusion of a fibre into a fibre bundle.
  
  More precisely, let $BG = EG /G$ be a classifying space for $G$, and $ BF = EF/F$ be a classifying space for $F$:
  since the group is free, $BF$ shall be a  bouquet of circles. Moreover, since $R$ is a subgroup of $F$, 
  its classifying space $BR $ is a covering of $BF$, indeed it is $ BR = EF /R$, and is a CW-complex of  dimension 1
  (whence, the theorem of Nielsen that a subgroup $R$ of a free group $F$ is also free).
  
 The surjection $ F \ra G$ yields an action of $F$ on $EG$, which restricts to a trivial action of $R$. 
 
 We use the construction of equivariant cohomology to obtain a new classifying space
 $$  B'F : = (EF \times EG) / F  .$$ 
 
 The projection of the product  $(EF \times EG) \ra EG$ yields a fibre bundle 
   $(EF \times EG)/ F   \ra EG/ G $ with fibre $ EF/R$.
   
   Hence an associated sequence of continuous maps 
   $$ BR = EF/R \ra  (EF \times EG)/ F = B'F   \ra EG/ G = BG,$$ 
   associated to the fibre bundle $ B'F \ra BG$ and with 1-dimensional fibre $BR$.
   
   We obtain associated maps of first homology groups 
  
   $$  H_1(BR, \ZZ) = R / [R,R]  \ra   H_1(BF, \ZZ) = F / [F,F] \ra   H_1(BG, \ZZ) = G / [G,G],$$
   and since  $BR, B'F$ have vanishing i-th  homology groups for $i \geq 2$,
   application of the Serre spectral sequence for fibre bundles (\cite{serre})
   shows that we have the exact sequence
    $$  0 \ra H_2(G, \ZZ) \ra H_1(BR, \ZZ)_F = (R / [R,R])_F  \ra  F / [F,F] \ra   G / [G,G] \ra 0.$$
    
    In the above exact sequence appears the group of co-invariants $(R / [R,R])_F$, which is the quotient of $R / [R,R]$ by the relations of the form
    $ frf^{-1} = r, \ \forall f \in F, \ r \in R$, or, equivalently, by the subgroup $[F,R]$ generated by the commutators $ frf^{-1}  r^{-1}$.
    In other words,
    $$ (R / [R,R])_F = R / [F,R]. $$ 
    
    Now, from the exact sequences
  $$ (*) \   0 \ra H_2(G, \ZZ) \ra R / [F,R] \ra \sK \ra 0$$
  and
  $$ 0 \ra \sK   \ra  F / [F,F] \ra   G / [G,G] \ra 0,$$ 
  since $ \ker (F \ra  G / [G,G] ) = R [F,F]$, we infer that
  $$  \sK  =    (R [F,F]) / [F,F] = R / ( R \cap [F,F]),$$
  where the last equality follows from the third isomorphism theorem for groups. 
   
   Since  $ [F,R] \subset  ( R \cap [F,F])$, Hopf's theorem follows from (*) above.
   
   \smallskip
   
     {\bf Proof II (direct argument):} follows from the following lemma, since the kernel of $R / [F,R] \ra   F / [F,F]$ 
     clearly equals $( R \cap [F,F]) / [F,R]$.
  \qed
  
 The relation between algebra and geometry in Hopf's theorem  is further explained in lemma 2.5 of \cite{CLP3}.
 
 \begin{lemma}
 Let $BG$ be a CW-complex which is a  classifying space for the group $G = F / R$,  such that 
 
 i)  its 1-skeleton $BG^1$  has all the  1-cells in  bijective  correspondence with  a set of generators of $F$,
 
 ii)   its 2-skeleton $BG^2$  has all the  2-cells in  bijective  correspondence with a set of generators of $R$ (hence $R / [R,R]$ is isomorphic to
 the relative homology group   $H_2(BG^2, BG^1, \ZZ)$).
 
 Then the following   exact sequence of relative homology 
 $$ 0 \ra H_2(BG, \ZZ) \ra H_2(BG, BG^1, \ZZ) \ra H_1(BG^1, \ZZ)  \ra H_1(BG, \ZZ) \ra 0 $$
  is isomorphic to 
  $$  0 \ra  ( [F,F] \cap R) / [F,R] \ra  R / [F,R] \ra   F / [F,F]    \ra  G / [G,G]  \ra 0.$$
 \end{lemma}

\Proof {\em (Sketch)}
 The main point  is to show  that the obvious surjection  $$H_2(BG^2, BG^1, \ZZ) \ra H_2(BG, BG^1, \ZZ)$$
 reads out  algebraically as the surjection $R / [R,R] \ra R / [F,R] $.
 
   In order to explain this, assume that $\{r_j | j \in J\}$ is a system 
 of free generators of $ R$. Now, $BG^3$ is obtained attaching 3-cells in order to kill the second homotopy group of $BG^2$ which, in
 turn, by Hurewicz' theorem, is the second homology group of its universal cover $\tilde{BG^2}$. But it is then easy to see that generators for
 $H_2(\tilde{BG^2}, \ZZ)$ are obtained applying the covering transformations of $G = F/R$ to the 2-cells corresponding to the elements 
  $\{r_j | j \in J\}$. One can then show  that the boundaries of these 3-cells yield relations on $R / [R,R]$  which are exactly those of
  the form $[a,r_j]$, for $ a \in F$.
  
  \qed
  
  \begin{ex}
  For a finite cyclic group $G = \ZZ / m$, $F= \ZZ$, hence $ [F,F] = 0$, $H_2 ( \ZZ / m, \ZZ) =0 $
  and the result fits with previous calculations.
  \end{ex}

  The spectral sequence argument is however more powerful, and yields a much more general result (see \cite{weibel} and references therein).

\begin{theo} {\bf (Lyndon-Hochshild-Serre spectral sequence)} 
Let $H$ be any group, and $R$ a normal subgroup, and let $ G = H / R$ be the factor group.

Assume that $A$ is an $H$-module, so that $A_H$ and $A^H$ are $G$-modules.
Then there are converging first quadrant spectral sequences

$$ E^2_{p,q} = H_p (G ; H_q (R, A)) \Rightarrow H_{p+q}(H, A) $$ 
$$ E_2^{p,q} = H^p (G ; H^q (R, A)) \Rightarrow H^{p+q}(H, A) $$ 
\end{theo}

\begin{ex}
Assume that the group $H$ is a metacyclic group, i.e., we have an exact sequence
$$ 1 \ra R \cong \ZZ/n \ra H \ra G \cong \ZZ/m \ra 1. $$

As previously remarked, $H_2 (\ZZ/n, \ZZ) =0$, $H_1 (\ZZ/n, \ZZ) =\ZZ/n$, while $H_0 (\ZZ/n, \ZZ) $ is the trivial
$\ZZ/m$-module $\ZZ$. Hence the spectral sequence has many zeros, and reveals an exact sequence
$$  \ZZ/m  \ra H_1 ( \ZZ/m, \ZZ/n) \ra H_2 (H, \ZZ) \ra 0. $$
Let $x$ be a generator of $ R \cong \ZZ/n$, and $y$ a generator of $   G \cong \ZZ/m $, which we shall identify to 
one of its lifts to $H$ (the order of $y$ may be assumed to be  equal to $m$ if and only if the exact sequence splits,
or, as one says $H$ is {\bf split metacyclic}).

Now, conjugation by $y$ sends $x$ to $ x^a$, where  $(a,n) = 1$, and we have also $ a^m \equiv 1 \ ( mod \ n)$.
Hence  the $G= \ZZ/m$-structure of $ R \cong \ZZ/n$,  is described by $y (b) = ab \in   \ZZ/n$.

Tensoring $R$ with the bar-complex of $G$ we find
  $$ \oplus_{g_1, g_2  \in G}  (\ZZ/n (g_1, g_2))  \ra \oplus_{g \in G} (\ZZ/n) g \ra \ZZ /n  $$ 
where however the last homomorphism is no longer the zero map, and its image is just 
the image of multiplication by $(a-1)$, i.e., $$(a-1) (\ZZ /n) =
(n, a-1) \ZZ / n \ZZ. $$

In fact
$ \ZZ/ n  \cong \ZZ[G] \otimes_{ \ZZ[G]  } \ZZ/n $, with  $ y^i \cdot 1 \equiv a^i \in \ZZ/n$. Therefore we obtain 
that $H_1 ( G, \ZZ/n)$ is the quotient of the kernel of multiplication by $(a-1)$, inside the direct sum
$ \oplus_{g \in G} (\ZZ/n)  g$,  by the subgroup generated by $g_1 g_2 \equiv g_1 + g_1(1) g_2$.

Writing the elements in $G$ multiplicatively, $g_1 = y^i, g_2 = y^j$, we get the relations
 $$  y^{i+j}  \equiv y^i  + a^i y^j. $$
 For $i=0$ we get thus $y^0 = 0$, and then inductively we find the equivalent relations
 $$ y^i = ( 1 + \dots +  a^{i-1} )y , \ \forall 1 \leq i \leq m-1, \ 0 = y^m =  ( 1 + \dots +  a^{m-1} )y.$$

 Hence $H_1 ( G, \ZZ/n) \subset  \ZZ/(n, r)$, where we set $ r : =  ( 1 + \dots +  a^{m-1} ) = (a^m -1)/(a-1)$,
 is the kernel of multiplication by $(a-1)$, hence it is cyclic of order $$d: = \frac{1}{n}(n,a-1)(n,r).$$
 
 The conclusion is that $H_2 (H, \ZZ)$ is a cokernel of $ \ZZ/m \ra \ZZ/d$.
 
 Edmonds shows in \cite{edmonds2}, with a smart trick, that the map  $ \ZZ/m \ra \ZZ/d$ is the zero map when we have a split metacyclic
 extension, hence $H_2 (H, \ZZ) \cong  \ZZ/d$ in this case. We shall show how to make the computation directly
 using Hopf's theorem: this has the advantage of writing an explicit generator inside $ R / [F,R]$, a fact which shall be proven
 useful in the sequel.
\end{ex}

\subsection{Calculating $H_2(G, \ZZ)$ via combinatorial group theory}

We treat again the case of a split metacyclic group $H$.

$H = F /R $, where $F$ is a free group on two generators $x,y$, and the subgroup $R$ of relations is normally generated by
$$  \xi : = x^n, \eta: = y^m , \zeta : = x^a y x^{-1} y^{-1} ,$$
with $(a,n)=1$, $a^m \equiv 1 ( {\rm mod}\  n)$.

Since $R$ is the fundamental group of the Cayley graph of $H$, whose vertices corresponding to the elements of 
$$H = \{ x^i y^j | 0 \leq i \leq n-1, 0 \leq j \leq m-1 \},$$
and $R$ is a free group on $mn+1$ generators. 

We use now the Reidemeister-Schreier algorithm (see \cite{Magnusks}), which shows that a basis of $R$ is given by:
$$ \eta_i : = x^i y^m x^{-i}, \   \xi : = x^n, \  \xi_{i,j} : = x^i y^j x  y^{-j } x ^{ - [i + a^j] }  ,\  0 \leq i \leq n-1, 1 \leq j \leq m-1,$$
and where $[b] \in \{ 0, \dots, n-1\}$ is the unique positive representative of the residue class of $b$ in $\ZZ / n$.

We want to calculate  $ H_2(G, \ZZ) $ as the kernel of $ R / [F,R] \ra \sK \ra 0$,
  keeping in mind that, since $F^{ab} = F / [F,F] $ is the free abelian group on generators $X,Y$, then 
  $  \sK   = ker (  F^{ab} \ra   G ^{ab} ) $ is generated by the images of  $  \xi , \eta  , \zeta  $ in $F^{ab}$,
  i.e., $ n X, m Y, (a-1) X$.
  
  We conclude that $  \sK  $ is a free abelian group with basis  $ (n, a-1) X, m Y,$
  and contains the  free abelian group with basis  $ n X, m Y.$
  
  Our trick is then to observe that $ \xi \mapsto nX, \ \eta = \eta_0  \mapsto mY$, hence $H_2(G, \ZZ)$ is the kernel of the surjection
 $$  Q: = (R / [F,R])/ (\ZZ \xi + \ZZ \eta))  \ra \sK / ( \ZZ (nX) \oplus \ZZ (mY)) = (n, a-1) \ZZ/ n \ZZ  \ra 0.$$
 
 This shall simplify our calculations considerably. 
 
 Observe that $  R / [F,R]$, since $  [R,R] \subset [F,R]$, is a quotient of $R^{ab}$.  $R^{ab}$ is a free abelian group on generators
 $\Xi, \Xi_{i,j}, \Theta_i$, which are the respective images of the elements  $\xi, \xi_{i,j}, \eta_i$.
 
 Thus $  R / [F,R]$, is the quotient of $R^{ab}$  by the relations of the form $$   r \sim f r f^{-1}, \ r \in R , f \in F.$$
 
 Down to earth, each generator is put to be equivalent to its conjugate by $x$, respectively  its conjugate by $y$.
 
 Conjugation by $x$ clearly sends $\eta_i $ to $\eta_{i+1}$, for all $ i < n-1$, while $ x \eta_{n-1} x^{-1} =  x^n y^m x^{-n}=
 \xi \eta \xi^{-1}$, thus in the abelianization we get $\Theta_i = \Theta , \forall i$, and in our further quotient $Q$ the
 classes of the $\Theta_i $'s are all zero. 
 
 Conjugation of $\xi_{i,j}$ by $x$ makes $$\xi_{i,j} \sim x x^i y^j x y^{-j} x^ { - [i + a^j]} x^{-1} = x^{i+1} y^j x y^{-j} x^ { - 1- [i + a^j]} = \xi_{i+1,j} \ {\rm or}\ \xi_{i+1,j}  \xi^{-1},$$
 hence the classes of $\Xi_{i,j}$ and of $\Xi_{i+1,j}$ are the same. Therefore the class of $\Xi_{i,j}$ is always equal to the one of $\Xi_{0,j}$.
 
 Conjugation by $y^j$ of $\xi = x^n$ yields that the class of $y^j x^n y^{-j}$ is trivial. But, observing that  the class of $x^n = \xi$ and all its conjugates is trivial, we can consider the exponents of  powers of $x$ as just residue
 classes modulo $n$, and avoid to use the square brackets. Then

 $$y^j x^n y^{-j} = (y^j x y^{-j})^n \sim (y^j x y^{-j} x^{- a^j} ) ( x^{a^j}y^j x y^{-j} x^{- 2 a^j} ) \dots ( x^{(n-1)a^j}y^j x y^{-j} x^{- n a^j} ) ,$$
 which means that the class of  $ n \Xi_{0,j} = 0$. We denote by $Z_j$ the class of $\Xi_{0,j} $. Hence $ n Z_j = 0$,
 and $Q$ is generated by the classes $Z_j$.
 
 In the calculation of conjugation by $y$ on $\Xi_{0,j}$ we observe that 
 $$ \xi_{i,j} \sim y  y^j x y^{-j} x^ { - [ a^j]} y^{-1}  = (y^{j+1} x  y^{-j-1} x^ { - [a^{j+1}]}) x^ {  [a^{j+1}]} y x^ { - [ a^j]} y^{-1} .$$
 
 Hence 
 
 $$ Z_j \sim y Z_j y^{-1} = Z_{j+1}  x^ {  a^{j+1}} y x^ { - a^j} y^{-1} = Z_{j+1}  x^ {  a^{j+1}} y x y ^{-1} x^ { -  a^j -a}   x^ {   a^j +a} y x^ { -  a^j - 1} y^{-1} =$$
  $$  = Z_{j+1}  Z_1^k x^ {  a^{j+1}+ ka } y  x^ { -  a^j -k} y ^{-1} .$$
  
  For $k = - a^j$ we obtained the desired relation 
  $$  Z_j =  Z_{j+1}  Z_1^{ - a^j}  \Leftrightarrow Z_{j+1} =  Z_j  Z_1^{ a^j},$$
  which can be written additively  as  
  $$  y^{i+1}  \equiv y^i  + a^i y. $$
  
  The conclusion is that we have a cyclic group whose order is  $d: = \frac{1}{n}(n, a-1)(n,r)$ and with generator 
  $  \frac{n}{(n, a-1)} Z_1  $, where  $ Z_1  $ is  the class of $\xi_{0,1} =  y x  y^{-1 } x ^{ - a}$.
  
 \subsection{Sheaves and cohomology on quotients, linearizations}
 
 The method of the bar-complex apparently settles the problem of calculating the cohomology groups of a $ K(G,1)$ space, and also to calculate the
 cohomology groups $ H^n (X, \QQ)$ for a rational $ K(G,1)$ space $X$, except that one must have a good knowledge of the group $G$,
 for instance a finite  presentation does not always suffice.
 
 More generally, Grothendieck in \cite{tohoku} approached the question of calculating the cohomology of a quotient 
 $ Y = X /  G$ where $X$ is not necessarily contractible. Using sheaf theory, one can reach a higher generality,
 considering $ G$-linearized sheaves, i.e. sheaves $\sF$ on $X$ with a lifting of the  action of $G$ to
 $\sF$ (this is the analogue to taking cohomology groups on non trivial $G$-modules).
 
  Since this concept is very useful in algebraic geometry, and lies at the heart of some important calculations,
  we shall try now to briefly explain it, following Mumford's treatment (\cite{abvar}, appendix to section 2, pages 22-24,
  section 7, pages 65-74, section 12, pages 108-122).
  
  The following is the basic example we have in mind, where we consider two formally different cases:
  
  \begin{enumerate}
  \item
  $X$ is a projective variety, $G$ is a finite group acting on $X$ with the property that for each point $x$ the orbit $ G x$ is contained in an affine open subset
  of $X$ (this property holds if there is an ample divisor $H$ whose linear equivalence class is fixed by $G$), or
  \item
  $X$ is a complex space, $G$ acts properly discontinuously on $X$ (in case where the quotient $ X/G$ is a projective variety, we can apply the G.A.G.A. principle
  by which there is an equivalence of categories between coherent algebraic sheaves and coherent holomorphic sheaves).
  
  \end{enumerate}
  
  Assume that we have a vector bundle $V$ on a quotient variety $ X/G$. Denote by $ p : X \ra X/G$
  the quotient map  and  take the fibre product $p^*(V) : =  V \times_{X/G} X$.
  We consider the action of $G$ on $V \times X$, which is trivial on the factor $V$, and the given one on $X$.
  In this way the action of $G$ extends to the pull back $ p^*(V)  =V \times_{X/G} X$ of $V$, and the sections of $V$ on $X/G$
  are exactly seen to be the $G$-invariant sections of $ p^*(V)$.
  
  Now, a main goal is to construct interesting varieties  as quotient varieties $ X/G$, and then study line bundles on them;  
  the following result is quite useful for this purpose.
  
  \begin{prop}
  Let $ Y = X/G$ be a quotient algebraic variety, $p : X \ra Y$ the quotient map.
  
  (1) Then there is a functor  between (the category of) line bundles $\sL'$ on $Y$ 
  and (the category of) $G$-linearized line bundles  $\sL$ on $X$, associating to $\sL'$ its pull back $p^* (\sL')$. The functor   $\sL \mapsto p_*(\sL)^G$ is a right inverse
  to the previous one, and $p_*(\sL)^G$ is invertible if the action is free, or if  both $X$ and $Y$ are  smooth \footnote{One can indeed relax the hypothesis, requiring only that $X$ 
  is normal and $Y$ is smooth.}.
  
  (2) Given a line bundle  $\sL$ on $X$, it admits a $G$-linearization if and only if there is a Cartier divisor $D$ on $X$ which is 
  $G$-invariant and such that $\sL \cong \hol_X(D)$ (observe that  $ \hol_X(D)= \{ f \in \CC(X) | div(f) + D \geq 0\}$ has an obvious $G$-linearization). 
  
  (3) A necessary condition for the existence of a linearization is that
  $$ \forall g \in G, \ g^* ( \sL ) \cong \sL .$$ If this condition holds, defining the Thetagroup of $\sL$ as:
  $$ \Theta (\sL, G) : =  \{ (\psi , g ) | g \in G, \ \psi : g^* ( \sL ) \cong \sL \ {\rm isomorphism} \},$$
  there is an exact sequence
  $$ 1 \ra \CC^* \ra  \Theta (\sL, G) \ra G \ra 1.$$
  
  The splittings of the above sequence correspond to the    $G$-linearizations of $\sL$.
  
  If the sequence splits, the  linearizations are a principal homogeneous space over the dual group $ \Hom (G, \CC^*) = : G^*$ of $G$
(namely, each linearization  is obtained from a fixed one  by multiplying with an arbitrary element in $ G^*$).
  
   \end{prop}
   
   \Proof
   In (1), the case where $G$ acts freely is trivial in the holomorphic context, and taken care of by proposition 2, page 70 of \cite{abvar}
   in the algebraic context.
   
   For (2) and (3) we refer to \cite{GIT}, section 3, chapter 1, and \cite{cat3}, lemmas 4-6.
   
   We only provide an argument for the last statement in (1), for which we do not know of a precise reference.
   The question is local, hence we may assume that $P$ is a smooth point of $X$, and that $G$ is equal to the stabilizer of $P$, hence it is a finite group.
   The invertible sheaf is locally isomorphic to $\hol_X$, except that this isomorphism does not respect the linearization.
   The $G$-linearization of $\sL$ corresponds then locally to a character $\chi : G \ra\CC^*$. In particular, there exists an integer
   $m$ such that the induced linearization on $\sL^{\otimes m}$ corresponds to the trivial character on $\hol_X$.
   
   Let $R$ be the ramification divisor of the quotient map $ p : X \ra Y = X/G$, and let $X^0 : =  X \setminus Sing(R)$,
   which is $G$-invariant. Set $Y^0 : = p (X^0 )$. 
   
   Step I: $p_*(\sL)^G$ is invertible on $Y^0$. 
   
   This is clear on the open set  which is the complement of the branch set $ p(R)$. At the points of $R$ which are smooth points,
   then the $G$ action can be linearized locally, and it  is a pseudoreflection;  so in this case the action of $G$ involves only the last variable,
   hence the result follows from the dimension 1 case, where every torsion free module is locally free.
   
   Step II. Let $ \sM : = p_*(\sL)^G$. We claim   that $\sM$ is invertible  at every point.
   
   Pick a point $y$ not in $Y^0$ and a local holomorphic Stein  neighbourhood  $U$ of $y$, say biholomorphic to a ball.
   Then $U^0 : = U \cap Y^0$ has vanishing homology groups $H^i(U^0 , \ZZ)$ for $i=1,2$.
   By the exponential sequence $H^1(U^0 , \hol^*_Y) \cong H^1(U^0 , \hol_Y)$, which is a vector space hence has no torsion elements
   except the trivial one.  Since $\sM^{\otimes m}$ is trivial on $U^0 $,
   we conclude that $\sM$ is trivial on $U^0 $.  By Hartogs' theorem  there is then an isomorphism between $G$-invariant sections 
   of $\sL$ on $p^{-1} (U)$ and sections of $\hol_Y$ on $U$.

   \qed
   
   Thus, the question of the existence of a linearization is reduced to the algebraic question of the splitting of the central extension 
   given by the Theta group. This question is addressed by group cohomology theory, as follows (see \cite{BAII}).
   
   \begin{prop}
  An exact sequence of groups 
  $$ 1 \ra M \ra \Ga \ra G \ra 1,$$
  where $M$ is abelian, is called an extension. Conjugation by lifts of elements of $G$ makes $M$ a $G$-module.
 Each choice of a lift $e_{g} \in \Ga$ for every element $ g \in G$ determines a 2-cocycle 
 $$ \psi (g_1,g_2) : = e_{g_1}   e_{g_2}  e_{g_1g_2}^{-1}.$$
 
The cohomology class $ [\psi] \in H^2(G, M)$ is independent of the choice of lifts, and in this way $H^2(G, M)$  is in bijection with
the set of  strict
isomorphism classes of extensions $ 1 \ra M \ra \Ga' \ra G \ra 1$ (i.e., one takes isomorphisms inducing the identity on $M$ and $G$).
Whereas isomorphism classes of extensions $ 1 \ra M' \ra \Ga' \ra G \ra 1$ (they should only induce the identity on $G$) are classified by 
isomorphism classes of pairs $(M', \psi')$.
   \end{prop}
   
   In particular, we have that the extension is a direct product $ M \times G$ if and only if $M$ is a trivial $G$-module (equivalenty,
   $M$ is in the centre of $\Ga$) and the class $ [\psi] \in H^2(G, M)$ is trivial.
   
    \begin{cor}
   Let $\sL$ be an invertible sheaf whose class in $Pic(X)$ is $G$-invariant. 
   Then the  necessary and sufficient condition for the existence of a linearization    is  the triviality  of the extension  class 
  $ [\psi] \in H^2(G, \CC^*)$ of the  Thetagroup $\Theta (G, \sL)$.
  \end{cor}
  
  The group $H^2(G, \CC^*)$ is the group of Schur multipliers (see again \cite{BAII}, page 369) , which for a finite group is, as we saw in theorem  \ref{Schur}, another incarnation of
  the second homology group $H_2(G, \ZZ)$. It occurs naturally when we have a projective representation of a group $G$.
  Since, from a homomorphism $ G \ra \PP GL(r, \CC)$  one can pull back the central extension
  $$ 1 \ra \CC^* \ra    GL(r, \CC)  \ra  \PP GL(r, \CC) \ra 1  \Rightarrow  1 \ra \CC^* \ra    \hat{G}  \ra G\ra 1,$$
  and the extension class $[\psi] \in   H^2(G, \CC^*)$ is the obstruction to lifting the projective representation
  to a linear representation $ G \ra  GL(r, \CC)$.
  
  It is an important remark that, if the group $G$ is finite, and $ n  = ord (G)$, then the cocycles take value in the group 
  of roots of unity $\mu_n : = \{ z \in \CC^* | z^n = 1 \}$.

   \begin{ex}
  Let $E$ be an elliptic curve, with the point  $O$ as the origin,and let $G$ be the group of 2-torsion points $ G : = E[2]$ acting by translation
  on $E$. The divisor class of $ 2 O$ is never represented by a $G$-invariant divisor, since all  the $G$-orbits consist of $4$ points,
  and the degree of $2O$ is not divisible by $4$. Hence, $\sL : = \hol_E (2O)$ does not admit a $G$-linearization.
  However, we have a projective representation on $ \PP^1 = \PP (H^0( \hol_E (2O))$, where each non zero element $\eta_1$ of the group
  fixes 2 divisors: the sum of the two points corresponding to $\pm \eta_1/2$, and its translate by another point $\eta_2 \in E[2]$.
  
  The two group generators  yield two linear transformations, which act on $ V = \CC x_0 \oplus \CC x_1$ as follows:
  $$ \eta_1(x_0) = x_1,  \eta_1(x_1) = x_0,  \ \eta_2 (x_j) = (-1)^j x_j.$$ 
  The linear group generated is however $D_4$, since $$\eta_1\eta_2 (x_0) = x_1, \ \eta_1\eta_2 (x_1) = - x_0.$$
  \end{ex}
  
  \begin{ex}
The previous example is indeed a special case of  the Heisenberg extension, and $V$ generalizes to the
 Stone von Neumann  representation associated   to an abelian group $G$, which is nothing else
 than the space $ V : = L^2 (G, \CC)$ of square integrable functions on $G$ (see \cite{igusa},\cite{abvar}).
 
 $G$ acts on  $ V =  L^2 (G, \CC)$ by translation $ f (x) \mapsto f (x - g)$, $G^*$ acts by multiplication with the given character
 $ f (x) \mapsto f (x) \cdot \chi ( x))$, and the commutator $[g, \chi ]$ acts by the scalar multiplication with the constant $ \chi (g)$.
 
 The Heisenberg group is the group of  automorphisms of $V$  generated by $G, G^*$ and by $\CC^*$ acting by  scalar multiplication,
 and there is a central extension
 $$ 1 \ra \CC^* \ra Heis(G) \ra G \times G^* \ra 1 ,$$
 whose class is classified by the $\CC^*$-valued  bilinear form  $(g, \chi) \mapsto  \chi (g)$, an element of 
 $$ \wedge^2 ( Hom ( G \times G^*, \CC^*)) \subset  H^2 ( G \times G^*, \CC^*).$$
 
 The relation with Abelian varieties $ A = V / \Lam$ is through the Thetagroup associated to an ample divisor $L$.
 
 In fact, since by the theorem of Frobenius the alternating form $ c_1(L) \in H^2(A, \ZZ) \cong \wedge^2 ( Hom (\Lam, \ZZ))$ admits,
 in a suitable basis for $\Lam$, the normal form
 \[ 
 D: =   \begin{pmatrix}
    0  &    D' \\
  - D' &    0
  \end{pmatrix}\\
  D' : = \diag (d_1, d_2, \dots , d_g), d_1| d_2 | \dots  | d_g.
\] 
The key property is that, if one sets $ G : = \ZZ^g/ D' \ZZ^g$, then $L$ is invariant exactly for the  translations in $( G \times G^*) \cong ( G \times G ) \subset A$, and  the Thetagroup of $L$  is just 
isomorphic to the Heisenberg group $ Heis(G)$. 

  \end{ex}
  
  The nice part of the story is the following useful result, which was used by Atiyah in the case of elliptic curves,  to study vector bundles over these (\cite{atiyah}).
  
  \begin{prop}\label{Heisenberg}
  Let $G$ be a finite abelian group, and let $ V : = L^2 (G, \CC)$ be the Stone-von Neumann representation.
  Then $V \otimes V^{\vee}$ is a representation of $( G \times G^*)$ and splits as the direct sum of all the characters of  $( G \times G^*)$,
  taken with multiplicity one.
  \end{prop}
  \Proof
  Since the centre $\CC^*$ of the Heisenberg group $ Heis(G)$ acts trivially on $V \otimes V^{\vee}$, 
    $V \otimes V^{\vee}$ is a representation of $( G \times G^*)$.
  Observe that $( G \times G^*)$ is equal to its  group of characters, and its cardinality equals to the dimension of  $V \otimes V^{\vee}$,
  hence it shall suffice, and it will be useful  for applications, to write for each character an explicit eigenvector.
  
  We shall use the notation $g, h, k$ for elements of $G$, and the notation $\chi, \eta, \xi$ for elements in the dual group.
  Observe that $V$ has two bases, one given by $\{ g \in G\}$, and the other given by the characters $\{ \chi \in G^* \}$.
  In fact, the Fourier transform yields an isomorphism of the vector spaces $ V : = L^2 (G, \CC)$ and  $ W : = L^2 (G^*, \CC)$
  $$ \sF (f) := \hat{f}, \  \hat{f}(\chi) : = \int  f(g) ( \chi, g) dg.$$
  The action of $h \in G$ on $V$ sends $ f(g) \mapsto f (g-h)$, hence for the characteristic functions in $\CC[G]$, $ g \mapsto g + h$.
  While   $\eta \in G^*$ sends $ f \mapsto  f \cdot \eta$, hence $ \chi \mapsto \chi + \eta$, where we use the additive notation also for the group of characters. 
  
Restricting $V$ to the finite Heisenberg group which is a central extension of $ G \times G^*$ by $\mu_n$, we get a unitary representation,
hence we identify $V^{\vee}$  with $\bar{V}$. This said, a basis of $V \otimes \bar{V}$ is given by the set $\{ g \otimes \bar{\chi}\}$.

Given a vector $ \sum_{g, \chi} a_{g, \chi}  (g \otimes \bar{\chi} )$ the action by $h \in G$ sends it to
$$ \sum_{g, \chi}  (\chi, h) a_{g -h, \chi}  (g \otimes \bar{\chi} ),$$ while the one by $ \eta \in G^*$ sends it to 
$$ \sum_{g, \chi}  (\eta, g) a_{g, \chi  - \eta}  (g \otimes \bar{\chi} ).$$

Hence one verifies right away that 
$$  F_{k, \xi} : =  \sum_{g, \chi}  ( \chi - \xi, g - k)  (g \otimes \bar{\chi} ) $$ 
is an eigenvector with eigenvalue $ (\xi, h) (\eta, k)$ for $ (h, \eta) \in ( G \times G^*)$.
  
  \qed
  
 \begin{rem}
 The explicit calculation reproduced above is based on the fact that the Fourier transform does not commute with the action of the Heisenberg group:
 the action of $G$ on $V$ corresponds to the action of $G = (G^*)^*$ on $W$, while the action of $G^*$ on $W$ corresponds to the multiplication
 of functions $ f(g) \in V$ with the conjugate $ \bar{\chi}$ of the character $\chi$.
 \end{rem} 

   In the next subsection we shall give an explicit example where the Heisenberg group is used for a geometrical construction.
   Now let us return to  the general discussion of cohomology of $G$-linearized sheaves. 
 
  In the special case  where the action is free there is an isomorphism between the category 
  of  $G$-linearized sheaves $\sF$ on $X$, and sheaves $\sF^G$ on the quotient $Y$; we have then
  
  \begin{theo} {\bf (Grothendieck, \cite{tohoku})}\label{tohoku}
  Let $ Y = X /  G$, where $G$ acts freely on $X$, and let $\sF$ be a $G$-linearized sheaf  on $X$, let $\sF^G$ be  the $G$-invariant direct image
  on $Y$.
  
  Then there is a spectral sequence converging to a suitable graded quotient of $H^{p+q}(Y, \sF^G)$
  and with $E_1$ term equal to $H^p (G, H^q(X, \sF)) $; we write then as a shorthand notation:
  $$ H^p (G, H^q(X, \sF)) \Rightarrow H^{p+q}(Y, \sF^G).$$
  \end{theo}
  
  The underlying  idea is simple, once one knows (see \cite{cartaneil}, \cite{h-s}) that there is a spectral sequence
  for the derived  functors of a composition of two functors: in fact one takes here  the functor $ \sF \mapsto  H^0(X, \sF)^G$,
  which is manifestly the composition of two functors $  \sF \mapsto  H^0(X, \sF)$ and $A \mapsto A^G$,
  but also the composition of two other functors, since $H^0(X, \sF)^G = H^0(Y, \sF^G)$.
  
  The main trouble is that the functor  $\sF \mapsto \sF^G$ is no longer exact if we drop the hypothesis that the action is free.
  
  A particular case is of course the one where $\sF$ is the constant sheaf $\ZZ$: in the case where the action is not free
  one has in the calculation to keep track of stabilizers $G_x$ of points $ x \in X$, or, even better, to keep track of fixed loci
  of subgroups $H \subset G$. 
 
 Since this bookkeeping is in some situations too difficult or sometimes also of not so much use, Borel (\cite{borelSTG})
 proposed  (in a way which is similar in spirit  to the one of Grothendieck) to consider a generalization of the
 cohomology of the quotient (see \cite{fultonLectures} for a nice account of the theory and its applications in algebraic geometry).
 
 \begin{defin}
 Let $G$ be a Lie group acting on a space $X$, and let $BG = EG / G$ be a classifying space for $G$.
 Here we do not need to specify whether $G$ acts on the left or on the right, since 
 for each left action 
 $(g,x) \mapsto gx$ there is the mirror action $xg : = g^{-1} x$. 
  
   Letting  $R$ be an arbitrary ring of coefficients,
 the
 {\bf equivariant cohomology} of $X$ with respect to the action of $G$ is defined in one of the following equivalent ways.

 1) Assume that $G$ acts on both spaces $X$ and $EG$ from the left: then we have the diagonal action $ g (e,x) : = (ge, gx)$
 and  $$ H^i_G(X, R): =   H^i ((EG \times X)/ G, R) .$$
 
  2) (as in \cite{fultonLectures}) Assume without loss of generality that $G$ acts on $X$ from the left, and on $EG$ from the right:
then 
 $$ H^i_G(X, R): =   H^i (EG \times^G X, R) ,$$
  where the space $EG \times^G X $ is defined to be the quotient of the Cartesian product $ EG \times X$ by the relation
 $$ (e \cdot g, x) \sim  ( e, g \cdot x).$$ 

 \end{defin}
 
\begin{rem}
Observe that, in the special case where $X$ is a point, and $G$ is discrete, we reobtain the cohomology group $H^i (G , R)$.
\end{rem}

\begin{ex}
 The consideration of arbitrary Lie groups is rather natural, for instance if $G = \CC^*$, then $BG = \PP^{\infty}_{\CC} =
( \CC^{\infty} \setminus \{0\} )/ \CC^*.$

Similarly the infinite Grassmannian $\sG r (n, \CC^{\infty} ) $ is a classifying space for $\GL (n, \CC)$,
it is a quotient (cf. \cite{milnorCC}) of the Stiefel manifold of frames  $$\sS t (n, \CC^{\infty} ) : = 
\{ (v_1, \dots, v_n ) | v_i \in  \CC^{\infty} , \rank \langle v_1, \dots, v_n  \rangle = n\}.$$

Clearly the universal vector subbundle on the Grassmannian is obtained as the invariant direct image  of
a $\GL (n, \CC)$- linearised vector bundle on the Stiefel manifold.
 \end{ex}
 
 The names change a little bit: instead of a linearized vector bundle one talks here of $G$-equivariant vector bundles,
 and there is a theory of equivariant characteristic classes, i.e., these classes are equivariant for morphisms of $G$-spaces
 (spaces with a $G$-action).
 
 Rather than dwelling more on abstract definitions and general properties, I prefer at this stage to mention the most
 important issues of $G$-linearized bundles for geometric rational $K (\pi, 1)$'s.
 
 \subsection{Hodge Bundles of weight  {\bf = 1.}}

{\bf 1) Hodge Bundle for families of Abelian varieties}.

 Consider  the family of  Abelian Varieties of dimension $g$, and with
a polarization of type $ D' \sim (d_1, d_2, \dots , d_g)$ over Siegel's upper half space
$$  \sH_g : = \{ \tau \in \Mat(g,g,\CC)| \tau = ^t \tau, \im (\tau) > 0 \}.$$ 

The family is given by the quotient of the trivial bundle $H : = \sH_g \times \CC^g$
by the action of $ \la \in \ZZ^{2g}$ acting by

$$  (\tau, z) \mapsto ( \tau, z + (D', - \tau) \la ).$$ 

The local system corresponding to the first cohomology groups of the fibres is the $  \Sp (D, \ZZ) $-linearized
local system (bundle with fibre  $\ZZ^{2g}$) on $H_{\ZZ}  : = \sH_g \times \ZZ^{2g}$, with the  obvious action of $  \Sp (D, \ZZ) $

$$  M \in \Sp (D, \ZZ) , \ M ( \tau, \la) = ( M  (\tau ) , M \la), $$
 \[  
 M: =   \begin{pmatrix}
    \al  &    \be \\
   \ga &    \de
  \end{pmatrix} , \ \
 M ( \tau) =  - D' (D' \al - \tau \ga)^{-1} ( D' \be - \tau \de).
\] 

Obviously $ H_{\ZZ} \otimes _{\ZZ} \CC$ splits as $ H^{1,0} \oplus H^{0,1}$ and the above vector bundle $H $ equals $ H^{1,0}$
and is called the Hodge bundle.

Stacks are  nowadays the language in order to be able to treat the above local system and the Hodge bundle as living over the 
 quotient moduli space $\sA_{g,D'} : =  \sH_g /  \Sp (D, \ZZ) $.
 
 {\bf 2) Hodge Bundle for families of curves}.
 
 Over Teichm\"uller space we have a universal family of curves, which we shall denote by $p_g :  \sC_g  \ra \sT_g$.
 
 The reason is the following: letting $M$ be here a compact oriented Riemann surface of genus $g$, and $ \sC (M)$  the space of complex structures
on $M$, it is clear that there is a universal tautological family of complex structures
parametrized by $ \sC (M)$, and with total space
$$\mathfrak U_{ \sC (M)} : =  M   \times  \sC (M) ,$$
on which the group   $\sD iff ^+ (M)$ naturally acts, in particular  $\sD iff ^0(M)$.

A rather simple consequence of the  Lefschetz' lemma \ref{LefLemma}  is that $\sD iff ^0(M)$ acts freely on $ \sC (M)$: in fact Lefschetz' lemma \ref{LefLemma}   implies that  for each complex structure
$C$ on $M$ the group of biholomorphisms $Aut(C)$ contains no automorphism (other than the identity) which is homotopic to the identity, hence a fortiori 
no one that is differentiably isotopic
to the identity.

We have the local system $ H_{\ZZ} : = \sR ^1 (p_g)_* (\ZZ)$, and again $$ H_{\ZZ} \otimes _{\ZZ} \CC =  H^{1,0} \oplus H^{0,1}, \ 
\ H^{1,0} : = (p_g)_* (\Omega^1_{ \sC_g |T_g}).$$
\
$H : = H^{1,0} $ is called the Hodge bundle, and it is indeed the pull back of the Hodge bundle on $\sA _g$ by the Torelli map, or period map,
which associates to a complex structure $C$ the Jacobian variety $ \Jac (C) : = H^{1,0}(C) ^{\vee} / H_1 (M, \ZZ)$.

The i-th Chern class of  the Hodge bundle  yield some cohomology classes which   Mumford (\cite{EnsMath} and \cite{mumshaf})
denoted as the $\la_i$-class (while the notation  $\la$-class is reserved for $\la_1$). 

Other classes, which play a crucial role in Mumford's conjecture (\cite{mumshaf},  see also theorem \ref{m-w}), are the classes
$\sK_i$, defined as 
$$ \sK_i : = (p_g)_*  (K^{i+1}), \ \  K : = c_1 (\Omega^1_{ \sC_g |\sT_g}).$$

\subsection{A surface in a Bagnera-De Franchis threefold. }
 
We want to describe here a construction given in joint work with Ingrid Bauer and Davide Frapporti (\cite{bcf}) of a surface
$S$ with an ample canonical divisor $K_S$, and with $K_S^2 = 6$, $p_g : = h^0(\hol_S(K_S)) = 1, \ q : =  h^0(\Omega^1_S)=1.$

Let $A_1$ be an elliptic curve, and let $A_2$ be an Abelian surface with a line bundle $L_2$ yielding a polarization of type $(1,2)$
(i.e., the elementary divisors for the Chern class of $L_2$ are $ d_1 = 1, d_2 = 2$). Take as $L_1$ the line bundle $\hol_{A_1} ( 2 O)$,
and let $L$ be the line bundle on $A' : = A_1 \times A_2$ obtained as the exterior tensor product of $L_1$ and $L_2$,
so that
$$ H^0 (A', L) =  H^0 (A_1, L_1)  \otimes  H^0 (A_2, L_2) .$$

Moreover, we choose the origin in $A_2$ so that the space of sections $H^0 (A_2, L_2) $ consists only of {\em even} sections
(hence, we shall no longer be free to further change the origin by an arbitrary  translation).

We want to take a Bagnera-De Franchis threefold $ X: = A/ G$, where $A = (A_1 \times A_2) / T$, and $ G \cong T \cong \ZZ/2$,
and have a surface $S \subset X$ which is the quotient of a $G\times T$ invariant $D \in |L|$, so that $S^2 = \frac{1}{4} D^2 = 6$.

We write as usual $ A_1 = \CC / \ZZ \oplus \ZZ \tau$, and we let $A_2 = \CC^2 / \Lam_2$, with $\la_1,\la_2, \la_3, \la_4$
a basis of $ \Lam_2$ on which the Chern class of $L_2$ is in Frobenius normal form. We let then 
\begin{equation}\label{BCF1}  G : = \{ Id, g \}, \ \ g (a_1 + a_2 ) : = a_1 + \tau /2 - a_2 + \la_2 / 2, \\\forall a_1 \in A_1, a_2 \in A_2\\
  \end{equation} 
  \begin{equation}\label{BCF2}  
  T : = ( \ZZ/2 ) ( 1/2 + \la_4 / 2) \subset A = (A_1 \times A_2) .
  \end{equation} 

Now, $G \times T$ surjects onto the group of two torsion points $A_1[2]$ of the elliptic curve (by the homomorphism associating to an affine transformation its translation vector)
and also on the subgroup $ ( \ZZ/2 )  ( \la_2 / 2) \oplus ( \ZZ/2 )  ( \la_4 / 2) \subset A_1[2]$,
and both $H^0 (A_1, L_1) $ and $  H^0 (A_2, L_2)$ are the Stone-von Neumann representation
of the Heisenberg group which is a central $\ZZ/2$ extension of $ G \times T$.

By proposition \ref{Heisenberg}, since in this case (recall the notation  $A = V / \Lam$), we have  $ V \cong \bar{V}$ (the only roots of unity occurring are just $\pm 1$),
we conclude that there are exactly 4 divisors in $|L|$, invariant by $ a_1 + a_2 \mapsto a_1 - a_2$,
and by $a_1 + a_2 \mapsto  a_1 + \tau /2 + a_2 + \la_2 / 2$, and  $a_1 + a_2 \mapsto  a_1 + 1 /2 + a_2 + \la_4 / 2$.

Hence these four divisors descend to give four surfaces $S \subset X$.

This construction is used in \cite{bcf} to prove the following.

\begin{theo}\label{BCF}
Let $S$ be a surface of general type with invariants $K_S^2 = 6$, $p_g =q = 1$ such that there exists an unramified double cover
$ \hat{S} \ra S$ with $ q (  \hat{S} ) = 3$, and such that the Albanese morphism $ \hat{\al} :  \hat{S}  \ra A$ is birational onto it image $Z$,
a divisor in $A$ with  $ Z^3 = 12$.

Then the canonical model of $\hat{S}$ is isomorphic to $Z$, and the canonical model of $S$ is isomorphic to $Y = Z / (\ZZ/2)$, 
which a divisor in a Bagnera-De Franchis threefold $ X: = A/ G$, where $A = (A_1 \times A_2) / T$, $ G \cong T \cong \ZZ/2$,
and where the action is as in (\ref{BCF1}, \ref{BCF2}).

These surfaces exist, have an irreducible four dimensional moduli space, and their Albanese map $\al : S \ra A_1 = A_1/ A_1[2]$ has 
general fibre a non hyperelliptic curve of genus $g=3$.

\end{theo}

\Proof
By assumption the Albanese map  $ \hat{\al} :  \hat{S}  \ra A$ is birational onto $Z$, and we have $ K_{ \hat{S}}^2 = 12 = K_Z^2$,
since by adjunction $\hol_Z(Z)$ is the dualizing sheaf of $Z$ (so $Z$ restricts to the canonical divisor $K_Z$ of $Z$).

We argue similarly to \cite{bc-inoue}, Step 4 of theorem 0.5, page 31. Denote by $W$ the canonical model of $\hat{S}$, and observe that by adjunction (see loc. cit.)
we have $ K_W =  \hat{\al}^* (K_Z ) - \mathfrak A$, where $\mathfrak A$ is an effective \QQ-Cartier divisor.

We observe now that $K_Z$ is ample (by adjunction, since $Z$ is ample) and   $K_W$ is also ample, hence  we have an inequality,
$$ 12 = K_W^2 =  (\hat{\al}^* (K_Z ) - \mathfrak A)^2 = K_Z^2 - (\hat{\al}^* (K_Z ) \cdot  \mathfrak A) - (K_W \cdot  \mathfrak A) \geq K_Z^2 = 12,$$
and since both terms are equal to $12$, we conclude that $\mathfrak A= 0$, which means that $K_Z$ pulls back to $K_W$ hence
$W$ is isomorphic to $Z$. We have  a covering involution $ \iota : \hat{S} \ra \hat{S}$, such that $ S = \hat{S} / \iota$. 
Since the action of $\ZZ/2$ is free on $\hat{S}$, $\ZZ/2$  also acts freely on $Z$.

Since $Z^3 = 12$, $Z$  is a divisor of type $(1,1,2)$ in $A$.  The covering involution $ \iota : \hat{S} \ra \hat{S}$
can be lifted to an involution $g$ of $A$, which we write as an affine transformation  $ g (a) = \al a + \be$. 

 We have Abelian subvarieties $A_1 = \ker (\al - Id)$, $A_2 = \ker (\al + Id)$,
and since the irregularity of $S$ equals $1$, $A_1$ has dimension $1$, and $A_2$ has dimension $2$.

We observe preliminarly  that $g$ is fixed point free: since otherwise the fixed point locus would be non empty 
 of  dimension one, so it would intersect the ample divisor $Z$, contradicting that $ \iota : Z \ra Z$ acts freely.
 
 Therefore $Y = Z / \iota $ is a divisor in the Bagnera- De Franchis threefold $ X = A / G$, where $G$ is the group of order two generated by $g$.
 
 We can then write (using the notation introduced in prop. \ref{BdF}) the Abelian threefold $A$ as $ A' / T$, and since $ \be_1 \notin T_1$ we have only  two possible cases.
 
 Case 0) : $ T = 0$.
 
 Case 1) : $ T \cong \ZZ/2$.
 
 We further observe that since  the divisor $ Z$ is $g$-invariant, its polarization is $\al$ invariant,
in particular its Chern class $c \in \wedge ^2 ( \Hom (\Lam, \ZZ))$, where $A = V / \Lam$.

Since $ T =  \Lam / ( \Lam_1 \oplus \Lam_2)$, $c$ pulls back to 
$$c'   \in \wedge ^2 ( \Hom ( \Lam_1 \oplus \Lam_2, \ZZ)) = \wedge ^2 ( \Lam_1^{\vee})  \oplus \wedge ^2 ( \Lam_2^{\vee}) \oplus (\Lam_1^{\vee})\otimes (\Lam_2^{\vee}) ,$$
and by invariance $c'  = (c'_1 \oplus c'_2 ) \in \wedge ^2 ( \Lam_1^{\vee})  \oplus \wedge ^2 ( \Lam_2^{\vee}) .$

So Case 0) bifurcates in the cases: 

Case 0-I) $c'_1$ is of type $(1)$, $c'_2$ is of type $(1,2)$.

Case 0-II) $c'_1$ is of type $(2)$, $c'_2$ is of type $(1,1)$.

Both cases can be discarded, since they lead to the same contradiction. Set  $D: = Z$: then $D$ is  the  divisor of zeros on 
$A = A_1 \times A_2$ of a section of a line bundle $L$ which is an  exterior tensor product of $L_1$ and $L_2$.
Since
$$ H^0 (A, L) =  H^0 (A_1, L_1)  \otimes  H^0 (A_2, L_2) ,$$
and $H^0 (A_1, L_1)$ has dimension one in case 0-I), while  $H^0 (A_2, L_2)$ has dimension one in case 0-II), 
we conclude that $D$ is a reducible divisor, a contradiction, since $D$ is smooth and connected.

In case 1), we denote $A' : = A_1 \times A_2$, and we let $D$ be the inverse image of $Z$ inside $A'$.
Again $D$ is smooth and connected, since $\pi_1(\hat{S})$ surjects onto $\Lam$. Now $ D^2 = 24$, so the Pfaffian
of $c'$ equals $4$, and there are a priori several  possibilities.

Case 1-I) $c'_1$ is of type $(1)$.

 Case 1-II) $c'_2$ is of type $(1,1)$.
 
 Case 1-III) $c'_1$ is of type $(2)$, $c'_2$ is of type $(1,2)$.
 
 Cases 1-I) and 1-II) can be excluded as case 0), since then $D$ would be reducible. 
 
 We are then left only with case 1-III), and we may, without loss of generality, assume that 
 $H^0 (A_1, L_1)=  H^0 (A_1, \hol_{A_1} (2 O))$,
 and assume that we have chosen the origin so that all the sections of
 $H^0 (A_2, L_2)$ are even.
 
We have $ A = A' / T $, and we may write the generator of $T$ as $t_1 \oplus t_2$,
and write $ g (a_1 + a_2 ) = (a_1 + \be_1) + ( a_2 - \be_2)$.

By the description of Bagnera-De Franchis varieties given in subsection 4.1  we have that $t_1$ and $\be_1$ are  a basis of the group
of $2$ torsion points of the elliptic curve $A_1$.

Now our divisor $D$, since all sections of $L_2$ are even, is $ G \times T$ invariant if and only if it
is invariant by $T$ and by translation by $\be$. 

This condition however implies that translation by $\be_2$ of $L_2$ is isomorphic to $L_2$, 
and similarly for $t_2$. 
It follows that $\be_2, t_2$ are  basis of the kernel $K_2$ of the map $ \phi_{L_2}: A_2 \ra Pic^0 (A_2)$,
associating to $y$ the tensor product of the translation of $L_2$ by $y$ with the inverse of $L_1$.

The isomorphism of $G \times T$ with both $K_1 : = A_1 [2]$ and $K_2$ allows to identify
both $ H^0 (A_1, L_1) $ and $ H^0 (A_2, L_2)$ with the Stone von Neumann representation $L^2 ( T)$:
observe in fact that there is only one alternating function $(G \times T)  \ra \ZZ/2$,  independent of the chosen basis.

Therefore, there are exactly $4$ invariant divisors in the linear system $|L|$. 

Explicitly, if  $ H^0 (A_1, L_1) $ has basis $x_0, x_1$ and $ H^0 (A_2, L_2) $ has basis $y_0, y_1$,
then the invariant divisors correspond to the four eigenvectors
$$  x_0 y_0 + x_1  y_1, x_0 y_0 - x_1  y_1, x_0 y_1 + x_1  y_0, x_0 y_1 - x_1  y_0. $$

To prove irreducibility of the above family of surfaces, it suffices to show that all the four 
invariant divisors occur in the same connected family.

To this purpose, we just observe that the monodromy of the family of elliptic curves $E_{\tau} : = \CC / ( \ZZ \oplus \ZZ \tau)$
on the upperhalf plane has the effect that a transformation in $SL ( 2 , \ZZ)$ acts on the subgroup $E_{\tau}  [2]$ of points
of $2$-torsion by its image matrix in $GL ( 2 , \ZZ /2)$, and in turn the effect on the Stone von Neumann representation is the one
of twisting it by a character of $E_{\tau}  [2]$.

This concludes the proof that the moduli space is irreducible of dimension $4$, since  the moduli space of elliptic curves,
respectively  the moduli space of Abelian surfaces with a  polarization of  type $(1,2)$, are irreducible, 
of respective dimensions $1,3$.
The final assertion is a consequence of the fact that $Alb(S) = A_1 / (T_1 + \La \be_1 \Ra)$,
so that the fibres of the Albanese map are just divisors in $A_2$ of type $(1,2)$. Their self intersection 
equals $4 = 2 (g-1)$, hence $g=3$.

In order to establish that the general fibre of the Albanese map is non hyperelliptic, it suffices to prove the following lemma.

\begin{lemma}
Let $A_2$ be an Abelian surface, endowed with a divisor $L$ of type $(1,2)$, so that there is an isogeny
of degree two $f : A_2 \ra A'$ onto a principally polarised Abelian surface, with kernel generated by a point $t$ of 2-torsion, and such that $L = f^*(\Theta)$.
Then the only curves $C \in |L|$ that are hyperelliptic are contained in  the pull backs of a translate of
$\Theta$ by a point of order $2$ for a suitable such isogeny $f : A_2 \ra A'$.
In particular, the general curve $C \in |L|$ is not hyperelliptic.
\end{lemma}
\Proof
Let $C \in |L|$, and consider $D : = f_* (C) \in | 2 \Theta|$. There are two cases.

Case I):  $ C + t = C$.  Then $D = 2 B$, where
$B$ has genus $2$, so that $ C = f^* (B)$, hence, since $ 2B \equiv 2 \Theta$, $B$ is a translate of $\Theta$ by a point of order $2$.
There are exactly two such curves, and for them $ C \ra B$ is \'etale.

Case II): $ C + t \neq  C$.  Then $C \ra D $ is birational, $f^* (D) = C \cup (C + t)$. Now, $C+t$ is also  linearly equivalent to $ L$,
hence $C \cap (C+t)$ meet in the $4$ base points of the pencil $|L|$.  Hence $D$ has two double  points and geometric genus 
equal to $3$. These double points are the intersection points of $\Theta$ and a translate by a point of order $2$,
and are points of $2$-torsion.

The sections of $H^0(\hol_{A'} (2 \Theta))$ are all even 
and $ | 2 \Theta|$ is the pull-back of the space of hyperplane sections of the Kummer surface
$\sK \subset \PP^3$, the quotient $\sK = A' / \{\pm 1\}$.

Therefore the image $E'$ of each such curve $D$ lies in the pencil of planes through $2$ nodes of $\sK$.

$E'$ is a plane quartic, hence $E'$ has geometric genus $1$, and we conclude that $C$ admits an involution
$\s$ with quotient an elliptic curve $E$ (normalization of $E'$), and the double cover is branched in $4$ points.

Assume that $C$ is hyperelliptic, and denote by $h$ the hyperelliptic involution, which lies in the centre of $ Aut(C)$.
Hence we have $(\ZZ/2)^2$ acting on $C$, with quotient $\PP^1$. We easily see that there are exactly six
branch points, two being the branch points of $ C/h \ra \PP^1$, four being the branch points of $E \ra \PP^1$.
It follows that there is an \'etale quotient $ C \ra B$ , where $B$ is the genus $2$ curve, double cover of $\PP^1$
branched on the six points.
 
Now, the inclusion $ C \subset A_2$ and the degree $2$ map $ C \ra B$ induces a degree two isogeny $ A_2 \ra J(B)$,
and $C$ is the pull back of the Theta divisor of $J(B)$, thus it cannot be a general curve.

\hfill{\it QED for the lemma.}

\qed

\begin{defin}
Let us call a surface $S$ as in theorem \ref{BCF} a {\bf Sicilian  surface with $q=p_g=1$}.

Observe that the fundamental group of $S$ is isomorphic to the  fundamental group $\Ga$ of $X$, and that $\Ga$, fitting into
the exact sequence 
$$  1 \ra \Lam \ra  \Ga \ra G = \ZZ/2 \ra 1,$$ 
is generated by the union of the set $\{ g , t\},$ where
$$ g ( v_1 +  v_2 ) =  v_1 + \tau / 2  - v_2 + \la_2 / 2 $$ 
$$ t ( v_1 +  v_2 ) =  v_1 + 1 / 2  + v_2 + \la_4 / 2 $$
with the set of translations by the elements of a basis $\la_1, \la_2, \la_3, \la_4$ of $\Lam_2$.  

It is therefore a semidirect product of $\ZZ^5 = \Lam_2 \oplus \ZZ t$ with the infinite cyclic group generated by $g$:
conjugation by $g$ acts as $-1$ on $\Lam_2$, and it sends $t \mapsto t - \la_4$ (hence $ 2 t - \la_4$ is an eigenvector for the eigenvalue $1$).
\end{defin}

We shall now give a topological characterization of Sicilian surfaces with $q=p_g=1$, following the lines
of \cite{inoue}.

Observe in this respect that $X$ is a $K(\Ga, 1)$ space, so that its cohomology and homology are just 
group cohomology, respectively homology, of the group $\Ga$. 

\begin{cor}
A Sicilian surface $S$ with $q=p_g=1$ is characterized by the following properties:

\begin{enumerate}
\item
$K_S^2 = 6$
\item
$ \chi(S) = 1$
\item
$\pi_1(S) \cong \Ga$, where $\Ga$ is as above,
\item
the classifying map $f : S \ra X$, where $X$ is the Bagnera-De Franchis threefold which is a classifying space for $\Ga$,
has the property that $ f_* [S] : = Y$ satisfies $ Y^3=6$.

\end{enumerate}

In particular, any surface homotopically equivalent to a Sicilian surface is a Sicilian surface, and we get a connected component of the moduli space of surfaces
of general type which is stable under the action of the absolute Galois group \footnote{As we shall see in the last sections of the article,
the absolute Galois group permutes the connected components of the moduli space of surfaces
of general type, and this action is indeed  faithful. }.

\end{cor}

\Proof
Since $\pi_1(S) \cong \Ga$, first of all $ q(S) =1$, hence also $p_g(S) = 1$.
By the same token
there is a double \'etale cover $\hat{S} \ra S$ such that $ q(\hat{S} ) = 3$,
and the Albanese image of  $\hat{S} $, counted with multiplicity, is the inverse image $Z$ of $Y$,
therefore $Z^3 = 12$. From this, it follows that $\hat{S} \ra Z$ is birational, since the class of $Z$ is
indivisible.

We may now apply the previous theorem in order to obtain the classification.

Observe finally that the condition $(\hat{\al}_* \hat{S})^3 = 12$ is not only a topological condition,
it is also invariant under Galois autorphisms.

\qed

\section{Regularity of classifying maps and fundamental groups of projective varieties}

\subsection{Harmonic maps}

Given a continuous map $f : M \ra N$ of differentiable manifolds, we can approximate it, as already partly explained, by a differentiable one,
homotopic to the previous one.
Indeed, as we already explained, we may assume that $ N \subset \RR^n, M \subset \RR^m$ and, by a partition of unity argument, that $M$ is an open set in $\RR^h$.
Convolution approximates then f by a differentiable function $F_1$ with values in a tubular neighbourhood  $T (N) $ of $N$,
and then the implicit function theorem applied to the normal bundle provides a differentiable retraction $ r : T(N) \ra N$.
Then $ F : = r \circ F_1$ is the required approximation, and the same retraction provides a homotopy between $f$ and $F$
(the homotopy between $f$ and $F_1$ being obvious).

If however $M,N$ are algebraic varieties, and algebraic topology tells us about the existence of a continuous map $f$ as above,
we would wish for more regularity, possibly holomorphicity of the homotopic map $F$.

Now, Wirtinger's theorem (see \cite{mumErgebnisse}) characterises complex submanifolds as area minimizing ones, so the first idea is to try to
deform a differentiable mapping $f$ until it minimizes some functional.

We may take the Riemannian structure inherited form the chosen embedding, and assume that $(M, g_M), (N, g_N)$
are Riemannian manifolds.

If we assume that $M$ is compact, then one defines the {\bf Energy} $\sE (f)$ of the map as the integral:
$$  \sE (f) : = 1/2  \int_M |Df|^2 d\mu_M,$$
where  $Df$ is the derivative of the differentiable map $f$, $d \mu_M$ is the volume element on $M$,
and $|Df|$  is just its  norm as a differentiable section of a bundle
endowed with a metric: $$Df \in H^0 (M, \sC^{\infty}(TM^{\vee} \otimes f^* (TN))).$$ 

\begin{rem}
(1) in more concrete terms, the integral in local coordinates  has the form

$$  \sE (f) : = 1/2 \int_M   \Sigma_{\al,\be, i,j} [(g_N)_{\al,\be} \frac{ \partial f_{\al}}{\partial x_i}  \frac{ \partial f_{\be}}{\partial x_j} (g_M)^{-1}_{i,j}] \sqrt{det (g_{ij})}.$$

(2) Linear algebra shows that, once we identify $TM, TN$ with their dual bundles via the Riemannian metrics,
$|Df|^2 = Tr ( (Df)^{\vee} \circ Df )$, hence we integrate over $M$  the sum of the eigenvalues of the endomorphism $(Df)^{\vee} \circ Df  : TM \ra TM$.

(3) the first variation of the energy function vanishes precisely when $f$ is a {\bf harmonic} map, i.e., 
$\Delta (f) = 0$, where $\Delta (f) : = Tr (\nabla  (Df)),$ $\nabla$ being the connection on $TM^{\vee} \otimes f^* (TN)$ 
induced by the Levi-Civita connections on $M$ and $N$.

(4) the energy functional enters also in the study of geodesics and Morse theory (see \cite{milnorMT})
\end{rem}

These notions were introduced by Eells and Sampson in the seminal paper \cite{eells-sampson}, which used the {\bf heat flow}
$$  \frac{\partial f_t} {\partial t}  = \Delta (f) $$
in order to find extremals for the energy functional. These curves in the space of maps are (as explained in  \cite{eells-sampson}) the analogue of gradient lines
in Morse theory, and the energy functional decreases on these lines.

The obvious advantage of the flow method with respect to discrete convergence procedures (`direct methods of the calculus of variations')  is that here it is clear that all the maps
are homotopic to each other! \footnote{ The flow method made then its way further through the work of Hamilton (\cite{hamilton}), 
 Perelman and others (\cite{perelman1}, \cite{perelman2}, \cite{perelman3}) , leading to the solution
of the three dimensional Poincar\'e conjecture (see for example \cite{morgan-tian} for an exposition). }

The next theorem  is one of the most important results, first  obtained in \cite{eells-sampson}

\begin{theo}{\bf (Eells-Sampson)}\label{harmonic}
Let $M, N$ be compact Riemannian manifolds, and assume that the sectional curvature $\sK _N $ of $N$ is 
semi-negative ($\sK _N \leq 0$): then every continuous map $f_0 : M \ra N$ is homotopic to a harmonic map
$f : M \ra N$. Moreover the equation $\Delta (f) = 0$ implies, in case where $M,N$ are real analytic manifolds,
the real analyticity of $f$. 
\end{theo}

\begin{rem}
Observe that, in the case where $N$ has strictly negative sectional curvature, Hartmann (\cite{hartmann}) proved the unicity 
of the harmonic map in each homotopy class.
\end{rem}

Not only the condition about the curvature is necessary for the existence of a harmonic representative in each homotopy class,
but moreover it constitutes the main source of connections with the concept of classifying spaces, in view of the classical (see \cite{milnorMT},
\cite{cartanEspacesRiemann}) theorem of Cartan-Hadamard establishing a deep link between curvature and topology.

\begin{theo} {\bf (Cartan-Hadamard)}
Suppose that $N$ is a complete Riemannian manifold, with semi-negative ($\sK _N \leq 0$) sectional curvature: then 
the universal covering $\tilde{N}$ is diffeomorphic to an Euclidean space, more precisely given any two points
there is a unique geodesic joining them.
\end{theo}

\begin{rem}
The reader will notice that the hypotheses of theorem \ref{harmonic} and of Hartmann's theorem apply naturally to two projective curves
$M, N$ of the same genus $g \geq 2$,
taken with the metric  of constant curvature $-1$ provided by the uniformization theorem: then one may take
for $f_0$ a diffeomorphism, and apply the result, obtaining a unique harmonic map $f $ which Samson in \cite{sampson}
shows to be also  a diffeomorphism. The result obtained is that to the harmonic map $f $ one associates a 
quadratic differential $ \eta_f \in H^0 ( \Omega_M^{\otimes 2})= H^0 (\hol_M ( 2 K))$, and that $\eta_f$
determines the isomorphism class of $N$. This result constitutes  another approach to Teichm\"uller space $\sT _g$.

\end{rem}

Thus  in complex dimension $1$ one cannot hope for a stronger result, to have a  holomorphic map rather than 
just a harmonic one. The surprise comes from the fact that, with suitable assumptions, the hope can be realized in higher dimensions,
with a small proviso: given a complex manifold $X$, one can define the conjugate manifold $\bar{X}$ as the same differentiable
manifold, but where in the decomposition $ TX \otimes_{\RR} \CC = T^{(1,0)} \oplus T^{(0,1)} $ the roles  of
$T^{(1,0)} $ and $ T^{(0,1)}$ are interchanged (this amounts, in case where $X$ is an algebraic variety, to replacing  the defining 
polynomial equations by polynomials obtained from the previous ones by applying complex conjugation to the  coefficients,
i.e., replacing each $P(x_0, \dots, x_N)$ by $ \overline{P(\overline{x_0}, \dots, \overline{x_N})}$).

In this case the identity map, viewed as a map $\iota:  X \ra \bar{X}$ is no longer holomorphic, but antiholomorphic.
Assume now  that we have a harmonic map $ f : Y \ra X$: then also  $\iota \circ f$ shall be harmonic, but a theorem
implying that $f$ must be holomorphic then necessarily implies that there is a complex isomorphism between $X $ and $\bar{X}$.
Unfortunately, this is not  the case, as one sees, already in the case of elliptic curves;  but then one may restrict the hope to proving that $f$ is either  holomorphic or 
antiholomorphic.

A breakthrough in this direction was obtained by Siu (\cite{siuannals}) who proved several results, that we shall discuss in the next sections. 

\subsection{K\"ahler manifolds and some archetypal theorem}

The assumption that a complex manifold $X$ is a K\"ahler manifold is that there exists a Hermitian metric on the tangent bundle $T^{(1,0)}$ 
whose associated $(1,1)$ form $\xi$  is closed. In local coordinates the metric is given by
$$ h = \Sigma_{i,j} g_{i,j} d z_i d\bar{z_j }, \ {\rm with} \ d \xi = 0, \xi : = ( \Sigma_{i,j} g_{i,j} d z_i \wedge d\bar{z_j } ). $$

 Hodge theory shows that the cohomology of a compact K\"ahler manifold $X$ has a Hodge-K\"ahler decomposition, where $H^{p,q} $ is the space of harmonic forms of type $(p,q)$,
 which are in particular $d$-closed (and $d^*$-closed):
$$ H^m (X, \CC) = \oplus_{p,q \geq 0, p+q = m}  H^{p,q} , \ H^{q,p}  = \overline{H^{p,q} } , H^{p,q}  \cong H^q (X, \Omega^p_X). $$ 
We give just an elementary application of the above theorem, a characterization of complex tori (see \cite{catAV}, \cite{cat04}, \cite{cop} for other characterizations)

\begin{theo}\label{tori}
Let $X$ be a cKM, i.e., a compact K\"ahler manifold $X$, of dimension $n$. Then $X$ is a complex torus if and only if it has the same integral cohomology algebra of 
a complex torus, i.e. $H^*(X, \ZZ) \cong \wedge ^* H^1(X, \ZZ)$. Equivalently, if and only if $H^*(X, \CC) \cong \wedge ^* H^1(X, \CC)$ and $ H^{2n}(X, \ZZ) \cong \wedge ^{2n}H^1(X, \ZZ)$
\end{theo}

\Proof
Since $ H^{2n}(X, \ZZ) \cong \ZZ$, it follows that $ H^1(X, \ZZ)$ is free of rank equal to $2n$, therefore $dim_{\CC} (H^{1,0}) = n$.
We consider then, chosen a base point $x_0 \in X$, the Albanese map 
$$ a_X : X \ra \Alb(X) : = H^0(\Omega^1_X)^{\vee} / \Hom(H^1(X, \ZZ), \ZZ), \ x \mapsto \int_{x_0}^x.$$
Therefore we have a map between $X$ and the complex torus $T : = \Alb(X)$, which induces an isomorphism of first cohomology groups,
and  has degree 1, in view of the isomorphism 
$$ H^{2n}(X, \ZZ) \cong \Lam^{2n}(H^1(X, \ZZ))  \cong H^{2n}(T, \ZZ).$$

In view of the normality of $X$, it suffices to show that $a_X$ is finite. Let $Y$ be a subvariety of $X$ of dimension $m >0$ mapping to a point: then the cohomology (or homology class, in view of Poincar\'e duality)
class of $Y$ is trivial, since the cohomology algebra of $X$ and $T$ are isomorphic. But since $X$ is K\"ahler, if $\xi$ is the K\"ahler form,
$ \int_Y \xi^m >0$, a contradiction, since this integral depends only (by the closedness of $\xi$) on the homology class of $Y$.

\qed

One can conjecture that a stronger theorem holds, namely

\begin{conj}\label{Q-tori}
Let $X$ be a cKM, i.e., a compact K\"ahler manifold $X$, of dimension $n$. Then $X$ is a complex torus if and only if it has the same rational cohomology algebra of 
a complex torus, i.e. $H^*(X, \QQ) \cong \wedge ^* H^1(X, \QQ)$. Equivalently, if and only if $H^*(X, \CC) \cong \wedge ^* H^1(X, \CC)$.
\end{conj}

\begin{rem}

Observe that $H^*(X, \QQ) \cong \wedge ^* H^1(X, \QQ) \Leftrightarrow H^*(X, \CC) \cong \wedge ^* H^1(X, \CC)$
by virtue of the universal coefficients theorem.

 The same argument of the previous theorem \ref{tori} yields that, 
since $ H^{2n}(X, \ZZ) \cong \ZZ$,  $dim_{\CC} (H^{1,0}) = n$ and the Albanese map 
$$  a_X : X \ra A: = \Alb(X) : = H^0(\Omega^1_X)^{\vee} / \Hom(H^1(X, \ZZ), \ZZ), \ x \mapsto \int_{x_0}^x$$
is finite, and it suffices to show that it is unramified (\'etale), since we have an isomorphism $ H^1(X, \ZZ) \cong H^1(A, \ZZ)$.

One sets therefore $R$ to be the ramification divisor of $a_X$, and $ B = {a_X}_* (R)$ the branch divisor. 

There are two cases: 

Case i) $B$ is an ample divisor, and $A$ and $X$  are projective.

Case ii) : $B$ is non ample, thus, by a result of Ueno (\cite{keno}) there is a subtorus $A_0 \subset A$ such that $B$ is the pull-back of an ample divisor 
on $A/A_0$.

Case ii) can be reduced, via Ueno's fibration to case i), which is the crucial one.

Since we have an isomorphism of rational cohomology groups, and Poincar\'e duality holds, we get that ${a_X}_*$ and $a_X^*$ are, in each degree, isomorphisms
with rational coefficients. 

The first question is to show that the canonical divisor, which is just the ramification divisor $R$, is an ample divisor.

Then to observe that the vanishing of the topological Euler Poincar\'e characteristic $ e(X)$ and of the  Euler Poincar\'e characteristics $\chi (\Omega^i_X) = 0$ for all $i$
implies (using Riemann Roch?) that the first Chern class of $X$ is zero, thus obtaining a contradiction to the ampleness of $K_X$.

This step works at least in dimension $n=2$, since then we have the Noether formula $ e(X) + K_X^2 = 12 \chi (\hol_X)$, hence we obtain $K_X^2  = 0$.

\end{rem}

\begin{rem}
Of course the hypothesis that $X$ is K\"ahler is crucial: there are several examples, due to Blanchard, Calabi, and Sommese  \cite{bla53,bla54, bla56,cal58, somm75}, of complex manifolds which
are diffeomorphic to a complex torus but are not complex tori: indeed $K_X $ is not linearly equivalent to a trivial divisor (see \cite{nankai} for references
to the cited papers , and \cite{cop} for partial results on the question whether
a compact complex manifold with trivial canonical divisor, which is diffeomorphic to a torus, is indeed biholomorphic to a complex torus).

\end{rem}

The previous theorem \ref{tori} allows a simple  generalization which illustrates well the use of topological methods in moduli theory.

\begin{theo}
Let $X = A/G$ be a Generalized Hyperelliptic manifold of complex dimension $n$ and assume that $Y$ is a compact K\"ahler Manifold satisfying the
following properties:
\begin{enumerate}
\item
$\pi_1(Y) \cong \pi_1(X) \cong \Ga$
\item
$H^{2n}(Y, \ZZ) \cong H^{2n}(\Ga, \ZZ)$ via the natural homomorphism of cohomology groups induced by the continuous map of $Y$ to the classifying space $X = B \Ga$,
associated to the homomorphism $\pi_1(Y) \cong \pi_1(X) \cong \Ga$.
\item
$H^{i}(Y, \CC) \cong H^{i}(\Ga, \CC) \ \forall i$ via the natural homomorphism analogous to the one defined in (2).
\end{enumerate}
Then 
\begin{enumerate}
\item $Y$ is a Generalized Hyperelliptic manifold $ Y = A' / G$, 
\item
if $G$ is abelian, then $Y$ is (real) affinely equivalent to $X$, and moreover
\item
$Y$ is a complex deformation of $X$ if and only if $X$ and $Y$ have the same Hodge type. 
\end{enumerate}

\end{theo}
\Proof

Without loss of generality, we may assume that $G$ contains no translations.

Since $\Ga = \pi_1(X)$, we have an exact sequence $ 1\ra  \Lam \ra \Ga \ra G \ra 0$,
and we let $W$ be the unramified covering of $Y$ corresponding to the surjection $\pi_1(Y) \cong \Ga \ra G$.

Then the Albanese variety $A'$ of $W$ has dimension $n$, and the group $G$ acts on $A'$. 

Main claim: $G$ acts freely on $A'$. 

This follows since, otherwise, there is an element $g' \in G \setminus{1_G}$ and a lift $g$ of $g'$ in $\Ga$ such that
a positive power $g^m$ of $g$ equals the neutral element $1_{\Ga}$ (see Step I of the proof of proposition \ref{affinetype}).
 But then $g'$ would not act freely
on $A$.

Now the Albanese map 
$ W \ra A'$ induces a holomorphic map $ f : Y \ra A'/G$, and we know that  
$A'/G$ is a classifying space with 
 $\pi_1(A'/G) \cong \Ga$.  By our hypothesis $f$ induces an
isomorphism of  cohomology groups and has degree equal to 1. 
The same argument as in theorem \ref{tori}
shows that
$f$ is an isomorphism. Therefore  $Y$ 
is also a Generalized Hyperelliptic manifold of complex dimension $n$.

For the second assertion, we simply apply proposition \ref{affinetype},
stating that the exact sequence $ 1\ra  \Lam \ra \Ga \ra G \ra 0$
determines, when $G$ is abelian,  the real affine type of the action of $\Ga$.

The last assertion is a direct consequence of remark \ref{Hodgetype}.

\qed

The following is one more characterization of complex tori and Abelian varieties.

\begin{cor}\label{abelianKp1}
Let $Z$ be a projective variety. Then $Z$ is a projective $K(\pi, 1)$ with $\pi$ an abelian group $\Leftrightarrow$ $Z$ is an Abelian variety.
Similarly, if $Z$ is a cKM,  $Z$ is a  $K(\pi, 1)$ with $\pi$ abelian  $\Leftrightarrow$ $Z$ is a complex torus.
\end{cor}
\Proof
Since $\pi$ is finitely generated, we can write $\pi = \ZZ^m \oplus T$, where $T$ is a finite abelian group. The subgroup $\ZZ^m $ 
is the fundamental group of a Galois cover $ W \ra Z$ with group $T$. In this case $m$ must be even $m=2n$, by the K\"ahler assumption,
and since $H^* (W, \ZZ) = H^* (\ZZ^{2n}, \ZZ)$, we obtain that $ n = dim_{\CC} (W) = dim_{\CC} (Z)$ and we can apply theorem
\ref{tori} to infer that $W$ is a complex torus.

Now, $ Z = W / T$, where $T$ is abelian, and $T$ acts trivially by conjugation on $\pi_1(W) \cong \ZZ^{2n}$. Hence $T$ is a group of translations, and
we obtain that $Z$ is a complex torus. Moreover, since the fundamental group of a torus is torsion free, we can actually conclude that $T=0$.

\qed

\subsection{Siu's results on harmonic maps}

The result by Siu that is the simplest to state is  the following

\begin{theo}\label{siu}
(I) Assume that $f : M \ra N$ is a harmonic map between two compact K\"ahler manifolds and that the curvature tensor of $N$ is strongly
negative. Assume further that the real rank of the derivative $Df$ is at least $4$ in some point of $M$. Then $f$ is either holomorphic or antiholomorphic.

(II) In particular, if $\dim_{\CC} (N) \geq 2$ and $M$ is homotopy equivalent to $N$, then $M$ is either biholomorphic or antibiholomorphic to $N$.
\end{theo}

From  theorem \ref{siu} follows an important consequence.

\begin{cor}\label{Siu-cor}
Assume that $f : M \ra N$ is a continuous map between two compact K\"ahler manifolds and that the curvature tensor of $N$ is strongly
negative. Assume further that there is a $j \geq 4 $ such that $ H_j(f, \ZZ) \neq 0$: then $f$ is homotopic to a map $F$ which is
either holomorphic or antiholomorphic.
\end{cor}
\Proof
$f$ is homotopic to a harmonic map $F$. One needs to show that at some point  the real rank of the derivative $DF$ is at least $4$.
If it were not so, then by Sard's lemma the image $Y : =  F(M)$ would be a subanalytic compact set of  Hausdorff dimension at most 3.

\begin{lemma}  Let $Y $  be a subanalytic compact set, or just a compact set  of  Hausdorff dimension at most h. Then $ H_i(Y, \ZZ) = 0$ for $ i \geq  h+1$.
\end{lemma}
The lemma then implies  $ H_i(Y, \ZZ) = 0, \ \forall  i \geq  4$, contradicting  the existence of an $i \geq 4$ such that $ H_i(f, \ZZ) \neq   0$ .

{\em First proof of the lemma:}  $Y$ is a finite union of locally closed submanifolds of dimension $\leq h$, by the results of  \cite{b-m}.
Using the exact sequence of Borel-Moore homology (see \cite{B-M}) relating the homology of $Y = U \cup F$, where $Y$ is the union of a closed set $F$
with an open subset $U$, and induction on the number of such locally closed sets, we obtain that $ H_i(Y, \ZZ) = 0$ for $ i \geq  h + 1$.

{\em Second proof of the lemma:}  $Y$ is a compact of Hausdorff dimension at most $h$. By theorem VII 3 of \cite{H-W}, page 104,
follows that the dimension of $Y$ is at most $h$. In turn, using theorem 2' , page 362 of \cite{alexandrov}, it follows (take $A$ ibidem to be a point $P$ )
that  $ H_i(Y, \ZZ) = H_i(Y,  P , \ZZ) = 0$ for $ i \geq  h+1$.

\qed

\qed

\begin{rem}
(1)
The hypothesis that the rank of the differential should  be at least 4 at some point is needed for instance to avoid that the image of $F$  is a complex
projective curve $C$  contained inside  a projective variety $N$ with such a strongly negative K\"ahler metric. Because in this case
$f$ could be  the composition of   a holomorphic map $ g : M \ra C'$ to a curve $C'$ of the same genus as $C$, but not isomorphic to $C$,
composed with a diffeomorphism between $C'$ and $C$.  

(2) Part II of  theorem \ref{siu} follows  clearly from part I: because, as we have seen,  the homotopy equivalence between $M$ and $N$
is realised by a harmonic map $f : M \ra N$. From part I of the theorem one knows that $f$ is holomorphic or antiholomorphic
(in short, one says that $f$ is {\bf dianalytic}). W.l.o.g. let us assume that $f$ is holomorphic (replacing possibly $N$ by $\bar{N}$).
There remains only to prove that $f$ is biholomorphic. The argument is almost the same as the one used in the archetypal theorem \ref{tori}.

{\bf Step 1:} $f$ is a finite map; otherwise, since $f$ is proper, there would be a complex curve $C$ such that $f(C)$ is a point. But in a K\"ahler manifold
the homology class in dimension $2m$ of a subvariety of complex dimension $m$, here $C$, is never trivial. 
On the other hand, since $f$ is a homotopy equivalence, and $f(C)$ is a point, this class must be zero.

{\bf Step 2:} $f$ is a map of degree one since $f$ induces an isomorphism between the last non zero homology groups of $M,N$ respectively.
If $n = \dim_{\CC} N, m = \dim_{\CC} M$, then these groups are $H_{2m}(M, \ZZ)$, respectively  $H_{2n}(N, \ZZ)$; hence $n=m$
and $$H_{2m} (f) : H_{2m}(M, \ZZ) \cong \ZZ[M] \ra H_{2m}(N, \ZZ)  \cong \ZZ [N]$$
is an isomorphism.

{\bf Step 3:} $f$ is holomorphic, finite and of degree 1.Therefore there are open subsets whose complements are Zariski closed
such that $f : U \ra V \subset N$ is an isomorphism. Then the inverse of $f$ is defined on $V $, the complement of a complex analytic set,
and by the Riemann extension theorem (normality of smooth varieties) the inverse extends to $N$, showing that $f$ is biholomorphic.

(3) We have sketched the above argument since it appears over and over in the application of homotopy equivalence to proving
isomorphism of complex projective (or just K\"ahler) manifolds. Instead, the proof of part I is based on the Bochner-Nakano
formula, which was later further generalised by Siu in \cite{siuJDG}, and is too technical to be fully discussed here.

\end{rem}

Let us  try however to describe precisely the main hypothesis of strong negativity of the curvature, which is a stronger condition than the
strict negativity of the sectional curvature.

As we already mentioned, the assumption that $N$ is a K\"ahler manifold is that there exists a Hermitian metric on the tangent bundle $T^{(1,0)}$ 
whose associated $(1,1)$ form is closed. In local coordinates the metric is given by
$$ \Sigma_{i,j} g_{i,j} d z_i d\bar{z_j }, \ {\rm with} \ d ( \Sigma_{i,j} g_{i,j} d z_i \wedge d\bar{z_j } ) = 0. $$
The curvature tensor is a $(1,1)$ form with values in $(T^{(1,0)}) ^{\vee} \otimes T^{(1,0)}$, and using the Hermitian metric
to identify $(T^{(1,0)}) ^{\vee}  \cong  \overline{T^{(1,0)}} = T^{(0,1)}$, and their conjugates 
$( (T^{(0,1)}) ^{\vee} = \overline{(T^{(0,1)})  } \cong  T^{(1,0)}$ ) we write as usual the curvature tensor as a section $R$ of 
$$ (T^{(1,0)}) ^{\vee} \otimes  (T^{(0,1)}) ^{\vee}  \otimes  (T^{(1,0)}) ^{\vee} \otimes  (T^{(0,1)}) ^{\vee}.$$

Then seminegativity of the sectional curvature is equivalent to
$$ - R ( \xi \wedge \bar{\eta} - \eta \wedge  \bar{\xi},  \overline {\xi \wedge \bar{\eta} - \eta \wedge  \bar{\xi}}) \leq 0 ,$$
for all pairs of complex tangent vectors $\xi, \eta$ (here one uses the isomorphism $T^{(1,0)} \cong TN$,
and one sees that the expression depends only on the real span of the two vectors $\xi, \eta$).

Strong negativity means instead that
$$ - R ( \xi \wedge \bar{\eta} - \zeta \wedge  \bar{\theta}, \overline{ \xi \wedge \bar{\eta} - \zeta \wedge  \bar{\theta}}) <  0 ,$$
for all 4-tuples of complex tangent vectors $\xi, \eta, \zeta, \theta$.

 The geometrical meaning is the following (see \cite{5book}, page 71): the sectional curvature is a quadratic form on $\wedge^2 (TN)$,
 and as such it extends to the complexified  bundle  $\wedge^2 (TN) \otimes \CC$ as a Hermitian form. Then strong negativity
 in the sense of Siu is also called negativity of the {\bf Hermitian sectional curvature} $ R (v,w, \bar{v}, \bar{w})$ for all vectors
  $v,w \in (TN) \otimes \CC$.
  
  Then a reformulation of the result of Siu (\cite{siuannals}) and Sampson (\cite{sampsonK}) is the following:
  
  \begin{theo}
  Let $M$ be a compact K\"ahler manifold, and $N$ a Riemannian manifold with semi-negative Hermitian sectional curvature. Then every harmonic map $f : M \ra N$ is pluri-harmonic.
  \end{theo}

Now, examples of varieties $N$ with a strongly negative curvature are the balls in $\CC^n$, i.e., the BSD of type $I_{n,1}$; Siu finds out that 
 (\cite{siuannals},  see also \cite{Cal-Ves}) for the  irreducible bounded symmetric domains of type 
$$I_{p,q}, {\rm for} \  p q \geq 2, \  II_n, \forall n \geq 3, III_n, \forall n \geq 2,IV_n, \forall n \geq 3,$$
the metric is not strongly negative, but just very strongly seminegative, where very strong negativity simply means
negativity of the curvature as a Hermitian form on $T^{1,0} \otimes T^{0,1}= T^{1,0} \otimes \overline{T^{0,1}}$.

Indeed, the bulk of the calculations is to see that there is an upper bound for the nullity of the  Hermitian sectional curvature, i.e. for the rank of  the real subbundles of $TM$ where the 
Hermitian sectional curvature restricts identically to zero (in Siu's notation, then one considers always the case where very strongly seminegativity holds,
and 2-negative means strongly negative, adequately negative means that the nullity cannot be maximal.).

Hence Siu derives several results, 

\begin{theo} {\bf (Siu)}
Suppose that $M, N$ are compact K\"ahler manifolds and the curvature tensor of $N$ is negative of order $k$.
Assume that $f : M \ra N$ is a harmonic map such that, $\exists i \geq 2k$ such that   $ H_{i}(f , \ZZ) \neq 0$. Then $f$ is either holomorphic or
antiholomorphic.

\end{theo}

\begin{theo} {\bf (Siu)}
Suppose that $M$ is a compact K\"ahler manifold and $N$ is a locally symmetric manifold $\sD /\Ga$, where $\sD$ is an irreducible 
bounded symmetric domain of one of the following types: 
$$I_{p,q}, {\rm for} \  pq \geq 2, \  II_n, \forall n \geq 3, III_n, \forall n \geq 2,IV_n, \forall n \geq 3.$$

Assume that $f : M \ra N$ is a harmonic map such that,  $n : = dim_{\CC} N$,   $ H_{2n}(f , \ZZ) \neq 0$. Then $f$ is either holomorphic or
antiholomorphic.

In particular, if $M$ is a compact K\"ahler manifold homotopically equivalent to  $N$ as above, then either $ M \cong N$ or $ M \cong \bar{N}$.
\end{theo}

The most general result is (\cite{5book}, page 80, theorems 6.13-15):

\begin{theo} {\bf (Siu)}\label{cite}
Suppose that $M$ is a  compact K\"ahler manifold and  $N$ is a locally Hermitian symmetric space $\sD / \Ga = G/K $, where the irreducible decomposition
of $\sD$ contains no dimension 1 factors. 

(1) Assume that $f : M \ra N$ is a harmonic map such that, $n : = dim_{\CC} N$,   $ H_{2n}(f , \ZZ) \neq 0$. Then $f$ is  holomorphic for some invariant complex structure
on $G/K$.

(2) Assume that $M$ is homotopically equivalent to $N$: then $f$ is  biholomorphic to $G/K$ for some invariant complex structure
on $G/K$.

\end{theo}

\begin{rem}
A few words are needed to explain the formulation `for some invariant complex structure
on $G/K$'. In the case where $\sD$ is irreducible, we can just take the conjugate complex structure. But if  $\sD = \sD_1 \times \dots \times \sD_l$,
one can just take the conjugate complex structure on a subset $\sI$ of the indices $ j \in \{ 1, \dots, l\}$.
\end{rem}

\subsection{Hodge theory and existence of maps to curves}

Siu also used harmonic theory in order to construct holomorphic maps from K\"ahler manifolds to projective curves. The first result in this direction
was the  theorem of \cite{siucurves}, also obtained by Jost and Yau
(see \cite{JostYau}  and also \cite{J-Y83} for other results).

\begin{theo}\label{siucurves}{\bf (Siu)}
Assume that a compact K\"ahler manifold $X$ is such that there is a surjection  $\phi : \pi_1(X) \ra \pi_g$, where $ g\geq 2$ and, as usual, $\pi_g$ is the fundamental group of a projective curve of genus $g$.
Then there is a projective curve $C$ of genus $g' \geq g$ and a fibration $f : X \ra C$ (i.e., the fibres of $f$ are connected)  such that $\phi$ factors through $\pi_1 (f)$.
\end{theo}

In this case the homomorphism leads to a harmonic map to a curve, and one has to show that the Stein factorization yields a map to some Riemann surface
which is holomorphic for some complex structure on the target.

In this case it can be seen more directly  how the K\"ahler assumption, which boils down to K\"ahler identities, is used.

Recall that Hodge theory shows that the cohomology of a compact K\"ahler manifold $X$ has a Hodge-K\"ahler decomposition, where $H^{p,q} $ is the space of harmonic forms of type $(p,q)$:

$$ H^m (X, \CC) = \oplus_{p,q \geq 0, p+q = m}  H^{p,q} , \ H^{q,p}  = \overline{H^{p,q} } , H^{p,q}  \cong H^q (X, \Omega^p_X). $$

The Hodge-K\"ahler decomposition theorem has a long story,  revived by Griffiths in \cite{griff1, griff2}, and was proven by Picard in special cases. It entails the following consequence:

{\bf Holomorphic forms are closed, i.e., $ \eta \in H^0 (X, \Omega^p_X) \Rightarrow d \eta = 0$.}

At the turn of last century this fact was then used by Castelnuovo and de Franchis (\cite{CdF}, \cite{deFranchis}):

\begin{theo} {\bf (Castelnuovo-de Franchis)}
Assume that $X$ is a compact K\"ahler manifold, $ \eta_1 , \eta_2  \in H^0 (X, \Omega^1_X)$ are $\CC$-linearly independent, and  the wedge product $ \eta_1 \wedge \eta_2$ is d-exact.
Then $ \eta_1 \wedge \eta_2 \equiv 0 $ and there exists a fibration $f : X \ra C$ such that $\eta_1, \eta_2 \in f^* H^0 ( C, \Omega^1_C)$. In particular,
$C$ has genus $g \geq 2$.
\end{theo}

Even if the proof is well known, let us point out that the first assertion follows from the Hodge-K\"ahler decomposition, while $ \eta_1 \wedge \eta_2 \equiv 0 $ implies the existence
of a non constant rational function $\fie$ such that $\eta_2 = \fie \eta_1$. This shows that the foliation defined by the two holomorphic forms has Zariski closed leaves,
and the rest follows then rather directly taking the Stein factorization of the rational map $\fie : X \ra \PP^1$.

Now, the above result, which is holomorphic in nature, combined with the Hodge decomposition, produces results which are topological in nature
(they actually only depend on the cohomology algebra structure of $H^*(X, \CC)$). 

To explain this in the most elementary case, we start from the  following simple observation.  If two linear independent vectors in the first cohomology group $H^1 (X, \CC)$ of a K\"ahler manifold have wedge product which is trivial
in cohomology, and  we represent them as $ \eta_1 + \overline{\om_1},  \eta_2 + \overline{\om_2},$ for $ \eta_1 , \eta_2 , \om_1, \om_2  \in H^0 (X, \Omega^1_X)$,
then by the Hodge decomposition and the first assertion of the theorem of  Castelnuovo-de Franchis 
$$ (\eta_1 + \overline{\om_1}) \wedge (  \eta_2 + \overline{\om_2}) = 0 \in H^2 (X, \CC)$$
implies $$ \eta_1  \wedge   \eta_2 \equiv 0 ,  \  \om_1  \wedge   \om_2 \equiv 0.$$

We can apply Castelnuovo-de Franchis unless $\eta_1  ,   \eta_2$   are $\CC$-linearly dependent, and similarly $\om_1  ,   \om_2$.
W.l.o.g. we may assume $\eta_2  \equiv 0$ and $\om_1 \equiv 0$. But then $\eta_1 \wedge  \overline{\om_2} = 0$ implies that the Hodge norm
 $$ \int_X (\eta_1 \wedge  \overline{\om_2}) \wedge \overline{(\eta_1 \wedge  \overline{\om_2}) } \wedge \xi^{n-2}= 0 ,$$
 where $\xi$ is here the K\"ahler form. A simple trick is to observe that
  $$ 0 = \int_X (\eta_1 \wedge  \overline{\om_2} ) \wedge \overline{(\eta_1 \wedge  \overline{\om_2}) } \wedge \xi^{n-2}=   - \int_X (\eta_1 \wedge  \om_2) \wedge \overline{(\eta_1 \wedge  \om_2}) \wedge \xi^{n-2} ,$$ therefore   the same integral yields that the Hodge norm of $\eta_1 \wedge \om_2$ is zero,
 hence $\eta_1 \wedge \om_2 \equiv 0;$ the final conclusion is that  we can in any case apply Castelnuovo-de Franchis and find a map to a projective curve $C$
 of genus $ g \geq 2$.
 
 More precisely, one gets  the following theorem (\cite{albanese}):
 
 \begin{theo}{\bf (Isotropic subspace theorem)}
 On a compact K\"ahler manifold $X$ there is a bijection between isomorphism classes of fibrations $ f : X \ra C$ to a projective curve of genus $g\geq 2$,
 and real  subspaces $V \subset H^1 (X, \CC)$ (`real' means that $V$ is self conjugate, $\overline{V} = V$) which have dimension $2g$ and are of the form $ V = U \oplus \bar{U}$, where $U$ is a maximal isotropic subspace for 
 the wedge product $$H^1 (X, \CC) \times H^1 (X, \CC) \ra H^2 (X, \CC).$$
 \end{theo}
 
 It is interesting that the above result implies the  following theorem of Gromov (\cite{Gromov}),
 which  in turn obviously implies theorem \ref{siucurves} of Siu  (see \cite{bauer,GromovL2}  for related results). 
 
 \begin{theo}{\bf (Gromov's few relations theorem)}
 Let $X$ be a compact K\"ahler manifold and assume that there exists a surjection of its fundamental group 
 $$ \Ga : = \pi_1 (X) \ra G = \langle x_1, \dots, x_n | R_1(x), \dots, R_m (x) \rangle ,$$
 onto a finitely presented group that  has `few relations', more precisely  where $ n \geq m-2$.
 Then there exists a fibration $f : X \ra C$, onto a projective curve $C$ of genus $g \geq \frac{1}{2} (n-m)$.
 
 If moreover $ \pi_1 (X)  \cong G$, then the first cohomology group $H^1 (X, \CC )$ equals $f^* H^1 (C, \CC )$.
 \end{theo}
 
 {\em Proof:}
We saw in a  section  5.6  that if the fundamental group of $X$, $\Ga : = \pi_1(X)$ admits a surjection onto $G$, then the induced classifying
continuous map $\phi : X \ra   B G$  has the properties 
that its induced action on first cohomology
$$ H^1 (\phi) : H^1 (B G, \CC) \ra H^1 (X, \CC)$$ is injective 
and 
the image $W$ of $ H^1 (\phi) $ is such that each element $ w \in W$ is contained in an  isotropic subspace of rank $\geq n-m$.

Hence it follows immediately that there is an isotropic subspace of dimension $\geq n-m$ and a fibration onto a curve of genus $g \geq n-m$.

For the second assertion, let $U$ be a maximal isotropic subspace contained in $W = H^1 (X, \CC )$: 
then there exists a fibration $f : X \ra C$, onto a projective curve $C$ of genus $g \geq 2$ such that $f^* H^1 (C, \CC ) = U \oplus \bar{U}$.
We are done unless $U \oplus \bar{U} \neq H^1 (X, \CC )$.

But in this case we know that $W = H^1 (X, \CC )$ is the union of such proper subspaces $U \oplus \bar{U}$.

These however are of the form $f^* H^1 (C, \ZZ ) \otimes \CC$, hence they are a countable number; by Baire's theorem their union cannot be the whole
$\CC$-vector space $W = H^1 (X, \CC )$.

 \qed
 
 Not only one sees clearly how the K\"ahler hypothesis is used, but indeed Kato (\cite{Kato}) and Pontecorvo {\cite{pontecorvo}) showed how the
 results are indeed false without the  K\"ahler assumption, using twistor spaces of algebraic surfaces  which are $\PP^1$-bundles over a 
 projective curve $C$ of genus $g \geq 2$. A simpler example was then found by Kotschick (\cite{5book}, ex. 2.16, page 28): a primary Kodaira surface
 (an elliptic bundle over an elliptic curve, with $b_1(X) = 3, b_2(X) = 4$).
 
 We do not mention in detail generalisations of the Castelnuovo-de Franchis theory to higher dimensional targets
 (see  \cite{G-L1},  \cite{albanese}, \cite{G-L2}, \cite{simpsonTors}), since these shall not be used in the sequel.
 
 We want however to mention another result (\cite{cime}, see also \cite{cko} for a weaker result) which again, like the isotropic subspace theorem, determines explicitly the genus of the target curve
 (a result which is clearly useful for classification and moduli problems).
 
 \begin{theo}\label{orb-fibr}
 Let $X$ be a compact K\"ahler manifold, and let $f : X \ra C$ be a fibration   onto a projective  curve $C$, of genus $g$,
 and assume that there are exactly $r$ fibres which are multiple with multiplicities $ m_1, \dots m_r \geq 2$. Then $f$ induces
an orbifold
fundamental group exact sequence
$$ \pi_1 (F) \rightarrow \pi_1 (X) \rightarrow  \pi_1 (g; m_1, \dots m_r) 
\rightarrow 0,$$  where $F$ is a smooth fibre of $f$, and 
$$  \pi_1 (g; m_1, \dots m_r) : = $$
$$:=  \ \langle  \al_1, \be_1, \dots, \al_g, \be_g, \ga_1, \dots \ga_r | \ \  \Pi_1^g [\al_j, \be_j] \Pi_1^r \ga_i  = \ga_1^{m_1}= \dots = \ga_r^{m_r} = 1\rangle.
$$ 
 Conversely, let $X$ be a compact K\"ahler manifold and let $(g, m_1, \dots m_r)$
be a hyperbolic type,  i.e., assume that $  2g-2 + \Sigma_i (1 - \frac{1}{m_i} ) > 0.$ 

Then
 each  epimorphism
$\phi : \pi_1 (X) \rightarrow  \pi_1 (g;
m_1, \dots m_r)$ with finitely generated kernel is obtained from a fibration $f : X \rightarrow C$ of
type $(g; m_1, \dots m_r)$.
 \end{theo}

\subsection{Restrictions on fundamental groups of projective varieties}

A very interesting (and still largely unanswered) question posed by J.P. Serre in \cite{ser5} is:

\begin{question} {\bf (Serre)}
1) Which are the {\bf projective groups,} i.e., the groups $\pi$ which occur as fundamental groups $\pi = \pi_1(X)$ of a complex projective manifold?

2) Which are the  {\bf K\"ahler groups,} i.e., the groups $\pi$ which occur as fundamental groups $\pi = \pi_1(X)$ of a compact K\"ahler  manifold ?

\end{question}

\begin{rem}
i) Serre himself proved (see \cite{shaf}) that the answer to the first question  is positive for every finite group.
In this chapter we shall only limit ourselves to mention some examples and results to give a general idea,
especially about the use of harmonic maps, referring the reader to the book \cite{5book},
entirely  dedicated to this subject.

ii) For question 1), in view of the Lefschetz hyperplane section theorem \ref{hyperplanesection}, the class of projective groups $\pi$ 
 is exactly the class of groups which occur as fundamental groups $\pi = \pi_1(X)$ of a complex projective  smooth surface.

iii) if   $\pi $ and $\pi'$ are projective (resp. K\"ahler), the same is true  for the Cartesian product (take $X \times X'$).
\end{rem}

The first obvious restriction for a group $\Ga$ to be a K\"ahler group (only a priori the class of K\"ahler groups is a larger class than the one of projective groups)
is that their first Betti number  $b_1$ (rank of the abelianization $\Ga^{ab} = \Ga / [ \Ga, \Ga]$) is even; since if $\Ga = \pi_1(X)$, then $\Ga^{ab}  = H_1(X, \ZZ)$, and $H_1(X, \ZZ) \otimes \CC $
has even dimension by the Hodge-K\"ahler decomposition.

A more general restriction is that the fundamental group $\Ga$ of a compact differentiable manifold $M$ must be finitely presented: since by Morse theory
(see \cite{milnorMT})  $M$ is homotopically equivalent to a finite CW-complex ($M$  is obtained attaching finitely many cells of dimension $\leq \dim_{\RR} (M)$).

Conversely (see \cite{ST}, page 180), 

\begin{theo}\label{ST}
Any finitely presented group $\Ga$ is the fundamental group of a compact oriented 4-manifold. 

\end{theo}

Recall for this the

\begin{defin} {\bf (Connected sum)}
Given two differentiable manifolds $M_1, M_2$ of the same dimension $n$, take respective points $ P_i \in M_i, i=1,2$ and respective open neighbourhoods $B_i$ which are diffeomorphic to balls, and with smooth boundary $ \partial B_i \cong S^{n-1}$. Glueing together the two manifolds with boundary $ M_i \setminus B_i, \partial B_i $
we obtain a manifold $ M_1 \sharp  M_2$, which is denoted by the connected sum of  $M_1, M_2$ and whose diffeomorphism class is independent of the choices 
made in the construction. Also, if $M_1, M_2$ are oriented, the same holds for $ M_1 \sharp  M_2$, provided the diffeomorphisms $ \partial B_i \cong S^{n-1}$ are chosen to be
orientation preserving.
\end{defin}

\begin{defin} {\bf (Free product)}
If $\Ga, \Ga'$ are finitely presented groups, 
 $$ \Ga = \langle x_1, \dots, x_n | R_1(x), \dots, R_m (x) \rangle, \Ga' = \langle y_1,  \dots, y_h | R'_1(y), \dots, R'_k (y) \rangle,$$
their free product is the finitely presented group
$$ \Ga * \Ga'   = \langle x_1, \dots, x_n ,  y_1, \dots, y_h | R_1(x), \dots, R_m (x),  R'_1(y), \dots, R'_k (y) \rangle.$$

\end{defin}

{\em Idea for the proof of theorem \ref{ST}:} By the van Kampen theorem, the fundamental group of the union $ X = X_1 \cup X_2$ of two open sets with connected intersection 
$X_1 \cap X_2$ is the quotient of the free product $\pi_1 (X_1 ) * \pi_1 (X_2) $ by the relations $ i_1 (\ga) = i_2(\ga), \forall \ga \in \pi_1 (X_1 \cap X_2)$,
where $i_j :  \pi_1 (X_1 \cap X_2) \ra  \pi_1 (X_j)$ is induced by the natural inclusion $ (X_1 \cap X_2) \subset  X_j  $.

Now, assume that $\Ga$ has a finite presentation  $$ \Ga = \langle x_1, \dots, x_n | R_1(x), \dots, R_m (x) \rangle.$$
Then one considers the connected sum
of $n$ copies of $S^1 \times S^3$, which has then fundamental group equal to a free group $\F_n  =  \La x_1, \dots, x_n \Ra$.

Realizing the  relations $R_j(x)$ as  loops connecting the base point to $m$ disjoint embedded circles (i.e., diffeomorphic to $S^1$), one can perform 
the so called {\bf surgery} replacing tubular neighbourhoods of these circles, which are diffeomorphic to $ S^1 \times B^3$, and have boundary $S^1 \times S^2$,
each by the manifold with boundary $B^2 \times S^2$. These manifolds are simply connected, hence by van Kampen we introduce the relation $ R_j(x) = 1$,
and finally one obtains a $M$ with fundamental group $\cong \Ga$.

\qed

There are no more restrictions, other  than finite presentability, if one requires that $\Ga$ be the fundamental group of a compact complex manifold,
as shown by Taubes (\cite{Taubes}), or if one requires that $\Ga$ be the fundamental group of a compact symplectic 4-manifold,
as shown by Gompf (\cite{Gompf}). 
Notice however that, by the main results of Kodaira's surface classification (see \cite{bpv}), the fundamental groups of complex non projective surfaces 
form a very restricted class, in particular  either their Betti number $b_1$ is equal to $1$, or they sit in an exact sequence of the form 
$$ 1 \ra \ZZ \ra \Ga \ra \pi_g \ra 1. $$ 

In fact Taubes builds on the method of Seifert and Threlfall in order to construct a compact complex manifold $X$ with $\pi_1 (X) = \Ga$.
He takes a (differentiable) 4-manifold $M$ with $\pi_1(M) = \Ga$ and then, taking the  connected sum with  a suitable number of copies of
the complex projective plane with opposed orientation, he achieves that $M$ allows a metric with anti self dual Weyl tensor.
This condition on the metric of $M$, by a theorem of Atiyah, Hitchin and Singer (\cite{AHS})  makes  the twistor space $Tw(M)$  a complex manifold  (not just an almost complex manifold).
Note that (see \cite{LeBrunAMS}, especially  page 366) the twistor space $Tw(M)$ of an oriented Riemannian  4-manifold 
is  an $S^2$-bundle over $M$, such that the fibre over $P\in M$
equals the sphere bundle of the rank 3 vector subbundle $ \Lam^+ \subset \Lam^2 (TM)$.  Here  $ \Lam^2 (TM) = \Lam^+ \oplus \Lam^-$
is the eigenspace decomposition for the $*$-operator $ * : \Lam^2 (TM) \ra \Lam^2 (TM) $, such that $*^{2} = 1$.

The fact that  $Tw(M)$  is  an $S^2$-bundle over $M$ is of course responsible of the isomorphism $ \pi_1 (Tw(M)) \cong \pi_1(M)$.

Returning to Serre's question, a first result was given by Johnson and Rees (\cite{J-R}}, which was later extended by other authors (\cite{Abr} and \cite{Gromov}):

\begin{theo}{\bf (Johnson-Rees)}

Let $\Ga_1 , \Ga_2$ be two finitely presented groups admitting some non trivial  finite quotient. Then the free product $\Ga_1 * \Ga_2$
cannot be the fundamental group of a normal projective variety. More generally, this holds for any direct product $H \times (\Ga_1 * \Ga_2)$.
\end{theo}

The idea of proof is to use the fact  that  the   first cohomology group $H^1(X, \CC)$  carries a nondegenerate skew-symmetric form, obtained simply  
from the cup product in 1-dimensional cohomology multiplied  with the $(n-1)$-th power of the K\"ahler class.

An obvious observation is that the crucial hypothesis is  projectivity. In fact, if we take a quasi-projective variety, then the theorem does not hold true: it suffices
to take as  $X$ the complement of $d$ points in a projective curve of genus $g$, and then $\pi_1(X) \cong \F_{2g + d -1}$ is a free group.
This is not just  a curve phenomenon: already Zariski showed that if $C$ is a plane sextic with equation $ Q_2(x)^3 - G_3(X)^2=0$, where $Q_2, G_3$ are generic forms of respective degrees $2,3$, then $\pi_1 (\PP^2 \setminus C) \cong \ZZ/2 * \ZZ/3$.

One important ingredient for the theorem of Johnson and Rees is the Kurosh subgroup theorem, according to which $\Ga_1 * \Ga_2$
would have a finite index subgroup of the form $ \ZZ * K$, which would therefore also be the fundamental group of a projective variety. 
Another proof is based on the construction of a flat bundle $E$ corresponding to a homomorphism $ \ZZ * K \ra \{\pm 1\}$ with
some cohomology group $H^i(X, E)$ of odd dimension, contradicting the extension of Hodge theory to flat rank 1 bundles (see \cite{go-mi}).

Arapura, Bressler and Ramachandran answer in particular one question raised by Johnson and Rees, namely they show (\cite{Abr}):

\begin{theo}{\bf (Arapura, Bressler and Ramachandran)}
If $X$ is a compact K\"ahler manifold, then its universal cover $\tilde{X}$ has only one end. And $\pi_1(X)$ cannot be a free product amalgamated by a finite subgroup.
\end{theo}
Here, the ends of a non compact topological space $X'$ are just the limit, as the compact subset $K$ gets larger,  of the connected components of the complement set $ X' \setminus K$.

The theorem of Johnson and Rees admits the following consequence

\begin{cor}{\bf (Johnson-Rees)}\label{J-R}
Let $X_1, X_2$ be smooth projective  manifolds of dimension $n \geq 2$ with non zero first Betti number, or more generally with the property that $\pi_1(X_j)$ has
a non trivial finite quotient.
 Then $X_1 \sharp X_2$ cannot be homeomorphic to a projective manifold.
\end{cor}

It is interesting to observe that the previous result was extended by Donaldson (\cite{don}) also to the case where  $X_1, X_2$ are simply connected
smooth projective surfaces.  One could believe that the combination of the two results implies tout court that the connected sum
 $X_1 \sharp X_2$ of two projective manifolds  of dimension at least two cannot be homeomorphic to a projective manifold. This is unfortunately not clear (at least to the author),
 because of a deep result by Toledo (\cite{toledo}), giving  a negative answer to another question posed by Serre
 
 \begin{theo}{\bf (Toledo)}
 There exist projective manifolds $X$ whose fundamental group $G : = \pi_1(X) $ is not residually finite.
 
This means  that  the natural homomorphism of $G : = \pi_1(X) $ 
 into its profinite completion $ \hat{G}$  is not injective ($\widehat{\pi_1(X) }$ is also called the algebraic fundamental group and denoted by $\pi_1 (X)^{alg} $).
 
The profinite completion is defined as the following inverse  limit
 $$  \hat{G} : = \lim_{K < G, } G/K ,$$
 where  $K$ runs through the set of normal subgroups  of  finite   index.  
 
 \end{theo}
 
 Even the following question does not have yet a positive answer.
 
 \begin{question}
 Does there exist a projective manifold $X$ such that $G =  \pi_1 (X) \neq 0$, while  $ \hat{G} = : \pi_1 (X)^{alg}  =  0$?
 ( $\pi_1(X) $ would then be an infinite group).
 \end{question}
 
 If the answer were negative, then the hypothesis in corollary  \ref{J-R} would be that $X_1, X_2$ are not simply connected.
 
 Concluding this section, it is clear how the results of Siu (and Mostow, \cite{mostow}), can be used to obtain restrictions
 for  K\"ahler groups. 
 
 For instance, using the above mentioned techniques of classifying spaces and harmonic maps, plus some Lie theoretic arguments, Carlson and Toledo (\cite{C-T}),
 while giving alternative proofs of the results of Siu, prove:
 
 \begin{theo}{\bf (Carlson-Toledo)}
 If $X$ is a compact K\"ahler manifold, then its fundamental group $\pi_1(X)$ cannot be isomorphic to a discrete subgroup $ \Ga < SO(1,n) $
 such that $\sD / \Ga$ be compact, where $\sD$ is the hyperbolic space $SO(1,n) / SO(n)$.

 \end{theo}
 
 Similar in spirit is the following theorem of Simpson (\cite{simpsonHiggs}), for which we recall that a lattice $ \Ga  $
contained in a  Lie group $G$  is a discrete subgroup such that $ G / \Ga$ has finite volume, and that a reductive Lie group is of Hodge type if it has a compact Cartan subgroup. 
 
  \begin{theo}{\bf (Simpson)}
 If $X$ is a compact K\"ahler manifold, its fundamental group $\pi_1(X)$ cannot be isomorphic to a lattice $ \Ga  $
contained in a simple Lie group not of Hodge type.

 \end{theo}
 
 For more information on K\"ahler groups (fundamental group of a cKM  $X$), we refer to the book \cite{5book}, to the survey article by Campana \cite{campanaFG}, and to \cite{campanaJAG}.
 
 \subsection{K\"ahler versus projective, Kodaira's problem and Voisin's negative answer.}
 
 Many of the restrictions for a group $\Ga$  to be a projective group (fundamental group of a projective variety $X$) are indeed restrictions 
 to be a K\"ahler group. For long time, it was not clear which was the topological difference between projective smooth varieties 
 and compact K\"ahler manifolds.  Even more, there was a question by Kodaira whether any compact K\"ahler manifold $X$ would
 be a deformation (or even a direct  deformation) of 
 a projective manifold $Y$, according to the following well known definition.
 
 \begin{defin} \label{deformation}
(1) Given two  compact complex manifolds $Y, X$, $Y$  is a {\bf direct deformation } of $X$ if  there is a smooth proper and connected family $p: \sX \ra B$, where $B$ is a smooth connected complex curve, 
such that  $X$ and $Y$ isomorphic to some fibres of $p$.

(2) We say  instead that $Y$ is {\bf  deformation equivalent} to $X$ if $Y$ is equivalent to $X$ for the
 equivalence relation generated by the (symmetric)  relation of direct deformation.

\end{defin}

Indeed, Kodaira proved:

\begin{theo} {\bf (Kodaira, \cite{kod1})}
Every compact K\"ahler surface is a deformation of a projective manifold.
\end{theo}

It was for a long time suspected that Kodaira's question, although true in many important cases (tori,  surfaces),
would be in general false. The counterexamples  by Claire Voisin were based on topological ideas,
namely on the integrality of some multilinear algebra structures on the cohomology of projective varieties.

To explain this, recall that a K\"ahler form $\xi$ on a K\"ahler manifold $X$ determines the Lefschetz operator on forms of type $(p,q)$:
$$ L : \sA^{p,q} \ra  \sA^{p+1,q+1}, L (\phi) : = \xi \wedge \phi  $$ 
whose adjoint 
$$ \Lam : = L^* : \sA^{p,q} \ra  \sA^{p-1,q-1} $$

is given by $ \Lambda = w * L *$
where 
$$w := \sum_{p,q} (-1)^{p-q} \pi_{p,q},\hspace{1cm} (Weil),$$
$$c:=  \sum_{p,q} (i)^{p-q} \pi_{p,q},\hspace{1cm} (Chern-Weil ),$$
are the Weil   and Chern-Weil operators respectively.
Here $\pi_{p,q}$ is the projector onto the space of  forms of type $(p,q)$.

These operators satisfy the commutation relations:
\begin{enumerate}
\item
 $ [ L, w ] = [L, c] = [\Lambda, w] = [\Lambda, c ]  = 0$
 \item
$ \ B := [L, \Lambda] = \sum (p+q -n)\pi_{p,q},$
\item
$ [L, B] = -2L .$
\end{enumerate}

In this way the bigraded algebra of differential forms is a representation of the Lie Algebra  $sl(2, {\CC})$, which has 
the  basis:
$$ b = \left( \begin{array}{c c}
-1 & 0\\
0 &  1
\end{array} \right),$$
$$l = \left( \begin{array}{c c}
0 & 0\\
1 &  0
\end{array} \right),$$
$$\lambda = \left( \begin{array}{c c}
0 & 1\\
0 &  0
\end{array} \right),$$
and with  Lie bracket given by

 $$b =[l, \lambda],$$
$$-2l = [l,b],$$
$$2 \lambda = [\lambda, b].$$

For finite dimensional representations of $sl(2, {\CC})$ one has the following structural result.

\begin{prop}
\label{lefschetz}
Let $\rho: sl(2, {\CC}) \longrightarrow End(W)$ be a representation of
$sl(2, {\CC})$, $dim(W) < \infty$, and set $L := \rho(l)$, $\Lambda
:= \rho(\lambda)$, $B := \rho(b)$. \\
Then we have:
\begin{enumerate}
\item $W = \oplus_{\nu \in \ZZ } W_{\nu}$ is a finite direct 
sum, where $W_{\nu}$ is the eigenspace of  $B$ with eigenvalue $\nu$.
\item $L(W_{\nu}) \subset W_{\nu +2}$, $\Lambda(W_{\nu}) \subset W_{\nu -2}$.
\item Let $P = \{w \ | \ \Lambda w =0\}$, be the space of primitive vectors, then we have a direct sum decomposition
$$W = \oplus_{r \in {\NN}} L^r (P).$$
\item
Moreover, the irreducible representations of $sl(2, {\CC})$ are isomorphic to $S^m({\CC}^2)$,
where
$S^m({\CC}^2)$ is the 
$m$-th symmetric power of the natural representation of $SL(2,
{\CC})$ on ${\CC}^2$  (i.e., the space of homogeneneous polynomials of degree  $m$).
\item Let $P_{\mu} := P \cap W_{\mu}$: then $P_{\mu} =0$ for $\mu >0$ and
$$P = \oplus_{\mu \in \ZZ , \ \mu \leq 0} P_{\mu}.$$
\item $L^r : P_{-m} \longrightarrow W_{-m +2r}$ is injective for $r
\leq m$ and zero for $r >m$.
\item $L^r: W_{-r} \longrightarrow W_r$ is an  isomorphism $\forall r$, and
$$W_{\mu} = \oplus_{r \in {\NN}, \ r \geq \mu} L^r(P_{\mu - 2r}).$$
\end{enumerate}
\end{prop}

On a compact K\"ahler manifold the projectors $\pi_{p,q}$ commute with the Laplace operator, so that
the components $\psi_{p,q}$ of a harmonic form $\psi$ are again harmonic,
hence the operators $L, \Lam,\dots $ preserve the finite dimensional subspace $\sH$
of harmonic forms, and one can apply  the above proposition \ref{lefschetz}. The following theorem goes often after the name of Lefschetz decomposition,
and is essential in order to prove the third Lefschetz theorem (see \ref{hyperplanesection}).

\begin{theo}\label{Lefdec}

i) Let $(X, \xi)$ be a compact K\"ahler manifold. Then $\Delta$ commutes with the
operators $*$, $\partial$, $\overline{\partial}$, $\partial^*$,
$\overline{\partial}^*$, $L$, $\Lambda.$

ii) In particular, a k- form $\eta = \sum_{p+q = k} \eta_{p,q}$ 
is harmonic if and only if all the   $\eta_{p,q}$'s are  harmonic. 

iii) Hence there is the so-called  {\bf Hodge Decomposition} of the  de Rham cohomology:
$$ \oplus_k H^k_{DR} (X, \CC) \cong  \oplus_k{\mathcal H}^k(X, {\CC})
= \oplus_{p.q} {\mathcal H}^{p,q},$$ where  ${\mathcal H}^{p,q}$ is the space of the
 harmonic
Forms of Type $(p,q)$, 

iv) 
{\bf Hodge Symmetry} holds true, i.e., we have  ${\mathcal H}^{p,q}= \overline{ {\mathcal H}^{q,p}}$.

v) There is a canonical isomorphism ${\mathcal H}^{p,q} 
\cong H^q ( X, \Omega_X^p)$.

vi) The Lefschetz operators $L, \Lam$ make $\sH$ a representation of 
$\mathfrak s \mathfrak l (2,
{\CC})$ \footnote{hence also of $SL(2,
{\CC})$}, with eigenspaces $H_{\mu} : = \oplus_{ p+q = \mu + n} {\mathcal H}^{p,q} $,
and the space $P^{p,q}$ of primitive forms of type $(p,q)$  allows a decomposition 

$$ {\mathcal H}^{p,q}  = \oplus_{r \in {\NN}, \ r \geq \mu= p + q - n} L^r(P_{p-r, q-r}),$$
where $P^{p,q}= 0$ for $ p + q > n$.

Moreover, the Hermitian product on $P^{p,q}$ given by 
$$ \langle \eta, \eta\rangle : = i^{n-2pq} L^{(n-p-q)}(\eta \wedge \overline{\eta})$$
is positive definite.
\end{theo}

Now, a projective manifold is a compact K\"ahler manifold when endowed with the Fubini-Study metric
$$ - \frac{1}{2 \pi i} \partial \overline{\partial} log
(\sum | z_i|^2) .$$
This is the
Chern class $ c_1 ({\mathcal L})$ where $\sL$ is the line bundle $\hol_X(1)$.

Kodaira's embedding theorem characterizes projective manifolds as those compact K\"ahler manifolds which admit
a K\"ahler metric whose associated form is integral, i.e., it is the first Chern class of a positive line bundle. 

\begin{defin}
Let $M$ be a compact differentiable manifold of real even dimension $dim_{\RR}(M) = 2n$.

1)  $M$ is said to admit a real polarized Hodge structure if there exists a real closed two form $\xi$
and a Hodge decomposition of its cohomology algebra (i.e., as in iii) of theorem \ref{Lefdec})
such that the operators  $ L (\psi) : = \xi \wedge \psi$, and $\Lam :  w * L *$ satisfy  properties iv), v) and vi) of theorem \ref{Lefdec}.

2)  $M$ is said to admit an integral  polarized Hodge structure if one can take $\xi$ as above
in $H^2(M, \ZZ)$.

\end{defin}

Thus the differentiable manifolds underlying a complex K\"ahler manifolds admit a real polarized Hodge structure,
while those underlying a complex projective manifold admit an integral polarized Hodge structure.

Voisin showed that these properties can distinguish between K\"ahler and projective manifolds.
We need here perhaps to recall once more that deformation equivalent complex manifolds
are diffeomorphic.

\begin{theo}{( Voisin, \cite{voisin1})}
In any complex dimension $n \geq 4$ there exist compact K\"ahler manifolds which are not homotopically
equivalent to a complex projective manifold.

\end{theo}

The construction is not so complicated: taking a complex torus $T$ with an appropriate endomorphism $\varphi$,
and blowing up $ T \times T$ along four subtori
$$ T \times \{0\},    \{0\}  \times T, \Delta_T,  \Ga_{\varphi},$$
where $\Delta_T$ is the diagonal, and $  \Ga_{\varphi}$ is the graph of $\varphi$,
one obtains the desired K\"ahler manifold $X$.

The crucial property is the following lemma

\begin{lemma}{\bf (Voisin, \cite{voisin1})}
Assume that the characteristic polynomial $f$ of $\varphi$,  a monic integral polynomial, has all eigenvalues 
of multiplicity 1, none of them is real, and moreover the Galois group of its splittiing field is the symmetric group
$\sS_{2n}$.

Then $T$ is not an Abelian variety.
\end{lemma}

We refer to the original paper for details of the construction and proof. However, since our $X$ is obviously bimeromorphic
to the torus $ T \times T$ , which is a deformation of a projective manifold, later Voisin went on to prove a stronger result
in \cite{voisin2}.

\begin{theo}{( Voisin, \cite{voisin2})}
In any complex dimension $n \geq 10$ there exist compact K\"ahler manifolds such that no compact smooth bimeromorphic
model $X'$ of $X$ is  homotopically
equivalent to a complex projective manifold. Indeed, no such model $X'$ has the same rational cohomology ring
of a projective manifold.

\end{theo}

That the rational cohomology ring was the essential topological invariant of K\"ahler manifolds, had been found already some 30 years
before by Deligne, Griffiths, Morgan and Sullivan (in \cite{DGMS}), using Sullivan's theory of rational homotopy theory
and minimal models of graded differential algebras (see \cite{GriffithsMorgan} for a thourough introduction,
and for some of the notation that we shall introduce here without explanation).

\begin{theo}\label{DGMS}
Let $X,Y$ be  compact complex   manifolds,  which are either K\"ahler, or such that   the $d d^c$ lemma holds for them
(as for instance it happens when $X$ is bimeromorphic to a cKM).

1) Then the real homotopy type of $X$ is determined by its real cohomology ring $H^*(X, \RR)$ , and similarly the effect of
every holomorphic map $ f : X \ra Y$ on real homotopy types  is a formal consequence of the 
induced homomorphism of real cohomology rings $f^* : H^*(Y, \RR) \ra H^*(X, \RR)$.

2) If moreover $X$ is simply connected, then the graded Lie algebra of real homotopy groups $\pi_*(X) \otimes_{\ZZ}\RR$
depends only on the cohomology ring $H^*(X, \RR)$. In particular, all Massey products of any order are zero over $\QQ$.

3)  If instead $X$ is not  simply connected, then the real form of the canonical tower of nilpotent quotients of $\pi_1(X)$
(the real form is obtained by taking $\otimes_{\ZZ}\RR$ of graded pieces and extension maps) is completely
determined by  $H^1(X, \RR)$ and the cup product map $$H^1(X, \RR) \times H^1(X, \RR)\ra H^2(X, \RR).$$ 

\end{theo}

The article \cite{DGMS} contains several proofs, according to the taste of the several authors.

A basic  property of compact K\"ahler  manifolds which plays here a key role is the following lemma, also called
{\em principle of the two types} in \cite{gh}.

\begin{lemma} {\bf ($ d d^c$-Lemma)}
Let $\phi$ be a differential form such that

1) $\partial \phi = \overline{\partial} \phi = 0$ (equivalently, $ d \phi = d^c \phi = 0$, where $d^c$ is the real operator 
 $ -i ( \partial - \overline{\partial} )$),
 
 2) $ d (\phi) = 0 $ (or  $ d^c (\phi) = 0$),
 
 then
 
 3) there exists $\psi$ such that $ \phi =\partial  \overline{\partial} \psi $  (equivalently,  there exists $\psi$ such that $ \phi = d d^c (\psi)$).
\end{lemma}

One can see  the relation of the above lemma with the vanishing of Massey products,
which we now define in the simplest case.

\begin{defin}
Let $X$ be a topological space, and let $a,b,c \in H^i(X , R)$ be cohomology classes such that
$$ a \cup b = b \cup c = 0,$$
where $R$ is any ring of coefficients.

Then the triple Massey product $<a,b,c>$ is defined as follows:
take cocycle representatives $a', b', c'$ for $a,b,c,$ respectively,
and cochains $ x,y \in C^{2i-1}(X, R)$ such that
$$ dx = a' \cup b', dy = b' \cup c' . $$ 

Then  $<a,b,c>$ is the class of 
$$a' \cup y + (-1)^{i+1} x \cup c'  $$
inside $$ H^{3i -1}(X, R) / (a \cup  H^{2i -1}(X, R)+  H^{2i -1}(X, R)\cup c).$$
\end{defin}

Concerning the choice of the ring  $R$ of coefficients: if we take $X$ a cKM, and $R = \RR$, or $R= \QQ$, then these Massey products
vanish, as stated in \cite{DGMS}.

Instead things are different for torsion coefficients, as shown by Torsten Ekedahl in \cite{teke}.

\begin{theo}
i) There exist a smooth complex projective surface $X$ and $a,b,c \in H^1( X, \ZZ/l)$ such that $ a \cup b = b \cup c = 0,$
but  $<a,b,c> \neq 0.$

ii) Let $X \ra\PP^6$ be an embedding and let   $Y$ be the blow up of $\PP^6$ along $X$. 

Then there exist 
 $a'',b'',c'' \in H^3( X, \ZZ/l)$ such that $ a'' \cup b'' = b'' \cup c'' = 0,$ but  $<a'',b'',c''> \neq 0.$

\end{theo}

\begin{rem}
1) As already mentioned,  theorem \ref{Lefdec} is an essential ingredient in the proof of the third Lefschetz theorem.
Indeed, for a smooth hyperplane section $W = H \cap X$ of a projective variety $X$, the cup product with the hyperplane class
$h \in H^2(X, \ZZ)$ corresponds at the level of forms to the operator $ L$ given by wedge product with $\xi$, the first Chern form
of $\hol_X(H)$.

2) Most K\"ahler manifolds, for instance a general complex torus, do not admit any nontrivial analytic subvariety 
(nontrivial means: different from $X$ or from a point). Whereas for projective varieties $X$ of dimension $n \geq 2$
one can take  hyperplane sections successively and obtain a surface $S$ with $\pi_1(S) \cong \pi_1(X)$.
Hence the well known fact that the set of  fundamental groups of smooth projective varieties is just the set of
fundamental groups of smooth algebraic surfaces.

3) As explained in the book by Shafarevich (\cite{shaf}), the same idea was used by J.P. Serre to show that
any finite group $G$ occurs as the fundamental group of a smooth  projective surface $S$.
Serre considers the regular representation on the Cartesian product of $Y$ indexed by $G$, which we denote as usual by $Y^G$, and where for instance $Y = \PP^m$ is a simply connected projective variety
of dimension $m \geq 3$. The quotient $X: = (Y^G)/G$ is singular on the image of the big diagonal, which however has codimension 
at least $m$. Cutting $X$ with an appropriate number of general hyperplanes one obtains a surface $S$ with $\pi_1(S) \cong G$.

4) Serre's result is used in Ekedahl's theorem: Ekedahl shows the existence of a group $G$ of order $l^5$ such that there
is a principal  fibration $K(G,1)  \ra K((\ZZ/l)^3,1)$ with fibre a $ K((\ZZ/l)^2,1)$; moreover he deduces the existence of elements
$a,b,c$ as desired from the fact that there is some CW complex with non zero triple Massey product.

Then the blow up of $\PP^6$ along a surface $X$ with $\pi_1(X) \cong G$ is used to obtain an example with classes in the third cohomology
group $H^3(Y, \ZZ/l)$.

\end{rem}

 \subsection{The Shafarevich conjecture}

 One of the many reasons of the beauty of the theory of curves is given by the uniformization theorem,
stating  that any complex manifold $C$ of
dimension 1 which is not  of special type (i.e., not $\PP^1$, $\CC$,
$\CC^*$,
or an elliptic curve) has as universal covering the unit disk $\BB_1
= \{ z  \in \CC | |z| < 1 \}$, which is biholomorphic to the upper
half plane
$\HH  = \{ z  \in \CC |Im (z) > 0 \}$.

In other terms, the universal cover is either $\PP^1$, $\CC$, or $\HH $ according to the sign of 
the curvature of a metric with constant curvature (positive, zero or negative).

In higher dimensions there are simply connected projective varieties that have a positive dimensional 
moduli space: already in the case of  surfaces, e.g., smooth surfaces
in $\PP^3$,  of degree at least 3, there is an uncountable family of pairwise non isomorphic varieties.

So, if there is some analogy, it must  be a weaker one, and the first possible direction is to relate somehow
Kodaira dimension with curvature.

In the case where the canonical divisor $K_X$ is ample, there is the theorem of 
Aubin and Yau (see
\cite{Yau78}, \cite{Aubin}) showing the existence, on a projective
manifold with ample canonical divisor $K_X$, of a K\"ahler - Einstein
metric,
i.e. a  K\"ahler  metric
$\omega$ such that
$$ Ric (\omega) = - \omega . $$

This theorem is  partly the right substitute for the uniformization theorem in dimension  $n >1$,
but the K\"ahler - Einstein
condition, 
i.e. the existence of a  K\"ahler  metric
$\omega$ such that
$$ Ric (\omega) = c \omega ,   $$
forces $K_X$ ample if $c <0$, $K_X$ trivial if $c=0$, and $- K_X$ ample if $c >0$.

In the case where $K_X$ is trivial, Yau showed the existence of a K\"ahler - Einstein metric
in \cite{Yau78}, while the existence of such a metric on Fano manifolds (those with ample anticanonical divisor
$- K_X$) has only recently been settled, under a stability assumption (see \cite{donKE}, \cite{donKEfull}, \cite{tianKE}).

In the general case, one is looking for suitable metrics on singular varieties, but we shall not dwell
on this here, since we are focusing on topological aspects: we refer for instance to \cite{EGZ} (see also \cite{guenancia}).
Once more, however,  curvature influences topology:
Yau showed in fact  (\cite{yau}) that, for a projective manifold with
ample canonical divisor $K_X$,   the famous  Yau inequality is valid
$$ K_X^n \leq \frac{2(n + 1)}{n}  K_X^{n-2} c_2 (X),$$
  equality
holding if and only if the universal cover $ \widetilde{X}$ is the unit
ball $\BB_n$ in
$\CC^n$.

The second possible direction is to investigate properties of the universal covering $\tilde{X}$ 
of a projective variety. For instance, one analogue of projective curves of genus $ g \geq 2$,
whose universal cover is the unit ball $\BB_1 \subset \CC$, 
is given by the  compact complex manifolds $X$ whose universal
covering
$\tilde{X}$ is biholomorphic to a bounded domain $\Omega\subset \CC^n$.

Necessarily such a manifold $X$ is projective and has ample
canonical divisor $K_X$ (see  \cite{kodemb}, \cite{Kodaira-Morrow},
Theorem 8.4
page 144, where the Bergman metric is used, while the method of
Poincar\'e series is used in
\cite{sieg},Theorem 3  page 117 , see also
\cite{Kollar-Shaf}, Chapter 5).

Moreover, it is known, 
by a   theorem of Siegel (\cite{Siegel}, cf. also
\cite{Kobayashi},
Theorem 6.2),  that $\Omega$ must be 
holomorphically convex, indeed  a Stein domain. We recall these concepts (see for instance \cite{AG}, \cite{g-r}).

\begin{defin}\label{Stein}
A complex space $X$ is said to be a Stein complex space if and only if one of the following equivalent conditions hold:
\begin{itemize}
\item
1S) $X$ is a closed analytic subspace of $\CC^N$, for some $N$,
\item
2S) for any coherent analytic sheaf $\sF$, we have $H^j(X, \sF) = 0$, for all $j \geq 1$,
\item
3S) $X$ coincides with the maximal spectrum $Spec^m(Hol_X)$ of its algebra of global holomorphic functions $Hol_X : = H^0(X, \hol_X)$, where 
$Spec^m(Hol_X)$ is the subspace of $Spec(Hol_X)$ consisting of the maximal ideals.
\item
4S) Global holomorphic functions separate points and $X$ possesses a $\sC^{\infty}$ strictly plurisubharmonic exhaustion function,
i.e., there is a real valued proper
function $ f : X \ra \RR$ such that the Levi form $\sL (f) : = i \partial \overline{\partial} f$ is  strictly positive definite on the Zariski tangent space of each point.
\end{itemize}
A complex space is said to be holomorphically convex if and only if one of the following equivalent conditions hold:
\begin{itemize}
\item
1HC) For each compact $K$, its envelope of holomorphy  $K' : = \{ x \in X | |f(x)|  \leq max_K |f|\}$ is also compact
\item
2HC) $X$ admits a proper holomorphic map $ s : X \ra \Sigma$ to a complex Stein space (this map is called Steinification, or Cartan-Remmert reduction),   which induces an isomorphism
of function algebras $ s^* : Hol_{\Sigma} \ra Hol_X$.
\item
3HC) $X$ admits a real valued twice differentiable proper function $ f : X \ra \RR$ and a compact $K \subset X$ such that

3i)  there is a number $c \in \RR$ with $ f^{-1} (( - \infty, c]) \subset K$

3ii) the Levi form $\sL (f) : = i \partial \overline{\partial} f$ is semipositive definite and strictly positive definite outside $K$.
\end{itemize}

Finally, a complex space  $X$ is  Stein manifold iff it is  holomorphically convex and global holomorphic functions
separate points.

\end{defin}

Now, a compact complex manifold (or space) is by definition holomorphically convex, so this notion captures two extreme behaviours
of the universal covers we just described: bounded domains in $\CC^n$, and compact manifolds.
In his book \cite{shaf} Shafarevich answered the following

{\bf Shafarevich's question:} is the universal covering $\tilde{X}$ of a projective variety $X$, or of a compact K\"ahler manifold $X$,
 holomorphically convex?

\begin{rem}
The answer is negative for a compact complex manifold that is not K\"ahler, since for instance Hopf surfaces $X = (\CC^2 \setminus \{0\})/ \ZZ$
have $\CC^2 \setminus \{0\}$ as universal cover, which is not Stein since it has the same algebra of holomorphic
functions as the larger affine space $\CC^2$.

\end{rem}

Observe that all the $K(\pi, 1)$ projective manifolds we have considered so far satisfy the Shafarevich property,
since their universal covering is  either $\CC^n$ or a bounded domain in $\CC^n$ (this holds also for Kodaira fibred surfaces,
by Bers' simultaneous uniformization, \cite{Bers}).

An interesting question concerns  projective varieties whose universal cover is a bounded domain $\sD$  in $\CC^n$.
In this case the group $Aut (\sD)$ contains an infinite cocompact subgroup, so  it is natural to look first at domains which have  a
big group of automorphisms, especially at {\bf bounded homogeneous
domains},
i.e., bounded domains such that the group $\mathrm{Aut}(\Omega)$ of
biholomorphisms of $\Omega$ acts transitively.

But, as already mentioned, 
 a classical result of J. Hano (see \cite{Hano} Theorem IV, page 886, and Lemma 6.2, page 317 of  \cite{milnorcurv}) asserts that a bounded homogeneous domain that covers a compact complex manifold is a bounded symmetric  domain.
 
These naturally occur as such universal covers: Borel proved in  \cite{Bo63} that for each bounded symmetric domain
there exists  a compact free
quotient $ X
= \Omega / \Ga$, called a Hermitian locally symmetric projective manifold (these were also called  compact Clifford-Klein forms  of the
symmetric domain $\Omega$).

By taking finite ramified coverings of such locally symmetric varieties, and blow ups of points of  the latter, we easily
obtain many examples which are holomorphically convex but not Stein. But if we blow up some subvariety $Y \subset X$ of positive dimension some care has
to be taken: what is the inverse image of $Y$ in the universal cover $\tilde{X}$? It depends on the image of $\pi_1(Y) $ inside $\pi_1(X)$:
if the image is finite, then we obtain a disjoint union of compact varieties, else the connected components of the inverse image are not compact.

A more general discussion leads to the following definition.

\begin{defin}
Assume that the Shafarevich property holds for a compact K\"ahler manifold:  then the fundamental group $\pi_1(X)$ acts properly discontinuously on
$\tilde{X}$ and on its Steinification $\Sigma$, hence one has a quotient complex space of $\Sigma$ which is holomorphically dominated by $X$:
$$ Shaf(X) : = \Sigma / \pi_1(X) , shaf_X: X=  \tilde{X} / \pi_1(X) \ra Shaf(X) = \Sigma / \pi_1(X).$$

$ Shaf(X) $ is then called the Shafarevich variety of $X$, and $ shaf_X$ is called the Shafarevich morphism. They are characterized by the following universal property: 

{\bf a subvariety  $Y \subset X$ is mapped to a point in $ Shaf(X) $ if and only if, letting $Y^n$ be the normalization of $Y$,   the image of $\pi_1(Y^n) \ra \pi_1(X)$ is finite.}

\end{defin}
So, a first attempt towards the question was in the 90's to verify the existence of the Shafarevich morphism $ shaf_X$.

A weaker result was shown by Koll\'ar (\cite{Kol-Shaf-art}) and by Campana (\cite{campanaShaf}): we borrow here the version of K\'ollar, which
we find more clearly formulated.

\begin{theo}
Let $X$ be a normal projective variety: then $X$ admits a rational Shafarevich map $shaf_X: X  \dashrightarrow Shaf(X)$,
with the following properties:

1) $shaf_X$ has connected fibres

2) there are countably many subvarieties $D_i \subset X$ such that for every subvariety $Y \subset X$, $ Y \nsubseteq \cup_i D_i$,
$shaf_X (Y)$ is a point if and only if, letting $Y^n$ be the normalization of $Y$,   the image of $\pi_1(Y^n) \ra \pi_1(X)$ is finite.

Moreover, the above properties determine $Shaf(X)$ up to birational equivalence.
\end{theo}

\begin{rem}\label{contradiction}
1) A subtle point concerning the definition of Shafarevich morphism is the one of passing to the normalization of 
a subvariety $Y$. In fact, the normalization $Y^n$ has a holomorphic map to $X$ since $X$ is normal, and the fibre product
$Y^n \times_X \tilde{X}$ has a holomorphic map to $\tilde{X}$. If the image of $\pi_1(Y^n) \ra \pi_1(X)$ is finite,
this fibre product consists of a (possibly infinite) Galois  unramified covering of $Y^n$ whose components are compact.

If  $\tilde{X}$ is holomorphically convex, then these components map to points in the Steinification $\Sigma$.

2) Consider the following special situation: $X$ is a smooth rational surface, and $Y$ is a curve of geometric genus zero
with a node. Then the normalization is $\PP^1$ and the Shafarevich map should contract $Y^n$.

Assume however that $\pi_1(Y) \cong \ZZ $ injects into $\pi_1(X)$: then the 
 inverse image of $Y$ in $\tilde{X}$ consists of an infinite chain of $\PP^1$'s, which should be  contracted by the Steinification morphism, hence 
 this map would not be proper. The existence of a curve satisfying these two properties would then give a negative answer to
 the Shafarevich question.
 
3) In other words, the Shafarevich property implies that:
 
  (***)  if  $\pi_1(Y^n) \ra \pi_1(X)$ has  finite image , then also $\pi_1(Y) \ra \pi_1(X)$ has  finite image. 

\end{rem}

So, we can define a subvariety to be {\bf Shafarevich bad} if 

  (SB)    $\pi_1(Y^n) \ra \pi_1(X)$ has  finite image , but  $\pi_1(Y) \ra \pi_1(X)$ has  infinite image. 

The rough philosophy of the existence of a rational Shafarevich map is thus that Shafarevich bad subvarieties do not move in families
which fill the whole $X$, but they are contained in a countable  union of subvarieties $D_i$ (one has indeed to take into account
also finite union of subvarieties, for instance two $\PP^1$'s crossing transversally in 2 points: these are called by Koll\'ar `normal cycles'). 

But still, could they exist, giving a negative answer to the Shafarevich question?

Potential counterexampes were  proposed by Bogomolov and Katzarkov (\cite{bog-katz}), considering fibred surfaces
$f : X \ra B$. Indeed, for each projective surface, after blowing up a finite number of points, we can always obtain such a fibration.

Assume that  $f : X \ra B$ is a fibration   of hyperbolic type: then, passing to a finite unramified covering
of $X$ and passing to a finite covering of the base $B$, we may assume that the fibration does not have
 multiple fibres.
 
 Then (see theorem \ref{orb-fibr}) the surjection $\pi_1(X) \ra \pi_1(B)$ induces an infinite unramified covering
 $$ \hat{f} :  \hat{X}  \ra  \hat{B}  : = \HH = \BB_1,$$ and we shall say that this fibration is obtained by {\em opening the base}.
 
 If instead  $f : X \ra B$ is a fibration   of parabolic type, we get an infinite unramified covering
 $$ \hat{f} :  \hat{X}  \ra \hat{B}  : = \CC.$$
 
 In the elliptic case, $B= \PP^1$ and we set $\hat{X} : = X, \hat{B}  : = B$.

Now, the Shafarevich question has a positive answer if the fundamental group of $\hat{X}  $, the image of the fundamental group of a general fibre $F$
 inside $\pi_1(X)$, is finite. 
 
 Assume instead that this image is infinite and look at components of the singular fibres: we are interested to see whether we find some Shafarevich bad cycles. For this, it is important to describe the fundamental group $\pi_1(\hat{X} )$  as a quotient of $\pi_g = \pi_1 (F)$, where $F$ is a general fibre,
 and $g$ is its genus.
 
Observe preliminarily that we may assume that there are singular fibres, else either we have a Kodaira fibration, or an isotrivial fibration (see \cite{miglio})
so that, denoting by $B^*$ the complement of the critical locus of $ \hat{f} $, the fundamental group of $B^*$ is a free group.

For each singular fibre $F_t, t \in B$, a neighbourhood of $F_t$ is retractible to $F_t$, and because of this we have a continuous map
$ F \ra F_t$.

 The  {\bf group of local vanishing cycles} is defined as the kernel of  $ \pi_1(F )\ra  \pi_1(F_t)$,and denoted by $Van_t$. 
 
 Then we obtain a description of $ \pi_1(\hat{X} ) $ as 

$$ \pi_1(\hat{X} ) = \pi_g / \La  \cup_t  \{ Van_t  \} \Ra  ,$$
where $   \La  M \Ra$ denotes the subgroup normally generated by $M$.

Bogomolov and Katzarkov consider the situation where the fibre singularities are exactly nodes, and then for each node $p$
there is a vanishing cycle $van_p$ on the nearby smooth fibre. In this case we divide by the subgroup normally generated by
the vanishing cycles $van_p$. 

Their first trick is now to replace the original fibration by the pull back via a map $\varphi : B' \ra B$
which is ramified at each critical value $t$ of multiplicity $N$.

In this way they obtain a surface $X'$ which is singular, with singular points $q$ of type $A_{N-1}$, i.e., with local equation $z^N = xy$. Since these 
are  quotient singularities $\CC^2 / (\ZZ/N)$, they have  local fundamental group $\pi_{1, loc} (X'\setminus \{q\}) = \ZZ/N$.

As a second step, they construct another surface $S_N$ such that the image of $\pi_1(S_N)$ has finite index inside $\pi_1 (X' \setminus Sing(X'))$.

They show (lemma 2.7 and theorem 2.3 of loc. cit.) the following.

\begin{prop}
The Bogomolov -Katzarkov procedure constructs a new fibred surface $S_N$ such that, if we open the base of the fibration, we get a surface  $\hat{S_N}$ such that
$$ \pi_1(\hat{S_N} ) = \pi_g / \La  \cup_p  \{ van_p^N  \} \Ra  .$$

\end{prop}

 In this way they propose to construct counterexamples 
to the Shafarevich question (with non residually finite fundamental group), provide certain group theoretic questions have an affirmative answer. 

\begin{rem}
We want to point out a topological consequence of the Shafarevich property, for simplicity we consider
only the case of a projective surface $X$ with infinite fundamental group.

Assume that $\tilde{X}$ is holomorphically convex, and let  $ s : \tilde{X} \ra \Sigma$ the Cartan-Remmert reduction morphism,
where $\Sigma$ is Stein and simply connected. If $ dim (\Sigma) = 1,$ then $\Sigma $ is contractible, otherwise
we know (\cite{AF1}) that $\Sigma$ is homotopy equivalent to a CW complex of real dimension $2$.

In the first case we have a fibration with compact fibres of complex dimension $1$, in the second case we have 
a discrete set such that the fibre has complex dimension $1$, and the conclusion is that:

{\bf  if a projective surface $X$ satisfies the Shafarevich property, then $\tilde{X}$ is  homotopy equivalent to a CW complex of real dimension  at most $2$.} 

To my knowledge, even this topological corollary of the Shafarevich property is yet unproven.

\end{rem}

On the other hand, the main assertion of the Shafarevich property is that the quotient of $\tilde{X}$  by the equivalence relation which
contracts to  points the compact analytic subspaces of $\tilde{X}$ is not only a complex space, but a Stein space.
This means that one has to produce a lot of holomorphic functions on $\tilde{X}$, in order to embed the quotient 
as a closed analytic subspace in some $\CC^N$. 

There are positive results, which answer the Shafarevich question in affirmative provided the fundamental group of $X$
satisfies certain properties,  related somehow to some of the  themes treated in this article, which is the
existence of  certain homomorphisms
to fundamental groups of classifying spaces, and to the theory of harmonic maps.  
To give a very very simple idea, an easy  result in this directions is that $\tilde{X}$ is a Stein manifold
if the Albanese morphism $  \al : X \ra A : = Alb(X)$  is a finite covering.

Now, since there has been a  series of results in this direction,we refer to the introduction and bibliography in 
the most recent one  (\cite{ekpr}) for some history of the problem and  for more information  concerning  previous results,
especially the first ones due to Jost and Zuo (for instance, \cite{j-z1}, \cite{j-z2}, \cite{j-z3},\cite{G-S}) and other ones.

We want however to directly cite some  previous  results due to Katzarkov-Ramachandran (\cite{Katz-Ram}), 
respectively to Eyssidieux (\cite{eyss-inv}) (see also \cite{eyssidieuxLN} for a general introduction) which are  simpler to state. We  only want to recall that 
a reductive representation is one whose image group is reductive, i.e., all of its representations are semi-simple, or completely reducible
(any invariant subspace $W \subset V$has an invariant complement).

\begin{theo}{\bf (Katzarkov-Ramachandran)}

Let $X$ be a normal K\"ahler compact  surface, and $ X' \ra X$ an unramified cover with Galois group $\Ga$ (so $X = X' / \Ga$).
If $\Ga$ does not contain $\ZZ$ as a finite index subgroup and it admits an  almost faithful (i.e., with finite kernel)
Zariski dense representation in a connected reductive complex Lie group, then $ X'$  is holomorphically convex. 
\end{theo}

\begin{theo}{\bf (Eyssidieux)}

Let $X$ be a smooth projective variety,  consider a homomorphism
$\rho : \pi_1(X) \ra G$, and let $\tilde{X}_{\rho} : = \tilde{X}/ ker (\rho)$ be the connected unramified covering of $X$ associated to $ker (\rho)$. 

1) If $G = \GL (n, \CC)$ and $\rho $ is a reductive representation, then there is a relative Shafarevich morphism
$$ shaf_{\rho} : X \ra  Shaf_{\rho}(X)$$
to a normal projective variety, satisfying the universal property that for each normal variety $Z$ mapping to $X$, its image in $ Shaf_{\rho}(X)$
is a point if and only if the image of $\pi_1(Z) \ra \pi_1(X)$ is finite. Moreover, the connected components of the  fibres 
of $\tilde{X}_{\rho}  \ra Shaf_{\rho}(X)$ are compact.

2) Let $M$ be  a quasi compact and absolute constructible (i.e., it remains constructible after the action of each $\s \in Aut(\CC)$)  set of conjugacy classes of reductive linear representations of $\pi_1(X)$,
and let $H_M$ be the intersection of the respective kernels. Then $ \tilde{X}_M : =  \tilde{X}/ H_M$ is holomorphically convex.
\end{theo}

The main result of \cite{ekpr} consists of two parts, the first one being  the following:

\begin{theo}{\bf (Eyssidieux-Katzarkov-Pantev-Ramachandran)}

Let $X$ be a smooth projective variety,  and let $G$ be a reductive algebraic group defined over $\QQ$. Consider the Betti character scheme 
$ M : = M_B (X, B)$ such that, for a $\CC$-algebra $A$ of finite type, the set of $A$-valued points $M(A) = M( Spec(A)) $ parametrizes the representations
$$ (*) \ \  \rho : \pi_1(X) \ra G(A) .$$

1) Denote by $\tilde{H}_M^{\infty}$ the kernel of all such representations: then the  associated Galois  connected unramified covering  $\tilde{X}_M^{\infty}$
is holomorphically convex.

2) there is a non-increasing family of normal subgroups $\tilde{H}_M^{k}$ (obviously containing $\tilde{H}_M^{\infty}$),  which 
correspond to  homomorphisms
$\rho : \pi_1(X) \ra G(A)$, with the property that  $A$ is an Artinian local $\CC$-algebra, and  that the Zariski closure of the image has $k$-step unipotent radical.

Then the associated Galois  connected unramified covering  $\tilde{X}_M^{k}$ is holomorphically convex.

\end{theo}

\subsection{Strong and weak rigidity for projective $K(\pi, 1)$ manifolds}

We want here to see how the moduli problem is solved for  many projective $K(\pi, 1)$ manifolds, which were mentioned in section 3.
However, we want here not to be entangled in the technical discussion whether a variety of moduli exists, and, if so, which properties
does it have. So, we use the notion of  deformation equivalence introduced by Kodaira and Spencer.

\begin{defin} {\bf (Rigidity)}
\begin{enumerate}
\item 
 Let $X$ be a projective manifold: then we say that $X$ is  {\bf strongly rigid} if, for each  projective manifold $Y$ homotopy equivalent
to $X$, then $Y$ is isomorphic to $X$ or to the complex conjugate variety $\bar{X} $ ($X \cong Y$ or $\bar{X} \cong Y$).
\item
We say instead that $X$ is  {\bf weakly rigid} if, for each other projective manifold $Y$ homotopy equivalent
to $X$, then either $Y$ is {\bf direct deformation } of $X$, or $Y $  is {\bf direct deformation } of $\bar{X}$ (see \ref{deformation}).

\item
We say instead that $X$ is  {\bf quasi  rigid} if, for each other projective manifold $Y$ homotopy equivalent
to $X$, then either $Y$ is {\bf  deformation equivalent } to $X$, or $Y$  is {\bf deformation equivalent  } to $\bar{X}$
(recall  that deformation equivalence is the equivalence relation generated by the relation of direct deformation).
\end{enumerate}
\end{defin}

\begin{rem}
(a) This is the intuitive meaning of the above definition: assume that there is a moduli space $\mathfrak M$, whose points correspond to isomorphism classes
of certain varieties, and such that for each flat family $p: \sX \ra B$ the natural map $ B \ra \mathfrak M$ associating to $b \in B$ the isomorphism class of
the fibre $ X_b : = p^{-1} (b)$ is holomorphic. Then to be a direct deformation of each other means to belong to the same irreducible
component of the moduli space $\mathfrak M$, while being deformation equivalent means to belong to the same connected 
component of the moduli space $\mathfrak M$.

(b) The same definition can be given in the category of compact K\"ahler manifolds, or in the category of compact complex manifolds.

(c) Of course, if $p: \sX \ra B$ is a smooth proper family with base $B$ a  smooth connected complex manifold, then all the fibres $X_b$
are direct deformation equivalent.

(d)  We shall need in the sequel some more technical generalization of these notions, which amount, 
given  a product of varieties $X_1 \times X_2 \times \dots \times X_h$, to take the complex conjugate of a certain number of factors.
\end{rem}

\begin{defin}
 Let $X$ be a projective manifold: then we say that $X$ is  {\bf strongly *  rigid} if, for each other projective manifold $Y$ homotopy equivalent
to $X$, then $Y$ is isomorphic to $X$ or to another projective variety $Z$ obtained by $X$ as follows: there is a Galois  unramified covering $X'$
of $X$ ( thus $X = X' / G$)  which splits as a product $$X' = X'_1 \times X'_2 \times \dots \times X'_h.$$ 
Set then, after having fixed a subset $\sJ \subset \{1, \dots, h\}$: 
$Z' _j  : =  \overline { X'_j }$ for $j \in \sJ$, and $Z'_j  : =   X'_j $ for $j \notin \sJ$.  Consider  then the action of $G$ on $Z' = (Z'_1 \times Z'_2 \times \dots \times Z'_h)$, and let 
$Z : = Z' / G$.

Replacing the word: `isomorphic' by `direct deformation equivalent', resp. `deformation equivalent', we define the notion 
of {\bf weakly *  rigid}, resp. {\bf quasi  *  rigid}.

\end{defin}

\begin{enumerate}
\item
The first example, the one of projective curves, was essentially already discussed: we have the universal family $$p_g:  \sC_g \ra \sT_g,$$
hence, according to the above definition, {\bf projective curves are weakly rigid} (and in a stronger way, since  we do not need to allow for complex conjugation, 
as the complex conjugate $\bar{C}$ is  a direct deformation of $C$).
\item
{\bf Complex tori are weakly rigid in the category of compact K\"ahler manifolds}, since any manifold with the same integral 
cohomology of a complex torus is a complex torus
(theorem \ref{tori}).

Moreover, complex tori are 
 parametrized  by an open set $\mathcal T_n$ of the complex
Grassmann Manifold $Gr(n,2n)$, image of the open set of matrices
$\{ \Omega \in Mat(2n,n; \CC) \ | \ i^n det  (\Omega \overline
{\Omega}) > 0 \},$ as follows:
we consider
a fixed lattice $ \Ga \cong \ZZ^{2n}$, and 
to each matrix $ \Omega $ as above we associate the subspace 
$$ V =  \Omega \CC^{n},$$ so that
$ V \in Gr(n,2n)$ and $\Ga \otimes \CC \cong V \oplus \bar{V}.$

Finally, to $ \Omega $ we associate the torus $Y_V : = V / p_V (\Ga)$,
$p_V : V \oplus
\bar{V} \ra V$ being the projection onto the first addendum.

As  it was shown in \cite{cat02} (cf. also
\cite{cat04}) $\mathcal T_n$ is
a connected component of  Teichm\"uller space.

\item
Complex tori are not quasi rigid in the category of compact complex manifolds.
Sommese generalized some construction by Blanchard and Calabi, obtaining  (\cite{somm75}) that  the space of complex
structures on a  six dimensional real torus is not connected.

\item
{\bf Abelian varieties are quasi rigid}, but not weakly rigid. In fact, we saw that all the Abelian varieties of dimension $g$ admitting a polarization of type $D'$
are contained in a family over $\sH_g$. Moreover, products $E_1 \times \dots \times E_g$ of elliptic curves admit polarizations of each possible type.

\item
Locally symmetric manifolds $\sD / \Ga$ where  $\sD$ is  irreducible   of dimension $ > 1 $ are {\bf strongly rigid}
by Siu's theorem \ref{cite}.

\item
Locally symmetric manifolds $\sD / \Ga$ with the property that $\sD$ does not have  any irreducible  factor of dimension 1 are strongly * rigid
by Siu's theorem \ref{cite}.
\item
Varieties isogenous to a product (VIP)  are {\bf weakly * rigid} in all dimensions, according to  theorem \ref{isog},
that we are going to state soon (see below). They  are {\bf weakly rigid} in dimension $n=2$ only if we require the homotopy equivalence to be orientation preserving.
\item
Among the VIP's, the  {\bf strongly * rigid} are exactly  the quotients $ X = (C_1 \times C_2 \times \dots \times C_n)/ G$ where $G$ not only acts freely, but
satisfies the following property.
Denote by   $G^0 \subset G$ (see \cite{isogenous} for more details)  the subgroup which does not permute the factors, and observe that 
$$ G^0 \subset  \Aut(C_1) \times  \Aut(C_2 )\times \dots \times  \Aut(C_n) .$$ Then $G^0$ operates on each curve $C_i$, and the required condition is that this action is rigid,
more precisely we want : $ C_i / G^0 \cong \PP^1$ and the quotient map $ p_i : C_i  \ra C_i / G^0 \cong \PP^1$ is branched in three points (we shall also say that we have then a {\bf triangle curve}).

Again, in dimension $n=2$ we get some strongly rigid surfaces, which have been called (ibidem) {\bf Beauville surfaces}. 
We propose therefore to call {\bf Beauville varieties}  the strongly * rigid VIP's. 
\item
 Hyperelliptic surfaces are weakly rigid, as classically known, see below for a more general result.
\item
 It is unclear, as we shall  see, whether Kodaira fibred surfaces are weakly *  rigid, however any surface homotopically equivalent to a Kodaira fibred surface
 is also a Kodaira fibred surface.  Indeed a stronger result, with a similar method to the one of \cite{isogenous} was shown by Kotschick
(see \cite{kotschick} and theorem  \ref{kot} below), after some partial results by Jost and Yau. 
\item
Generalized hyperelliptic varieties $X$ are  a class of varieties for which a weaker property holds.
Namely,  if $Y$ is a compact K\"ahler manifold which is homotopy equivalent to $X$ (one can relax this assumption a bit, obtaining stronger results),
then $Y$ is the quotient of $\CC^n$ by an affine action of $\Ga : = \pi_1 (X) \cong \pi_1 (Y) $ which,  using proposition \ref{affine}, can be shown 
to have  the same real affine type
as the action yielding $X$ as a quotient. But the Hodge type could be different, except in special cases where weak rigidity holds.
\end{enumerate}

The following (see \cite{isogenous} and \cite{cat03}) is the main result concerning surfaces isogenous to a product, and is a stronger result than weak *  rigidity.
\begin{theo}\label{isogenous}

a) A
     projective smooth surface $S$ is isogenous   to a product of two curves of respective genera $g_1, g_2 \geq 2$ ,  if and only if
the following two conditions are satisfied:

1) there is an exact sequence
$$
1 \rightarrow \pi_{g_1} \times \pi_{g_2} \rightarrow \pi = \pi_1(S)
\rightarrow G \rightarrow 1,
$$
where $G$ is a finite group and where $\pi_{g_i}$ denotes the fundamental
group of a projective curve of genus $g_i \geq 2$;

2) $e(S) (= c_2(S)) = \frac{4}{|G|} (g_1-1)(g_2-1)$.

\noindent
b) Write $ S = (C_1 \times C_2) / G$. Any surface $X$ with the
same topological Euler number and the same fundamental group as $S$
is diffeomorphic to $S$ and is
also isogenous to a product. There is a smooth proper family with connected smooth base manifold $T$, $ p : \sX \ra T$
having two fibres respectively isomorphic to $X$, and $Y$, where $Y$ is one of the 4 surfaces $ S = (C_1 \times C_2) / G$,
$ S_{+-} :  = (\overline{C_1} \times C_2) / G$, $ \bar{S }= (\overline{C_1}  \times \overline{C_2} ) / G$,
$ S_{-+} :  = (C_1 \times \overline{C_2} ) / G = \overline{ S_{+-}   }$.

c) The corresponding subset of the moduli space of surfaces of general type
$\mathfrak{M}^{top}_S = \mathfrak{M}^{diff}_S$, corresponding to
surfaces orientedly homeomorphic,
resp. orientedly diffeomorphic to $S$, is either
irreducible and connected or it contains
two connected components which are exchanged by complex
conjugation.

In particular, if $S'$ is orientedly diffeomorphic to $S$, then $S'$ is
deformation equivalent to $S$ or to $\bar{S}$.

\end{theo}

{\em Idea of the proof}
 $\Ga : = \pi_1 (S) $  admits a subgroup $\Ga'$ of index $d$ such that 
$ \Ga' \cong  (\pi_{g_1} \times \pi_{g_2})$.  Let $S'$ be the associated unramified covering of $S$. Then application of the isotropic subspace theorem or of theorem \ref{orb-fibr} yields 
a pair of holomorphic maps $f_j : S' \ra C_j$, hence a holomorphic map $$ F : =  f_1 \times f_2 : S' \ra C'_1 \times C'_2.$$

 Then the fibres of $f_1$ have genus $h_2 \geq g_2$,
hence by the Zeuthen Segre formula (\ref{ZS}) $ e(S') \geq  4 (h_2 -1 ) (g_1 - 1)$, equality holding if and only if all the fibres are smooth.

But   $ e(S') = 4 (g_1-1)(g_2-1) \leq 4 (h_2 -1 ) (g_1 - 1)$, so $h_2= g_2$, all the fibres are smooth hence isomorphic to $C'_2$; therefore $F$ is an isomorphism.

The second assertion follows from the refined Nielsen realization theorem, theorem  \ref{refinedNR}, considering  the action of the index two subgroup $G^0$ of $G$.
Notice that if a group, here $G^0$, acts on a curve $C$, then it also acts on the complex conjugate curve, thus $ G^0 \subset \Aut (C) \Rightarrow G^0 \subset \Aut (\bar{C}) $.
But since complex conjugation is orientation reversing, these two actions are conjugate by $ \Out (\pi_g)$, not necessarily by $ \sM ap_g = \Out^+ (\pi_g)$.

Recall now that the exact sequence 
$$  1 \ra  \pi_{g_1} \times \pi_{g_2}  \ra \Ga^0 \ra G^0 \ra 1$$
yields two injective homomorphisms $\rho_j :  G^0  \ra Out (\pi_{g_j} )$ for $j=1,2$. 
Now, we choose the isomorphism $ \pi_1 (C_j) \cong \pi_{g_j} $ in such a way that it is orientation preserving. The isomorphism $ \pi_1(X) \cong \pi_1(S)$,
and the realization  $ X = (C'_1 \times C'_2) / G$ then gives an isomorphism $\phi_j :  \pi_1 (C'_j) \cong \pi_{g_j} $. If this isomorphism is orientation preserving,
then the two actions of $G^0$ on $C_j, C'_j$ are deformation of each other by the refined Nielsen realization. If instead this isomorphism is orientation reversing,
then we replace $C_j$ by $ \overline {C_j}$ and now the two actions are conjugate by $\sM ap_g$, so we can apply the refined Nielsen realisation theorem
to  $ \overline {C_j}, C'_j$.
Finally, in the case where the homeomorphism of $X$ with $S$ is orientation preserving, then either both $\phi_1, \phi_2$ are orientation preserving, or they are both orientation reversing.
Then $X$ is a deformation either of $S$ or of $\bar{S}$.

\qed

The first part of the following theorem is instead a small improvement of theorem 7.1 of \cite{isogenous}, the second part is theorem 7.7 ibidem, and relies on the results of Mok
(\cite{Mok1}. \cite{Mok2}).

\begin{theo}\label{isog}

(1) Let  $Y$ be a projective variety of dimension $n$ with $K_Y$ ample and that  $\Ga : = \pi_1 (Y) $  admits a subgroup $\Ga'$ of index $d$ such that 
$$ \Ga' \cong  \pi_{g_1} \times \pi_{g_2} \times \dots \times \pi_{g_n},  \ g_i \geq 2  \ \forall i,$$ 
and moreover $H^{2n} (\Ga, \ZZ) \ra H^{2n} ( Y, \ZZ) $ is an isomorphism.

Then $Y$ is a variety isogenous to a product.

(2) Let $X$ be projective variety with universal covering the polydisk $\sH^n$, and let $Y$ be a projective variety of dimension $n$ with $K_Y$ ample and
such  that  $ \pi_1 (Y) \cong  \Ga : = \pi_1 (X)$, and assume that $H^{2n} (\Ga, \ZZ) \ra H^{2n} ( Y, \ZZ) $ is an isomorphism. Then also $Y$ has   $\sH^n$ as universal covering
(hence we have a representation of $Y$ as a quotient $Y = \sH^n / \Ga$). 

(3) In particular $X$ is weakly * rigid, and strongly * rigid if the action of $\Ga$ is irreducible (this means that there is no finite index subgroup $\Ga' < \Ga$ and an isomorphism 
$\sH^n \cong \sH^m \times \sH^{n-m}$ such that $$ \Ga'  =  \Ga' _1 \times \Ga' _2, \ \   \Ga' _1 \subset  \Aut ( \sH^m ), \  \Ga' _2 \subset \Aut ( \sH^{n-m}).$$

\end{theo}

{\em Proof of (1):}
As in the previous theorem we take the associated unramified covering $Y'$ of $Y$ associated to $\Ga'$, and with the same argument we produce a holomorphic
map  $$ F : Y' \ra  Z : = C_1 \times C_2 \times \dots \times C_n.$$

We claim that $F$ is an isomorphism. Indeed, by the assumption $H^{2n} (\Ga, \ZZ) \cong H^{2n} ( Y, \ZZ) $ follows that $H^{2n} (\Ga', \ZZ) \cong H^{2n} ( Y', \ZZ) $,
hence $F$ has degree 1 and is a birational morphism.  Let $R$ be the ramification divisor of $F$, so that $ K_{Y'} = F^* ( K_Z) + R$.
Now, since $F$ is birational we have $ H^0( m K_{Y'} ) \cong H^0 (m K_Z)$, for all $ m \geq 0$.
This implies for the corresponding   linear systems of divisors: $$   | m K_{Y'} | =  F^* (|m K_Z| ) + m R.$$  For $ m >> 0$ this shows that $K_{Y'} $ cannot be ample,
contradicting the ampleness of $K_Y$.

\qed 

\begin{rem}
Part (2) of the above theorem is shown in \cite{isogenous}, while part (3) follows by the same argument given in the proof of theorem \ref{isogenous} in the case where $X$ is isogenous to a product.
We omit the proof of (3) in general.
\end{rem}

\begin{theo}\label{kot}{\bf (Kotschick)}
Assume that $S$ is a compact K\"ahler surface, and that

(i)  its fundamental group sits into an exact  sequence, where $ g,b \geq 2$:

$$ 1 \ra \pi_g \ra \pi_1(S) \ra \pi_b \ra 1   $$

(ii) $ e(S) = 4 (b-1) ( g-1)$.

Then $S$  has a smooth holomorphic fibration $ f : S \ra B$, where $B$ is a projective curve of genus $b$, and where all the fibres are smooth projective curves
of genus $g$.  $f$ is a Kodaira fibration if and only if the associated homomorphism $ \rho : \pi_b \ra \sM ap_g$ has image of infinite order, else it
is a surface isogenous to a product of unmixed type and where the action on the first curve is free.
\end{theo}

\Proof
 By theorem \ref{orb-fibr} the above exact sequence yields a fibration $ f : S \ra B$ such that there is a surjection $\pi_1(F) \ra \pi_g$, where $F$ is a smooth fibre.
 Hence, denoting by $h$ the genus of $F$, we conclude that $ h \geq g$, and again we can use the Zeuthen-Segre formula to conclude that $h = g$ and that all fibres
 are smooth. So $F$ is a smooth fibration. Let $C' \ra C$ be the unramified covering associated to $ \ker (\rho)$: then the pull back family $ S' \ra C'$ has
 a topological trivialization, hence is a pull back of the universal family $\sC_g \ra \sT_g$ for an appropriate holomorphic map $\fie : C' \ra \sT_g$.
 
 If $ \ker (\rho)$ has finite index, then $C'$ is compact and, since Teichm\"uller space is a bounded domain in $\CC^{3g-3}$, the holomorphic map is constant.
 Therefore $S'$ is a product $ C' \times C_2$ and, denoting by $G : = \im (\rho)$, $S = (C' \times C_2)$, and we get exactly the surfaces 
 isogenous to a product such that  the action of $G$ on the curve $  C'$ is free. 
 
 If instead $G : = \im (\rho)$ is infinite, then the map of $C'$ into Teichm\"uller space is not constant, since the isotropy group of a point corresponding to a curve $F$ is, 
 as we saw, equal to the group
 of automorphisms of $F$ (which is finite). Therefore, in this case, we have a Kodaira fibration.
 
 \qed
 
 \begin{rem}
 Jost and Yau (\cite{J-Y83}) proved a weaker result,  that a deformation of the original examples by Kodaira of Kodaira fibrations are again Kodaira  fibrations. 
 In our joint paper with Rollenske \cite{CR} we gave examples of Kodaira fibred surfaces which are rigid, and we also described the irreducible connected components
 of the moduli space corresponding to the subclass of those surfaces which admit two different Kodaira  fibrations. 
 \end{rem}
 
 \begin{question}
Are Kodaira fibred surfaces weakly * rigid?  
 \end{question}
 
 Let us explain why one can ask this question: a Kodaira fibration $ f : S \ra B$ yields a holomorphic map $\fie : B \ra \mathfrak M_g$.
 The unramified covering $B'$ associated to $\ker (\rho) $, $\rho : \pi_1 (B) \ra \sM ap_g$ admits then a 
 holomorphic map $\fie' : B' \ra \sT_g$ which lifts $\fie$.
 If $\mathfrak M_g$ were a classifying space the homomorphism $\rho$ would determine the homotopy class of $\fie$.
 A generalization of the cited theorems of Eells and Samson and of Hartmann (\cite{eells-sampson}, \cite{hartmann}, see theorem  \ref{harmonic})would then show that in each class there is a unique
 harmonic representative, thus proving that, fixed the complex structure on $B$, there is a unique holomorphic representative, if any.
 This might help to describe the locus thus obtained in $\sT_b$ (the set of maps  $\fie' : B' \ra \sT_g$ that are $\rho$-equivariant).
 
 But indeed, our knowledge of Kodaira fibrations is still scanty, for instance the following questions are  still open.
 
  \begin{question}
1) Given an exact sequence  $ 1 \ra  \pi_g \ra \pi \ra \pi_b \ra 1$ where $g \geq 3 $, $b \geq 2$, does there exist a Kodaira fibred surface
$S$ with fundamental group $\pi$?

2) (Le Brun's question, see \cite{CR}) : are there Kodaira fibred surfaces with slope $ c_1^2 / c_2  > 8/3$ ?

3) Are there surfaces admitting three different Kodaira fibrations?
 \end{question}
 
 A brief comment on question 1) above. The given exact sequence determines, via conjugation of lifts, a
 homomorphism $\rho : \pi_b \ra \Out^+ (\pi_g) = Map_g$. If the image $G$ of $\rho$ is finite, then
 $\pi$ is the fundamental group of a surface isogenous to a product, and the answer is positive, as we previously explained.
 
 Case II) : $\rho$ is injective, hence (as  $\pi_b $ has no nontrivial elements of finite order) $G = Im (\rho)$ acts freely on Teichm\"uller space $\sT_g$, and $M : = \sT_g/G$ is a classifying space 
 for $\pi_b$. There is a differentiable map of a curve of genus $b$ into $M$, and the question is whether one can deform it to obtain
 a harmonic or holomorphic map.
 
 Case III) :  $\rho$ is not  injective. This case is the most frequent one, since, given a Kodaira fibred surface $ f : S \ra B$,
 for each surjection $ F : B' \ra B$, the pull-back $ S' : = S \times_B B'$ is again a Kodaira fibration and $\rho '$ factors as the composition of
 the surjection $\pi_1(B') \ra \pi_1(B)$ with $\rho : \pi_1(B) \ra Map_g$.
 
 It is interesting to observe, via an elementary calculation, that the slope $ K^2_S / e(S) > 2$ can only decrease for a branched covering
 $B' \ra B$, and that it tends to $2$ if moreover the weight of the branch divisor $ \beta: = \sum (1 - 1/m_i)$ tends to $\infty$.
 
 Similar phenomena ($ker (\rho) \neq 0$, and the slope decreases) occur if we take   general hypersurface sections of large
 degree   $ B \subset \mathfrak M_g^*$  in the Satake compactification of the moduli space
 $\mathfrak M_g$  which    intersect neither   the boundary
 $ \partial \  \mathfrak M_g^* = \mathfrak M_g^* \setminus  \mathfrak M_g$ nor the locus of curves with automorphisms (both have codimension
 at least 2).
 
\subsection{Can we work with locally symmetric varieties?}
 
As we shall see in the sequel, Abelian varieties and curves are simpler objects to work with, because explicit constructions may be performed
via bilinear algebra in the former case, and via  ramified coverings in the latter. For locally symmetric varieties, the constructions are more difficult,
hence rigidity results are not enough to decide in concrete cases whether a  given variety is locally symmetric. 

 Recently, as a consequence of the theorem of Aubin and Yau on the existence of a K\"ahler-Einstein metric on projective varieties with ample
 canonical divisor $K_X$, there have been several simple explicit criteria which guarantee that a  projective variety with ample
 canonical divisor $K_X$ is locally symmetric (see \cite{yau2}, \cite{V-Z}, \cite{C-DS1}, \cite{C-DS2}). 
 
A difficult question  which  remains open is the one   of studying actions of a finite group $G$  on them, and 
especially of describing the $G$-invariant divisors.

\section{Inoue type varieties}

While a couple of hundreds   examples are known today of families of minimal surfaces of general type with geometric genus 
$p_g (S) : = \dim H^0(\hol_S(K_S)) = 0$  (observe that for these surfaces $ 1 \leq K_S^2 \leq 9$),
for  the value $K_S^2= 7$ there are  only two examples known (cf. \cite{inoue}
and \cite{yifan}), and for a long time
only one   family of
such surfaces was  known,  the one constructed by Masahisa Inoue
(cf. \cite{inoue}).

The attempt to prove that  Inoue surfaces form a connected component of the moduli space of surfaces of general type
proved to be successful (\cite{bc-inoue}), and was based on a weak rigidity result: the topological type of an Inoue surface determines an irreducible
connected component of the moduli space (a phenomenon similar to the
one which was  observed   in several papers, as \cite{keumnaie},
\cite{burniat1} \cite{coughlanchan}, \cite{bc-CMP}).

  The starting point was the calculation of  the fundamental group of an Inoue
surface with $p_g = 0$ and $K_S^2 =7$: it  sits
in an extension ($\pi_g$ being as usual  the fundamental group of a projective
curve of genus $g$):
$$ 1 \rightarrow \pi_{5} \times \mathbb{Z}^4 \rightarrow \pi_1(S)
\rightarrow (\mathbb{Z}/2\mathbb{Z})^5 \rightarrow 1.
$$

This extension is given geometrically, i.e., stems from the observation (\cite{bc-inoue})  that
an Inoue surface $S$ admits an unramified
$(\ZZ / 2\ZZ)^5$ - Galois covering $\hat{S}$ which is an
ample divisor in $E_1 \times E_2 \times D$, where
$E_1, E_2$ are elliptic curves and $D$ is a projective  curve of genus
$5$; while Inoue described $\hat{S}$ as  a complete intersection of two non ample divisors
in the product   $E_1 \times E_2 \times E_3 \times E_4 $ of four elliptic curves.

It turned out that the ideas needed to treat this special family of Inoue surfaces could be put in a rather general framework,
valid in all dimensions,  setting then  the
stage for  the investigation and search for a new class of varieties, 
which we proposed to call
Inoue-type varieties. 

\begin{defin} {\bf (\cite{bc-inoue})}
Define a complex projective manifold $X$ to be an {\bf Inoue-type manifold} if
\begin{enumerate}
\item
$ dim (X) \geq 2$;
\item
there is a finite group $G$ and an   unramified $G$-covering $
\hat{X} \ra X$,
(hence $ X = \hat{X} / G$) such  that
\item
$ \hat{X}$ is an ample divisor inside a $K(\Ga, 1)$-projective manifold
$Z$, (hence by the theorems of  Lefschetz, see theorem \ref{hyperplanesection}, $\pi_1 ( \hat{X}) \cong \pi_1 (Z)
\cong \Ga$) and moreover
\item
the action of $G$ on $ \hat{X}$   yields a faithful action on $\pi_1
( \hat{X}) \cong \Ga$:
in other words the exact sequence
$$ 1 \ra \Ga  \cong \pi_1 ( \hat{X}) \ra \pi_1 ( X) \ra G \ra 1$$
gives an injection $ G \ra \Out (\Ga)$, defined by conjugation by lifts of elements of $G$;
\item
the action of $G$ on $ \hat{X}$   is induced by an action of $G$ on $Z$.
\end{enumerate}

Similarly one defines the notion of an  {\bf Inoue-type variety}, by
requiring the same properties
for a variety $X$ with canonical singularities.
\end{defin}

\begin{ex}
Indeed, the examples of Inoue, which also allow him to find another description of  the surfaces constructed by Burniat (\cite{burniat}),
are based on products of curves, and there  the group $G$  is  a $\ZZ/2$-vector space.

In fact things can be done more algebraically (as in \cite{bc-CMP}). If we take a  $(\ZZ/2)^3$-covering of $\PP^1$ branched on 4 points,
it has equations in $\PP^3$:
$$ x_1^2 +   x_2^2 +  x_3^2 = 0, \  x_0^2  - a_2   x_2^2  - a_3  x_3^2 = 0, $$ 
the group $G \cong (\ZZ/2)^3$ is the group of transformations $ x_i \mapsto \pm x_i$,  the quotient is 
$$\PP^1 = \{ y \in \PP^3 |  y_1+   y_2 +  y_3 = 0, \  y_0  - a_2   y_2  - a_3  y_3 = 0\}, $$
and the branch points are the 4 points $ y_i = 0, i = 0,1,2,3$.

The group $G$, if we see the elliptic curve as $ \CC /(\ZZ + \ZZ \tau)$, is the group of affine  transformations
$$ [z] \mapsto  \pm [z]  + \frac{1}{2} ( a + b \tau), \ a,b \in \{ 0,1\},$$
with linear coefficient $\pm 1$ and translation vector a point of 2-torsion.

Algebra becomes easier than theta functions if one takes several  square (or cubic) roots: for instance 3 square roots define
a curve of genus 5  in $\PP^4$:
$$ x_1^2 +   x_2^2 +  x_3^2 = 0, \  x_0^2  - a_2   x_2^2  - a_3  x_3^2 = 0,  x_4^2  - b_2   x_2^2  - b_3  x_3^2 = 0.$$ 
One can take more generally a $(\ZZ/2)^n$ covering of $\PP^1$ branched on $n+1$  points, which is a curve of genus 
$ g = 2^{n-2} (n-3) + 1$, or curves corresponding to similar Kummer coverings ($(\ZZ/m)^n$ coverings of $\PP^1$ branched on $n+1$  points).

\end{ex}

The above   definition of Inoue type manifold, although imposing a strong restriction on $X$, is too general,
and in order to get 
weak rigidity type results  it is convenient to impose restrictions on the fundamental group $\Ga$ of $Z$,
for instance the most interesting case
is the one where $Z$ is a product of Abelian varieties, curves, and other locally symmetric varieties with
ample canonical bundle.

\begin{defin}\label{SIT}
   We shall say  that an  Inoue-type manifold $X$  is 

\begin{enumerate}
\item
  a {\bf special Inoue type manifold}
if  moreover
$$ Z = (A_1 \times \dots \times A_r) \times  (C_1 \times \dots \times
C_h) \times (M_1 \times \dots \times
M_s)$$ where each
$A_i$ is an Abelian variety,  each $C_j$ is a curve of genus $ g_j \geq 2$,
and $M_i$ is  a
compact quotient of an irreducible bounded symmetric domain of 
dimension at least 2 by a
torsion free subgroup;
\item
  a {\bf 
classical Inoue type manifold}
if  moreover

$ Z = (A_1 \times \dots \times A_r) \times  (C_1 \times \dots \times
C_h) $ where each
$A_i$ is an Abelian variety,  each $C_j$ is a curve of genus $ g_j \geq 2$;
\item
a special Inoue type manifold  is said to be
{\bf diagonal }
if  moreover:
\begin{itemize}
\item[(I)] the action of $G$ on $ \hat{X}$   is induced by a  diagonal 
action on $Z$, i.e.,
\begin{equation}
  G \subset    \prod_{i=1}^r \Aut(A_i)
\times
\prod_{j=1}^h\Aut(C_j) \times \prod_{l=1}^s\Aut(M_l)
\end{equation}
and furthermore:
\item[(II)] 
the  faithful action on $\pi_1 ( \hat{X}) \cong \Ga$,
induced by conjugation by lifts of elements of $G$ in  the exact sequence
\begin{equation}\label{prodFG}
 1 \ra \Ga= \Pi_{i=1}^r (\Lambda_i) \times   \Pi_{j=1}^h (\pi_{g_j}) 
\times \Pi_{l=1}^s
(\pi_1 (M_l))
\ra \pi_1 ( X) \ra G \ra 1
\end{equation}
(observe that each factor $\Lambda_i$, resp. $\pi_{g_j}, \pi_1 (M_l)$ is a
normal subgroup), satisfies  the {\bf Schur property}
$$(SP) \ \ \ \ \Hom (V_i, V_j)^{G} = 0,\forall i \neq j.   $$
Here $V_j : = \Lambda_j\otimes \QQ$  and, in order that the Schur property holds,  it  suffices for instance to verify  that 
for each $\Lambda_i$ there is a subgroup $H_i$ of $G$ for which $ \Hom (V_i, V_j)^{H_i} = 0, \forall j \neq i  $.
\end{itemize}

\end{enumerate}

\end{defin}

The Schur property (SP) plays an important role in order to show that an
Abelian variety  with such a $G$-action on its fundamental group
must split as a product.

 Before stating the main general result of \cite{bc-inoue} we need the following definition, which was already used in \ref{tori} for the characterization of complex tori 
 among K\"ahler manifolds.

\begin{defin}
Let $Y$, $Y'$ be two projective manifolds with isomorphic fundamental groups. We identify the respective fundamental groups $\pi_1(Y) = \pi_1(Y') = \Gamma$. Then we say that the condition {\bf (SAME HOMOLOGY)} is satisfied for $Y$ and $Y'$ if there is an isomorphism $\Psi : H_*(Y', \ZZ) \cong
H_*(Y, \ZZ)  $
of homology groups which is compatible with the homomorphisms
$$ u \colon H_*(Y, \ZZ)  \ra   H_*(\Ga , \ZZ)   ,  u'   \colon
 H_*(Y', \ZZ)  \ra   H_*(\Ga , \ZZ)  ,   $$
i.e., $\Psi$ satisfies $u \circ \Psi = u'$.

\end{defin}

 We can now state the following

\begin{theorem}\label{special diagonal}\
  Let $X$ be a diagonal special  Inoue type manifold, and let $X'$ be a
projective manifold with $K_{X'}$ nef and
with the same fundamental group as $X$, which moreover either

(A)
is homotopically equivalent to $X$;

\noindent
or satisfies the following weaker property:

(B) let $\hat{X'}$ be the corresponding unramified covering of $X'$. Then 
$\hat{X}$ and $\hat{X}'$ satisfy the condition
{\bf (SAME HOMOLOGY)}.

\noindent
Setting $ W : = \hat{X'}$, we have that

\begin{enumerate}
\item
$X'=  W / G$ where $W$ admits a generically finite
morphism 
$f : W \ra Z'$, and where
$Z'$ is also a $K(\Ga, 1)$ projective manifold,
of the form $ Z' = (A'_1 \times \dots \times A'_r) \times
(C'_1 \times \dots \times  C'_h) \times (M'_1 \times \dots \times
M'_s)$.

Moreover here $M'_i$ is either $M_i$ or its complex conjugate,
and the product decomposition corresponds to the product decomposition
(\ref{prodFG}) of the fundamental group of $Z$.

The  image cohomology class  $f_*([W])$ corresponds, up to sign,
to the cohomology class of $\hat{X}$.
\item
The morphism $f$ is
finite if $n = \dim X$ is odd, and  it is
generically injective if 

(**) the cohomology class of $\hat{X}$ (in $H^* (Z, \ZZ)$) is indivisible,
or if every strictly submultiple  cohomology class cannot be represented by an effective
$G$-invariant divisor  on any pair  $(Z', G)$ homotopically equivalent to $(Z,G)$.
\item
 $f$ is an embedding if moreover 
 $K_{X'}$ is ample, 
 
(*)
 every such divisor $W$ of $Z'$  is ample, and

(***)  $\ \ \ K_{X'}^n = K_{X}^n$.\footnote{ This last property for
algebraic surfaces follows automatically from homotopy invariance.}
\end{enumerate}
In particular, if $K_{X'}$ is ample and (*), (**) and (***) hold, also $X'$ is a diagonal SIT (special  Inoue type) manifold.

A similar  conclusion holds under the alternative assumption that  the homotopy equivalence sends the canonical class of $W$
to that of $\hat{X}$: then  $X'$ is a minimal resolution of a diagonal SIT (special  Inoue type) variety.

\end{theorem}

Hypothesis  (A) in  theorem \ref{special diagonal} is used  to derive the
conclusion that also
$W : = \hat{X'} $ admits a holomorphic map $f'$ to a complex manifold $Z'$
with the same structure as $Z$, while
 hypotheses  (B) and following ensure  that the morphism is birational onto its image,
and the class of the image
divisor $f' (\hat{X'} )$ corresponds to $\pm$ that of $\hat{X} $  under
the identification
$$ H_*(Z' , \ZZ) \cong  H_*(\Ga , \ZZ) \cong   H_*(Z , \ZZ).$$

Since $K_{X'}$ is ample, one uses (***) to conclude that $f'$ is an isomorphism with its image.

The next question which the theorem leaves open  is weak * rigidity, for which 
  several ingredients should come into play: the Hodge type,
a fine analysis of
the structure of the action of $G$ on $Z$, the problem of  existence of hypersurfaces on which $G$ acts freely and the study 
of  the family of such  invariant effective divisors, in particular whether the family has a  connected base.

We have  restricted  ourselves to special Inoue type manifolds in order to be able to use the regularity results for classifying maps
discussed in the previous section
(the diagonality assumption is only a simplifying assumption).

Let us  now sketch the proof of  theorem \ref{special diagonal}.

\noindent
\Proof {\em of Theorem \ref{special diagonal}.}

\noindent
{ \bf Step 1}

The first step consists in showing that  $W : = \hat{X'}$ admits a holomorphic
mapping to a manifold $Z'$ of the above type $ Z' = (A'_1 \times 
\dots \times A'_r) \times
(C'_1 \times \dots \times  C'_h) \times (M'_1 \times \dots \times
M'_s)$, where $M'_i$ is either $M_i$ or its complex conjugate.

First of all, by the  results of Siu and others (\cite{siuannals}, 
\cite{siu2},\cite{fibred},  \cite{cime}, Theorem 5.14) cited in section 6, $W$ admits a 
holomorphic map to
a product manifold of the desired type  $$Z_2'\times Z'_3 = (C'_1 \times \dots \times 
C'_h) \times (M'_1 \times
\dots
\times  M'_s).$$

Look now at the  Albanese variety $\Alb (W)$ of the K\"ahler manifold 
$W$, whose fundamental
group is the quotient of the Abelianization of $\Ga = \pi_1 (Z)$ by its torsion subgroup. 
Write the fundamental group of
$\Alb (W)$ as the first homology group of $A  \times  Z_2  \times Z_3$,
i.e., as 
$$ H_1 ( \Alb (W)) = \Lambda \oplus H_1 ( Z_2, \ZZ) \oplus (H_1 ( Z_3, \ZZ)/ Torsion),$$ 
($\Alb (Z_2)$ is the product of Jacobians
$J : = (\Jac(C_1) \times \dots \times  \Jac(C_h))  $).

Since however, by the universal property,  $\Alb (W)$ has a holomorphic map to 
$$B' : = \Alb (Z'_2) \times \Alb (Z'_3),$$
inducing a splitting of the lattice  $H_1 (\Alb (W), \ZZ) = \Lambda \oplus 
H_1 (B', \ZZ)$,
it follows that $\Alb(W)$ splits as $A' \times B'$.

Now, we want to show that the Abelian variety $A'$ ($W$ is assumed to 
be a projective
manifold) splits as desired.  This is in turn a consequence of assumption (3) in definition \ref{SIT}.
In fact, the group $G$ acts on the Abelian variety $A'$ as a group of 
biholomorphisms,
hence it acts on $\Lambda \otimes \RR$ commuting with multiplication by 
$\sqrt{-1}$.
Hence multiplication by $\sqrt{-1}$ is an isomorphism of $G$ representations, and then (3)
 implies  that $\Lambda_i\otimes \RR$ is stable by  multiplication by 
$\sqrt{-1}$;
whence $\Lambda_i\otimes \RR$ generates a subtorus $A'_i$. Finally, $A'$ splits because $\Lambda$ is
the direct sum of the sublattices $\Lambda_i$. We are through with the 
proof of step 1.

\noindent
{ \bf Step 2}

Consider now the holomorphic map $ f \colon W \ra Z'$. We shall show 
that  the image
$ W' : = f (W)$ is indeed a divisor in $Z'$.
 For this we use  the assumption  (SAME HOMOLOGY) and, in fact, the claim is an immediate consequence of the following lemma.

\begin{lemma}\label{cohalg}
Assume that $W$ is a K\"ahler manifold, such that
\begin{itemize}
\item[i)] there is an isomorphism of
fundamental groups $\pi_1 (W) = \pi_1 (\hat{X})= \Ga$, where $ 
\hat{X}$ is a smooth ample
divisor in a $K(\Ga, 1)$ complex projective manifold $Z$;
\item[ii)] there exists a holomorphic map  $ f : W \ra Z'$, where $Z'$ is another 
$K(\Ga, 1)$ complex  manifold,  such that $f_* : \pi_1 (W) \ra  \pi_1 (Z') = \Ga$ is an 
isomorphism, and moreover
\item[iii)]  {\bf (SAME HOMOLOGY)} there is an isomorphism $\Psi : H_*( W, \ZZ) \cong
H_*(\hat{X}, \ZZ)  $
of homology groups which is compatible with the homomorphisms
$$ u \colon H_*(\hat{X}, \ZZ)  \ra   H_*(\Ga , \ZZ)   ,  u'   \colon
 H_*( W, \ZZ)  \ra   H_*(\Ga , \ZZ)  ,   $$
i.e., we have $u \circ \Psi = u'$.

\end{itemize}
Then $f$ is a generically finite  morphism of $W$ into $Z'$, and 
the  cohomology class  $f_* ([W])$  in
$$ H^*(Z', \ZZ) = H^*(Z, \ZZ)
=H^*(\Ga, \ZZ)  $$ corresponds to $\pm 1$  the one  of $ \hat{X}$.
\end{lemma}

{\it Proof of the Lemma.}
 We can identify the homology groups of $W$ and $ \hat{X}$ under $\Psi : H_*(W, \ZZ) \cong
H_*(\hat{X}, \ZZ) $, 
and then the image
in the homology groups of
$ H_*(Z', \ZZ) = H_*(Z, \ZZ) =H_*(\Ga , \ZZ)  $ is the same.

We apply the above consideration to the fundamental classes of the oriented manifolds $W$ and $ \hat{X}$,
which are generators of the infinite cyclic top degree homology groups $H_{2n}(W, \ZZ)$,
respectively $H_{2n}( \hat{X}, \ZZ)$.

This implies a fortiori that $ f \colon W \ra Z'$ is generically finite:
since then the homology class  $f_* ([W])$ (which we identify to a cohomology class  by virtue of Poincar\'e duality)
equals  the class of $ \hat{X}$, up to sign.

\qed

\noindent
{\bf Step 3}

We claim that  $ f \colon W \ra Z'$ is generically 1-1 onto its image $W'$.

\noindent
Let  $d$ be the degree of $ f \colon W \ra W'$.

\noindent
Then   $f_* ([W]) = d [W']$, hence if the class of 
$\hat{X}$ is indivisible, then obviously $d=1$.

Otherwise, observe that the divisor $W'$ is an effective $G$-invariant divisor and use our assumption (**).

\noindent
{\bf Step 4}

 Here  we are going to prove that $f$ is an embedding under the additional hypotheses that $K_X^n =K_{X'}^n$ and that $W'$ is ample,
 as well as $K_{X'}$.
 
We  established that  $f$ is birational onto its image $W'$, hence it is a 
desingularization of $W'$.

We  use now adjunction. We claim that, since (by our assumption on $K_{X'}$) $K_W$ is nef, there exists an effective 
divisor $\mathfrak A$,
called the adjunction divisor, such that
$$K_W = f^* (K_{Z'} + W') - \mathfrak A.$$

This can be shown by taking the Stein factorization
$$ p \circ h \colon W \ra W^N \ra W', $$
where $W^N$ is the normalization of $W'$.

Let $\sC $ be the conductor ideal  $\mathcal{H}om (p_*\hol_{W^N}, \hol_{W'} )$
viewed as an ideal $\sC \subset \hol_{W^N}$; then the Zariski  canonical divisor 
of $W^N$ satisfies 
$$K_{W^N} =   p^*   (K_{W'}) - C =   p^*   (K_{Z'} + W') - C$$
where $C$ is the Weil divisor associated to the conductor ideal (the equality on the Gorenstein locus
of $W^N$ is shown for instance in \cite{bucharest}, then it suffices to take the direct image from the open set to
the whole of $W^N$).

In turn, we would have in general 
$K_W = h^* (K_{W^N } ) - \mathfrak B $, with $\mathfrak B $ not necessarily effective; but, by Lemma 2.5 of \cite{alessio},
see also Lemma 3.39 of \cite{KollarMori},
and since $- \mathfrak B $ is h-nef, we conclude that $\mathfrak B $ is effective.
We establish the claim by setting $\mathfrak A : = \mathfrak B + h^* C$.

Observe that, under the isomorphism of homology groups, $f^* (K_{Z'} + W')$
corresponds to $(K_{Z} + \hat{X})|_{\hat{X}}= K_{\hat{X}} $, in 
particular we have
$$ K_{\hat{X}}^n = f^* (K_{Z'} + W')^n = (K_W + \mathfrak A)^n . $$
If we assume that $K_W$ is ample, then $ (K_W + \mathfrak A)^n \geq (K_W)^n $,
equality holding if and only if $  \mathfrak A = 0$.

Under assumption (**), it follows that  $$K_{\hat{X}}^n = |G| K_X^n = 
|G| K_{X'}^n= K_W^n,$$
hence  $  \mathfrak A = 0$. Since however $K_W$ is ample, it follows 
that $f$ is an embedding.

If instead we assume that $K_W$ has the same class as $ f^* (K_{Z'} + W')$,
we conclude first that necessarily $\mathfrak B =0$, and then we get that $C= 0$.

Hence $W'$ is normal and  has canonical singularities.

\noindent
{\bf Step 5}

Finally, the group $G$ acts on $W$, preserving the direct summands of 
its fundamental group
$\Ga$. Hence, $G$ acts on the curve-factors, and on  the locally symmetric factors.

By assumption, moreover, it sends the summand $\Lambda_i$ to itself, 
hence we get a well defined linear action  on each
Abelian variety $A'_i$, so that we have a diagonal linear action of 
$G$ on $A'$.

Since however the image of $W$ generates $A'$, we can extend the 
action of $G$ on $W$ to a compatible  affine
action on $A'$.
In the case where $G$ is abelian\footnote{in \cite{bc-inoue} and \cite{bcf} we made this assertion for any finite group $G$, see however remark \ref{erratum}}, we can show that the {\em real affine type} of the action on $A'$ is uniquely determined. This is taken care of by the following lemma.

\begin{lemma}\label{affine}
Given a  diagonal special Inoue type manifold, if $G$ is abelian,  the real affine type of the action of $G$ on the
Abelian variety $ A = (A_1 \times \dots \times A_r )$ is determined by the fundamental group
exact sequence

$$ 1 \ra \Ga= \Pi_{i=1}^r (\Lambda_i) \times   \Pi_{j=1}^h (\pi_{g_j}) 
\times \Pi_{l=1}^s
(\pi_1 (M_l))
\ra \pi_1 ( X) \ra G \ra 1.$$ 
\end{lemma}
{\em Proof.}
Define as before $\Lambda : = \Pi_{i=1}^r \Lambda_i = \pi_1 (A)$; moreover, since all the summands on
the left hand side are normal in $ \pi_1 ( X)$, set 
$$\overline{G} : =  \pi_1 ( X) / (( \Pi_{j=1}^h \pi_{g_j} )
\times (\Pi_{l=1}^s
\pi_1 (M_l) )).$$ 

Observe now that $ X$ is the quotient of its universal covering 
$$\tilde{X} = \CC^m \times \prod_{j=1}^h\mathbb{H}_j \times \prod_{l=1}^s\mathcal{D}_l$$ by its fundamental group,
acting diagonally (here $\mathbb{H}_j$ is a copy of Poincar\'e's upper half plane, while $ \mathcal{D}_l$ is an irreducible  bounded symmetric domain
of dimension at least two), hence we  obtain that $\overline{G}$ acts on $\CC^m$ as a group of affine transformations.

This action yields  a homomorphism $$ \alpha : \overline{G} \twoheadrightarrow \im (\alpha) =: \hat{G} \subset \Aff (m, \CC):$$
let  $\sK$  the kernel of $\alpha$, and let $$\overline{G}_1 : = \ker (  \alpha_L : \overline{G} \ra \GL  (m, \CC)).$$

$\overline{G}_1 $ is obviously Abelian, and contains $\Lambda$, and maps onto a lattice $\Lambda' \subset \hat{G}$.

Since $\Lambda$ injects into $\Lambda'$, $\Lambda \cap \sK = 0 $, whence $\sK$ injects into $G$, therefore $\sK$ is a torsion subgroup; since $\Lambda'$ is free, we obtain
$$\overline{G}_1 = \Lambda'  \oplus \sK,$$
and we finally get $$  \sK = \Tors (\overline{G}_1 ), \ \  \hat{G} =  \overline{G}  / \Tors (\overline{G}_1 ).$$

Since our action is diagonal, we can write $\Lambda'  = \oplus_{i=1}^r \Lambda'_i $, and the linear action of
the group $G_2 : = G / \sK$ preserves the summands. 

Since $ \hat{G} \subset \Aff (\Lambda' \otimes \RR)$, we can now apply proposition \ref{affinetype} to it. 

This shows that the affine group  $ \hat{G}$ is uniquely determined.

Finally, using the image groups $G_{2,i}$ of $G_2$ inside $ \GL ( \Lambda'_i)$,
we can  define uniquely groups of affine transformations of $A_i$ which fully determine
the diagonal action of $G$ on $A$ (up to real affine automorphisms of each $A_i$).

\qed

  The proof of  Theorem \ref{special diagonal} is now completed.

\qed 
\bigskip

In order to obtain weak rigidity results, one has to 
use, as an   invariant for group actions on tori and Abelian varieties,  the Hodge type, introduced in definition \ref{HodgeT}
(see also remark \ref{Hodgetype}).

\begin{rem}
In the previous theorem special assumptions are needed in order to guarantee that for each manifod $X'$ homotopy
equivalent to $X$ the classifying holomorphic map $f : \hat{X' } \ra Z'$ be birational onto its image, and indeed an embedding.

However, there is the possibility that an Inoue type variety deforms to one for which $f$ is a covering of finite degree.
This situation should be analyzed and the singularities of the image of $f$ described in detail, so as to lead to a
 generalization of the   theory of    Inoue type varieties, including   `multiple'  Inoue type varieties (those for which $f$ has
 degree at least two). 

\end{rem}

The study of moduli spaces of Inoue type varieties, and their connected and irreducible components, relies very much
on  the study of moduli spaces of varieties $X$ endowed with the action of a finite group $G$: and it  is for us a strong motivation 
to pursue this line of research. 

This topic will occupy a central role in the following sections, first in general, and then in the special case of algebraic curves.

\section{ Moduli spaces of  surfaces and higher dimensional varieties}

 Teichm\"uller theory works out quite fine in the case of projective curves, as well as other approaches,
like Geometric Invariant Theory (see \cite{GIT},  \cite{EnsMath}, \cite{Giescime}, \cite{Giestata}), which provides a quasi-projective
moduli space $\mathfrak M_g$ endowed with a natural
compactification $\overline{\mathfrak M_g}$ (this is a coarse moduli space for the so-called moduli stable projective
curves: these are the  reduced curves with at most
nodes as singularities,  such that their group of automorphisms is finite).

In higher dimensions one has a fully satisfactory theory of `local moduli' for compact complex manifolds
or spaces,  but there are difficulties  with the global theory.

So, let us start from the local theory, developed by Kodaira and Spencer, and culminating in the results of Kuranishi and Grauert.

\subsection{Kodaira-Spencer-Kuranishi theory}

While describing complex structures as integrable almost complex structures (see theorem \ref{NN}),
it is convenient to view an almost  complex structure as a differentiable $(0,1)$ - form with values in the dual of the
cotangent bundle $(TY^{1,0})^{\vee}$.

This representation leads to the  Kodaira-Spencer-Kuranishi
theory of local deformations, addressing precisely  the study of the {\bf small
deformations} of a complex manifold $ Y = (M, J_0)$.

In this theory, complex structures correspond to closed such forms which,  by Dolbeault's theorem,  determine a cohomology
class  in $H^1 ( \Theta_Y)$, where $ \Theta_Y$ is the sheaf of holomorphic 
sections of the holomorphic tangent bundle  $TY^{1,0}$.

We shall use here unambiguously the double notation $TM^{0,1} = TY^{0,1} $, $ TM^{1,0} = TY^{1,0} $ to refer  to the 
splitting of the complexified tangent bundle determined by the complex structure $J_0$. 

$J_0$ is a point in $\sC (M)$, and  a neighbourhood in the space of almost complex structures
corresponds to a distribution of subspaces which are globally defined as  graphs of an endomorphism
$$ \phi :  TM^{0,1} \ra TM^{1,0},$$
called a {\bf small variation of almost complex structure}, since one then defines 
 $$TM^{0,1}_{\phi} : = \{ (u, \phi (u))| \ u \in TM^{0,1} \} \subset TM^{0,1} \oplus  TM^{1,0}.$$
 
 In terms of the old $\bar{\partial}$ operator, the new one is simply obtained by considering
 $$\bar{\partial}_{\phi} : =  \bar{\partial} + \phi, $$
 and the integrability condition is given by the Maurer-Cartan equation
 
$$ (MC) \ \ \bar{\partial} (\phi) + \frac{1}{2} [ \phi, \phi ] = 0. $$

Observe that, since our original complex structure $J_0$ corresponds to $\phi = 0$, the derivative
of the above equation at $\phi = 0$ is simply
$$ \bar{\partial} (\phi)= 0,$$
hence the tangent space to the space of complex structures consists of the space of $ \bar{\partial} $-closed
forms of type $(0,1)$ and with values in the bundle $ TM^{1,0}$.

One  can restrict oneself (see e.g. \cite{handbook})  to consider only the class of such forms $\phi$ in the
Dolbeault cohomology group  $$H^1 (\Theta_Y): =  ker ( \bar{\partial} ) / Im ( \bar{\partial})  ,$$
by looking at the action of the group of diffeomorphisms which are exponentials of global vector fields on $M$.

Representing these cohomology classes by harmonic forms, the integrability equation becomes easier to solve
via the following Kuranishi equation.

Let $\eta_1, \dots , \eta_m \in H^1 (\Theta_Y)$ be a basis for the space of harmonic (0,1)-forms with values in 
$TM^{1,0}$, and set $t : =(t_1, \dots, t_m) \in \CC^m$, so that  $ t \mapsto \sum_i t_i \eta_i$ 
establishes an isomorphism $\CC^m \cong H^1 (\Theta_Y)$. 

Then the {\em Kuranishi slice}  is obtained
by associating to $t$ the unique power series solution of the following equation:
$$ \phi (t) =   \sum_i t_i \eta_i + \frac{1}{2} \bar{\partial}^* G [ \phi (t), \phi (t)] , $$
satisfying moreover  $\phi (t) =   \sum_i t_i \eta_i + $ higher order terms
($G$ denotes here  the Green operator).  

The upshot is that for these forms $ \phi (t)$
the integrability equation simplifies drastically; the result is summarized in the following definition.

\begin{defin}
The  Kuranishi space $\frak B (Y)$ is defined as the germ of complex subspace of $H^1 (\Theta_Y)$
defined by
$\{ t \in \CC^m | \  H  [ \phi (t), \phi (t)] = 0 \} $, where $H$ is the harmonic projector onto the space 
$ H^2 (\Theta_Y)$ of harmonic forms of type $(0,2)$ and with values in $TM^{1,0}$.

 Kuranishi space $\frak B (Y)$ parametrizes  the set of small variations of almost complex structure $ \phi (t)$
which are integrable. Hence over $\frak B (Y)$ we have a family of complex
structures which deform the complex structure of $Y$.

\end{defin}

It follows then  that the  Kuranishi space $\frak B (Y)$
surjects onto the germ of
the  Teichm\"uller space at the point corresponding to the  given
complex structure $Y = (M, J_0)$.

It fails badly in general to be a homeomorphism (see \cite{montecatini}, \cite{handbook}): for instance Teichm\"uller space, 
in the case of the Hirzebruch Segre surface $\FF_{2n}$ (blow up of the cone onto the rational normal curve of degree $2n$),
 consists of $n+1$ points $p_0, p_2,  \dots ,p_{2n}$, such that the open sets are exactly the sets $\{ p_{2i} | i \leq k\}$.

But a  first consequence is that Teichm\"uller space  is locally connected by holomorphic arcs,
hence the determination of the connected components of $\sC (M)$, respectively 
 of $\sT (M)$, can be done using the original definition of 
deformation equivalence, given by Kodaira and Spencer in \cite{k-s58} (see definition \ref{deformation}).

One can define deformations not only for complex manifolds, but also for complex spaces.

\begin{defin}
1) A {\bf deformation}  of a compact complex space $X$ is a pair
consisting of

1.1) a flat proper morphism $ \pi : \X \ra T$ between connected complex spaces (i.e.,
$\pi^* : \hol_{T,t} \ra \hol_{\X,x}$ is a flat ring extension for each
$ x$ with $ \pi (x) = t$)

1.2) an isomorphism $ \psi : X \cong \pi^{-1}(t_0) : = X_0$ of $X$ with a fibre
$X_0$ of $\pi$.

2.1) A {\bf small deformation} is the germ $ \pi : (\X, X_0) \ra (T, t_0)$ of
a deformation.

2.2) Given a deformation  $ \pi : \X \ra T$ and a morphism $ f : T' \ra T$
with $ f (t'_0) = t_0$, the {\bf pull-back} $ f^* (\X)$
is the fibre product $ \X': = \X \times_T T'$ endowed with the projection onto
the second factor $T'$ (then $ X \cong X'_0$).

3.1) A small deformation  $ \pi : \X \ra T$ is said to be {\bf versal} or {\bf complete} if every
small deformation $ \pi : \X' \ra T'$ is obtained from it via pull back; it is said
to be {\bf semi-universal} if the differential of $ f : T' \ra T$
 at $t'_0$ is uniquely determined, and {\bf universal}  if the morphism
 $f$ is uniquely determined.

4) Two compact complex manifolds $X,Y$ are said to be {\bf directly deformation
equivalent} if there are 

4i) a deformation $ \pi : \X \ra T$ of $X$ with $T$
  irreducible and where all the fibres are smooth,
and 

4ii) an isomorphism $ \psi' : Y \cong \pi^{-1}(t_1)  = : X_1$ of $Y$ with a
fibre $X_1$ of $\pi$.
\end{defin}

\begin{rem} 
The technical assumption of flatness replaces, for families of spaces,   the condition on $\pi$ to be a submersion,
necessary in order that the fibres be smooth manifolds.
\end{rem}

Let's however come back to the case of complex manifolds, observing that in a small deformation
of a compact complex manifold one can shrink the base $T$ and assume that all the fibres are smooth.

We can now state the results of Kuranishi and Wavrik (\cite{kur1}, \cite{kur2}, \cite{Wav}) in the language of deformation theory.

\begin{theo}
{\bf (Kuranishi).} Let $Y$
be a compact
complex manifold: then

I) the Kuranishi family   $ \pi : (\Y, Y_0) \ra (\frak B (Y), 0)$ of
$Y$ is semiuniversal.

II) $(\frak B (Y), 0)$ is unique up to (non canonical) isomorphism, and is a germ of analytic
subspace of the vector space
$ H^1 (Y, \Theta_Y)$, inverse image of the origin under a local
holomorphic map (called Kuranishi map) 
$$ k : H^1 (Y, \Theta_Y) \ra H^2 (Y, \Theta_Y) $$ whose differential vanishes
at the origin. 

Moreover  the  quadratic term in the Taylor development of
$k$ is  given by  the bilinear map $ H^1 (Y, \Theta_Y) \times
   H^1 (Y, \Theta_Y) \ra  H^2 (Y, \Theta_Y)$, called Schouten bracket,
   which is the composition of  cup product  followed by Lie bracket of vector fields.

III) The Kuranishi family is a versal deformation of $Y_t$ for $t \in  \frak B (Y)$.

IV) The Kuranishi family is universal if $H^0 (Y, \Theta_Y)= 0.$

V) {\bf (Wavrik)} The Kuranishi family is universal if $ \frak B (Y)$ is reduced and $h^0 (Y_t, \Theta_{Y_t}) : = \rm{dim}\ H^0 (Y_t, \Theta_{Y_t}) $
is constant  for $t \in   \frak B (Y)$ in a suitable neighbourhood of $0$.
\end{theo}

The Kodaira Spencer map, defined soon below, is to be thought as the derivative of a family of complex structures.

\begin{defin}
The  Kodaira Spencer map of a family $ \pi : (\Y, Y_0) \ra (T, t_0)$
of complex manifolds having a smooth base $T$ is defined as follows: consider the cotangent bundle sequence of the fibration
$$ 0 \ra \pi^* (\Omega^1_T) \ra  \Omega^1_{\Y} \ra  \Omega^1_{\Y| T} \ra 0,$$
and the direct image sequence of the  dual sequence of bundles,
$$ 0 \ra \pi_* (\Theta_{\Y| T}) \ra   \pi_* (\Theta_{\Y}) \ra   \Theta_T  \ra  \sR^1  \pi_* (\Theta_{\Y| T}) .$$

Evaluation at the point $t_0$ yields a map $\rho$ of the tangent space to $T$ at $t_0$ into $ H^1 (Y_0, \Theta_{Y_0})$,
which is the derivative of the variation of complex structure. 

\end{defin}

The Kodaira Spencer map and the implicit functions theorem allow  to determine 
 the Kuranishi space and the Kuranishi family in many cases.

\begin{cor}
Let $Y$
be a compact
complex manifold and assume that we have a family $ \pi : (\Y, Y_0) \ra (T, t_0)$ with smooth base $T$, such that $Y \cong Y_0$,
and such that the Kodaira Spencer map $\rho_{t_0}$ surjects onto $ H^1 (Y, \Theta_Y)$.

Then the Kuranishi space $\frak B (Y)$ is smooth and there is a submanifold $T' \subset T$ which maps isomorphically to $\frak B (Y)$;
hence the Kuranishi family is the restriction of $\pi$ to $T'$. 
\end{cor}

The key point is that, by versality  of the Kuranishi family, there is a morphism $f : T \ra \frak B (Y)$ inducing $\pi$ as a pull back,
and $\rho$ is the derivative of $f$: then one uses the implicit functions theorem. 

This approach clearly works only if $Y$ is {\bf unobstructed}, which simply means
that $\frak B (Y)$ is smooth. In general it is difficult to describe the obstruction map, and even calculating the quadratic term
is nontrivial (see \cite{quintics} for an interesting  example).

Even if  it is difficult to calculate the obstruction map,  Kuranishi theory  gives a lower bound for the `number of moduli'
of $Y$, since it shows that $\frak B (Y)$ has dimension  $\geq  h^1 (Y, \Theta_Y) - h^2 (Y, \Theta_Y)$. In the case of curves
$ H^2 (Y, \Theta_Y) = 0$, hence curves are unobstructed; in the case of a surface $S$ 
$$ \rm{dim} \frak B (S ) \geq  h^1 ( \Theta_S) - h^2 ( \Theta_S) = - \chi ( \Theta_S) + h^0 ( \Theta_S) = 10 \chi (\hol_S) - 2 K^2_S + h^0 ( \Theta_S).$$

The above is the Enriques inequality (\cite{enr},  see also \cite{enr96}, and \cite{clemens} for an improvement: observe that Max Noether postulated equality),  proved by Kuranishi in all cases and also for non algebraic surfaces.

\subsection{Kuranishi and Teichm\"uller}

Ideally, we would like to have  that Teichm\"uller space, up to now only defined as a topological space,
is indeed a complex space, locally isomorphic to Kuranishi space.

In fact, we already remarked that there is a locally surjective continuous map of $ \frak B (Y)$ to the germ
 $\sT (M)_Y$ of $\sT (M)$ at the point corresponding to the complex structure yielding $Y$. 
 For curves this map is a local homeomorphism, and this fact provides a complex structure on Teichm\"uller space.
 
 Whether this holds in general is related to the following definition.

\begin{defin}
A compact complex manifold $Y$  is said to be {\bf rigidified}  if $Aut(Y) \cap \sD iff ^0(Y) = \{ Id_Y \}.$

A compact complex manifold $Y$  is said to be cohomologically rigidified if $Aut(Y) \ra Aut (H^* (Y, \ZZ))$
is injective, and rationally cohomologically rigidified if  $Aut(Y) \ra Aut (H^* (Y, \QQ))$ is injective.
\end{defin}

\begin{rem}
We refer to \cite{cai-wenfei} for the proof that surfaces of general type with $ q \geq 3$ are cohomologically rigidified,
and for examples showing that there are  rigidified surfaces of general type which are not cohomologically rigidified.
\end{rem}

In fact,   it is clear that there is a universal tautological family of complex structures
parametrized by $ \sC (M)$ (the closed subspace $ \sC (M)$  of $\sA \sC (M)$ consisting of the set of complex structures
on $M$), and with total space
$$\mathfrak U_{ \sC (M)} : =  M   \times  \sC (M) ,$$
on which the group   $\sD iff ^+ (M)$ naturally acts, in particular  $\sD iff ^0(M)$.

The main observation is that $\sD iff ^0(M)$ acts freely on $ \sC (M)$ if and only if for each complex structure
$Y$ on $M$ the group of biholomorphisms $Aut(Y)$ contains no nontrivial automorphism which is differentiably isotopic
to the identity. Thus, the condition of being rigidified implies that the tautological family of complex structures
descends to a universal family of complex structures on Teichm\"uller space:

$$\mathfrak U_{ \sT (M)} : =   (M   \times  \sC (M)) / \sD iff ^0(M) \ra  \sC (M)) / \sD iff ^0(M) = \sT (M) ,$$
on which the mapping class group acts.

Fix now a complex structure yielding a compact complex manifold $Y$, and compare with the Kuranishi family
$$\sY \ra \frak B (Y).$$

Now, we already remarked that there is a locally surjective continuous map of $ \frak B (Y)$ to the germ
 $\sT (M)_Y$ of $\sT (M)$ at the point corresponding to the complex structure yielding $Y$. 

 The following was observed in \cite{handbook}.
 
\begin{rem}
If

1) the Kuranishi family is universal at any point

2)   $ \frak B (Y) \ra \sT (M)_Y$ is injective (it is then a local homeomorphism at every point)

then Teichm\"uller space has a natural structure of complex space.

\end{rem} 

In  many cases (for instance, complex tori)  Kuranishi and Teichm\"uller space coincide, in spite of the fact that
 the manifolds are not rigidified. For instance we showed in \cite{handbook}:

\begin{prop}\label{kur=teich}
 1)  The continuous map $\pi \colon  \frak B (Y) \ra \sT (M)_Y$ is a local homeomorphism
  between Kuranishi space and Teichm\"uller space if there is an
  injective continuous map
  $f \colon  \frak B (Y) \ra Z$, where $Z$ is Hausdorff, which factors through $\pi$.
  
 2)  Assume that $Y$ is a compact K\"ahler manifold and that the local period map $f$ is injective:
  then $\pi \colon  \frak B (Y) \ra \sT (M)_Y$ is a local homeomorphism.
  
  3) In particular, this holds if $Y$ is K\"ahler with trivial or torsion canonical divisor.

\end{prop}

\begin{rem}
1) The condition of being rigidified implies the condition $H^0 (\Theta_Y)=0$ (else there is
a positive dimensional Lie group of biholomorphic self maps), and is obviously implied by the condition of being cohomologically rigidified.

2) By the cited Lefschetz' lemma  compact curves of genus $ g \geq 2$ are rationally cohomologically rigidified,
and it is an interesting question whether compact complex manifolds of general type are rigidified
(see   \cite{cai}, and \cite{cai-wenfei} for recent progresses).

3) One can dispense with many assumptions, at the cost of having a more complicated result. For instance,
Meersseman shows in \cite{meerssemann}, explaining the situation via  classical examples  of \cite{catrav}, that one can take a quotient of the several Kuranishi families
(using their semi-universality), obtaining, as analogues of Teichm\"uller spaces, respectively moduli spaces, an analytic Teichm\"uller groupoid, 
respectively a Riemann groupoid, both independent up to analytic Morita equivalence
of the chosen countable disjoint union of Kuranishi spaces. These lead to some stacks called by the author `analytic Artin stack'
 (see \cite{artinstacks} and \cite{fantechi} as a reference on Artin stacks).
\end{rem}

\subsection{Varieties with singularities}

For higher dimensional varieties moduli theory works better if one considers varieties with moderate singularities,
rather than smooth ones (see \cite{Kolhand}), as we shall illustrate in the next sections for the case of algebraic surfaces.

The Kuranishi theory does indeed extend perfectly to all compact complex spaces, and Kuranishi's theorem 
was generalized by  Grauert 
 (see  \cite{grauert1, grauert}, see also \cite{sernesi} for the algebraic analogue).

\begin{theo}
{\bf Grauert's Kuranishi type theorem for complex spaces.} Let $X$
be a compact
complex space: then

I) there is a semiuniversal deformation $ \pi : (\X, X_0) \ra (T, t_0)$ of
$X$, i.e., a deformation such that every  small deformation
$ \pi' : (\X', X'_0) \ra (T', t'_0)$ is the pull-back of $\pi$ for
an appropriate morphism $f : (T', t'_0) \ra (T, t_0)$ whose
differential at $t'_0$ is uniquely determined.

II) $(T, t_0)$ is unique up to isomorphism, and is a germ of analytic
subspace of the
vector space
$\T ^1$ of first order deformations. 

$(T, t_0)$ is the  inverse image of the origin under a local
holomorphic map (called Kuranishi map) 
$$ k :
\T ^1  \ra \T ^2 $$  to the finite dimensional vector space $\T ^2$ (called {\bf obstruction space}), and  whose differential vanishes
at the origin (the point corresponding to the point $t_0$).

If $X$ is reduced, or if the singularities of $X$ are local complete intersection singularities, then $\T ^1 =  {\rm Ext }^1 (\Omega^1_X, \hol_X ).$ 

If the singularities of $X$ are local complete intersection singularities, then $ \T ^2 = {\rm Ext }^2 (\Omega^1_X, \hol_X) $ .
\end{theo}

 \bigskip
 
 Indeed, the singularities which occur for the canonical models of
 varieties of general type are called canonical singularities and are somehow tractable (see \cite{youngguide});
 nevertheless the study  of their deformations may present highly  nontrivial problems.

For this reason, we restrict now ourselves to the case of complex dimension two, where these singularities are easier to describe.

\section{Moduli spaces of surfaces of general type}

\subsection{Canonical models of surfaces of general type.}

The classification theory of algebraic varieties proposes to classify the birational equivalence classes of projective varieties.

Now, in  the birational class of a non ruled projective surface there is, by the theorem of Castelnuovo (see e.g. \cite{arcata}),
a unique (up to isomorphism) minimal model $S$ (concretely, minimal means that $S$ contains no (-1)-curves, i.e., curves
$E$ such that $E \cong \PP^1$, and with $E^2=-1$).

We shall assume from now on  that $S$ is a smooth minimal (projective) surface of general
type: this is equivalent (see \cite{arcata}) to the two conditions:

(*) $K_S^2 > 0$ and
$K_S$
   is nef, where as well known,  a divisor $D$ is said to be
nef if, for each irreducible curve $C$, we have $ D \cdot C \geq 0$.

It is very important that, as shown by Kodaira in \cite{kod-1}, the class of non  minimal surfaces is stable by small deformation;
on the other hand, a small deformation of a minimal algebraic surface of general type is again minimal (see prop. 5.5 of \cite{bpv}).
Therefore, the class of minimal algebraic surfaces of general type is stable by deformation in the large.

Even if the canonical divisor $K_S$ is nef, it does not however need to  be an ample divisor, indeed 

{\em  The canonical divisor $K_S$ of a minimal surface of general type $S$
  is ample iff
there does not exist an irreducible curve $C$  on $S$ with $K\cdot
C=0$ $\Leftrightarrow $ there is no (-2)-curve $C$ on $S$, i.e., a curve such  that $C\cong \PP^1$, and $C^2=-2$.}

The number of (-2)-curves  is bounded by  the rank of the Neron Severi lattice 
$NS(S) \subset H^2(S, \ZZ)/ Torsion $ of $S$ ($NS(S)$ is the image of $Pic(S) = H^1 (\hol^*_S)$ inside  $H^2(S, \ZZ)/ Torsion $,
and, by Lefschetz' (1,1) theorem, it is the intersection with the Hodge summand $H^{1,1}$); the (-2)-curves can be contracted by a contraction $\pi \colon S \ra X$, where $X$ is  a normal surface which is
called the
{\bf canonical model} of $S$.

The singularities of $X$ are called Rational Double Points (also called  Du Val or Kleinian singularities, see \cite{duval}), and $X$ is a Gorenstein variety,
i.e. (see \cite{hart}) the dualizing sheaf $\omega_X$ is invertible, and the associated Cartier divisor $K_X$, called again 
canonical divisor,  is   ample and such that $ \pi^* ( K_X) = K_S$.

$X$  is also   the projective spectrum (set of homogeneous
prime ideals) of the canonical ring $$\sR (S) : = \sR(S,K_S) : =  \bigoplus_{m \geq 0}H^0 (\hol_S (mK_S).$$

This definition generalizes to any dimension, since for a variety of general type
 (one for which there is a pluricanonical map which is birational onto its image) the canonical ring 
 is a finitely generated graded ring, as was proven by Birkar, Cascini, Hacon and McKernan in \cite{bchm}.
 And one defines a canonical model as a  variety with canonical singularities and with ample canonical divisor $K_X$.

More concretely, the  canonical model of a surface of general type
is  directly obtained as the image, for $m \geq 5$,  of  the  $m$-th pluricanonical map of $S$ (associated to the sections in $H^0 (\hol_S (mK_S))$)
as shown by Bombieri  in \cite{bom}  (extended to any characteristic by Ekedahl \cite{ekedahl}, and  in \cite{cf, 4auth}).

\begin{theo}{\bf (Bombieri)}
   Let $S$ be a minimal surface of general type,
and consider the m-th pluricanonical map $\varphi_m$ of $S$ (associated to the linear system $|mK_S|$) for $m\ge 5$, or for $m=4$
when $K^2_S \ge 2$.

Then $\varphi_m$ is a  birational morphism  whose image  is isomorphic to the  canonical model $X$, embedded by its
 m-th pluricanonical map.
 
 \end{theo}
 
 \begin{cor}{\bf (Bombieri), \cite{arcata}}
 Minimal surfaces $S$ of general type with given  topological invariants $e(S), b^+(S) $ (here, $b^+(S) $ is the positivity
 of the intersetion form on $H^2(S, \ZZ)$, and the pair of  topological invariants is equivalent to the pair of 
  holomorphic invariants  $K^2_S,   \chi (S) := \chi (\hol_S)$) are `bounded',
  i.e., they belong to a finite number of  families having a  connected base.
  
  In particular, for fixed  Euler number and positivity $e(S), b^+(S) $ we have a finite number of differentiable and topological types.
 \end{cor}
 
 In fact, one can  deduce  from here effective upper bounds  for the number of these families, hence for the types
 (see \cite{Cat-Chow} and \cite{paris}); see \cite{lp} for recent work on lower bounds and references to other recent and less recent results.
 
Results in the style of Bombieri have been recently obtained by Hacon and Mckernan in \cite{hm-pluric}, and from these one
obtains (non explicit) boundedness results for 
 the varieties of general type with fixed invariants $K^n_X, \chi(X)$ .

 \subsection{ The Gieseker moduli space}

When one deals with  projective varieties or projective subschemes
the most natural parametrization, from the point of view of deformation theory,
is given by the Hilbert scheme, introduced by Grothendieck (\cite{groth}).

For instance, in  the case  of surfaces of general type with fixed invariants  $\chi  (S) = a $ and $K^2_S = b$,
 their 5-canonical models $X_5$ are surfaces with
Rational Double Points as singularities and of degree $ 25 b$ in a fixed
projective space $ \PP^N$, where $ N + 1 = P_5: = h^0 (5 K_S) =
\chi (S)  + 10 K^2_S = a + 10 b$.

The  Hilbert polynomial of $X_5$ is the polynomial 
$$ P (m) :=  h^0 (5 m K_S) =
a + \frac{1}{2}( 5 m -1)  5 m  b .$$

Grothendieck (\cite{groth})  showed that, given a Hilbert polynomial (see \cite{hart}),  there is

i) an integer $d$ and

ii) a subscheme $\HHH = \HHH_P$ of the
Grassmannian of codimension $P(d)$- subspaces of $ H^0 (\PP^N, \hol_{\PP^N} (d))$,
called Hilbert scheme, such that

iii) $\HHH $  parametrizes the degree $d$ graded pieces $ H^0
(\sI_{\Sigma}(d))$ of the homogeneous ideals of all the subschemes $\Sigma
\subset \PP^N$
having the given Hilbert polynomial $P$.

The Hilbert point of $\Sigma$ is the Pl\"ucker point 
$$ \Lambda^{P(d)} (r_{\Sigma}^{\vee}) \in \PP (H^0 (\PP^N, \hol_{\PP^N} (d))^{\vee}   )$$

where  $ r_{\Sigma} $ is  the restriction homomorphism (surjective for $d$ large)
$$ r_{\Sigma} :  H^0 (\PP^N, \hol_{\PP^N} (d)) \ra  H^0 (\Sigma, \hol_{\Sigma} (d)).$$

Inside $\HHH$ one has the open sets
$$  \HHH ^* =  \{ \Sigma | \Sigma {\rm \ is \ smooth}  \} \subset \HHH ^0 ,$$
where
$$ \HHH ^0 : = \{ \Sigma | \Sigma {\rm \ is\ reduced \ with \ only \
rational \ Gorenstein  \
singularities } \}. $$

Gieseker showed in \cite{gieseker},  replacing 
the 5-canonical embeddding by an $m$-canonical embedding with
much higher $m$,  the following

\begin{theo}{\bf (Gieseker)}
The moduli space    of canonical models of surfaces of general type
with invariants $\chi, K^2$ exists as a quasi-projective  scheme $$ \mathfrak M_{\chi, K^2}^{can}$$
which  is called the {\bf Gieseker moduli space}.
\end{theo}

Similar theorems hold also for higher dimensional varieties, see \cite{lastopus}, and \cite{Kolhand} for the most recent developments,
but we do not state here the results, which are more  complicated and technical.

\subsection{Components of moduli spaces and deformation equivalence}

We mentioned previously that the relation of deformation equivalence is a good general substitute for
the condition that two manifolds belong to the same connected (resp.: irreducible) component of the moduli space.

Things get rather complicated and sometimes pathological  in higher dimension, especially since, even for varieties of general type,
 there can be three models, smooth, terminal, and canonical model, and the latter are singular.
 
In the case of surfaces of general type things still work out fine (up to a certain extent),
since the main issue is to compare the deformations of 
minimal models versus the ones of canonical models.

We have then  the following theorem.

\begin{theo}\label{can=min}

Given two minimal surfaces of general type $S, S'$ and their respective
canonical models $X, X'$, then

$S$ and $S'$ are deformation equivalent $\Leftrightarrow$ $X$ and $X'$ are
deformation equivalent  $\Leftrightarrow$ $X$ and $X'$ yield two points in the same connected component of the Gieseker moduli space.
\end{theo}

One  idea behind the proof is the  observation that, in order to analyse
deformation equivalence, one may restrict oneself to the case of families parametrized by a base $T$
with $ dim (T) = 1$: since two points in a complex space $ T \subset \CC^n$ (or in an algebraic variety)
belong to the same irreducible component of $T$ if and only if they
belong to an irreducible curve $ T ' \subset T$.
And one may further reduce to the case where  $T$ is smooth simply by taking the
normalization $ T^0 \ra T_{red} \ra T$ of the reduction $T_{red}$
of $T$, and taking the pull-back of the family  to $T^0$.

A less trivial  result which is used is the 
 so-called simultaneous resolution of singularities (cf. \cite{tju},\cite{brieskorn},
\cite{brieskorn2}, \cite{nice})

\begin{theo}\label{simultaneous}
{\bf (Simultaneous resolution according to Brieskorn and Tjurina).}
Let $T : = \CC^{\tau}$ be the basis of the semiuniversal deformation
of a Rational Double Point $ (X,0)$. Then there exists a ramified Galois
cover $ T' \ra T$, with $T'$ smooth $T' \cong  \CC^{\tau}$ such that the pull-back $ \X ' : = \X \times _T T'$
admits a simultaneous resolution of singularities
$ p : \SSS' \ra \X'$ (i.e., $p$ is bimeromorphic,
all the fibres of the composition $  \SSS' \ra \X' \ra T'$
are smooth, and the fibre over  $t' = 0$ is equal  to the minimal resolution
of singularities of $ (X,0)$).
\end{theo}

Another  important observation is  that  the local analytic structure of the Gieseker moduli space is determined by the action of
the group of automorphisms of $X$ on the Kuranishi space of $X$.

\begin{rem}
Let $X$ be the canonical model of a minimal surface of general type $S$ with invariants 
$\chi, K^2$. The isomorphism class of $X$ defines a point  $[X] \in \mathfrak M_{\chi, K^2}^{can}$.

Then the germ of complex space $(\mathfrak M_{\chi, K^2}^{can},{[X]})$ is analytically
isomorphic to the quotient $\mathfrak B (X) / Aut (X)$ 
of the Kuranishi space of $X$ by the finite group $ Aut (X) = Aut (S)$.
\end{rem}

Let $S$ be a minimal surface of general type and let $X$ be its
canonical model. To avoid confusion between the corresponding Kuranishi spaces,
denote by $\Def(S)$ the Kuranishi space for $S$, 
respectively  $\Def(X)$ the Kuranishi space of $X$.

Burns and Wahl  \cite{b-w}, inspired by \cite{atiyahDP} explained the relation holding between $\Def(S)$ and $\Def(X)$.

\begin{theo}{\bf  (Burns - Wahl)}
  Assume that $K_S$ is not ample and let $\pi :S \ra X$ be the
canonical morphism.

   Denote by $\mathcal{L}_X$ the space of local deformations of the
singularities of $X$ (Cartesian product of the corresponding Kuranishi  spaces) and by
$\mathcal{L}_S$ the space of deformations of a neighbourhood of the
exceptional locus  $Exc ( \pi)$ of $\pi$. Then
$\Def(S)$ is realized as the fibre product associated to the Cartesian diagram

\begin{equation*}
\xymatrix{
\Def(S) \ar[d]\ar[r] & \Def (S_{Exc(\pi)})= :  \mathcal{L}_S \cong \CC^{\nu}, \ar[d]^{\lambda} \\
\Def(X) \ar[r] & \Def (X_{Sing X}) = : \mathcal{L}_X \cong \CC^{\nu} ,}
\end{equation*}
where $\nu$ is the number of rational $(-2)$-curves in $S$, and
$\lambda$ is a Galois covering
with Galois group $W := \oplus_{i=1}^r W_i$, the direct sum of the
Weyl groups $W_i$ of the singular points of $X$ (these are generated by reflections, hence yield
a smooth quotient, see \cite{chevalley}).
\end{theo}

An immediate consequence is the following

\begin{cor}{\bf  (Burns - Wahl)}
  1) $\psi:\Def(S) \ra \Def(X)$ is a finite morphism, in
particular, $\psi$ is surjective.

\noindent
2) If the derivative of $\Def(X) \ra \mathcal{L}_X$ is not surjective (i.e., the
singularities of $X$ cannot be independently  
   smoothened  by the first order infinitesimal deformations of $X$),
   then $\Def(S)$ is singular.
\end{cor}

In the next section we shall see  the role played by  automorphisms.

\subsection{Automorphisms and canonical models}

Assume  that $G$ is a group with a faithful action on a complex manifold $Y$:  then $G$ acts naturally on the
sheaves associated to $\Omega^1$, hence we have a linear representation of $G$ on the vector spaces 
that are the cohomology groups of such sheaves: for instance $G$ acts linearly on 
 $H^q (\Omega^p_Y)$, and also on the vector spaces
$ H^0 ( (\Omega^n_Y)^{\otimes m}) = H^0 (\hol_Y (mK_Y))$, hence on the canonical ring 
 $$\sR (Y) : =\sR(Y,K_Y) : =  \bigoplus_{m \geq 0}H^0 (\hol_Y (mK_Y)).$$

 If $Y$ is a variety of general type, then  $G$ acts linearly on the vector space $H^0 (\hol_Y (mK_Y))$,
 hence linearly on the m-th pluricanonical image $Y_m$, which is an algebraic variety bimeromorphic to $Y$.
 Hence $G$ is isomorphic to a subgroup of   the algebraic group $ Aut (Y_m)$.  Matsumura (\cite{matsumura}) used the structure theorem for linear algebraic groups 
 (see \cite{hum})  to show that, if $G$ were infinite,  then $ Aut (Y_m)$ would contain a non trivial Cartan subgroup ($\CC$ or $\CC^*$), hence $Y$ would be uniruled: a contradiction. 
 Hence we have the (already mentioned) theorem:

 \begin{theo}{\bf (Matsumura)}
 The automorphism group of a variety $Y$ of general type is finite.
 
 \end{theo}
 
The above considerations apply now to the $m$-th pluricanonical image $X$ of   a variety  $Y$ of general type.

I.e., we have an embedded variety $X \subset \PP(V)$ and a linear representation of a finite group $G$ on the vector space $V$, such that $X$
is $G$-invariant (for $Y$ of general type one has
 $V := V_m : = H^0 (\hol_Y (mK_Y))$).
 
Now, since we work over $\CC$, the vector space $V$ splits uniquely, up to permutation of the summands,   as a direct sum
of irreducible representations
$$ (**) \ \ V_m = \bigoplus_{\rho \in Irr (G)} W_{\rho}^{n(\rho)}.$$

We come now to the basic notion of a family of $G$-automorphisms (this notion shall be further explained in  the next section). 

\begin{defin}
A family of $G$-automorphisms is a triple $$ ((p \colon \sX \ra T), G, \alpha )$$ where:
\begin{enumerate}
\item
$ (p \colon \sX \ra T)$ is a flat family in a given category (a smooth family for the case of minimal models of surfaces of general type)
\item
$G$ is a (finite) group
\item
$\alpha \colon G \times \sX \ra \sX $ yields a biregular  action  $ G \ra Aut (\sX)$, which is compatible with the projection $p$ and with the trivial action of $G$ on the base $T$ (i.e., $p(\alpha (g,x)) = p(x) , \ \forall g \in G, x \in \sX$).

\end{enumerate}

As a shorthand notation, one may also write $ g(x)$ instead of $\alpha (g,x)$, and by abuse of notation say that the family
of automorphisms is a deformation of the pair $(X_t,G)$ instead of the triple  $(X_t,G, \alpha_t)$.

\end{defin}

The following then holds.

\begin{prop}\label{actiontype}

1) A family of automorphisms of  manifolds of general type (not necessarily minimal models) induces
a family of automorphisms of canonical models.

2) A family of automorphisms of canonical models induces, if the basis $T$ is connected, 
a constant decomposition type $(**)$ for $V_m (t)$. 

3) A  family of automorphisms of  smooth complex manifolds  admits a differentiable trivialization,
i.e., in a neighbourhood of $ t_0 \in T$,  a diffeomorphism as a family with $ ( S_0 \times T, p_T, \alpha_0 \times Id_T)$;
in other words, with the trivial family for which $ g (y,t) = (g(y),t)$.

\end{prop}

We refer to \cite{handbook} for a sketchy  proof,  let us just observe that if we have a continuous family of finite dimensional representations,
the multiplicity of an irreducible summand is a locally constant function on the parameter space $T$ (being given by the scalar product of the
respective trace functions, it is an integer valued continuous function, hence locally constant). 

Let us then consider the case of a family of canonical models of varieties of general type: by 2) above, and shrinking the base in order to
make the addendum $ \sR (p)_m  = p_* (\hol_{\sS} (mK))$ free, we get an embedding of the family
$$ (\XX, G)  \hookrightarrow T \times (\PP ( V_m = \bigoplus_{\rho \in Irr (G)} W_{\rho}^{n(\rho)} ), G).$$
In other words, all the canonical models $X_t$ are contained in a fixed projective space, where also
the action of $G$ is fixed.

Now, the canonical model $X_t$ is left invariant by the action of $G$ if and only if its Hilbert point is
fixed by $G$. Hence, we get a closed subset $  \HHH _0 (\chi, K^2)^G$ of the Hilbert scheme
 $  \HHH ^0 (\chi, K^2)$

$$  \HHH _0 (\chi, K^2)^G: = \{ X | \omega_X \cong \hol_X(1),\ X {\it is \ normal,} \  \ga (X) = X \  \forall \ga \in G\}.$$

For instance in the case of surfaces one has the following theorem (see \cite{vieh}, \cite{lastopus}.
\cite{Kolhand} and references therein for results on moduli spaces of canonically polarized varieties in higher dimension).

\begin{theo}\label{can-aut=closed}
The surfaces of general type which admit an action of a given pluricanonical type $(**)$
i.e., with a fixed irreducible G- decomposition of their canonical ring, form a closed 
subvariety $( \mathfrak M_{\chi, K^2}^{can})^{G, (**)}$ of the moduli space $ \mathfrak M_{\chi, K^2}^{can}$.
\end{theo}

\begin{rem}
The situation for the minimal models of surfaces of general type is different, because then the subset of the moduli space
where one has a fixed differentiable type is not closed, as showed in \cite{burniat3}.

The puzzling phenomenon
which we discovered in  joint work with Ingrid Bauer (\cite{burniat2}, \cite{burniat3},
on the moduli spaces of Burniat surfaces) is that 
deformations of automorphisms differ for canonical and for minimal models.

\end{rem}

More precisely, let $S$ be a minimal surface of general type and let $X$ be its canonical model. Denote by $\Def(S)$, 
resp. $\Def(X)$, the base of the Kuranishi family of $S$, resp. of $X$.

Assume now that we have $1 \neq G \leq \Aut(S) = \Aut(X)$. 

Then we can consider the Kuranishi space of $G$-invariant  deformations of $S$, denoted by $\Def(S,G)$, and respectively consider
 $\Def(X,G)$; we have
a  natural map $\Def(S,G) \ra \Def(X,G)$. 

We indeed show in \cite{burniat3}  that 
this map needs not  be surjective, even if  surjectivity would seem plausible;
we have the following result:  

\begin{theo}\label{path}
 The deformations of nodal Burniat surfaces with $K^2_S =4$ to extended  Burniat surfaces with $K^2_S =4$ 
yield examples where $\Def(S,(\ZZ/2\ZZ)^2) \ra \Def(X,(\ZZ/2\ZZ)^2)$ is not surjective.

Moreover, whereas for the canonical model we have: $$\Def(X,(\ZZ/2\ZZ)^2) = \Def(X),$$ for the minimal models 
we have $$\Def(S,(\ZZ/2\ZZ)^2) \subsetneq \Def(S)$$
and indeed  the subset $\Def(S,(\ZZ/2\ZZ)^2)$ corresponds to the locus of nodal Burniat surfaces.

The moduli space of pairs $(S,G)$  of an  extended (or nodal) smooth projective and minimal Burniat surface $S$ with $K^2_S =4$, taken together with its canonical   $G: = (\ZZ/2\ZZ)^2$-action, 
is disconnected; but its image in the Gieseker moduli space
is a connected open set.
\end{theo}

\subsection{ Kuranishi subspaces for automorphisms of a fixed type}

Proposition \ref{actiontype} is quite useful when one analyses the deformations of a given $G$-action,
say on a compact complex manifold: it tells us that 
we have to look at the complex structures for which the given differentiable action is holomorphic.
Hence we derive (see \cite{montecatini}):

\begin{prop}
Consider a fixed action of a finite group $G$ on a compact complex manifold
 $X$.
Then we obtain a closed subset $Def(X,G)$ of the  Kuranishi space,
corresponding to deformations which preserve the given action, 
and yielding a maximal family of deformations of the $G$-action.

The subset $Def(X,G)$ is the intersection  $ Def(X) \cap H^1(\Theta_X)^G$.
\end{prop}

\begin{rem}
 1)  The proof is based on the well known and already cited Cartan's lemma (\cite{cartan}), that the action of a finite
group in an analytic  neighbourhood of a fixed point can be linearized.
The proof and the result extend also
to encompass compact complex spaces.

2) The above proposition is, as we shall now show in the example of projective curves, quite apt to estimate
the dimension of the subspace $Def(X,G) \subset Def(X)$. For instance, if $Def(X)$ is smooth,
then the dimension of $Def(X,G)$ is just equal to $ dim  (H^1(\Theta_X)^G)$.

\end{rem}

We want to show how to calculate the tangent space to the Kuranishi space of $G$-invariant deformations in some cases.

To do this, assume that $X$ is smooth, and that there is a normal crossing divisor $D = \cup_i D_i$ such that

\begin{enumerate}
\item
$D_i$ is a smooth divisor
\item
the $D_i$ 's intersect transversally
\item
there is a cyclic sugbroup $ G_i \subset G$ such that $D_i = Fix (G_i)$
\item
the stabilizer on the smooth locus of $D_{i_1} \cap D_{i_2} \cap \dots \cap D_{i_k} $ equals
$G_{i_1} \oplus G_{i_2} \oplus  \dots \oplus  G_{i_k} $.

\end{enumerate}

Then, setting $ Y : = X/G$, $Y$ is smooth and the branch divisor $B$ is a normal crossing divisor $B = \cup_i B_i$.

Define as usual the sheaf of logarithmic forms 
$\Omega^1_X ( \log D_1, \dots ,\log D_h)$, as the sheaf generated by $\Omega^1_X$, and by the logarithmic derivatives of the equations
$s_i$ of the divisors $D_i$,  define  $\Theta_Y (- \log D_1, \dots ,- \log D_h)$ as the dual sheaf (this sheaf can also be more generally defined as the sheaf 
of derivations carrying the ideal sheaf of $D$ to itself). We have then the following proposition.

\begin{prop}
Let $p : X \ra Y = X/G$ the quotient projection. Then $ p_* ( \Theta_X)  = \Theta_Y (- \log B_1, \dots ,- \log B_h)$.

In particular, $H^i (\Theta_X ) = H^i (\Theta_Y (- \log B_1, \dots ,- \log B_h)).$
\end{prop}

{\em Idea of proof}

The basic idea of the proof is the calculation one does in dimension $n=1$ (and is the same one to be done at the generic point of
the divisor $D_i$, which maps  to $B_i$). 

Namely, assume that the local quotient map is given by $ w = z^m$, the action being given by $ z \mapsto \e z$,
$\e$ being a primitive $m$-th root of unity. It follows that  $ dz \mapsto \e dz$, and, dually, $(\partial /  \partial z)  \mapsto \e^{-1} (\partial /  \partial z) $.

In particular, $$d log(z) = dz/z = (1/m)  dw/w \ {\it  and} \  z (\partial /  \partial z) =  m w (\partial /  \partial w)$$ are invariant.

Then a vector field 
$$ \theta = f(z) \partial /  \partial z =  f(z)/z (z \partial /  \partial z ) =  f(z)/z m w (\partial /  \partial w)$$ is invariant if and only if
$f(z)/z$ is invariant, hence if and only if $ f (z)/ z = F(w)$, where $F(w)$ is holomorphic.

\qed

For instance, in the case of algebraic curves, where $Def(X)$ is smooth, the dimension of $ Def (X,G)$ equals the dimension 
of $$ H^1 (\Theta_Y (- \log B_1, \dots ,- \log B_d))=  H^1 (\Theta_Y (- y_1, \dots ,- y_d)).$$ 

Let $g'$ be the genus of $Y$, and observe that  this vector space is Serre dual
to $H^0 ( 2 K_Y +  y_1+  \dots + y_d))$.  By  Riemann-Roch its  dimension equals  $3 g' -3 + d$ whenever $3 g' -3 + d \geq 0$
(in fact $ H^0 (\Theta_Y (- y_1, \dots ,- y_d))= 0$  when $2 - 2g' - d < 0$).

Calculations with the above  sheaves are the appropriate ones to calculate the deformations of ramified coverings, see for instance \cite{cat1},\cite{Pardini}, \cite{bidouble},  \cite{burniat2}, and especially  \cite{burniat3}. 

The main problem with the group of automorphisms on minimal models of surfaces (see \cite{handbook}) is that a limit of isomorphisms
need not be an isomorphism
\footnote{ this means:  there exists an algebraic  family of surfaces $ \sS = \cup_{t \in T} S_t $ over a smooth curve $T$, and an algebraic  family of maps $\varphi_t : S_t \ra S_t, t \in T$, 
i.e., such that 
the union of the graphs of the $\varphi_t$ 's  is a closed algebraic set in the fibre product $\sS \times_T \sS$, with the property that
$\varphi_t$ is an isomorphism exactly for $ t \neq t_0$.}, it can be for instance a Dehn twist on a vanishing cycle (see also  \cite{seidelcime}, \cite{cime}).

\section{Moduli spaces of symmetry marked varieties}
\subsection{Moduli marked varieties}

We give now the definition of a symmetry marked variety for projective varieties, but one can similarly give the same definition
for complex or K\"ahler manifolds; to understand the concept of a marking, it suffices to consider a cyclic group acting on a variety $X$.
A marking consists in this case of the choice of a generator for the group acting on $X$. The marking is very important when we
have several actions of a group $G$ on some projective varieties, and we want to consider the diagonal action of $G$
on their product.

\begin{defin}\label{marked}
\begin{enumerate}
\item
A {\em  $G$-marked (projective) variety} is a triple $(X, G, \eta)$ where $X$ is a projective variety, $G$ is a group and $\eta \colon  G \ra \Aut(X)$ is an injective
homomorphism;
\item
equivalently, a marked variety is a triple $(X, G, \alpha)$ where $\alpha \colon X \times G \ra X$ is  a faithful action of the group $G$ on $X$.
\item   Two marked  varieties  $(X,G, \alpha)$, $(X',G, \alpha')$ are said to be {\em isomorphic} if there is an  isomorphism $ f \colon X \ra X'$
 transporting the action $\alpha \colon X \times G \ra X$ into the action $\alpha' \colon X' \times G \ra X'$, 
 i.e., such that $$ f \circ \alpha = \alpha' \circ ( f \times \id) \Leftrightarrow   \eta '   = Ad(f) \circ \eta, \ \ 
 Ad (f) (\phi) : = f \phi f^{-1}.$$
  \item
If $G $ is a subset of $\Aut(X)$, then the natural marked variety is the triple $(X,G,i)$, where $ i \colon G \ra \Aut(X)$
  is the inclusion map, and it shall sometimes be denoted simply by the pair $(X,G)$.
\item  A marked curve  $(D,G, \eta)$ consisting of a smooth projective curve of genus $g$ and a  faithful  action of the group $G$ on $D$
is said to be a {\em marked triangle curve of genus $g$} if $ D / G \cong \PP^1$ and the quotient morphism
 $p \colon D \ra  D / G \cong \PP^1$ is branched in three points.

 \end{enumerate}
\end{defin} 

\begin{rem}
Observe that:

1)  we have a natural action of $\Aut(G)$ on $G$-marked varieties, namely, if $\psi \in Aut (G)$, $$ \psi (X,G,\eta) : = (X, G, \eta \circ \psi^{-1}).$$
 The corresponding equivalence class of a  $G$-marked variety is  defined to be  a {\em  $G$-(unmarked) variety.}

2) the action of the group $\Inn(G)$ of inner automorphisms does not change the isomorphism class of $(X,G,\eta)$ since, for $\ga \in G$,  we may set 
$ f : =  \eta(\ga)$, $\psi : = Ad (\ga)$, and then  $\eta \circ \psi = Ad(f) \circ \eta$, since 
$ \eta ( \psi (g)) = \eta ( \ga  g \ga^{-1} ) = \eta ( \ga) \eta( g ) (\eta(\ga)^{-1} ) = Ad(f) (\eta(g)) $.

3) In the case where $ G = \Aut(X)$, we see that $\Out(G)$ acts simply transitively on the isomorphism classes of
the $\Aut(G)$-orbit of $  (X,G,\eta)$.
\end{rem}

In the spirit of the above concept we have already given, in the previous section, the definition of a family of $G$-automorphisms (we shall also speak of a family of $G$-marked varieties).

The local deformation theory of a $G$-marked variety $X$, at least in the case where the group $G$ is finite,  is simply given by the fixed locus $ \Def(X)^G$ 
of  the natural $G$ -action on  the Kuranishi space $ \Def(X)$. 
As we mentioned previously (see \ref{path}), one  encounters difficulties, when $X$ has dimension at least $2$,  to use the Kuranishi approach for a global theory of moduli of $G$-marked 
 minimal models.

But we have  a  moduli space of $G$-marked varieties in the case of curves of genus $g \geq 2$ and in the case of canonical models of surfaces
of general type, and similarly also for canonical models
in higher dimension. 

In fact, one considers again a fixed linear representation space for the group $G$, the associated projective action,
$$(\PP ( V_m = \bigoplus_{\rho \in Irr (G)} W_{\rho}^{n(\rho)} ), G).$$
and a  locally closed subset of the Hilbert scheme

$  \HHH _{can} (\chi, K^2)^G \subset \HHH ^0 (\chi, K^2)$

$$  \HHH _{can}  (\chi, K^2)^G: = \{ X | \omega_X \cong \hol_X(1),\\  \ga (X) = X \  \forall \ga \in G, $$
$$ \ X {\it is \ normal, \ with \ canonical \ singularities\ }  \}.$$

In this case one divides the above subset  by the subgroup 
$$\sC (G) \subset GL (V_m)$$
which is the centralizer of $G \subset  GL (V_m)$.

\begin{rem}
1) Assume in fact that there is an isomorphism of the marked varieties $(X,G)$ and $(X', G)$
(here the marking is furnished by an m-th  pluricanonical embedding of $X, X' \hookrightarrow \PP ( V_m)$, and by  the fixed action 
$ \rho : G  \hookrightarrow  GL (V_m)$).

Then there is an isomorphism $ f : X \cong X'$ with $ f \circ \ga = \ga \circ f$ $\forall \ga \in G$. Now, $f$ induces a linear map of $V_m$
which we denote by the same symbol: hence we can write the previous condition as 
$$ f \circ \ga \circ f^{-1} = \ga , \ \forall  \ga \in G \Leftrightarrow f \in \sC (G) .$$

2) By Schur's lemma  $\sC (G) =  \Pi_{\rho \in Irr (G)} GL (n(\rho) , \CC)$ is a reductive group.

3) Observe that, since the dimension of $V_m$ is determined by the holomorphic invariants of $X$, then there is only a finite number of
possible representation types for the action of $G$ on $V_m$.

4) Assume  moreover  that $X, X'$ are isomorphic varieties: then there is such a linear isomorphism $ f : V_m \ra V_m$ sending $X$ to $X'$.

The variety $X'$ has therefore two $G$-markings, the marking $\eta' : G \ra Aut (X')$ provided by the pluricanonical embeddding
$X' \hookrightarrow \PP ( V_m)$, and another one given, in  an evident notation,  by $ Ad (f) \circ \eta$.

Therefore $Aut (X')$ contains two subgroups, $G : = Im (\eta')$, and $G ' : = Im ( Ad (f) \circ \eta)$. $G = G'$ if and only if $f$ lies in the normalizer $\sN_G$ of $G$.

5) This happens if $G$ is a full subgroup (this means that $G = Aut (X')$).  In general
$f$ will just belong to  the normalizer $\sN_{X'}$ of $Aut (X')$, and $\sN_G$ shall just be
 the finite index subgroup of $\sN_{X'}$ which stabilises, via conjugation,  the subgroup $G$.

6) In other words,  we define the marked moduli space (of a fixed representation type) as the quotient 

$$ \frak M[G] : =   \HHH _{can}  (\chi, K^2)^G / \sC_G .$$ 

 Then the forgetful map from the marked moduli space to the moduli space  $\frak M[G] \ra \frak M$, which factors through the quotient $ \frak M(G)$ by
$\sN_G / \sC_G$, yields  finite maps  $\frak M[G] \ra \frak M(G) \ra \frak M^G$.

We omit here technical details in higher dimension, which have been treated   by Binru Li   in his Bayreuth Ph. D. Thesis (in preparation).

Observe that the action of $ Out (G)$ does not need to respect the representation type.

\end{rem} 

Let us see now  how the picture works in the case of curves: this case is    already  very enlightening  and intriguing.

\subsection{Moduli of curves with automorphisms}

There are several   `moduli spaces' of curves with automorphisms.
First of all, given a finite group $G$, we define a   subset $\mathfrak M_{g,G}$ of 
the moduli space $\mathfrak M_g$ of smooth curves of genus $g>1$:   $\mathfrak M_{g,G}$ is the locus of the curves that admit an effective action 
by the group $G$. It turns out that  $\mathfrak M_{g,G}$ is a Zariski closed algebraic subset.

 In order to understand   the irreducible components of 
 $\mathfrak M_{g,G}$ we have seen that Teichm\"uller theory plays an important role: it shows  the connectedness, given an injective   homomorphism 
$\rho \colon G \to Map_g$, of the locus 
$$  \sT_{g, \rho}  : = Fix ( \rho (G)).$$ 

Its image  $\mathfrak M_{g, \rho}$ in   $\mathfrak M_{g,G}$ is a Zariski closed irreducible subset (as observed in \cite{CLP2}). 
Recall that to a curve $C$ of genus $g$ with an action by $G$ we can associate several discrete 
invariants that are constant under deformation.

The first is the above  \emph{topological type}
of the $G$-action: it is a homomorphism 
$\rho \colon G \to Map_g$, which is
well-defined up to inner conjugation (induced by
different choices of an isomorphism  $Map (C)  \cong  Map_g$).

We  immediately see that the locus  $\mathfrak M_{g, \rho}$ is first of all  determined by the subgroup  $\rho(G)$
and not by the marking. Moreover, this locus remains the same not only  if we change $\rho$  modulo the action by ${\rm Aut}(G)$, but also if we change
 $\rho$  by the adjoint action by ${\rm Map}_g$.
 
\smallskip

\begin{defin}

1) The moduli space of G-marked curves of a certain topological type $\rho$  is the quotient of the Teichm\"uller submanifold
$\sT_{g, \rho} $ by the centralizer  subgroup $\sC_{\rho(G)}$ of the subgroup $\rho(G)$ of the mapping class group. 
We get a normal complex space which we shall denote  $\frak M_g[\rho]$.
$\frak M_g[\rho] = \sT_{g, \rho} / \sC_{\rho(G)}$ is a  finite covering of a Zariski closed subset of  the usual moduli space (its image $\frak M_{g, \rho}$), therefore it is quasi-projective, by 
the theorem of Grauert and Remmert. 

2) Defining  $\frak M_g(\rho)$ as the quotient of $\sT_{g, \rho} $ by the normalizer  $\sN_{\rho(G)}$ of $\rho(G)$, we call it the 
moduli space of curves with a G-action of a given   topological type. It is again a normal quasi-projective variety.
\end{defin}

\begin{rem}
1) If we consider $G'  : = \rho(G)$ as a subgroup $ G' \subset Map_g$, then we get a natural $G'$-marking for any $C \in Fix(G') = \sT_{g, \rho} $.

2) As we said,  $ Fix(G') = \sT_{g, \rho} $ is independent of the chosen marking, moreover the projection $ Fix(G') = \sT_{g, \rho} \ra \frak M_{g, \rho}$
factors through a finite map $\frak M_g(\rho) \ra \frak M_{g, \rho}$. 
\end{rem}

The next question is whether $\frak M_g(\rho)$ maps 1-1 into the moduli space of curves. This is not the case, as we shall easily see. 
Hence we give the following definition.

\begin{defin}
Let  $ G \subset Map_g$ be a finite group, and let $C $ represent a point in $Fix(G)$. Then we have a natural inclusion $G \ra A_C : = Aut (C)$,
and $C$  is a fixed point for the subgroup $A_C \subset Map_g$:
$A_C$  is indeed  the stabilizer of the point $C$ in $Map_g$, so that locally (at the point of $ \frak M_g$ corresponding to $C$) we get  
a complex analytic isomorphism $ \frak M_g =  \sT_g / A_C. $

We define $H_G : = \cap_{ C \in Fix(G)} A_C$ and we shall  say that $G$ is a {\bf  full subgroup} if $ G = H_G$.
Equivalently, $H_G$ is the largest subgroup $H$ such that $ Fix(H) = Fix (G)$.

This implies that $H_G$ is a full subgroup.
\end{defin}

\begin{rem}
The above definition shows that the map  $ \sT_{g, \rho} \ra \frak M_{g, \rho}$ factors through   $\frak M_g (H_G) : = \sT_{g, \rho} / \sN_{H_G}$, hence we restrict our attention only to full subgroups.

\end{rem}

\begin{prop}
If $H$ is a full subgroup  $ H  \subset Map_g$, and $\rho : H \subset Map_g$ is the inclusion homomorphism,
then $ \frak M_g (\rho)  $ is the normalization of $\frak M_{g, \rho} $.

\end{prop}

\Proof
Recall that $ \frak M_g (H)  = Fix (H) / \sN_H$, and that for a general curve $C \in Fix (H)$ we have $ A_C = H$.

Assume then that $ C = \ga (C')$, and that both $C, C'  \in Fix (H)$. Then $ \ga A_C \ga^{-1} = A_C$, 
hence for $C$ general  $\ga \in \sN_H$.

This means that the finite surjective morphism $\varphi :  \frak M_g (H)  \ra \frak M_{g, H} $ is generically bijective; $ \frak M_g (H)  $
is a normal variety, being locally the quotient of a smooth variety by a finite group: hence  $\varphi $  is the normalization 
morphism and $ \frak M_g (H)  $ is the normalization of $\frak M_{g, H} $.

\qed

Next we investigate when the natural morphism $ \frak M_g (H)  \ra \frak M_{g, H} $ is not injective.
In order that this be the case, we have already seen that we must pick up a curve $C$ such that $ H \subset A_C \neq H$.
Now, let $C'$  be another point of $ Fix (H) $ which has the same image in the moduli space: this means that $C'$
 is in the orbit of $C$ so there is an element  $\ga \in Map_g$ carrying $C'$ to $C$.
 
 By conjugation $\ga$  sends $H$ to another subgroup $H'$ of $A_C$. We can assume that $ H \neq H'$, else
 $\ga$ lies in the normalizer $N_H$ and we have the same point in $ \frak M_g (H) $. Hence the points of $ \frak M_g (H) $  that have the same image as $C$ correspond to the subset of subgroups
  $H' \subset A_C$ which  are conjugate to $H$ by the action of some $ \ga \in Map_g$.

  \begin{ex}\label{notfull1}
  Consider the genus $2$ curve $C$ birational to the affine curve  with equation $ y^2 = (x^6-1)$. 
  Its canonical ring is generated by $x_0, x_1, y$ and is the quotient $\CC[x_0, x_1, y] / (y^2 - x_1^6 + x_0^6)$.
  
  Its group of automorphisms is generated (as a group of projective transformations)  by $ a (x_0, x_1, y) : = (x_0, \e x_1, y)$ where $\e$ is a primitive sixth root of $1$,
  and by $ b (x_0, x_1, y) : = (x_1,  x_0, i y)$. $a$ has order $6$, $b$ has order $4$, the square $b^2$ is the hyperelliptic involution  $h (x_0, x_1, y) : = (x_0,  x_1, -y)$.
  
  We have $$ b \circ a \circ b^{-1} = a^{-1}, $$ a formula which implies the (known) fact that the hyperelliptic involution lies in the centre of the group 
  $ A_C = Aut (C)$ .
  
  Taking as $H$ the cyclic subgroup generated by $b$, the space $ Fix (H)$ has dimension equal to $1$ since $ C \ra C/H = \PP^1$ is branched in $4$ points.
  The quotient of $ A_C = Aut (C)$ by $b^2$ is the dihedral group $D_6$, and we see that $H$ is conjugate to $6$ different subgroups of $A_C$.
  \end{ex}

 Of course an important question in order to understand the locus in $ \frak M_g $ of curves with automorphisms  is the determination of all the non full subgroups $G \neq H_G$.
 
 For instance Cornalba (\cite{cornalba}) answered this question for  cyclic groups of prime order, thereby obtaining a full determination of the irreducible components
 of $ Sing (\frak M_g) $. In fact, for $ g \geq 4$, the locus $ Sing (\frak M_g) $ is the locus of curves admitting a nontrivial automorphism, so this locus is the locus
 of curves admitting a nontrivial automorphism of prime order. In this case, as we shall see in the next subsection, the topological type is easily determined.
 We omit to state Cornalba's result in detail: it amounts in fact to a list  of all non full such subgroups of prime order. 
 We limit ourselves to  indicate the simplest example of such a situation.
 
   \begin{ex}\label{notfull2}
  Consider the genus  $ g = \frac{p-1}{2}$ curve $C$, birational to the affine curve  with equation $ z^p  = (x^2-1)$, where $p$ is an odd prime number. Letting $H$ be the cyclic group
  generated by the automorphisms $ (x,z) \mapsto (x, \e z)$, where $\e$ is a primitive $p$-th root of unity, we see that the quotient morphism
  $ C \ra C / H \cong \PP^1$ corresponds to the inclusion of fields $ \CC(x) \hookrightarrow \CC(C)$. The quotient map is branched on the three points $ x = 1, -1 , \infty$,
  so that $ Fix (H)$ consists of just a point (the above curve $C$).
  
  An easy inspection of the above equation shows that the curve $C$ is hyperelliptic, with  involution $ h : (x,z) \mapsto (-x,  z)$ which is hyperelliptic since
  $ \CC(C)^h = \CC (z)$. In this case the locus $ Fix(H)$ is contained in the hyperelliptic locus, and therefore it is not an irreducible component of $ Sing (\frak M_g) $.
  
   \end{ex}

\subsection{Numerical and homological invariants of group actions on curves}

\medskip

As  already mentioned, given an effective 
 action of a finite group $G$ on $C$, we set $C':=C/G$, $g':=g(C')$, and we have the quotient  morphism $p\colon C \to C/G=: C'$, a $G$-cover.
 
 The geometry of $p$
encodes several numerical invariants that are constant on $M_{g,\rho}(G)$: first of all 
the genus $g'$ of $C'$, then the  number $d$ of branch points $y_1, \dots , y_d \in C'$.

We call the set $B=\{ y_1, \ldots , y_d\}$ the branch locus, and for each $y_i$ we denote by $m_i$ the multiplicity
of $y_i$ (the greatest number dividing the divisor $p^{*} (y_i)$). We choose an ordering of $B$ such that 
$ m_1 \leq  \dots \leq  m_d$. 

These numerical invariants  $g', d, m_1 \leq  \dots \leq  m_d$  form the
so-called  {\bf primary numerical type}.

  $p \colon C\to C'$ is determined (Riemann's existence theorem)  by the monodromy, a surjective homomorphism:
$$
\mu \colon \pi_1(C'\setminus B) \to G \, .
$$
We have:
$$ \pi_1(C' \setminus B) \cong \Pi_{g', d}:= \langle \ga_1, \ldots , \ga_d, \al_1, \be_1 , \ldots , \al_{g'}, \be_{g'} | \prod_{i=1}^d \ga_i \prod_{j=1}^{g'} [\al_j, \be_j]=1 \rangle .$$

We set then  $c_i : = \mu (\ga_i)$, $a_j : = \mu (\al_j)$, $b_j : = \mu (\be_j)$, thus obtaining 
a Hurwitz generating vector, i.e. a vector 
$$
v : = (c_1, \ldots , c_d, a_1, b_1, \ldots , a_{g'}, b_{g'}) \in G^{d+2g'}
$$
s.t. 
\begin{itemize}
\item
$G$ is generated by the entries $c_1, \ldots , c_d, a_1, b_1, \ldots , a_{g'}, b_{g'}$,
\item
$c_i \not= 1_G$, $\forall i$, and 
\item
$\prod_{i=1}^d c_i \prod_{j=1}^{g'}[a_j,b_j]=1$.
\end{itemize}

We see that the monodromy $\mu$ is completely equivalent, once an isomorphism $ \pi_1(C' \setminus B) \cong \Pi_{g', d}$ is chosen,
to the datum of a Hurwitz generating vector (we also call  the sequence $c_1, \ldots ,c_d, a_1, b_1, \ldots , a_{g'}, b_{g'}$
of the vector's  coordinates a {\em Hurwitz generating system}).

A second numerical invariant of these components of $\mathfrak M_{g} (G)$ is obtained from the monodromy
{$\mu \colon \pi_1(C'\setminus\{ y_1, \dots , y_d\})\to G$}  of the restriction of $p$
to $p^{-1}(C'\setminus \{ y_1, \dots , y_d \})$,
and is called the \emph{$\nu$-type} or Nielsen function of the covering.

The Nielsen  function $\nu$ is a function defined on the set of conjugacy classes in $G$
which, for each conjugacy class $\mathcal{C}$ in $G$, counts the
number $\nu (\sC)$ of local monodromies $c_1, \dots, c_d$ which belong to $\mathcal{C}$ (observe that the numbers 
 $m_1 \leq  \dots \leq  m_d$ are just the orders of the local monodromies).  

Observe in fact that  the generators $\ga_j$ are well defined only  up to conjugation in the group $ \pi_1(C'\setminus\{ y_1, \dots , y_d\})$, hence
the local monodromies are well defined only  up to conjugation in the group $G$.

We have already observed that  the  irreducible closed algebraic sets $M_{g,\rho}(G)$ 
depend only upon what we call the `unmarked topological type', which is defined
as the conjugacy class of the subgroup  $\rho (G)$  inside $Map_g$. This concept remains however still mysterious,
due to the complicated nature of the group $Map_g$. Therefore one tries to use more geometry to get a grasp on the topological type.

The following is immediate
by Riemann's existence theorem and the irreducibility of
the moduli space $\mathfrak M_{g',d}$ of $d$-pointed curves of genus
$g'$.
Given $g'$ and $d$, the unmarked topological types
whose primary numerical type is  of the form $g', d, m_1, \dots , m_d$  are in bijection with
the quotient of the set of the corresponding  monodromies $\mu$ modulo 
the actions by $Aut(G)$ and by $Map(g',d)$. \\ 
Here $Map(g',d)$
is the full mapping class group of genus $g'$ and $d$ unordered  points.

Thus 
Riemann's existence theorem shows that  the components of the moduli space $$\mathfrak M (G): = \cup_g \mathfrak M_g (G)$$
 with numerical invariants  $g', d$ correspond to the 
following quotient set.

\begin{defin}
$$
\sA (g',d,G) : = {\rm Epi}(\Pi_{g',d}, G)/{\rm Map}_{g',d}\times {\rm Aut}(G) \, .
$$
 \end{defin}
 
 Thus  a first step toward the general problem 
consists in finding  a fine  invariant that distinguishes these orbits.

In the paper \cite{CLP2} we introduced a new homological invariant $\hat{\e}$ for $G$-actions on smooth curves
 (and showed that, in  the case where $G$ is the dihedral group $D_n$ of order $2n$,
$\hat{\e}$ is a fine invariant since it distinguishes the different  unmarked topological types).

This invariant generalizes the classical homological invariant in the unramified case.

\begin{defin}
Let $p\colon C \to C/G=: C'$ be unramified, so that $d=0$ and we have a monodromy $\mu \colon \pi_1(C')\to G$.

Since $C'$ is a classifying space for the group $\pi_{g'}$, we obtain a continuous map 
$$ m : C' \ra  BG, \  \pi_1(m) = \mu.$$

Moreover, $H_2(C', \ZZ)$ has a natural generator $[C']$, the fundamental class of $C'$ determined by the orientation
induced by the complex structure of $C'$. 

The homological invariant of the $G$-marked  action is then defined as:

$$\e : =  H_2(m) ([C']) \in H_2(BG, \ZZ) = H_2(G, \ZZ).$$

If we forget the marking we have to take $\e$ as an element in $ H_2(G, \ZZ)/ Aut(G).$

\end{defin}
 
\begin{prop}
Assume that a finite group $G$ has a fixed point free action on a curve of genus $g \geq3$.

 Let $p\colon C \to C/G=: C'$ be the quotient map and
pick an isomorphism $ \pi_1(C') \cong  \pi_{g'}$. Let $\mu : \pi_{g'} \ra G$ be the surjection corresponding to the monodromy of $p$,
and denote by $a_i : = \mu (\al_i), b_j : = \mu (\be_j), \ i,j = 1, \dots, g'$.

Assume that we have a realization $ G = F /R$ of the group $G$ as the quotient of a free group $F$, and that $\widehat{a_i}, \widehat{b_j}$
are lifts of $a_i, b_j$ to $F$.

Then the homological invariant $\e$ of the covering is the image of 
$$ \e' = \prod_1^{g'}[\widehat{a_i}, \widehat{b_i}]   \in [F,F] \cap R $$ 
into $ H_2(G, \ZZ) = ([F,F] \cap R) / [F,R]$.

\end{prop}

\Proof
$C'$ is obtained attaching a 2-cell $D$ to a bouquet of circles, with boundary $\partial D$ mapping to $\prod_1^{g'}[\al_i, \be_i]$.

Similarly the 2-skeleton $BG^2$ of $BG$ is obtained from a bouquet of circles, with fundamental group $\cong F$, attaching 2-cells
according to the relations in $R$.  Since $\prod_1^{g'}[{a_i}, {b_i}] = 1$ in $\pi_{g'}$, the relation $\prod_1^{g'}[\widehat{a_i}, \widehat{b_i}]$ is a product of elements of $R$.
In this way $m$ is defined on $D$, and the image of the fundamental class of $C'$ is the image of $D$, which is a sum of 2-cells whose 
boundary is exactly the loop $\e'$. If we write $\e$ as a sum of 2-cells, we get an element in $ H_2(BG^2, BG^1,\ZZ) $: but since $\e$
is a product of commutators, the boundary of the corresponding 2-chain is indeed zero, so we get an element in $ H_2(BG, \ZZ) = H_2(G, \ZZ) $.

\qed

\subsection{The refined homology invariant in the ramified case}

Assume now that $ p : C \ra C/G$ is  ramified. Then we define $H$ to be the 
minimal
normal subgroup of $G$ generated by the local monodromies $c_1, \dots, c_d$ $\Leftrightarrow$ $H$ is the (normal) subgroup 
generated by the transformations $\ga \in G$ which have some fixed point.

We have therefore a factorization of $p$ 

$$ C \ra C'' : = C/H \ra  C' : = C/G $$
where $ C'' \ra C'$ is an unramified $G/H$-cover, and where $ C \ra C''$ is totally ramified. 

The refined homology  invariant includes and extends two  invariants that have been studied in the literature, and were already mentioned:
the $\nu$-type  (or Nielsen type) of the cover  (also called shape in \cite{FV})
and
the class in the second homology group $H_2(G/H, \ZZ)$
 corresponding to the unramified cover
$p' \colon C''= C/H \to C'$.

The construction of the invariant $\hat{\e}$ is similar in spirit to the procedure used in the unramified case.
But we achieve a little less   in the `branched' case of a non-free action. In this case we are only able to associate, to
two given 
actions with the same $\nu$-type, an invariant in a quotient group of $H_2(G,\ZZ)$ which is the `difference' of the respective $\hat{\e}$- invariants. 
  Here is  the way we do it.

Let $\Sigma$ be the Riemann surface (with boundary) obtained from $C'$ after removing  open discs $\De_i$ around each of the branch points.

Take generators  $\ga_1, \dots , \ga_d$ formed by a simple path 
going from the base point  $y_0$ to a point  $z_i$ on the circle bounding the open discs $\De_i$, and by the circle $\partial \De_i$.

Fix then a CW-decomposition of $\Sigma$ as follows. The $0$-skeleton $\Sigma^0$ consists of  the  point $y_0$ and 
of the points $z_i$, $i=1, \dots , d$.  
The $1$-skeleton $\Sigma^1$ is given by  the geometric basis $\al_1, \dots, \be_{g'}, \ga_1, \dots , \ga_d$ and the $2$-skeleton $\Sigma^2$ consists of one cell.

The restriction $p_\Sigma$ of $p\colon C \to C'$ to $p^{-1}(\Sigma)$ is an unramified $G$-covering of $\Sigma$ and hence corresponds to a continuous map
$Bp_\Sigma \colon \Sigma \to BG$, well defined up to homotopy. Let $Bp_1 \colon \Sigma \to BG$ be a cellular approximation of $Bp_\Sigma$.
Since $Bp_1$ can be regarded as a map of pairs $Bp_1 \colon (\Sigma, \partial \Sigma)\to (BG,BG^1)$, the push-forward of the fundamental 
(orientation) class $[\Sigma, \partial \Sigma]$ gives an element
$$
{Bp_1}_*  [\Sigma, \partial \Sigma] \in H_2(BG,BG^1)=\frac{R}{[F,R]} \, .
$$
This element depends on the chosen cellular approximation $Bp_1$ of $Bp_\Sigma$, but its image in a quotient group which we denote by $K_\Gamma$
does not depend  on the chosen cellular approximation, as shown in \cite{CLP3}.
The main idea to define this quotient group  is that a homotopy between two 1-cellular approximations
can also be made cellular, hence the difference in relative homology is the sum of boundaries of cylinders; and if we have a cylinder $\sC$ with upper circle $a$,
lower circle $c$ (with the same orientation), and meridian $b$, then the boundary of the corresponding 2-cell is
$$\partial \sC =  a b c^{-1} b^{-1} .$$

\begin{defin}
1) For any finite group $G$, let $F$ be the free group
generated by the elements of $G$ and let $R\trianglelefteq F$ be the subgroup of relations, that is $G=F/R$.

Denote by $\hat{g}\in F$ the generator corresponding to $g\in G$.

For any $\Ga \subset G$, union of non trivial conjugacy classes, 
let $G_\Ga$ be the quotient group of $F$
by the  minimal normal subgroup $R_{\Ga}$ generated by $[F,R]$ and  by the elements
$\hat{a}\hat{b}\hat{c}^{-1}\hat{b}^{-1}\in F$, for any  $a, c \in \Ga$, $b\in G$,
such that  $b^{-1} ab=c$. 

2) Define instead  $$K_\Ga = R / R_{\Ga} \subset G_\Ga = F / R_{\Ga},$$
a central subgroup of $G_\Ga $.
We have thus a central extension

$$  1 \ra K_\Ga  \ra G_\Ga \ra G \ra 1.$$ 

3) Define
$$
H_{2,\Ga}(G)=\ker \left( G_\Ga \to G\times G_\Ga^{ab} \right) \, .
$$

Notice that 
$$
H_2(G,\ZZ)\cong \frac{R\cap [F,F]}{[F,R]} \cong \ker \left( \frac{F}{[F,R]} \to G\times G^{ab}_\emptyset \right) \, .
$$
\end{defin}

\begin{rem}
In particular, when $\Ga = \emptyset$, $H_{2,\Ga}(G)\cong H_2(G, \ZZ)$.

By \cite[Lemma 3.12]{CLP2} we have that the morphism
$$
R\cap [F,F] \to \frac{R}{R_\Ga} \, , \quad r \mapsto rR_\Ga
$$
induces a surjective group homomorphism 
$$ H_2(G,\ZZ) \to H_{2, \Ga}(G) .$$

\end{rem}

Now, to a given $G$-cover $p\colon C \to C'$ we associate the set $\Ga$ of { the local monodromies, i.e., of} the elements of $G$ which can be geometrically described as those
which  i) stabilize some point $x$ of $C$
 and ii) act on the tangent space at $x$ by a rotation  of angle
$\frac{2\pi}{m}$ where $m$ is the order of the stabilizer at $x$. 

In terms of the notation that we have previously  introduced  $\Ga =\Ga_v$ is simply  the union of the conjugacy classes of the $c_i$'s. 

\begin{defin}
The tautological lift $\hat{v}$ of $v$ is the vector:
$$(\widehat{c_1}, \dots , \widehat{c_d}; \widehat{a_1}, \widehat{b_1}, \dots , \widehat{a_{g'}}, \widehat{b_{g'}}).$$
Finally, define $\e(v)$ as the class in $K_{\Ga}$ of 
$$
\prod_1^d\widehat{c_j} \cdot \prod_1^{g'}[\widehat{a_i}, \widehat{b_i}] \, .
$$

\end{defin}

It turns out that  the image of $\e(v)$ in $K_{\Ga} /_{Inn(G)}$ is invariant under the action of $Map(g',d)$,
as shown in  \cite{CLP2}. 

Moreover the $\nu$-type of a Hurwitz monodromy vector $v$ can be deduced from $\e(v)$, as it is essentially the image of $\e(v)$ in the abelianized group 
$G_{\Ga}^{ab}$.

\begin{prop}
Let $v$  be a Hurwitz generating vector and let $\Ga_v\subset G$ be the union of the conjugacy classes of the $c_j$'s , $j\leq d$.
The abelianization $G_{\Ga_v}^{ab}$ of $G_{\Ga_v}$ can be described as follows:
$$
G_{\Ga_v}^{ab}\cong (\bigoplus_{\mathcal{C}\subset \Ga_v}\ZZ\langle \mathcal{C} \rangle) \oplus ( \bigoplus_{g\in G\setminus \Ga_v}\ZZ \langle g \rangle \, ),
$$
where $\mathcal{C}$ denotes a conjugacy class of $G$.  

Moreover the Nielsen function $\nu(v)$ coincides with the vector whose $\mathcal{C}$-components
are the corresponding components of  the image 
of $\e(v)\in G_{\Ga_v}$ in $G_{\Ga_v}^{ab}$.

 \end{prop}

\subsection{Genus stabilization of components of moduli spaces of curves with $G$-symmetry}

In order to take into account also the automorphism group $Aut(G)$, one has to consider 
$$
K^\cup : = {\coprod}_\Ga K_{\Ga}\, ,
$$
the disjoint union of all the $K_{\Ga}$'s. Now, the group $Aut(G)$ acts on $K^\cup$ and we get a map
$$
\hat{\e}\colon \sA (g',d,G)  \to (K^\cup)/_{Aut(G)} \, 
$$
which is induced by $v\mapsto \e(v)$. 

Next, one has to observe that the Nielsen functions of coverings have to satisfy a necessary condition,
consequence of the relation $$\prod_{i=1}^d c_i \prod_{j=1}^{g'}[a_j,b_j]=1.$$

\begin{defin}
An element 
$$
\nu = (n_{\mathcal{C}})_\mathcal{C} \in \bigoplus_{\mathcal{C}\not= \{ 1 \}}\NN\langle \mathcal{C}\rangle 
$$
is {\bf admissible} if the following equality holds in the $\ZZ$-module $G^{ab}$:
$$
\sum_\mathcal{C} n_{\mathcal{C}} \cdot [\mathcal{C}] =0 \, 
$$
(here $[\sC]$ denotes the image element of $\sC$ in the abelianization $G^{ab}$).
\end{defin}

The main result of \cite{CLP3} is the following `genus stabilization' theorem.

\begin{theo}\label{genusstab}
There is an integer $h$ such that for  $g' > h$ 
$$
\hat{\e}\colon \sA (g',d,G)  \to (K^\cup)/_{Aut(G)} \, 
$$
induces a bijection onto the set of admissible classes of refined homology invariants.

In particular, if $g' > h$, and we have two Hurwitz generating systems $v_1, v_2$ having the same Nielsen function,
they are equivalent if and only if the `difference'  $ \hat{\e} (v_1) \hat{\e} (v_2) ^{-1} \in H_{2, \Ga}(G) $ is trivial.

\end{theo}

The above result extends a nice theorem of Livingston, Dunfield and Thurston (\cite{Liv}, \cite{DT})  in the unramifed case, where also the statement is simpler.

\begin{theo}
For $g' > > 0$ $$
\hat{\e}\colon \sA (g',0, G)  \to H_2(G, \ZZ) /_{Aut(G)} \, 
$$
is a bijection.
\end{theo}

\begin{rem}
Unfortunately the integer $h$ in theorem \ref{genusstab}, which depends on the group $G$, is not explicit. 

\end{rem}

A key concept used in the proof is the concept of genus stabilization of a covering,
which we now briefly explain.

\begin{defin}

Consider a group action of $G$ on a projective curve $C$, and let $ C \ra C' = C/G$
the quotient morphism, with monodromy 
$$
\mu \colon \pi_1(C'\setminus B) \to G \, 
$$
(here $B$ is as usual the branch locus).
Then the first genus stabilization of the differentiable covering is defined geometrically by  simply adding
 a handle to the curve $C'$, on which the covering is trivial.
 
 Algebraically, given the monodromy homomorphism
$$ \mu : \pi_1(C' \setminus B) \cong \Pi_{g', d}:= \langle \ga_1, \ldots , \ga_d, \al_1, \be_1 , \ldots , \al_{g'}, \be_{g'} | \prod_{i=1}^d \ga_i \prod_{j=1}^{g'} [\al_j, \be_j]=1 \rangle \rightarrow G ,$$
we simply extend $\mu$ to $\mu^1 \colon  \Pi_{g'+1, d} \ra G $ setting
$$ \mu^1 (\al_{g' + 1 }) = \mu^1 ( \be_{g' + 1} ) = 1_G .$$ 

In terms of Hurwitz vectors and
 Hurwitz generating systems, we replace the vector 
$$
v : = (c_1, \ldots , c_d, a_1, b_1, \ldots , a_{g'}, b_{g'}) \in G^{d+2g'}
$$
by 
$$
v^1 : = (c_1, \ldots , c_d, a_1, b_1, \ldots , a_{g'}, b_{g'}, 1,1) \in G^{d+2g'+ 2}.
$$

The operation of first genus stabilization generates then an equivalence relation among monodromies (equivalently, Hurwitz generating systems),
called {\bf stable equivalence}.

\end{defin}

The most important step in the proof, the geometric understanding of the invariant $\e \in H_2(G, \ZZ)$ was obtained by Livingston \cite{Liv}.

\begin{theo}
Two monodromies $\mu_1$, $\mu_2$ are stably equivalent if and only if they have the same invariant $\e \in H_2(G, \ZZ)$.

\end{theo}

A purely algebraic proof of Livingston's theorem was given by Zimmermann in \cite{Zimmer}, while
a nice sketch of proof was given by Dunfield and Thurston in \cite{DT}.

{\em Idea of proof}

Since one direction is clear (the map to $BG$ is homotopically trivial on the handle that one adds to $C'$),
one has to show that two coverings are stably equivalent  if their invariant is the same.

The first idea is then to interpret second homology as bordism: given two maps of two curves $C'_1, C'_2 \ra BG$,
they have the same invariant in $H_2(G, \ZZ)$ (image of the fundamental classes $[C'_1],[ C'_2 ]$)
if and only if there is a 3-manifold $W$ with boundary $ \partial W = C'_1 - C'_2$, and a continuous map $ f \colon W \ra BG$
which extends the two maps defined on the boundary $ \partial W = C'_1 - C'_2$.

Assume now that there is relative Morse function for $(W, \partial W )$, and that one first adds all the 1-handles (for the critical points 
of negativity 1) and then all the 2-handles.

Assume that at level $t$ we have a curve $C'_t$, to which we add a 1-handle. Then the genus of $C'_t$ grows by 1,
and the monodromy is trivial on the meridian $\al$ (image of the generator of the fundamental group of the cylinder we are attaching).

Pick now a simple loop $\be'$ meeting $\al$ transversally in one point with intersection number 1: since the monodromy 
$\mu$ is surjective, there is a simple loop $\ga$ disjoint from $\al$ and $\be'$ with $\mu (\ga) = \mu (\be')$. Replacing $\be'$ by $\be := \be' \ga^{-1}$,
we obtain that the monodromy is trivial on the two simple loops $\al$ and $\be$. A suitable neighbourhood
of $\al \cup \be$ is then a handle that is added to $C'_t$, and with trivial monodromy. This shows that at each critical value of the
Morse function we pass from one monodromy to a stably equivalent one (for the case of a 2-handle, 
repeat the same argument replacing  the Morse function $F$
by its opposite $-F$).

\qed

\subsection{Classification results for certain concrete groups}

The first result in this direction was obtained by Nielsen (\cite{nielsen}) who proved that 
$\nu$ determines $\rho$ if $G$ is cyclic (in fact in this case $H_2(G, \ZZ) = 0$!).

In the cyclic case the Nielsen function for $G =  \ZZ/n$ is simply a function $\nu : (\ZZ/n)  \setminus \{0\} \ra \NN$,
and admissibility here simply means that 
$$  \sum_i   i  \cdot \nu(i) \equiv 0 \ ( mod \ n).$$ 

The class of $\nu$ is just the equivalence class  for  the equivalence relation  $\nu (i) \sim \nu_r ( i)$, $\forall r \in (\ZZ/n)^*$, 
 where $\nu_r (i) : =  \nu ( r i)$, $\forall i \in (\ZZ/n).$

From the refined  Nielsen realization theorem of  \cite{isogenous} (\ref{refinedNR})
it follows that the components of $\mathfrak M_g (\ZZ/n)$ are in bijection with the classes of Nielsen functions
(see also  \cite{singMg} for an elementary proof).

\begin{ex}\label{triple}

For instance, in the case $n=3$,  the components of $\mathfrak M_g (\ZZ/3)$ correspond to triples  of integers $ g', a , b \in \NN$
such that $ a \equiv b  \ ( mod \ 3),$ and $g-1 = 3 (g'-1) + a + b$ ($a : = \nu (1), b : = \nu(2)$).

For $g'=0$ we have, if we assume $ a \leq b$ (we can do this by changing the generator of $\ZZ/3$), $b = a + 3 r, r \in \NN$,
and $ 2a + 3r = g+2$, thus $ a \equiv -(g-1) \ (mod \ 3)$, and $ 2a \leq g+2$.

In the case where $g'=0$,  two such components $N_a$ ( labelled by $a$ as above) 
have images which do not intersect in $\mathfrak M_g$ as soon as $ g \geq 5$.

Otherwise we would  have two $\ZZ/3$-quotient morphisms $ p_1, p_2  : C \ra \PP^1$ and a birational map $C \ra C'' \subset \PP^1 \times \PP^1$.
Then $C$ would be the normalization of a curve of arithmetic genus $4$, so $ g \leq 4$.

\end{ex}

The genus $g'$ and the Nielsen class (which refine the primary numerical type), and the homological invariant 
 $h \in H_2 (G/H,\ZZ)$ (here $H$ is again the subgroup generated by the local monodromies)
determine the connected components of  $\mathfrak M_g (G)$ under some restrictions: for instance when 
$G$ is abelian  or when $G$ acts freely and is the semi-direct product of two finite cyclic groups 
(as it follows by combining results from  \cite{isogenous}, \cite{singMg},  \cite{Edm I} and  \cite{edmonds2}).

\begin{theo}{\bf (Edmonds)}
$\nu$ and $h \in H_2 (G/H,\ZZ)$ determine $\rho$ for $G$ abelian.
Moreover, if $G$ is split-metacyclic and the action is free, then $h$
determines $\rho$.
\end{theo}

However, in general, these invariants  are not enough to distinguish unmarked topological types,
as one can see already for non-free $D_n$-actions (see \cite{CLP2}). 
Already for dihedral groups, one needs the refined homological invariant $\hat{\e}$.

\begin{theo}[\cite{CLP2}]
For the dihedral group  $G=D_n$ the connected components of the moduli space $\mathfrak M_g (D_n)$
are in bijection, via the  map $\hat{\e}$, with the admissible classes of refined homology invariants.
\end{theo}

\begin{rem}
In this case, the classification is simple: two monodromies with the same Nielsen function differ by an element in $H_{2, \Ga}(D_n)$,
which is a quotient of $H_{2}(D_n, \ZZ)$. This last group is $0$ if $n$ is odd, and $\cong \ZZ/2$ for $n$ even.

More precisely: $H_{2,\Ga}(D_n)=\{0\}$ if and only if
\begin{itemize}
\item
$n$ is odd or
\item
$n$ is even and $\Ga $  contains some reflection or
\item
 $n$ is even and  $\Ga $  contains 
the non-trivial central element.
\end{itemize}
 In the remaining cases, $H_{2,\Ga}(D_n)=\ZZ/2\ZZ$.

\end{rem}

The above result completes the classification of the 
unmarked topological types for $G = D_n$, begun in \cite{CLP1}; moreover this result  entails the classification of the irreducible components of the loci $\mathfrak M_{g,D_n}$ (see the appendix to \cite{CLP2}). 

It is an interesting question: for which groups $G$ does the refined homology invariant  $\hat{\e}$ determine the connected components of 
$\mathfrak M_g (G)$?

In view of Edmonds' result in the unramified case, it is reasonable  to expect
 a positive answer for split metacyclic groups (work in progress by Sascha Weigl) or for some more general metacyclic or metaabelian groups. 

As mentioned in \cite{DT}, page 499, the group $G = \PP SL(2, \F_{13})$ shows that, for $g'=2$, in the unramified case there are different
components with trivial  homology
invariant $\e \in H_2(G,\ZZ)$: these topological types of coverings are therefore stably equivalent but not equivalent.
\medskip

 \subsection{$Sing(\mathfrak M_g)$ II:   loci  of curves with automorphisms in $\mathfrak M_g$}
 
 Several authors independently found restrictions in order that a finite subgroup $H $ of $Map_g$ be a not  full subgroup
 (Singerman, \cite{singer}, Ries \cite{ries}, and Magaard-Shaska-Shpectorov-V\"olklein,  \cite{mssv})
 
 We refer to lemma 4.1 of \cite{mssv} for the proof of the following result. 

\begin{theo}\label{MSSV} {\bf (MSSV)} 
 Suppose $H \subset  G \subset Map_g$ and assume $Z : = Fix (H) = Fix (G) \subset \sT_g$, with $H$ a proper subgroup of $G$, and let $C\in Z$.
Then
$$
\de : = dim (Z)  \leq 3.
$$

I) if  $\de = 3$, then $H$ has index $2$ in $G$, and $C \ra C/G$ is covering of $\PP^1$ branched on six points, $P_1, \dots, P_6$,
and with branching indices all equal to $2$.  Moreover the subgroup $H$ corresponds to the unique genus two double cover of $\PP^1$ 
branched on the six points, $P_1, \dots, P_6$ (by Galois theory, intermediate covers correspond to subgroups of $G$ bijectively).

II)  If $\de = 2$, then $H$ has index $2$ in $G$, and $C \ra C/G$ is covering of $\PP^1$ branched on five points, $P_1, \dots, P_5$,
and with branching indices  { $2,2,2,2,c_5$} (where obviously $c_5 \geq 2$).  Moreover the subgroup $H$ corresponds to a genus one double cover of $\PP^1$ 
branched on {  four  of the points  $P_1, \dots,  P_5$ which have branching index $2$}.

III)  If $\de = 1$, then there are three possibilities.

III-a) $H$ has index $2$ in $G$, and $C \ra C/G$ is covering of $\PP^1$ branched on four points $P_1, \dots, P_4$,
 with branching indices $2,2,2, 2 d_4$, where $d_4 > 1$.  Moreover the subgroup $H$ corresponds to the unique genus one double cover of $\PP^1$ 
branched on the four points, $P_1, \dots, P_4$.

III-b) $H$ has index $2$ in $G$, and $C \ra C/G$ is covering of $\PP^1$ branched on four points, $P_1, \dots, P_4$,
 with branching indices $2,2,c_3, c_4$, where $2 \leq c_3 \leq c_4 > 2$.  Moreover the subgroup $H$ corresponds to a genus zero double cover of $\PP^1$ 
branched on two points whose branching index equals 2.

III-c) $H$ is normal  in $G$, $ G/H \cong (\ZZ/2)^2$, moreover  $C \ra C/G$ is covering of $\PP^1$ branched on four points, $P_1, \dots, P_4$,
 with branching indices $2,2,2, c_4$, where $ c_4 > 2$.  Moreover the subgroup $H$ corresponds to the unique  genus zero  cover of $\PP^1$ 
with group $(\ZZ/2)^2$ branched on the three  points $P_1, P_2, P_3$ whose branching index equals 2.

\end{theo}

The main point in the theorem above is a calculation via the Hurwitz formula, showing that the normalizer of $H$ in $G$ is nontrivial; 
from this follows the classification of the singular locus of $\mathfrak M_g$, due to Cornalba, and which we do not reproduce here (see  \cite{cornalba},
and also \cite{singMg} for a slightly different proof).

\medskip

\subsection{Stable curves and their automorphisms, $Sing (\overline{\mathfrak M_g}$)}
The compactification of the moduli space of curves $\mathfrak M_g$ is given by the moduli space $\overline{\mathfrak M_g}$ of stable curves 
of genus $g$ (see \cite{d-m}, \cite{EnsMath}, \cite{Giescime}).

\begin{defin}
A stable curve of genus $g \geq 2$ is a reduced and connected projective curve $C$, not necessarily irreducible, whose singularities are only nodes,
such that

1) the dualizing sheaf $\omega_C$ has degree $2g-2$

2) its group of automorphisms $Aut(C)$ is finite (equivalently, each smooth irreducible component $D \subset C$  of genus zero intersects 
$ C \setminus D$ in at least 3 points). 

\end{defin}

Again the Kuranishi family $ Def(C)$ of a stable curve is smooth of dimension $3g-3$, and the locus of stable curves which admit a given
$G$-action forms a smooth submanifold of $Def(C)$. 

Hence again the singular locus of $\overline{\mathfrak M_g}$ corresponds to loci of curves with automorphisms which do not form a divisor.

\begin{ex}
We say that $C$ has an elliptic tail if we can write $ C = C_1 \cup E$, where $E$ is a smooth  curve of genus 1 intersecting $C_1$ in only one point $P$.

In this case $C$ has always an automorphism $\ga$ of order $2$: $\ga$ is the identity on $C_1$, and, if we choose $P$ as the origin of the elliptic 
curve $(E,P)$, $\ga$ is simply multiplication by $-1$.

An easy calculation shows that such curves with an elliptic tail form a divisor inside $\overline{\mathfrak M_g}$.

\end{ex}

We refer to \cite{singMg} (the second part) for the explicit determination of the loci of stable curves admitting an action by a cyclic group
of prime order, especially those contained in the boundary $\partial \overline{\mathfrak M_g}$. From this follows the description of the singular locus
$Sing (\overline{\mathfrak M_g})$, which we do not reproduce here.

An interesting question is to describe similar  loci also for other groups.

As mentioned, there are irreducible loci which are  contained in the boundary, while other ones  are the closures of irreducible  loci in $\mathfrak M_g$.
But it is possible that some loci which are disjoint in ${\mathfrak M_g}$ have closures that meet in $\overline{\mathfrak M_g}$.

We end this section with an explicit example, due to Antonio F. Costa, Milagros Izquierdo, and Hugo Parlier (cf. \cite{CMP}), but we give a different proof.
We refer to the notation of example \ref{triple} in the following theorem.

\begin{theo}
The strata $N_a$ and $N_{3+a}$ of $\mathfrak M_g(\ZZ/3)$ fulfill:

i)  $N_a \cap N_{3+a} = \emptyset$

ii)  $\overline{N_a} \cap \overline{N_{3+a}} \neq  \emptyset$ in  $\overline{\mathfrak M_g}$.
\end{theo}

\Proof
 By the fact that inside Kuranishi space the locus of curves with a given action of a group $G$ is a smooth manifold, hence locally irreducible,
one sees right away that these strata must intersect in a point corresponding to a curve with at least two different automorphisms of order 3. A natural choice is
to take a stable curve with an action of $(\ZZ/3)^2$.

The simplest choice is to take first  the Fermat cubic $$F : = \{ (x_0, x_1, x_2) | x_0^3 +  x_1^3 +  x_2^3 = 0\}$$ with the action of $(u,v) \in (\ZZ/3)^2$
$$ (x_0, x_1, x_2) \mapsto (x_0, \e^u  x_1, \e^v x_2).$$

We let then $C$ be a $(\ZZ/3)$-covering of the projective line, of genus $g-3$ and with $ \nu(1) = a $. 
Denote by $\psi$ the covering automorphism and  take a point $p \in C$ which is not a ramification point:
then  glue the three points $ \psi^v (p), \ v \in \ZZ/3,$ with the three points $ (1,0, - \e^v)$. In this way we get a stable curve $D$ of genus $g$,
and with an action of $ (\ZZ/3)^2$.

We have two actions of  $\ZZ/3$ on $D$,  the first induced by the action of $(0,1)$ on $F$ and of $\psi$ on $C$, and the second 
induced by the action of $(1,1)$ on $F$ and $\psi$ on $C$
(it is immediate to observe that the glueing is respected by the action).

Now, for the first action of the generator $\bar{1} \in \ZZ/3$ on $F$, the fixed points are the points $x_2=0$, with local coordinate $ t = \frac{x_2}{x_0}$, such that $ t \mapsto \e t$,
while for the second action  the fixed points are the points $x_0=0$, with local coordinate $ t = \frac{x_0}{x_1}$, such that $ t \mapsto \e^2  t$.

Hence the same curve has two automorphisms of order 3, the first with $\nu(1) = a + 3$, the second with $\nu(1)  = a $. As easily shown (cf. \cite{singMg}),
since, for any of these two actions of $\ZZ/3$,  the nodes are not fixed, then there exists a smoothing for both actions of $\ZZ/3$ on  the curve $D$.
This shows that the strata $N_a$ and $N_{a+3}$ have $D$ in their closure.

\qed

\subsection{Branch stabilization and relation with other approaches}

When $g'=0$ our $G_\Ga$ is related to the group $\widehat{G}$ defined in \cite{FV} (Appendix),
where the authors give a proof of a theorem by Conway and Parker. Roughly speaking the theorem says that:
if the Schur multiplier $M(G)$ (which is isomorphic to $H_2(G,\ZZ$))
 is generated by commutators, then the $\nu$-type is a fine stable invariant, when $g'=0$.\\

\begin{theo}{\bf (Conway -Parker, loc.cit.)}
In the  case $g'=0$, let $G = F / R$ where $F$ is free, and  assume that $  H_2(G,\ZZ) \cong \frac{[F,F]\cap R}{[F,R]}$ is generated by commutators. 
Then  there is an integer $N$ such that if the numerical function $\nu$ takes values $\geq N$,
then there is only one   equivalence class  with the given  numerical function $\nu$.
\end{theo}

Michael L\"onne, Fabio P{erroni and myself have been able to extend this result in the following way:

\begin{theo}
Assume that $g'=0$, or more generally  that the union $\Ga$ of nontrivial conjugacy classes of $G$ generates $G$.

Then there exists an integer $N$, depending on $\Ga$, such that if two Hurwitz generating systems $v,w$ satisfy 
$ \nu(v) \geq N \chi_{\Ga}, \nu(w) \geq N \chi_{\Ga}$, where  $\chi_{\Ga}$ is the characteristic function of the set of conjugacy classes
corresponding to $\Ga$, the Nielsen functions $\nu(v)$ and $\nu(w)$ are equivalent $\Leftrightarrow$ $v$ is equivalent to $w$ and they yield the same point in $\sA (g',d, G)$.
\end{theo}

It is still an open question how to extend the previous result without the condition that $\Ga$ generates $G$.

Similar results have been obtained by Viktor Kulikov and Slava Kharlamov (\cite{k-k2}), who use geometric arguments for the construction
of semigroups similar to the ones constructed by Fried and V\"olklein, while 
a purely  algebraic construction of a group similar to  our group $K_{\Ga}$ can be found in work of Moravec on unramified Brauer groups
(\cite{moravec}).

\bigskip

\subsection{Miller's description of the second homology of a group and developments}

Clair Miller found (\cite{miller}) another nice description of the second homology group $H_2(G, \ZZ)$, as follows.

\begin{defin}
Let $\langle G, G \rangle$ be the free group on all pairs $\langle x, y \rangle$ with $x,y \in G$.

Then there is a natural surjection of $\langle G, G \rangle$ onto the commutator subgroup $[G, G]$ sending 
$\langle x, y \rangle$ to the commutator $[x,y]$.  

Denote as in \cite{miller} by $ Z (G)$ the kernel of this surjection
(it might have been better to denote it by $ Z (\langle G, G \rangle)$), and then denote by $B (G)$
(it might have been  better to denote it by $ B (\langle G, G \rangle)$) the normal subgroup of $ Z (G)$ normally generated by the following elements

\begin{enumerate}
\item
$\langle x, y \rangle$ 
\item
$\langle x, y \rangle \langle y, x \rangle$ 

\item 
$\langle z, xy \rangle \langle y^x, z^x \rangle \langle x, z \rangle$, where $ y^x : = x   y x^{-1}$,
\item
$\langle z^x, y^x \rangle \langle x, [y,z] \rangle  \langle y, z \rangle$.

\end{enumerate}

\end{defin}

\begin{theo} {\bf (Miller)}
There is  a canonical isomorphism $ Z(G) / B(G) \cong H_2(G, \ZZ)$.

\end{theo}

Clearly we have then an exact sequence

$$ 0 \ra    H_2(G, \ZZ)  \cong Z(G) / B(G) \ra \langle G, G \rangle \ra [G, G] \ra 0$$
which corresponds to the previously seen
 $$ 0 \ra    H_2(G, \ZZ) = (R \cap [F,F] ) / [F,R]  \ra[F,F] / [F,R] \ra [G, G] \ra 0.$$
 
 The relations by Miller were later slightly modified by Moravec (\cite{moravec}) in a more symmetric fashion as follows:
 
 \begin{enumerate}
\item
$ \langle g,g  \rangle \sim 1$
\item
$ \langle g_1 g_2, h  \rangle \sim \langle g_2^{g_1}, h^{ g_1}   \rangle  \langle g_1 , h \rangle  $
\item
$ \langle g, h_1 h_2   \rangle \sim \langle g , h_1  \rangle \langle g^{h_1}, h_2^{ h_1}   \rangle  . $

\end{enumerate}

In the case where we consider also a finite union $\Ga$ of conjugacy classes, Moravec defined the following 
group 

\begin{defin} {\bf (Moravec)}
$$ G \wedge_{\Ga} G = \langle G, G \rangle / B_{\Ga}, $$
where $B_{\Ga}$ is defined by the previous relation (1), (2), (3) and by the further relation
$$ (4) \  \langle g,k  \rangle \sim 1, \ \forall k \in \Ga.$$

\end{defin}

The definition by Moravec and the one given in \cite{CLP3} are related, as the following easy proposition shows.

\begin{prop}
There is an exact sequence 

$$ 0 \ra    H_{2, \Ga} (G, \ZZ)  \ra G \wedge_{\Ga} G  \ra [G, G] \ra 0.$$

\end{prop}

Another  important source of interest for these quotient groups of $H_2(G, \ZZ)$ comes from rationality questions.

In general,  the stable cohomology of a finite group $G$, or more generally of an algebraic group $G$, is
obtained in the following way.

Let $V$ be a finite dimensional representation of $G$, and let $U$ be an open set of $V$ where $G$ acts freely:
then $U/G$ is considered as an algebraic approximation to a classifying space for $G$, and one can take
the limit over such representations ($V$ and $V'$ begin smaller than their direct sum) of the cohomology groups
of $U/G$. The same can of course be done also for homology groups and Chow groups (see  \cite{Totaro}).

The groups that are thus obtained are important to study the problem of rationality of the quotients $V/G$,
or of their stable rationality ($X$ is said to be stably rational if there is an integer $n$ such that $X \times \PP^n$ is rational).

We refer to  \cite{Bog}, \cite{moravec}, \cite{Bogstable}, \cite{BarPetr} and the literature cited there for more details.

\medskip

\section{Connected components of moduli spaces and the action of the absolute Galois group}

Let $X$ be a complex projective variety: let us quickly recall the 
notion of a conjugate variety.

\begin{rem}
1) $\phi \in Aut (\CC)$ acts on $ \CC [z_0, \dots z_n]$,
by sending $P (z) = \sum_{i =0}^n \ a_i z ^i
\mapsto  \phi (P) (z) : = \sum_{i =0}^n \ \phi (a_i) z ^i$.

2) Let $X$ be as above a projective variety
$$X  \subset  \PP^n_\CC,  X : = \{ z | f_i(z) = 0 \ \forall i \}.$$

The action of $\phi$ extends coordinatewise to $ \PP^n_\CC$,
and carries $X$ to another variety, denoted $X^{\phi}$,
and called the {\bf conjugate variety}. Since $f_i(z) = 0 $ implies
$\phi (f_i)(\phi (z) )= 0 $, we see that

$$  X^{\phi}  = \{ w | \phi (f_i)(w) = 0 \  \forall i \}.$$

\end{rem}

If $\phi$ is complex conjugation, then it is clear that the variety
$X^{\phi}$ that we obtain is diffeomorphic  to $X$; but, in general,
what happens when $\phi$ is not continuous ?

Observe that, by the theorem of Steiniz,  one has a surjection $ Aut (\CC) \ra Gal(\bar{\QQ} /\QQ)$,
and by specialization the heart of the question concerns the action of $Gal(\bar{\QQ} /\QQ)$
 on varieties $X$ defined over $\bar{\QQ}$.

For curves, since in general the dimensions of spaces of
differential forms of a fixed degree and without poles are the same
for $X^{\phi}$ and $X$, we shall obtain a curve of the same genus,
hence $X^{\phi}$ and $X$ are diffeomorphic.

\subsection{Galois conjugates of projective classifying spaces}

General questions of which the first  is answered in the positive in most concrete cases, 
and the second is answered in the negative in many cases, as we shall see, are the following.

\begin{question}
Assume that $X$ is a projective $ K (\pi, 1)$, and assume $\phi \in Aut (\CC)$.

A) Is then the conjugate variety $  X^{\phi} $ still a classifying space $ K (\pi', 1)$?

B) Is then $\pi_1 ( X^{\phi} ) \cong \pi \cong \pi_1 (X)$ ? 

\end{question}

Since $\phi$ is never continuous, there would be no reason to expect a positive answer to both questions A) and B),
except that Grothendieck   showed (\cite{sga1}, see also \cite{groth2}).

\begin{theo}
Conjugate varieties $ X,  X^{\phi} $ have isomorphic algebraic fundamental groups 
$$ \pi_1(X)^{alg} \cong  \pi_1(X^{\phi} )^{alg} ,$$
where $\pi_1(X)^{alg}$ is the profinite completion
of the topological fundamental group $\pi_1(X)$.

\end{theo}

We recall once more  that the profinite completion 
 of a group $G$ is the inverse limit
$$ {\hat G} = lim_{K\unlhd_f G}  (G/K),$$
of the factor groups $G/K$, $K$ being a normal subgroup of finite index in $G$; and since finite index subgroups of the fundamental group
correspond to finite unramified (\'etale) covers, Grothendieck  \cite{sga1} defined in this way the algebraic fundamental group
for varieties over other fields than the complex numbers, and also for more general schemes.

The main point of the proof of the above theorem  is that if we have $ f : Y \ra X$ which is \'etale, also the Galois conjugate
$f^{\phi} : Y^{\phi}  \ra X^{\phi} $ is \'etale ($f^{\phi} $ is just defined taking the Galois conjugate of the graph of $f$,
a subvariety of $ Y \times X$).

Since Galois conjugation gives an isomorphism of natural cohomology groups, which respects the cup product, 
as for instance the Dolbeault cohomology groups $ H^p(\Omega^q_X)$, we obtain interesting consequences in the direction of question A) above.
Recall the following definition.

\begin{defin}
Two varieties $X,Y$ are said to be isogenous if there exist a third variety $Z$, and \'etale finite morphisms $ f_X : Z \ra X$,
$ f_Y : Z \ra Y$.

\end{defin}

\begin{rem}
It is obvious that if $X$ is isogenous to $Y$, then $X^{\phi} $ is isogenous to $Y^{\phi} $.

\end{rem}

\begin{theo}
i) If $X$ is an Abelian variety, or isogenous to an Abelian variety, the same holds for any Galois conjugate  $X^{\phi} $.

ii) If $S$ is a Kodaira fibred surface, then any  Galois conjugate $S^{\phi} $ is also Kodaira fibred.

iii) If $X$ is isogenous to a product of curves, the same holds for any Galois conjugate  $X^{\phi} $.

\end{theo}

\Proof
i) $X$ is an Abelian variety if and only it is a projective variety and there is a morphism $X \times X \ra X$,
$(x,y) \mapsto (x \cdot y^{-1})$, which makes $X$ a group (see \cite{abvar}, it follows indeed that the group is commutative).
This property holds for $X$ if and only if it holds for $X^{\phi} $.

ii) The hypothesis is that there is $ f : S \ra B$ such that all the fibres are smooth and not all isomorphic: obviously the same property
holds, after Galois conjugation, for $f^{\phi} : S^{\phi}  \ra B^{\phi} $.

iii) It suffices to show  that the Galois conjugate of a product of curves is a product of curves. But since $X^{\phi}  \times  Y^{\phi} = (X \times Y) ^{\phi}  $
and the Galois conjugate of a curve $C$ of genus  $g$ is again a curve of the same genus $g$, the statement follows.

\qed

Proceeding with other projective $K(\pi, 1)$'s, the question becomes more subtle and we have to appeal to a famous theorem
by Kazhdan on arithmetic varieties (see \cite{Kazh70}, \cite{Kazh83},  \cite{milne}, \cite{C-DS1}, \cite{C-DS2}, \cite{V-Z}).

\begin{theo} Assume that $X$ is a projective manifold with
$K_X$  ample, and that
the universal covering $\tilde{X}$ is a bounded symmetric domain.

   Let $\tau \in
\mathrm{Aut}(\mathbb{C})$ be an automorphism of $\mathbb{C}$.

Then the conjugate variety  $X^{\tau}$ has universal covering
$\tilde{X^{\tau}} \cong \tilde{X}$.
\end{theo}

Simpler proofs follow from recent results  obtained together with Antonio Di Scala, and  based on  the Aubin-Yau theorem  and the results of Berger \cite{Berger}. These results yield a precise characterization of varieties possessing a  bounded symmetric
domain as universal cover, and can be rather useful in view of the fact that our knowledge and classification of these fundamental groups 
is not so explicit.

To state them in detail would require some space, hence  we just mention the simplest result (see \cite{C-DS1}).

\begin{theo}

Let $X$ be a compact complex manifold of dimension $n$ with  $K_X$ ample.

Then the
following two conditions (1) and (1'), resp. (2) and (2') are equivalent:

\begin{itemize}

\item[(1)] $X$ admits a slope zero   tensor $0 \neq  \psi   \in
H^0(S^{mn}(\Omega^1_X)(-m K_X) )$, (for some  positive
integer $m$);

\item[(1')]$X \cong \Omega / \Gamma$ , where
$\Omega$ is a bounded symmetric domain of tube type and $\Gamma$ is a
cocompact
discrete subgroup of
$\mathrm{Aut}(\Omega)$ acting freely.
\item[(2)] $X$ admits a semi special tensor $0 \neq \phi  \in
H^0(S^n(\Omega^1_X)(-K_X) \otimes \eta)$, where $\eta$ is a 2-torsion invertible sheaf,
 such that there is a  point
$p\in X$ for which the
corresponding hypersurface $F_p : = 
\{\phi_p = 0 \}\subset \PP (TX_p)$ is  reduced.
\item[(2')]
The universal cover of $X$ is a polydisk.
\end{itemize}

Moreover, in case (1),  the degrees and the multiplicities of the irreducible
factors of the polynomial  $\psi_p$ determine uniquely the universal
covering
$\widetilde{X}=\Omega$.

\end{theo}

\subsection{Connected components of  Gieseker's moduli space}
For the sake of simplicity we shall describe in this and the next subsection
 the action of the absolute Galois group on the set of connected components of the moduli space
of surfaces of general type. We first recall the situation concerning these components.

As we saw, all 5-canonical models of surfaces of general type with invariants $ K^2$, $\chi$ occur
in a big family parametrized by an open set of the Hilbert scheme $ \HHH ^0$ parametrizing 
subschemes with Hilbert polynomial $ P (m)= 
\chi  + \frac{1}{2}( 5 m -1)  5 m  K^2 ,$ namely the open set

$$ \HHH ^0 (\chi, K^2) : = \{ \Sigma | \Sigma {\rm \ is\ reduced \ with \ only \
R.D.P. 's \ as \
singularities \ }\}.
$$

 We shall however, for the sake of brevity, talk about connected components $\sN$ of the Gieseker moduli space $\frak M_{a,b}$
 even if these do not really parametrize families of canonical models.
 
We refer to \cite{perugia} for a more ample discussion of the basic ideas which we are going to sketch here.

$\frak M_{a,b}$ has a finite number of connected components, and 
these parametrize the 
deformation classes
of surfaces of general type with numerical invariants $\chi (S) = a, K^2_S = b$. By the classical theorem of Ehresmann 
(\cite{ehre}),
deformation equivalent varieties are diffeomorphic, and moreover, via a diffeomorphism carrying the canonical class 
to the canonical class.

Hence, fixed the two numerical invariants $\chi (S) = a, K^2_S = b$,
which are determined
by the topology of $S$ (indeed, by the  Betti numbers $b_i(S)$  of $S$ and by $b^+:=$ positivity of the intersection
form on $H^2(S, \RR)$),
we have a finite number of
differentiable types.

For some time the following question was open: whether two surfaces which are orientedly diffeomorphic would belong  to the same 
connected component of
the moduli space.

I conjectured (in \cite{katata}) that the answer should be negative, on the basis of some families of simply connected
surfaces of general type constructed in \cite{cat1}  and later investigated in \cite{cat2}, \cite{c-w}and \cite{clw}: these were shown to be homeomorphic by the results of Freedman (see \cite{free}, and \cite{f-q}),
and it was then relatively easy to show then (\cite{cat3}) that there were many connected components of the moduli space
corresponding to homeomorphic but non diffeomorphic surfaces. It looked like the situation should be similar
even if one would fix the diffeomorphism type.

Friedman and Morgan instead made 
the `speculation' that the answer to the DEF= DIFF question should be positive (1987) (see \cite{f-m1}), motivated by the new examples of homeomorphic 
but not diffeomorphic surfaces discovered by Donaldson (see \cite{don} and \cite{don5} for a survey on this topic).

The question was finally answered in the negative, and  in every
possible way (\cite{man4},\cite{k-k},\cite{cat03},\cite{c-w},\cite{bcg}).

\begin{theo}
(Manetti '98, Kharlamov -Kulikov 2001, C. 2001, C. - Wajnryb 2004,
Bauer- C. - Grunewald 2005 )

The Friedman- Morgan speculation
does not hold true and the  DEF = DIFF question has a negative answer.
\end{theo}

In my joint work with
Bronislaw Wajnryb (\cite{c-w}) the  DEF = DIFF question was  shown to have a negative answer also for simply connected surfaces
(indeed for some of the families of surfaces constructed in \cite{cat1}).

I refer to  \cite{cime} for a  rather comprehensive treatment of the above questions (and to \cite{a-k, abkp, cat02, cat09, clw} for the symplectic point of view,
\cite{Dolgachev, survey} for the special case of geometric genus $p_g=0$).

 \subsection{Arithmetic of moduli spaces and faithful actions of the absolute Galois group}

A basic remark is that all the schemes involved in the construction of the Gieseker moduli space are defined by equations involving only $\ZZ$-coefficients,
since the defining equation of the Hilbert scheme is a rank condition for a multiplication map
(see for instance \cite{green}), and similarly the condition  $ \omega_{\Sigma}^{
\otimes 5} \cong
   \hol_{\Sigma}(1) $ is also closed (see \cite{abvar}) and defined over $\ZZ$.
   
It follows that  the absolute Galois group  $ Gal (\overline{\QQ}, \QQ)$ acts on the Gieseker
moduli space $\frak M_{a,b}$. In particular, it acts on the set of its irreducible components, and
on the set of its connected components.

After an incomplete initial attempt in \cite{AbsoluteGalois} in joint work with Ingrid Bauer and Fritz Grunewald,
we were able in \cite{bcg2} to show:

\begin{theo}\label{faithful}

The absolute Galois group $Gal(\bar{\QQ} /\QQ)$
acts faithfully
on the  set of connected 
components of the Gieseker moduli space of  surfaces of general type,
$$ \frak M : = \cup_{x,y \in \NN, x,y  \geq 1} \frak M_{x,y}. $$

\end{theo}

Another  result in a similar direction had  been obtained by Easton and Vakil (\cite{east-vak}) using abelian coverings 
of the plane branched on union of lines.

\begin{theo}

The absolute Galois group $Gal(\bar{\QQ} /\QQ)$
acts faithfully
on the  set of irreducible 
components of the Gieseker moduli space of  surfaces of general type,
$$ \frak M : = \cup_{x,y \in \NN, x,y  \geq 1} \frak M_{x,y}. $$

\end{theo}

The main ingredients for the proof of theorem \ref{faithful} are the following ones.

\begin{enumerate}

\item
Define,
 for any complex number $ a \in \CC \setminus \{ -2g,0,1, \ldots , 2g-1 \} $, $C_a$  as the hyperelliptic
curve   of genus $g \geq 3$  which is the smooth complete model of the affine curve of equation
$$ w^2 = (z-a) (z + 2g) \Pi_{i=0}^{2g-1} (z-i) .$$ 

 Consider then two complex numbers  $a,b$ such that $a \in  \CC 
\setminus \QQ $: then $C_a
\cong C_b$ if and only if $a = b$.
\item
If $ a \in \bar{\QQ} $, then by Belyi's theorem (\cite{belyi}) there is a morphism $f_a : C_a \ra \PP^1$
which is branched only on three points, $0,1,\infty$.
\item
The normal closure $D_a$ of $f_a$ yields a triangle curve, i.e., a curve $D_a$ with the action of
a finite group $G_a$ such that $D_a / G_a \cong \PP^1$, and $ D_a \ra \PP^1$ is branched only on three points.
\item
Take surfaces isogenous to a product $ S = (D_a \times D')/G_a$ where the action of $G_a$ on $D'$ is free.
Denote by $\sN_a$ the union of connected components parametrizing such surfaces.
\item
Take all the possible twists of the $G_a$-action on $D_a \times D'$ via an automorphism $\psi \in Aut (G_a)$
(i.e., given  the action $(x,y) \mapsto (\ga x, \ga y)$, consider all the actions of the form
 $$(x,y) \mapsto (\ga x, \psi (\ga) y).$$
 One observes that, for each $\psi$ as above, we get more connected components in $\sN_a$.
\item 
Find by an explicit calculation (using (4) and (5))  that the subgroup of $Gal(\bar{\QQ} /\QQ)$ acting trivially on the set of connected components of the moduli space
would be a normal and abelian subgroup.
\item
Finally, this contradicts a known  theorem  (cf. \cite{F-J}).

\end{enumerate}
 
\begin{rem}
An interesting remark of a referee is that the meaning of the above theorem could be further elaborated as a parallel of the Drinfend theory of
Galois representation, built on the theorem of Belyi  (see \cite{drinfeld}). We leave this as a future task.
\end{rem}

\subsection{Change of fundamental group}

 Jean  Pierre Serre proved  in the 60's (\cite{serre})
 the existence of a field automorphism  $\phi \in
Gal(\bar{\QQ} /\QQ)$, and a variety $X$ defined over $\bar{\QQ}$ such that
$X$ and the Galois conjugate variety $X^{\phi}$  have
  non isomorphic fundamental groups.
  
  In \cite{bcg2} this phenomenon is vastly generalized, thus answering question B) in the negative.
  
  \begin{theo} \label{fundamentalgroup}  If $\sigma \in Gal(\bar{\QQ} /\QQ)$ is not in the conjugacy class of complex conjugation,
then there exists a surface isogenous to a product $X$ such that $X$ and the Galois conjugate surface $X^{\sigma}$  have non-isomorphic fundamental groups. 
\end{theo}

Since the argument for the above theorem is not constructive, let us observe that,
in  work   in collaboration with Ingrid Bauer and Fritz
Grunewald (\cite{almeria}, \cite{bcg2}), we discovered  wide classes of explicit algebraic surfaces
isogenous to a product 
for which the same phenomenon holds.

By the strong rigidity of locally symmetric spaces $ X = \sD / \Ga$ whose universal covering $\sD $ is
an  irreducible bounded symmetric domain of dimension $\geq 2$,  similar phenomena should also occur in this case.

\begin{rem}
Further developments have been announced  in \cite{gabino} by Gonzal\'ez-Diez and Jaikin-Zapirain: for instance the faithfulness of the action of the absolute Galois group
on the discrete set of the moduli space corresponding to Beauville surfaces, and the extension of theorem \ref{fundamentalgroup} to all 
automorphisms $\sigma$ different from
complex conjugation.

\end{rem}
\bigskip

\section{Stabilization results for the homology of moduli spaces of curves and Abelian varieties}

We have seen that the moduli space of curves is a rational $K(\pi,1)$,
since it can be written as 
a quotient of  the Teichm\"uller space $\mathcal{T}_g $ of a closed oriented real 2 -manifold $M$ of genus $g$
$$ 
\mathfrak{M}_g = \mathcal{T}_g / Map_g , \ \sT_g := \mathcal{C}S(M)/{\rm Diff}^0(M), \ 
$$
As a corollary, and as we saw,  the rational cohomology of the moduli space is calculated by group cohomology:
$$
H^*(\mathfrak{M}_g, \QQ) \cong H^*({\rm Map}_g, \QQ) \, .
$$

Harer showed, using the concept of genus stabilization that we have already introduced in section 10,
 that these cohomology groups stabilize with $g$; indeed, stabilization furnishes an inclusion of
 a space which is homotopically equivalent to $ \mathfrak{M}_g $ inside $\mathfrak{M}_{g+1} $
 (alternatively, one may say that $ Map_g \ra Map_{g+1}$, by letting the operation 
 be trivial on the added handle).

\begin{theo}[Harer \cite{harer}] 
Let ${\rm Map}_{g,r}^s$ be the mapping class group of an orientable surface $F$ of genus $g$
with $r$ boundary components and $s$ punctures. Then, for $g\geq 3k-1$, $H_k({\rm Map}_{g,r}^s , \ZZ)$ 
is independent of $g$ and $r$ as long as $r>0$, for $g\geq 3k$, $H_k({\rm Map}_{g,r}^s, \QQ)$ is
independent of $g$ and $r$ for every $r$, and for $g\geq 3k+1$,  $H_k({\rm Map}_{g,r}^s , \ZZ)$
is independent of $g$ and $r$ for every $r$.
\end{theo}

The ring structure of the cohomology of the ``stable mapping class group'' is described by a conjecture of Mumford (\cite{mumshaf}),
that has been proven by Madsen and Weiss (\cite{m-w}).

\begin{theo}{\bf (Mumford's conjecture)}\label{m-w}
The stable cohomology of the moduli space of curves is a polynomial algebra
$$ H^* ( Map_{\infty} ,  \QQ ) = \QQ [ \sK_1, \sK_2, \dots  ]$$
where the class $\sK_i$ is the direct image $(p_g)_* (K^{i+1})$ of the $(i+1)$-th power of the relative
canonical divisor of the universal family of curves.

\end{theo}

These results paralleled earlier results of Borel (\cite{borelcoh}) and Charney and Lee (\cite{CharneyLee})
on the cohomology of arithmetic varieties, such as the moduli space of Abelian varieties and some partial compactifications of them.

For instance, in the case of the moduli space of Abelian varieties we have the following theorem by Borel.

\begin{theo}{\bf (Borel)}\label{borel}
The stable cohomology of the moduli space of Abelian varieties  is a polynomial algebra
$$ H^* ( \sA_{\infty} , \QQ ) = \QQ [ \la_1, \la_3,  \la_5 \dots  ]$$
where the class $\la_i$ is the $i$-th Chern class of the universal bundle $H^{1,0}$.

\end{theo}

The theme of homology (cohomology-) stabilization is indeed a very general one, which has been recently
 revived through the work of several authors, also in other contexts ( see for instance \cite{dc-sa},  \cite{egh} , \cite{evw}, \cite{ght}).
 
 It would take too long to dwell here on this topic, which would deserve a whole survey article devoted to it.
 
 \subsection{Epilogue}
 
 There are many other  interesting topics which are very tightly related to the main theme of this article.
 
 For instance, there is a relation between symmetric differentials and the fundamental group (\cite{bog-deol1},\cite{bog-deol}, \cite{BKT}, \cite{klingler}).
 
 Brunebarbe, Klingler and Totaro  showed indeed  in \cite{BKT}  that some 'linearity' property of the fundamental group entails the existence of nontrivial
 symmetric differentials.
 
 \begin{theo}
 Let $X$ be a compact K\"ahler manifold, and let $k$ be any field.
 
 Assume that there is a representation 
 $$\rho : \pi_1(X) \ra  GL (r, k)$$
 with infinite image.
 
 Then the symmetric algebra of $X$ is nontrivial,
  $$ \oplus_{m\geq 0}  H^0 ( S^m (\Omega^1_X)) \neq \CC.$$
 
 \end{theo}

Another very interesting topic is Gromov's  h-principle,  for which we refer the reader to \cite{Eliashberg}, and  \cite{C-W}.  
 Perhaps not only the author is tired at this point.

\bigskip

\noindent
{\bf Acknowledgements.}
Thanks to Efim Zelmanov for the honour of  inviting me to write this article, and for his encouragement and support
during its preparation.

Thanks also to Maurizio Cornalba: some parts of this survey, which are directed to a rather general public, are
influenced by an unpublished draft which I wrote after
 a plenary lecture which I gave at the XIV Congress of the Italian Mathematical Union
 (Catania, 1991).

Thanks  to Ingrid Bauer, Michael L\"onne and Fabio Perroni, for the pleasure I got from  collaborating with them and 
for their  invaluable contributions to our exciting
joint research.

A substantial part of the paper was written when I was  visiting KIAS as research scholar: I am very grateful for the
excellent atmosphere and mathematical environment  I found at  KIAS. Lectures held in Bayreuth, KAIST Daejeon
and Centre Henri Lebesgue in Angers were also very useful.

Support of the Forschergruppe 790 `Classification of algebraic surfaces and compact complex manifolds' of the DFG,
and of the ERC Advanced grant n. 340258, `TADMICAMT'  is gratefully acknowledged.

Finally, thanks to the students and postdocs who listened to my lectures in Bayreuth for their remarks, 
and  to Ingrid Bauer, Binru Li, Sascha Weigl, and especially Wenfei Liu, for pointing out misprints and corrections to be made.



\begin{thebibliography}{Grif-SchmX}


\bibitem[Al75]{alexandrov}
Alexandrov, Paul
{\em Introduction \'a la Th\'eorie homologique de la dimension},
French translation of the original russian volume of 1975,
Editions MIR, Moscow (1977).



\bibitem[ABCKT96]{5book}
 Amor\'os, J.; Burger, Marc; Corlette, A.; Kotschick, D.; Toledo, D.
{\em Fundamental groups of compact K\"ahler manifolds. }
Mathematical Surveys and Monographs. 44. Providence, RI: American Mathematical Society (AMS). xi, 140 p. (1996). 


\bibitem[A-B-K-P-00]{abkp}
J. Amoros, F. Bogomolov, L. Katzarkov, T.Pantev
{\em Symplectic Lefschetz fibrations with arbitrary fundamental group},
With an appendix by Ivan
       Smith,
     {\bf J. Differential Geom. 54 }
     no. 3, (2000), 489--545.



\bibitem[A-F59]{AF1}
A. Andreotti, T.  Frankel,
{\em  The Lefschetz theorem on hyperplane sections.} Ann. of Math. (2) 69 (1959), 713--717. 




\bibitem[A-F69]{AF2}
A. Andreotti, T.  Frankel,
{\em  The second Lefschetz theorem on hyperplane
sections,}  in ' Global Analysis (Papers in Honor of K. Kodaira)' 
{\bf    Univ. Tokyo
Press, Tokyo } (1969), 1--20.
\bibitem[A-G62]{AG}
Andreotti, Aldo; Grauert, Hans,
{\em Th\'eor\'eme de finitude pour la cohomologie des espaces complexes.} Bull. Soc. Math. France 90 (1962) 193--259.

\bibitem[ABR92]{Abr}
Arapura, D.,  Bressler, P.,  Ramachandran, M. ,
{\em On the fundamental group of a compact K\"ahler manifold. }
Duke Math. J. 68 , no. 3, 477--488 (1992).

 \bibitem[ArCor09]{ar-cor}
 E. Arbarello, M.  Cornalba, {\em Teichm\"uller space via Kuranishi families.}
 {\bf  Ann. Sc. Norm. Super. Pisa Cl. Sci. (5) 8} (2009), no. 1, 89--116.
 
 
\bibitem[Art26]{art1}
E. Artin,
{\em Theorie der Z\"opfe,}
{\bf Hamburg Univ. Math. Seminar Abhandlungen 4-5} (1926), 47--72.

\bibitem[Art65]{art}
E. Artin,
{\em The collected papers of Emil Artin,}
edited by Serge Lang and John T. Tate,
{\bf Addison--Wesley Publishing Co., Inc.}, Reading, Mass.-London (1965).



 
\bibitem[ArtM74]{artinstacks}
M. Artin,
{\em Versal deformations and algebraic stacks.} Invent. Math. 27(1974), 165--189.

\bibitem[At57]{atiyah}
M.F. Atiyah,{\em Vector bundles over an elliptic curve. }
Proc. London Math. Soc. (3) 7 (1957) 414--452.

\bibitem[At58]{atiyahDP}
M.F. Atiyah,{\em
On analytic surfaces with double points.}
Proc. Roy. Soc. London. Ser. A 247 (1958)  237--244.

\bibitem[AHS78]{AHS}
Atiyah, M. F.; Hitchin, N. J.; Singer, I. M. 
{\em Self-duality in four-dimensional Riemannian geometry.}
 Proc. Roy. Soc. London Ser. A 362 , no. 1711, 425--461(1978).
 
\bibitem[A-B83]{a-b} 
 Atiyah, Michael F.; Bott, Raoul,
{\em The Yang-Mills equations over Riemann surfaces.}
Philos. Trans. R. Soc. Lond., A 308, 523--615 (1983).

\bibitem[Aub78]{Aubin} {\sc Aubin, T.:}
   {\it \' Equations du type Monge-Amp\'ere sur les vari\'et\'es
k\"ahl\'eriennes compactes.}
   Bull. Sci. Math. (2) 102 (1978), no. 1, 63-95.

\bibitem[A-K00]{a-k}
D. Auroux, L. Katzarkov ,
{\em Branched coverings of $\C \PP^2$ and invariants of symplectic
4-manifolds } {\bf Inv. Math. 142 } (2000), 631-673.

\bibitem[Bad04]{badescu}
 Badescu, Lucian 
 {\em Projective geometry and formal geometry. }
 Instytut Matematyczny Polskiej Akademii Nauk. Monografie Matematyczne (New Series), 65. Birkh\"auser Verlag, Basel, (2004). xiv+209 pp.

\bibitem[BdF08]{BdF}
Bagnera, G.; de Franchis, M.
{\em Le superficie algebriche le quali ammettono una rappresentazione parametrica mediante funzioni iperellittiche di due argomenti. }
Mem. di Mat. e di Fis. Soc. It. Sc. (3) 15, 253--343 (1908).

\bibitem[BarPetr07]{BarPetr}
Baranovsky, P. , Petrov, T. ,
{\em  Brauer groups and crepant resolutions,}
Adv. in Math. 209 (2007), 547--560.

\bibitem[BPHV]{bpv} W. Barth, C. Peters, A. Van de Ven, {\em Compact complex surfaces.}
 Ergebnisse der Mathematik und ihrer Grenzgebiete (3), 4. Springer-Verlag, Berlin, (1984); second edition by 
W. Barth, K. Hulek, C. Peters, A. Van de Ven,  Ergebnisse der Mathematik und ihrer Grenzgebiete. 3. Folge. A , 4. Springer-Verlag, Berlin, (2004). 

\bibitem[Barth70]{Barth1}
  Barth, Wolf {\em Transplanting cohomology classes in complex-projective space}, Amer. J. Math. 92 (1970), 951--967. 
  

\bibitem[B-L72]{bl} Barth ,  W.  and Larsen,M. E.  {\em  On the homotopy groups of complex projective algebraic manifolds}, Math. Scand. 30 (1972), 88--94. 

\bibitem[B-VdV74]{bvdv}
 Barth ,  W.  and Van de Ven, A.  {\em A decomposability criterion for algebraic 2-bundles on projective spaces}, Invent. Math. 25 (1974), 91--106.

\bibitem[Bau97]{bauer}
I. Bauer,
{\em  Irrational pencils on non-compact algebraic manifolds,}
{\bf   Internat. J. Math.
8 ,  no. 4} (1997), 441--450.

\bibitem[BCG05]{bcg}
I. Bauer, F. Catanese, F. Grunewald,
{\em Beauville surfaces without real structures.}
In: Geometric methods in algebra and number theory,
{\bf Progr. Math., 235}, Birkh\"auser (2005), 1--42.

\bibitem[BCG06]{almeria}
I.  Bauer, F. Catanese, F. Grunewald,
{\em Chebycheff and Belyi polynomials, dessins d'enfants, Beauville
    surfaces and group theory.} {\bf Mediterranean J. Math. 3,
no.2},  (2006) 119--143.


\bibitem[BCG07]{AbsoluteGalois} I.  Bauer, F. Catanese, F. Grunewald, {\em The absolute Galois group acts faithfully on
 the connected components of the moduli space of surfaces of general type},   arXiv:0706.1466 ,13 pages.

\bibitem[BCG14]{bcg2}
I.  Bauer, F. Catanese, F. Grunewald,
{\em Faithful actions of the absolute Galois group on connected components of moduli spaces},
Invent. Math. (2015) {\bf199}, 859-888
(published online in 2014,  
http://link.springer.com/article/10.1007/s00222-014-0531-2 

\bibitem[BC09a]{keumnaie} Bauer, I., Catanese, F., {\it The moduli
space of Keum-Naie surfaces}.
   Groups Geom. Dyn.  5  (2011),  no. 2, 231--250.


\bibitem[BC09b]{burniat1}
     Bauer, I., Catanese, F., {\it Burniat surfaces I: fundamental
groups and moduli of primary Burniat surfaces}.
   Faber, Carel (ed.) et al., Classification of algebraic varieties.
Based on the conference on classification of varieties,
Schiermonnikoog, Netherlands, May 2009. Z\"urich: European
Mathematical Society (EMS).
EMS Series of Congress Reports, 49--76 (2011).


\bibitem[BC10]{burniat2}
I. Bauer, F. Catanese,
{\em  Burniat surfaces. II. Secondary Burniat surfaces form three connected components of the moduli space}
 Invent. Math. 180 (2010), no. 3, 559--588. 
 
 \bibitem[BC10-b]{burniat3}
I. Bauer, F. Catanese,
{\em  Burniat surfaces III: deformations of automorphisms and extended Burniat surfaces},
Doc. Math. 18 (2013), 1089--1136..

\bibitem[BCP11]{survey} Bauer, I., Catanese, F.,  Pignatelli, R.
{\em Surfaces with geometric genus zero: a survey}.
Ebeling, Wolfgang (ed.) et al., Complex and differential geometry.
Conference held at Leibniz Universit\"at Hannover, Germany,
September 14--18, 2009. Proceedings. Berlin: Springer. Springer
Proceedings in Mathematics 8, 1-48 (2011).


\bibitem[B-C12]{bc-inoue}
Bauer, Ingrid, Catanese, Fabrizio 
{\em Inoue type manifolds and Inoue surfaces: a connected component of the moduli space of surfaces with $K^2=7, p_g=0$}. 
Geometry and arithmetic, 23--56, EMS Ser. Congr. Rep., Eur. Math. Soc., Z\"urich, (2012). 



\bibitem[B-C13]{bc-CMP}
Bauer, I.,  Catanese, F.
{\em  Burniat-type surfaces and a new family of surfaces with $p_g=0, K^2=3$}
 Rend. Circ. Mat. Palermo (2) 62 , no. 1, 37--60 (2013). 
 
 \bibitem[BCF14]{bcf}
 Bauer, I.,  Catanese, F., Frapporti, D.
{\em  Generalized Burniat type surfaces and Bagnera-de Franchis varieties}, J. Math. Sci. Univ. Tokyo {\bf 22} (2015), 55-111.

\bibitem[Belyi79]{belyi} G.V. Bely\u i,  {\em On Galois extensions of a maximal cyclotomic field,} {\bf Izv.
Akad. Nauk SSSR Ser. Mat. 43:2} (1979), 269--276. Translation in  Math. USSR- Izv. \textbf{14} (1980),
247--256.

\bibitem[Ber53]{Berger} {\sc Berger, M.:} {\it Sur les groupes
d'holonomie homog\'ene des vari\'et\'es \'a connexion affine et des
vari\'et\'es
riemanniennes.}  Bull. Soc. Math. France 83 (1953), 279-330

 \bibitem[Bers60]{Bers}
L. Bers,  {\em  Simultaneous uniformization. } {\bf Bull. Amer. Math. Soc. 66}  (1960),  94--97. 

\bibitem[BCS13]{bcs13}
Bestvina, Mladen; Church, Thomas; Souto, Juan 
{\em Some groups of mapping classes not realized by diffeomorphisms. }
Comment. Math. Helv. 88 (2013), no. 1, 205--220.

\bibitem[B--Mi98]{b-m}
Bierstone, Edward; Milman Pierre,
{\em  Semianalytic and subanalytic sets},
Publ. Math., Inst. Hautes \'Etud. Sci. 67, 5--42 (1988).

\bibitem[BF86]{bf} Biggers, R.; Fried, M.  { \it Irreducibility of
moduli spaces of cyclic unramified covers of
genus g curves.} Trans. Am. Math. Soc. 295, 59-70 (1986).

\bibitem[BCHM10]{bchm}
Birkar, Caucher; Cascini, Paolo; Hacon, Christopher D.; McKernan, James 
{\em Existence of minimal models for varieties of log general type. }
J. Amer. Math. Soc. 23 (2010), no. 2, 405--468.


\bibitem[Bir69]{Bir69} Birman, Joan S. {\it  Mapping class groups and their relationship to braid groups},
Comm. Pure Appl. Math. 22 (1969), 213--238.
 


\bibitem[Bir74]{birman}
J.S. Birman,
{\em Braids, links, and mapping class groups,}
 Annals of Mathematics Studies, {\bf No. 82},
Princeton University Press, Princeton, N.J.;
University of Tokyo Press, Tokyo (1974).

\bibitem[Blan53]{bla53}
A. Blanchard, {\em Recherche de structures analytiques complexes sur certaines vari\'et\'es. }
 C. R. Acad. Sci. Paris 236, (1953). 657--659.


\bibitem[Blan54]{bla54}
A. Blanchard, {\em Espaces fibr\'es k\"ahl\'eriens compacts.} C. R. Acad. Sci. Paris 238, (1954). 2281--2283. 



\bibitem[Blan56]{bla56}
A. Blanchard, {\em Sur les vari\'et\'es analytiques complexes. }  Ann. Sci. Ecole Norm. Sup. (3) 73 (1956), 157--202.

\bibitem[Bog78]{bogomolov}
 F. A. Bogomolov,
{\em Hamiltonian K\"ahlerian manifolds. }
{\bf Dokl. Akad. Nauk SSSR 243 , no. 5}, (1978) 1101--1104. 

\bibitem[Bog-Katz-98]{bog-katz}
 F. A. Bogomolov, L. Katzarkov, 
{\em Complex projective surfaces and infinite groups.} Geom. Funct. Anal. 8 (1998), no. 2, 243--272.




\bibitem[Bog78]{Bogstable}
 Bogomolov, F. A.,
{ \it Holomorphic tensors and vector bundles on projective manifolds. }
Izv. Akad. Nauk SSSR Ser. Mat. 42 (1978), no. 6, 1227 --1287, 1439.

\bibitem[Bog88]{Bog}
Bogomolov, F. A.,
{ \it The Brauer group of quotient spaces of linear representations. }
Math. USSR-Izv.  30 (3) (1988), 455--485.



\bibitem[Bog-DO11]{bog-deol1}
Bogomolov, Fedor; De Oliveira, Bruno 
{\em Symmetric differentials of rank 1 and holomorphic maps. }
Pure Appl. Math. Q. 7 (2011), no. 4, Special Issue: In memory of Eckart Viehweg, 1085--1103.

\bibitem[Bog-DO13]{bog-deol}
Fedor Bogomolov, Bruno De Oliveira,
{\em Closed symmetric 2-differentials of the 1st kind}
 Pure Appl. Math. Q. 9(4) (2013), 613--642 .

\bibitem[Bom73]{bom}
   E. Bombieri,  {\em Canonical models of surfaces of
general type}, {\bf Publ. Math. I.H.E.S., 42} (1973), 173-219.

\bibitem[B-H75]{arcata}
E. Bombieri, D. Husemoller,
{\em  Classification and embeddings of surfaces, }
in 'Algebraic geometry,
Humboldt State Univ., Arcata, Calif., 1974 ',
{ \bf  Proc. Sympos. Pure Math., Vol. 29, Amer. Math. Soc., Providence, R.I.}
     (1975), 329--420.

\bibitem[Bore60]{borelSTG}
 Borel, Armand
{\em Seminar on transformation groups.} With contributions by G. Bredon, E. E. Floyd, D. Montgomery and R. Palais. 
(Annals of Mathematics Studies, No. 46). Princeton, N. J.: Princeton University Press. 245 p. (1960).

\bibitem[Bore63]{Bo63} {\sc Borel, A.:} {\it Compact Clifford-Klein
forms of symmetric spaces}, Topology 2 1963 111 - 122.

\bibitem[Borel74]{borelcoh}
Borel, Armand
{\em Stable real cohomology of arithmetic groups. }
Ann. Sci. \'Ecole Norm. Sup. (4) 7 (1974), 235--272 (1975). 

\bibitem[B-M60]{B-M}
Borel, Armand; Moore, J.C.
{\em Homology theory for locally compact spaces.}
Mich. Math. J. 7, 137--159 (1960)

\bibitem[Briesk68-a]{brieskorn}
  E. Brieskorn,
{\em Rationale Singularit\"aten komplexer Fl\"achen,}
{\bf Invent. Math.
4 } ( 1967/1968 ), 336--358.

\bibitem[Briesk68-b]{brieskorn2}
  E. Brieskorn,
{\em  Die Aufl\"osung der rationalen Singularit\"aten holomorpher
Abbildungen,} {\bf  Math. Ann.  178} (  1968),  255--270.

\bibitem[Briesk71]{nice}
  E. Brieskorn,
{\em Singular elements of semi-simple algebraic groups,}  'Actes du Congr\'es
International des Math\'ematiciens' (Nice, 1970), Tome 2,
{\bf Gauthier-Villars}, Paris (1971), 279--284.

\bibitem[Bro]{Brown} Brown, K. S., {Cohomology of groups},
{Graduate Texts in Mathematics},
  {\bf 87},  {Springer-Verlag},  {New York}, (1982).
  
  \bibitem[BKT13]{BKT}
  Brunebarbe, Yohan; Klingler, Bruno; Totaro, Burt
{\em Symmetric differentials and the fundamental group. }
Duke Math. J. 162 (2013), no. 14, 2797--2813. 

\bibitem[Bur66] {burniat}
P. Burniat, {\em Sur les surfaces de genre $P_{12} > 1$}, Ann. Mat. Pura Appl. (4) 71 (1966), 1--24.
  
\bibitem[B-R75]{b-r}
D. Burns,M.  Rapoport,{\em On the Torelli problem for k\"ahlerian K3 surfaces.}
Ann. Sci. \'Ecole Norm. Sup. (4) 8 (1975), no. 2, 235--273.

\bibitem[B-W74]{b-w}
D. Burns, J.  Wahl,
{\em Local contributions to global deformations of surfaces,}
{\bf Invent. Math. 26}(1974), 67-88 .

\bibitem[Cai07]{cai}
 J.X.   Cai, {\em
Classification of fiber surfaces of genus 2 with automorphisms acting trivially in cohomology.}
Pacific J. Math. 232 (2007), no. 1, 43--59. 

\bibitem[C-L13]{cai-wenfei}
 Cai, Jin-Xing; Liu, Wenfei; Zhang, Lei 
 {\em Automorphisms of surfaces of general type with $ q \geq 2$ acting trivially in cohomology.} Compos. Math. 149 (2013), no. 10, 1667--1684.

\bibitem[Cal58]{cal58}
E. Calabi, {\em Construction and properties of some 6-dimensional almost complex manifolds. }
Trans. Amer. Math. Soc. 87 (1958),  407--438.

\bibitem[C-V60]{Cal-Ves}
Calabi, Eugenio; Vesentini, Edoardo
{\em On compact, locally symmetric K\"ahler manifolds. }
Ann. Math. (2) 71, 472--507 (1960).

\bibitem[Cam94]{campanaShaf}
Campana, Fr\'ed\'eric,
{\em
Remarques sur le rev\^etement universel des vari\'et\'es k\"ahl\'eriennes compactes.}
Bull. Soc. Math. France 122 (1994), no. 2, 255--284

\bibitem[Cam95A]{campanaFG}
Campana, Fr\'ed\'eric,
{\em Kodaira dimension and fundamental group of compact K\"ahler manifolds }
Lecture Notes Series, 7, Dipartimento di Matematica, Universit\"a di Trento (1995).

\bibitem[Cam95B]{campanaJAG}
Campana, Fr\'ed\'eric,
{\em Fundamental group and positivity of cotangent bundles of compact K\"ahler manifolds. }
J. Algebraic Geom. 4 (1995), no. 3, 487--502. 

\bibitem[C-W12]{C-W}
Fr\'ed\'eric Campana, J\"org Winkelmann,
{\em On h-principle and specialness for complex projective manifolds}
arXiv:1210.7369 .

\bibitem[C-T89]{C-T}
Carlson, James A.,  Toledo, Domingo,
{\em  Harmonic mappings of K\"ahler manifolds to locally symmetric spaces.}
 Inst. Hautes \'Etudes Sci. Publ. Math. No. 69, 173--201(1989).

 \bibitem[Car28]{cartanEspacesRiemann} 
Cartan, Elie
{\em Lecons sur la g\'eom\'etrie des espaces de Riemann. }
Paris: Gauthier-Villars (Cahiers scientifiques publi\'es sous la direction de G. Julia, 2). VI, 273 p. (1928).

\bibitem[Car35]{Cartan} 
{\sc Cartan, E.:} {\it Sur les domaines
born\'es homog\'enes de l'espace des $n$ variables complexes.}
Abhandl. Hamburg 11,
116--162 (1935).

\bibitem[Car57]{cartan}
H. Cartan,
{\em Quotient d'un espace analytique par un groupe
d'automorphismes,} in 'A symposium in honor of S. Lefschetz, 
Algebraic geometry
and topology' {\bf Princeton University Press}, Princeton, N. J. 
(1957),  90--102.

 \bibitem[CE56]{cartaneil}
 H. Cartan, S. Eilenberg {\em  Homological algebra.} Princeton University Press, Princeton, N. J., 1956. xv+390 pp.


\bibitem[Cast05]{CdF}
 Castelnuovo, G.,
{\em Sulle superficie aventi il genere aritmetico negativo.}
Palermo Rend. 20, 55--60 (1905).


\bibitem[Cat82]{catrav}
F. Catanese,
{\em  Moduli of surfaces of general type. } in  'Algebraic
geometry---open problems' (Ravello, 1982),
   {\bf Lecture Notes in Math., 997} Springer, Berlin-New York (1983),90--112.

\bibitem[Cat84]{cat1}
F. Catanese,
{\em On the Moduli Spaces of Surfaces of General Type},
{\bf J. Differential Geom. 19} (1984), 483--515.

\bibitem[Cat84b]{bucharest} Catanese, F. {\em  Commutative algebra methods and equations of regular surfaces. }
Algebraic geometry, Bucharest 1982 (Bucharest, 1982), 68--111, Lecture Notes in Math., 1056, Springer, Berlin, (1984).



\bibitem[Cat87]{cat2}
F. Catanese,
{\em  Automorphisms of Rational Double Points and
   Moduli Spaces  of Surfaces of General Type},
{\bf Comp. Math. 61 } (1987),  81-102.

\bibitem[Cat86]{cat3}
F. Catanese,
{\em Connected Components of  Moduli  Spaces},
{\bf J. Differential Geom. 24} (1986), 395-399.

\bibitem[Cat88]{montecatini}
F. Catanese,
{\em Moduli of algebraic surfaces.}  Theory of moduli (Montecatini
Terme, 1985),
   {\bf  Lecture
Notes in Math. 1337},  Springer, Berlin (1988), 1--83.

\bibitem[Cat89]{cat5}
F. Catanese,
{\em   Everywhere non reduced Moduli  Spaces },
{\bf Inv. Math. 98}, 293-310 (1989).

\bibitem[Cat91]{albanese}
F. Catanese,
   {\em Moduli and classification of irregular K\"ahler
manifolds (and algebraic
varieties) with Albanese general type fibrations.} {\bf Invent. math.
104},
    (1991), 263-289.
    
    \bibitem[Cat92]{Cat-Chow}
Catanese, Fabrizio,
   {\em  Chow varieties, Hilbert schemes and moduli spaces of surfaces of general type.}
    J. Algebraic Geom. 1 (1992), no. 4, 561--595.
    
  \bibitem[Caci93]{cacicetraro} 
    F. Catanese, C.  Ciliberto, {\em On the irregularity of cyclic coverings of algebraic surfaces.}
     Geometry of complex projective varieties (Cetraro, 1990), 89--115, Sem. Conf., 9, Mediterranean Press, Rende, (1993).
     
     \bibitem[Cat95]{catAV}
     Catanese, Fabrizio 
     {\em Compact complex manifolds bimeromorphic to tori. }, in `Abelian varieties' (Egloffstein, 1993), 55--62, de Gruyter, Berlin, (1995).

\bibitem[Cat96]{Levico}
F. Catanese,
   {\em Fundamental groups with few relations,}
in 'Higher dimensional complex varieties. Proceedings of the international
  conference, Trento, Italy, June 15--24, 1994'
  {\bf Walter de Gruyter, Berlin} (1996), 163-165.
  
  \bibitem[Cat94b]{paris}
   Catanese, Fabrizio {\em (Some) old and new results on algebraic surfaces.}
    First European Congress of Mathematics, Vol. I (Paris, 1992), 445--490, Progr. Math., 119, Birkh\"auser, Basel, (1994). 

\bibitem[C-F-96]{cf}
F. Catanese, M. Franciosi,
{\em  Divisors of small genus on algebraic surfaces and projective embeddings,}
in 'Proceedings of the Hirzebruch 65 Conference on Algebraic Geometry 
(Ramat Gan,
1993)',  {\bf  Israel Math. Conf. Proc., 9}, Bar-Ilan Univ., Ramat 
Gan, (1996),109--140.

\bibitem[CFHR99]{4auth}
F. Catanese, M. Franciosi, K. Hulek, M. Reid,
{\em  Embeddings of curves
and surfaces,}  {\bf Nagoya Math. J.  154 } (1999), 185--220.



\bibitem[Cat99]{bidouble}
F. Catanese, {\em Singular bidouble covers
and the construction of
interesting algebraic surfaces.}  In `Algebraic Geometry: Hirzebruch 70', Pragacz, Szurek, Wisniewski editors, {\bf AMS Cont. Math.  241} 
(1999), 97 - 120.


\bibitem[Cat00]{isogenous}
F. Catanese,
{\em Fibred surfaces, varieties isogenous to a product
and related moduli spaces.}, {\bf  Amer. J. Math. 122, no.1} (2000),
1--44.

\bibitem[C-F03]{c-f}
F. Catanese, P. Frediani, 
{\em Real hyperelliptic surfaces and
the orbifold fundamental group}   Journal of the Inst. Math. Jussieu {\bf 2} (2) ,(2003) 169--233.

\bibitem[Cat03]{cat03}
F. Catanese,
{\em   Moduli  Spaces of Surfaces and Real Structures},
{\bf Ann. Math.(2) \textbf{158}, no.2} (2003),577-592.

\bibitem[Cat02A]{nankai}
Catanese, Fabrizio 
{\em Deformation types of real and complex manifolds.}
 Contemporary trends in algebraic geometry and algebraic topology (Tianjin, 2000), 195�238, Nankai Tracts Math., 5, World Sci. Publ., River Edge, NJ,(2002).


\bibitem[Cat02]{cat02}
F. Catanese,
{\em  Symplectic structures of algebraic surfaces and deformation},
14 pages, math.AG/0207254.



\bibitem[Cat03b]{fibred}
F. Catanese,
{\em  Fibred K\"ahler and quasi projective groups},
Advances in Geometry, suppl. , Special Issue
dedicated to A. Barlotti's 80-th birthday (2003),
{\bf Adv. Geom. suppl.} (2003) S13--S27.

\bibitem[CKO03]{cko}
F. Catanese, J. Keum, K. Oguiso,
``{\em  Some remarks on the universal cover of an open K3 surface.}
{ \bf Math. Ann. 325, No.2,}  (2003), 279-286.


\bibitem[Cat04]{cat04}
F. Catanese,
{\em Deformation in the large of some complex manifolds, I} Ann. Mat. Pura Appl. (4)
183,  Volume in Memory of Fabio Bardelli, (2004), no. 3, 261--289.

\bibitem[CW04]{c-w}
F. Catanese, B. Wajnryb,
{\em Diffeomorphism of simply connected algebraic surfaces}
J. Differential Geom. 76 (2007), no. 2, 177--213.

\bibitem[Cat06]{localpi1}
Catanese, Fabrizio 
{\em Surface classification and local and global fundamental groups. I. }
Atti Accad. Naz. Lincei Cl. Sci. Fis. Mat. Natur. Rend. Lincei (9) Mat. Appl. 17 (2006), no. 2, 135--153.


\bibitem[Cat08]{cime} Catanese, F. {\it Differentiable and
deformation type of algebraic  surfaces, real and
symplectic structures.}
           Symplectic 4-manifolds and algebraic surfaces,  55--167,
Lecture Notes in Math., 1938, Springer, Berlin,
(2008).


\bibitem[Cat09]{cat09}
F. Catanese,
 {\em Canonical symplectic structures and deformations of algebraic surfaces},
Comm. in Contemp. Math. 11 (2009), n. 3, 481--493.

\bibitem[Cat09b]{perugia}
F. Catanese,
 {\em  Algebraic surfaces and their moduli spaces: real, differentiable and symplectic structures}
  Boll. Unione Mat. Ital. (9) 2 (2009), no. 3, 537--558.
  
  \bibitem[CatRol09]{CR}
F. Catanese, S. Rollenske, 
 {\em   Double Kodaira fibrations.} {\bf J. Reine Angew. Math. 628} (2009), 205--233. 
 
 \bibitem[COP10]{cop}
 Catanese, Fabrizio; Oguiso, Keiji; Peternell, Thomas 
 {\em On volume-preserving complex structures on real tori.}
  Kyoto J. Math. 50 (2010), no. 4, 753--775.
  
  \bibitem[Cat-Lo-Waj11]{clw}
F. Catanese, M. L\"onne, B. Wajnryb,
{\em
Moduli spaces and braid monodromy types of bidouble covers of the quadric,}
{\bf Geometry - Topology 15 } (2011) 351--396,
   URL: http://www.msp.warwick.ac.uk/gt/2011/15-01/p010.xhtml
   DOI: 10.2140/gt.2011.15.351.   


\bibitem[CLP11]{CLP1}
Catanese, Fabrizio, L\"onne, Michael,  Perroni, Fabio,
{\em  Irreducibility of the space of dihedral covers of the projective line of a given numerical type. }
Atti Accad. Naz. Lincei Cl. Sci. Fis. Mat. Natur. Rend. Lincei (9) Mat. Appl. 22 , no. 3, 291--309 (2011). 


\bibitem[Cat12]{singMg}
Catanese, Fabrizio,
{\em Irreducibility of the space of cyclic covers of algebraic curves of fixed numerical type and the irreducible components of $Sing (\bar{\mathfrak M_g})$}
in `Advances in Geometric Analysis', in honor of Shing-Tung Yau's 60th birthday, S. Janeczko, J. Li, D.H. Phong editors, Advanced Lectures in Mathematics 21,
International Press (Somerville, USA) and  Higher Education Press (Beijing, China) , 281--306 (2012).

\bibitem[Cat13]{handbook}
Catanese,Fabrizio,
{\em  A superficial working guide to deformations and moduli. }
Handbook of moduli. Vol. I, 161--215, Adv. Lect. Math. (ALM), 24, Int. Press, Somerville, MA, (2013).

\bibitem[CLP12]{CLP2}
Catanese, Fabrizio, L\"onne, Michael,  Perroni, Fabio,
{\em The irreducible components of the moduli space of dihedral covers of algebraic curves},
arXiv:1206.5498, to appear in Groups, Geometry and Dynamics.

\bibitem[CLP13]{CLP3}
Catanese, Fabrizio, L\"onne, Michael,  Perroni, Fabio,
{\em Genus stabilization for moduli of curves with symmetries}
arXiv:1301.4409, to appear in Algebraic Geometry. 

\bibitem[C-DS13]{C-DS1}
Catanese, Fabrizio; Di Scala, Antonio J. 
{\em A characterization of varieties whose universal cover is the polydisk or a tube domain.} Math. Ann. 356 (2013), no. 2, 419--438.

\bibitem[C-DS14]{C-DS2}
Catanese, Fabrizio; Di Scala, Antonio J. 
{\em A characterization of varieties whose universal cover is a bounded symmetric domain without ball factors.}
 Adv. Math. 257 (2014), 567--580.
 
\bibitem[ChCou10]{coughlanchan}
Chan Mario T.,  Coughlan 
S.,  {\it  Kulikov surfaces form a 
connected component of the 
moduli space},
  Nagoya Math. J. 210 (2013), 1--27. 

\bibitem[Ch-Le83]{CharneyLee}
Charney, Ruth; Lee, Ronnie 
{\em Cohomology of the Satake compactification. }
Topology 22 (1983), no. 4, 389--423.

\bibitem[CY13]{yifan}
Chen, Yifan
{\em A new family of surfaces of general type with $K^2= 7$ and $p_g = 0$. }
Math. Z. 275, No. 3-4, 1275--1286 (2013).


\bibitem[CDS14]{donKE}
Chen, Xiuxiong; Donaldson, Simon; Sun, Song 
{\em K\"ahler-Einstein metrics and stability.}
 Int. Math. Res. Not. IMRN (2014), no. 8, 2119--2125.

\bibitem[CDS12-3]{donKEfull}
Xiuxiong Chen, Simon Donaldson, Song Sun
 {\em K\"ahler-Einstein metrics on Fano manifolds. I: Approximation of metrics with cone singularities. II: Limits with cone angle less than $2\pi$. III: Limits as cone angle approaches $2\pi$ and completion of the main proof.}
arXiv:1211.4566,  arXiv:1212.4714,
arXiv:1302.0282 , J. Amer. Math. Soc. 28 (1) (2015), 183--197, 199--234, 235--278.


\bibitem[Chev55]{chevalley}
Chevalley, C.
{ \it Invariants of finite groups generated by reflections. }
Amer. J. Math. 77 (1955), 778--782.



\bibitem[Chmu92]{chmutov}
S.V. Chmutov, {\em Examples of projective surfaces with many singularities}
 J. Algebraic Geom. 1 (1992), no. 2, 191--196.
 
 
\bibitem[Cleb72]{Clebsch}
Clebsch A., {\it Zur Theorie der Riemann'schen 
Fl\"achen.} Math. Ann. 6, 216-230 (1872).
 



 \bibitem[Clem05]{clemens}
H. Clemens, {\em Geometry of formal Kuranishi theory.} Adv. Math. 198 (2005), no. 1, 311--365.

\bibitem[Cor87]{cornalba}
Cornalba, Maurizio,
{\em  On the locus of curves with automorphisms. } Ann. Mat. Pura Appl. (4) 149, 135--151  (1987).
Erratum:  Ann. Mat. Pura Appl. (4) 187, no. 1, 185--186  (2008).

\bibitem[Cor95]{alessio}
Corti, A.,
 {\it Factoring birational maps of threefolds after Sarkisov.}
 J. Algebraic Geom. 4 (1995), no. 2, 223--254.
 
 \bibitem[C-MI-P-13 ]{CMP}
 Antonio F. Costa, Milagros Izquierdo, Hugo Parlier,
{\em Connecting $p$-gonal loci in the compactification of moduli space}
arXiv:1305.0284.


\bibitem[dF05]{deFranchis}
de Franchis, M.
{\em Sulle superficie algebriche le quali contengono un fascio irrazionale di curve. }
Palermo Rend. 20, 49--54 (1905).

\bibitem[DC-S99]{dc-sa}
De Concini, Corrado; Salvetti, Mario
{\em Stability for the cohomology of Artin groups.} 
Adv. Math. 145 (1999), no. 2, 291--305. 

\bibitem[Dehn38]{dehn}
M. Dehn,
{\em Die Gruppe der Abbildungsklassen. (Das arithmetische Feld auf Fl\"achen.)}
{\bf Acta Math. 69} (1938), 135--206.

      \bibitem[Del81]{del}
 Deligne, Pierre, {\em  Le groupe fondamental du compl\'ement d'une courbe plane n'ayant que des points doubles ordinaires est ab\'elien (d'apr\`es W. Fulton).}  [The fundamental group of the complement of a plane curve having only ordinary double points is abelian (after W. Fulton)] Bourbaki Seminar, Vol. 1979/80, pp. 1--10, Lecture Notes in Math., 842, Springer, Berlin-New York, (1981).


 \bibitem[DGMS75]{DGMS}
 P. Deligne, P.  Griffiths, J, Morgan, D. Sullivan, {\em  Real homotopy theory of K\"ahler manifolds.}
  Invent. Math. 29 (1975), no. 3, 245--274.


\bibitem[Del-Most93]{DM}
P. Deligne, G.D. Mostow, 
{\em Commensurabilities among lattices in PU(1,n). }
Annals of Mathematics Studies, 132. Princeton University Press, Princeton, NJ, (1993) viii+183 pp..
      
      \bibitem[Del-Mum69]{d-m}
    P.    Deligne, D. Mumford {\em The irreducibility of the space of curves of given genus}
     Inst. Hautes \'Etudes Sci. Publ. Math. No. 36 (1969) 75--109.
     

     
     \bibitem[Dolga81]{Dolgachev} Dolgachev, I. {\it  Algebraic 
surfaces 
with $q=p_g=0$}.
Algebraic surfaces,  97-215, C.I.M.E. 
Summer
School 1977, 76, Liguori Editore, Napoli (1981), reedited by 
Springer, Heidelberg, (2010).

\bibitem[Dolg82]{dolg}
I. Dolgachev,
{\em  Weighted projective varieties,} in ' Group actions and vector fields
(Vancouver, B.C., 1981)', {\bf Lecture Notes in Math., 956, 
Springer}, Berlin, (1982),
34--71.



\bibitem[Don92]{don} S.K. Donaldson,
{\em Gauge theory and four-manifold topology.}
[CA] Joseph, A. (ed.) et al., First European congress of mathematics
(ECM),
   Paris, France, July 6-10, 1992. Volume I: Invited lectures (Part 1).
{\bf Basel:
Birkh\"auser, Prog. Math. 119}  (1994) 121-151 .

\bibitem[Don96]{don5} S.K. Donaldson,
{\em
The Seiberg-Witten Equations and 4-manifold topology. }
{\bf Bull. Am. Math. Soc., (N S) 33, 1}  (1996) 45-70.

\bibitem[Drin90]{drinfeld}
 Drinfel' d, V.G.
{\em On quasitriangular quasi-Hopf algebras and on a group that is closely connected with}
$Gal( \bar{\QQ}/ \QQ)$,  (Russian) Algebra i Analiz 2 (1990), no. 4, 149--181; translation in Leningrad Math. J. 2(1991), no. 4, 829--860.


\bibitem[Du-Th06]{DT}
Dunfield, N. M.; Thurston, W. P.
{\it Finite covers of random 3-manifolds.}
Invent. Math. 166 , no. 3, 457--521 (2006).

\bibitem[DuVal34]{duval}
P. Du Val, {\em
On singularities of surfaces which do not impose adjunction
conditions.}
Proc. Cambridge Philos. Soc. 30 (1934), 483--491.

\bibitem[EastVak07]{east-vak}
 R.W. Easton, R.  Vakil {\em Absolute Galois acts faithfully on the components of the moduli space of surfaces: 
 a Belyi-type theorem in higher dimension},  Int. Math. Res. Not. IMRN (2007), no. 20, Art. ID rnm080, 10 pp. 
 
 \bibitem[Eck-Fr51]{E-F} 
 B. Eckmann,A. Fr\"olicher, {\em Sur l'int\'egrabilit\'e des structures presque complexes.} {\bf C. R. Acad. Sci. Paris 232}, (1951) 2284--2286.
 
 \bibitem[Edm82]{Edm I} Edmonds, A. L., {\it Surface symmetry. {I}}, {Michigan Math. J.}, {\bf 29}  (1982), {n. 2}, {171--183}.
\bibitem[Edm83]{edmonds2} Edmonds, A. L., {\it Surface symmetry. {II}}, {Michigan Math. J.}, {\bf 30}  (1983), {n. 2}, {143--154}.


\bibitem[Eells-Sam64]{eells-sampson}
Eells, J., Sampson, J.H. 
{\em  Harmonic maps of Riemannian manifolds},
Amer. Jour. Math. 86 (1964), 109--160.

\bibitem[Ehr47]{ehre}
C. Ehresmann,
{\em Sur les espaces fibr\'es  diff\'erentiables.}
{\bf C.R. Acad. Sci. Paris 224} (1947), 1611-1612.

\bibitem[Ehr49]{ACS}
C. Ehresmann,
{\em Sur la th\'eorie des espaces fibr\'es. } Topologie alg\'ebrique, pp. 3--15. 
Colloques Internationaux du Centre National de la Recherche Scientifique, no. 12. 
Centre de la Recherche Scientifique, Paris, (1949). 

\bibitem[Eke86]{teke}
T. Ekedahl,
{\em Two examples of smooth projective varieties with nonzero Massey products. }
Algebra, algebraic topology and their interactions (Stockholm, 1983), 128--132, Lecture Notes in Math., 1183, Springer, Berlin, (1986).



\bibitem[Eke88]{ekedahl}
T. Ekedahl,
{\em Canonical models of surfaces of general type in positive
characteristic,} {\bf  Inst. Hautes \'Etudes Sci. Publ. Math.  No. 
67}  (1988), 97--144.


\bibitem[E-M02]{Eliashberg}
Eliashberg, Y.; Mishachev, N. {\em Introduction to the h-principle.}
 Graduate Studies in Mathematics, 48. American Mathematical Society, Providence, RI, (2002). xviii+206 pp.

\bibitem[Elk]{elkik}
R. Elkik, {\em Singularit\'es rationnelles et d\'eformations} Invent. Math. 47 (1978), no. 2, 139 --147.

\bibitem[EVW09]{evw}
Jordan S. Ellenberg, Akshay Venkatesh, Craig Westerland,
{\em Homological stability for Hurwitz spaces and the Cohen-Lenstra conjecture over function fields}
 arXiv:0912.0325 .

\bibitem[Ell07]{ellia}
Ellia, Philippe
{\em Codimension two subvarieties and related questions.}  Vector bundles and low codimensional subvarieties: state of the art and recent developments, 129-- 208, 
Quad. Mat., 21, Dept. Math., Seconda Univ. Napoli, Caserta, (2007). 

\bibitem[Enr96]{enr96} Enriques, F., 
{\em Introduzione alla geometria 
sopra le superficie algebriche}. 
Memorie della Societa'
Italiana delle Scienze (detta "dei XL"), s.3, 
to. X, (1896), 1--81.

\bibitem[ES09]{Enr-Sev}
Enriques, F.; Severi, F..
{\em M\'emoire sur les surfaces hyperelliptiques. }  
[Acta Math. 32, 283-- 392 (1909) and 
33, 321--403 (1910).

\bibitem[Enr49]{enr}
F. Enriques, {\em Le Superficie Algebriche.}
{\bf Zanichelli},
Bologna (1949) xv+464 pp. 

\bibitem[EGH10]{egh}
Erdenberger, C.; Grushevsky, S.; Hulek, K. 
{\em Some intersection numbers of divisors on toroidal compactifications of $\sA_g$}. J. Algebraic Geom. 19 (2010), no. 1, 99--132.

\bibitem[EV92]{esnaultviehweg}  Esnault, H.,  Viehweg, E., {\em Lectures on Vanishing Theorems}.
DMV Seminar, Band 20, Birkh\"auser Verlag Basel, 1992 .

 \bibitem[Eyss04]{eyss-inv}
 Eyssidieux, Philippe, {\em Sur la convexit\'e holomorphe des rev\^etements lin\'eaires r\'eductifs d'une vari\'et\'e projective alg\'ebrique complexe. }Invent. Math. 156 (2004), no. 3, 503--564.

\bibitem[EKPR12]{ekpr}
P. Eyssidieux, L. Katzarkov, T.  Pantev, M. Ramachandran, {\em Linear Shafarevich conjecture. } Ann. of Math. (2) 176 (2012), no. 3, 1545--1581.

\bibitem[Eyss11]{eyssidieuxLN}
P. Eyssidieux, {\em  Lectures on the Shafarevich conjecture on uniformization. }
Complex manifolds, foliations and uniformization, 101--148, Panor. Synth�ses, 34/35, Soc. Math. France, Paris, 2011.

\bibitem[EGZ09]{EGZ}
Eyssidieux, Philippe; Guedj, Vincent; Zeriahi, Ahmed {\em Singular K\"ahler-Einstein metrics.}
 J. Amer. Math. Soc. 22 (2009), no. 3, 607--639.
 
 \bibitem[Falt-C90]{FaltingsChai}
  Faltings, Gerd; Chai, Ching-Li 
  {\em Degeneration of abelian varieties. }
  With an appendix by David Mumford. Ergebnisse der Mathematik und ihrer Grenzgebiete (3) [Results in Mathematics and Related Areas (3)], 22. Springer-Verlag, Berlin, (1990). xii+316 pp.

 \bibitem[Fant00]{fantechi}
 B. Fantechi, {\em Stacks for everybody},
  European Congress of Mathematics, Vol. I (Barcelona, 2000), 349-359, Progr. Math., 201, Birkh\"auser, Basel, (2001).





\bibitem[FM12]{fm}
B.~Farb, D.~Margalit,
{\it A primer on mapping class groups}, volume~49 of {\em Princeton
  Mathematical Series},
Princeton University Press, Princeton, NJ, 2012.

\bibitem[Free82]{free}  M. Freedman, {\em The topology of
four-dimensional manifolds.},
{\bf J. Differential Geom. 17, n. 3} (1982), 357-454.

\bibitem[F-Q90]{f-q}  M. Freedman and F.Quinn, {\em Topology of 4- manifolds.},
{\bf Princeton Math. Series 39}  Princeton Univ. Press (1990).




\bibitem[F-M88]{f-m1}  R. Friedman and J.W.Morgan, {\em  Algebraic
   surfaces and four-manifolds: some conjectures and speculations.},
{\bf Bull. Amer. Math.Soc. 18} (1988),1--19.


  \bibitem[F-J08]{F-J}  M. D. Fried, M. Jarden {\em  Field arithmetic. } Third edition. Revised by Jarden. Ergebnisse der Mathematik und ihrer Grenzgebiete. 3. Folge. A Series of Modern Surveys in Mathematics, 11. Springer-Verlag, Berlin, 2008. xxiv+792 pp.


\bibitem[FV91]{FV}
Fried, M. D.; V{\"o}lklein, H.,  {\it The inverse {G}alois problem and rational points on moduli spaces},
  {Math. Ann.},  {\bf 290},  {1991},  {n. 4},  {771--800}.

\bibitem[Ful69]{Fulton}
W. Fulton: { \it Hurwitz schemes and 
irreducibility of moduli of algebraic curves.}
  Ann. of Math. (2) 90 ,542-575 (1969) .




\bibitem[Ful87]{fultonAMS}
 Fulton, William {\em On the topology of algebraic varieties. } Algebraic geometry, Bowdoin, 1985 (Brunswick, Maine, 1985), 15--46, Proc. Sympos. Pure Math., 46, Part 1, Amer. Math. Soc., Providence, RI, (1987).

\bibitem[Ful07]{fultonLectures}
 Fulton, William {\em Equivariant cohomology in algebraic geometry}, Lectures at Columbia University (2007).
 
\bibitem[FH79]{fh}
Fulton, William; Hansen, Johan {\em A connectedness theorem for projective varieties, with applications to intersections and singularities of mappings.}
 Ann. of Math. (2) 110 (1979), no. 1, 159--166.

\bibitem[FL81]{f-l}
Fulton, William; Lazarsfeld, Robert {\em Connectivity and its applications in algebraic geometry. } Algebraic geometry (Chicago, Ill., 1980), pp. 26--92, Lecture Notes in Math., 862, Springer, Berlin-New York, (1981).

\bibitem[Geer13]{gerardAV}
 van der Geer, Gerard 
 {\em The cohomology of the moduli space of abelian varieties.}
 Handbook of moduli. Vol. I,  415--457, Adv. Lect. Math. (ALM), 24, Int. Press, Somerville, MA, (2013).

\bibitem[Gie77]{gieseker}
D. Gieseker,
{\em Global moduli for surfaces of general type,}
{\bf Invent. Math. 43, no. 3}(1977), 233--282.

\bibitem[Gie82]{Giestata}
D. Gieseker,
{\em Lectures on moduli of curves. } Tata Institute of Fundamental Research Lectures on Mathematics and Physics, 69. Published for the Tata Institute of Fundamental Research, Bombay; Springer-Verlag, Berlin-New York, (1982). iii+99 pp.

\bibitem[Gie83]{Giescime}
D. Gieseker,
{\em Geometric invariant theory and applications to moduli problems.}
 Invariant theory (Montecatini, 1982), 45--73, Lecture Notes in Math., 996, Springer, Berlin, (1983). 

\bibitem[God35]{god} 
Godeaux, L., {\em  Les involutions cycliques 
appartenant \`a une 
surface alg\'ebrique}. Actual. Sci. Ind.,
{\bf  270}, Hermann, Paris, 
1935.

\bibitem[God58]{godement}
Godement, Roger {\em Topologie alg\'ebrique et th\'eorie des faisceaux. }  Actualit\'es Scientifiques et Industrielles 
No. 1252. Publ. Math. Univ. Strasbourg. No. 13 Hermann, Paris (1958) viii+283 pp. ; reprinted (1973).

 \bibitem[Go-Mi88]{go-mi}
 Goldman, William M.; Millson, John J. 
 {\em The deformation theory of representations of fundamental groups of compact K\"ahler manifolds. }
 Inst. Hautes \'Etudes Sci. Publ. Math. No. 67 (1988), 43--96.

\bibitem[Gompf95]{Gompf}
Gompf, Robert E., 
{\em A new construction of symplectic manifolds. }
Ann. of Math. (2) 142 , no. 3, 527--595 (1995).

\bibitem[GD-JZ14]{gabino}
Gonz\'alez-Diez, Gabino; Jaikin-Zapirain, Andrei,
{\em The absolute Galois group acts faithfully on regular dessins and on Beauville surfaces},
preprint 2014.

\bibitem[GoMP88]{g-mp}
Goresky, Mark; MacPherson, Robert {\em Stratified Morse theory. }
Ergebnisse der Mathematik und ihrer Grenzgebiete (3) 14. Springer-Verlag, Berlin, (1988) xiv+272 pp..



\bibitem[GR58]{g-r}
H. Grauert, R. Remmert,
{\em Komplexe R\"aume,}
{\bf Math. Ann. 136} (1958), 245--318.

\bibitem[Gra72]{grauert1}
H. Grauert,
{\em \"Uber die Deformation isolierter Singularit\"aten analytischer Mengen,}
{\bf Invent. Math. 15} (1972), 171-198.


\bibitem[Gra74]{grauert}
H. Grauert,
{\em Der Satz von Kuranishi f\"ur kompakte komplexe R\"aume,}
{\bf Invent. Math. 25}(1974), 107-142 .

\bibitem[Green88]{green}
M. Green, {\em  Restrictions of linear series to hyperplanes, and some results of Macaulay and Gotzmann},
 Algebraic curves and projective geometry (Trento, 1988), Lecture Notes in Math., 1389, Springer, Berlin, (1989),  76--86. 
 
 \bibitem[Greenb67]{greenberg}
Greenberg, Marvin J. , {\em Lectures on algebraic topology. } W. A. Benjamin, Inc., New York-Amsterdam x+235 pp., (1967 ). 

Second edition:
 Greenberg, Marvin J.; Harper, John R. {\em Algebraic topology. A first course. } Mathematics Lecture Note Series, 58. Benjamin/Cummings Publishing Co., Inc., Advanced Book Program, Reading, Mass., xi+311 pp. (1981). 
 
  \bibitem[G-L87]{G-L1}
  Green, Mark; Lazarsfeld, Robert 
  {\em Deformation theory, generic vanishing theorems, and some conjectures of Enriques, Catanese and Beauville.}
   Invent. Math. 90 , no. 2, 389--407 (1987). 
   
    \bibitem[G-L91]{G-L2}
    Green, Mark; Lazarsfeld, Robert 
    {\em Higher obstructions to deforming cohomology groups of line bundles.}
     J. Amer. Math. Soc. 4 , no. 1, 87--103 (1991).

\bibitem[Griff-68]{griff1}
P. Griffiths, {\em Periods of integrals on algebraic manifolds.I. Construction and properties of the modular varieties. -II. Local study of the period mapping.}  Amer. J. Math. 90 (1968) 568--626 and 805--865.

\bibitem[Griff-70]{griff2}
P. Griffiths, {\em  Periods of integrals on algebraic manifolds. III. Some global differential-geometric properties of the period mapping.} Inst. Hautes \'Etudes Sci. Publ. Math. No. 38 (1970) 125--180.

\bibitem[G-H78]{gh}
 Griffiths, Phillip; Harris, Joseph 
 {\em Principles of algebraic geometry.} Pure and Applied Mathematics. Wiley-Interscience [John Wiley and Sons], New York, (1978). xii+813 pp

\bibitem[G-M-81]{GriffithsMorgan}
P. Griffiths, J. Morgan, {\em   Rational homotopy theory and differential forms.}
 Progress in Mathematics, 16. Birkh\"auser, Boston, Mass., (1981). xi+242 pp.

\bibitem[Grom89]{Gromov}
M. Gromov,
{\em Sur le groupe fondamental d'une vari\'et\'e k\"ahl\'erienne,}
{\bf C. R. Acad. Sci.,
Paris, S\'er. I 308, No.3}(1989), 67-70.

\bibitem[Gro91]{GromovL2}
 Gromov, M. 
 {\em K\"ahler hyperbolicity and $L^2$ -Hodge theory. }
 J. Differential Geom. 33 , no. 1, 263--292 (1991).

 \bibitem[G-S92]{G-S}
 Gromov, Mikhail,  Schoen, Richard,
 {\em  Harmonic maps into singular spaces and p-adic superrigidity for lattices in groups of rank one. }
 Inst. Hautes \'Etudes Sci. Publ. Math. No. 76 , 165--246 (1992).

\bibitem[Groth57]{tohoku}
Grothendieck, A.
{\em Sur quelques points d'alg\`ebre homologique. }
Tohoku Math. J., II. Ser. 9, 119--221 (1957).

\bibitem[Groth60]{groth}
A. Grothendieck,
{\em Techniques de construction et th\'eoremes d'existence en g\'eom\'etrie
alg\'ebrique. IV, Les schemas de Hilbert}, Sem. Bourbaki, Vol. 13, 
(1960-61), 1--28.

\bibitem[SGA1]{sga1} A. Grothendieck (dirig\'e par), {\em Rev$\hat{e}$tements \'etales et groupe
fondamental,} S\'eminaire de g\'eom\'etrie alg\'ebrique du Bois Marie 1960-61 {\bf Springer Lecture Notes
in Math. 224}, (1971). Reedited by {\bf Soci\'et\'e Math\'ematique de France} in the series 'Documents
math\'ematiques (2003). 

\bibitem[Gro97]{groth2} A. Grothendieck,  {\em Esquisse d'un programme}. In: 'Geometric Galois actions,
1', {\bf London Math. Soc. Lecture Note Ser., 242. Cambridge Univ. Press, Cambridge} (1997), 5--48.

\bibitem[GHT13]{ght}
Samuel Grushevsky, Klaus Hulek, Orsola Tommasi,
{\em Stable cohomology of the perfect cone toroidal compactification of the moduli space of abelian varieties},
arXiv:1307.4646 .

\bibitem[Guena13]{guenancia}
H. Guenancia, {\em K\"ahler-Einstein metrics with cone singularities on klt pairs.}
 Internat. J. Math. 24 (2013), no. 5, 1350035, 19 pp.

\bibitem[Hac-McK06]{hm-pluric}
Hacon, Christopher D.; McKernan, James 
{\em Boundedness of pluricanonical maps of varieties of general type. }
Invent. Math. 166 (2006), no. 1, 1--25.

\bibitem[Ham82]{hamilton}
Hamilton, Richard S.
{\em Three-manifolds with positive Ricci curvature. }
J. Differ. Geom. 17, 255--306 (1982). 

 \bibitem[Hano57]{Hano} {\sc  Hano, J.:} {\it On Kaehlerian homogeneous spaces of unimodular Lie groups.}
Amer. J. Math. 79 (1957), 885-900.

\bibitem[Har85]{harer}
Harer, John L. 
{\em Stability of the homology of the mapping class groups of orientable surfaces.} Ann. of Math. (2) 121 (1985), no. 2, 215--249.

\bibitem[Hartm67]{hartmann}
Hartmann, P. ,
{\em  On homotopic harmonic maps,} Canad. J. Math. 19 , 673--687 (1967).

\bibitem[Harts74]{hartshorne}
Hartshorne, Robin
{\em Varieties of small codimension in projective space. }
Bull. Amer. Math. Soc. 80 (1974), 1017--1032. 

\bibitem[Hart77]{hart}
R. Hartshorne,
{\em  Algebraic geometry,}
{\bf Springer GTM 52} (1977).

\bibitem[HT80]{h-t}
A. Hatcher, W.Thurston,
{\em A presentation for the mapping class group of a closed orientable
    surface,}{\bf Topology 19, no. 3}(1980), 221--237.

\bibitem[Helga62]{Helgason} {\sc  Helgason, S.:} {\it Differential
geometry and symmetric spaces.}
   Pure and Applied Mathematics, Vol. XII. Academic Press, New
York-London (1962) xiv+486 pp.

   \bibitem[Helga78]{Helgason2} {\sc  Helgason, S.:} {\it  Differential
geometry, Lie groups, and symmetric spaces. } Pure and Applied
Mathematics,
80. Academic Press, Inc. [Harcourt Brace Jovanovich, Publishers], New
York-London, (1978) xv+628 pp.

\bibitem[H-S71]{h-s}
 P. J. Hilton,U.  Stammbach, {\em  A course in homological algebra.} Second edition. Graduate Texts in Mathematics,Vol. 4. Springer-Verlag, New York-Berlin, (1971), ix+338 pp.
 
  \bibitem[Hopf42]{Hopf}
Hopf, H.: {\em Fundamentalgruppe und zweite {B}ettische {G}ruppe}.
 {Comment. Math. Helv.}
 {14} {(1942)}, {257--309}.
 
\bibitem[Hor75]{quintics}
E. Horikawa,
{\em On deformations of quintic surfaces,}
{\bf Invent. Math. 31} (1975), 43--85.

\bibitem[Hub06]{hubbard}
J. H. Hubbard, {\em Teichm\"uller theory and applications to geometry, topology, and dynamics. Vol. 1. Teichm\"uller theory. } Matrix Editions, Ithaca, NY, (2006), xx+459 pp. 

\bibitem[Hum75]{hum}
J. E. Humphreys,
{\em Linear algebraic groups, }
{\bf Graduate Texts in Mathematics 21,  Springer-Verlag}
New
York - Heidelberg - Berlin,  (1975), XV, 247 p. .

\bibitem[H-W48]{H-W}
Hurewicz, Witold; Wallman, Henry
{\em Dimension theory. }
Princeton: Princeton University Press. 165 p. (1948).

\bibitem[Hur91]{Hurwitz}
Hurwitz, A.:{ \em Ueber Riemann'schen Fl\"achen 
mit gegebenen Verzweigungspunkten.} Math. Ann. 
39, 1--61 (1891).


\bibitem[Huy05]{huy}
D. Huybrechts,
{\em Complex geometry. An introduction.}
{\em  Universitext.
Springer-Verlag}, Berlin, (2005), xii+309 pp.

\bibitem[Ig72]{igusa}
Igusa, Jun-ichi 
{\em Theta functions.}
 Die Grundlehren der mathematischen Wissenschaften, Band 194. Springer-Verlag, New York-Heidelberg, (1972). x+232 pp.

\bibitem[In94]{inoue} Inoue, M.  {\em Some new 
surfaces of general 
type.} Tokyo J. Math. 17 (1994), no. 2, 
295--319.

 \bibitem[Jac80]{BAII}
 Jacobson, Nathan {\em Basic algebra. II. } W. H. Freeman and Co., San Francisco, Calif., xix+666 pp (1980()>

\bibitem[Kas-Schl72]{k-s}
A. Kas, M. Schlessinger, 
{\em On the versal deformation of a complex space with an isolated singularity.}
Math. Ann. 196 (1972), 23--29.


\bibitem[J-R87]{J-R}
Johnson, F. E. A.; Rees, E. G.,
{\em On the fundamental group of a complex algebraic manifold.} 
Bull. London Math. Soc. 19 , no. 5, 463--466 (1987). 


\bibitem[J-Y83]{J-Y83}
Jost, J\"urgen, Yau, Shing-Tung,
{\em Harmonic mappings and  K\"ahler manifolds. }
Math. Ann. 262, 145--166 (1983).

\bibitem[J-Y85]{J-Y85}
Jost, J\"urgen, Yau, Shing-Tung,
{\em A strong rigidity theorem for a certain class of compact complex analytic surfaces. }
  Math. Ann. 271, 143--152 (1985).
 
\bibitem[J-Y93]{JostYau}
Jost, J\"urgen, Yau, Shing-Tung,
{\em Applications of quasilinear PDE to algebraic geometry and arithmetic lattices.} Algebraic geometry and related topics (Inchon, 1992), 169--193,
 Conf. Proc. Lecture Notes Algebraic Geom., I, Int. Press, Cambridge, MA, (1993).

\bibitem[J-Z96]{j-z1}
Jost, J\"urgen; Zuo, Kang,
{\em  Harmonic maps and $Sl(r,\CC)$-representations of fundamental groups of quasiprojective manifolds.}  J. Algebraic Geom. 5 (1996), no. 1, 77--106.

\bibitem[J-Z97]{j-z2}
Jost, J\"urgen; Zuo, Kang,
{\em Harmonic maps of infinite energy and rigidity results for representations of fundamental groups of quasiprojective varieties.}
 J. Differential Geom. 47 (1997), no. 3, 469--503.
 
 
 \bibitem[J-Z00]{j-z3}
Jost, J\"urgen; Zuo, Kang,
{\em Harmonic maps into Bruhat-Tits buildings and factorizations of p-adically unbounded representations of $\pi_1$ of algebraic varieties. I. }
J. Algebraic Geom. 9 (2000), no. 1, 1--42.

\bibitem[Katata83]{katata}
{\em Open problems: Classification of algebraic and analytic manifolds.}
Classification of algebraic and analytic manifolds, Proc. Symp.
Katata/Jap. 1982.
   Edited by Kenji Ueno.
{\bf Progress in Mathematics, 39.
Birkh\"auser, Boston, Mass.} (1983), 591-630.

\bibitem[Kato90]{Kato}
Kato, Masahide,
{\em A non K\"ahler structure on an $S^2$-bundle over a ruled surface},
Unpublished manuscript, 15 pages, May 1992.

\bibitem[KR98]{Katz-Ram}
Katzarkov, L.; Ramachandran, M. {\em On the universal coverings of algebraic surfaces. }
Ann. Sci. \'Ecole Norm. Sup. (4) 31 (1998), no. 4, 525--535.

 
 \bibitem[Kazh70]{Kazh70} {\sc Kazhdan D.:} {\it Arithmetic varieties
and their fields of quasi-definition.}, Actes du Congr\`es
International des
Mathematiciens (Nice, 1970), Tome 2, pp. 321 --325. Gauthier-Villars,
Paris, (1971).

 \bibitem[Kazh83]{Kazh83}
Kazhdan, David
{\em On arithmetic varieties. II. }
Israel J. Math. 44 (1983), no. 2, 139--159. 

\bibitem[Ker83]{kerckhoff}
Kerckhoff, Steven P. {\em The Nielsen realization problem.}
 Ann. of Math. (2) 117, no. 2, 235--265 (1983). 
 
 \bibitem[Ke88]{keum} Keum, Y.H. {\it Some new surfaces of 
general 
type with $p_g = 0$.}  Unpublished manuscript 
(1988).

\bibitem[KK02]{k-k}
V. M. Kharlamov, V.S. Kulikov,
{\it On real structures of real surfaces}.
Izv. Ross. Akad. Nauk Ser. Mat. \textbf{66}, no. 1, 133--152 (2002);
translation in Izv. Math. \textbf{66}, no. 1, 133--150 (2002)

\bibitem[KK13]{k-k2}
  Kharlamov, V. M.; Kulikov, Vik. S., 
{\em The semigroups of coverings. }
 Izv. Ross. Akad. Nauk Ser. Mat. 77 (2013), no. 3, 163--198; translation in Izv. Math. 77 (2013), no. 3, 594--626.

\bibitem[Kling13]{klingler}
Klingler, Bruno
{\em Symmetric differentials, K\"ahler groups and ball quotients, }
Invent. Math. 192 (2013), no. 2, 257--286. 


\bibitem[Koba58]{Kobayashi} {\sc Kobayashi, S.:} {\it  Geometry of
bounded domains. }Trans. Amer. Math. Soc. 92 (1959) 267-290.

\bibitem[Koba80]{Koba80} {\sc Kobayashi, S.:}
  {\it First Chern class and holomorphic tensor fields.} Nagoya Math. J.
77  (1980), 5 - 11.

\bibitem[Koba80-2]{Koba80-2} {\sc Kobayashi, S.:} {\it The first Chern
class and holomorphic symmetric tensor fields.}, J. Math. Soc. Japan
32 (1980),
no. 2, 325-329.

\bibitem[KobNom63]{KobNom} {\sc
Kobayashi, S. and Nomizu, K.:} {\it Foundations of differential geometry.
Vol I. } Interscience Publishers, a division of John Wiley \& Sons, New York-Lond on (1963) xi+329 pp.



\bibitem[Kod-Mor72]{Kodaira-Morrow} {\sc  Kodaira,K.  and Morrow, J. :}
   {\it  Complex manifolds.} Holt, Rinehart and Winston, Inc., New
York-Montreal, Que.-London, (1971) vii+192 pp.




\bibitem[Kod54]{kodemb}
K. Kodaira, {\em On K\"ahler varieties of restricted type (an intrinsic characterization of algebraic varieties). }
{\bf Ann. of Math. (2) 60} (1954), 28--48.

\bibitem[K-S58]{k-s58}
K.  Kodaira, ; D. C. Spencer,
{\em On deformations of complex analytic structures. I,
     II,} {\bf  Ann. of Math. (2) 67 } (1958), 328 -466.
     
\bibitem[Kod60]{kod1}
K. Kodaira, {\em On compact complex analytic surfaces, I. }, Ann. Math. 71, 111--152 (1960).



\bibitem[Kod63-b]{kod-1}
K. Kodaira,
{\em On stability of compact submanifolds of complex manifolds,}
{\bf Am. J. Math. 85}(1963), 79-94 .

\bibitem[Kod67]{Kodsurf}
K. Kodaira, {\em  A certain type of irregular algebraic surfaces. } {\bf J. Analyse Math. 19} (1967), 207--215. 
 
\bibitem[Kod-Mor71]{k-m71}
 K. Kodaira, J. Morrow,
 {\em Complex manifolds.} Holt, Rinehart and Winston, Inc., New York-Montreal, Que.-London, (1971). vii+192 pp.
 
 \bibitem[Kod86]{Kodbook}
 K. Kodaira,{\em  Complex manifolds and deformation of complex structures. }
 Translated from the Japanese by Kazuo Akao. With an appendix by Daisuke Fujiwara. Grundlehren der Mathematischen Wissenschaften, 283. Springer-Verlag, New York, (1986) x+465 pp.

\bibitem[K-SB88]{k-sb}  J. Koll\'ar and N.I. Shepherd Barron, {\em
Threefolds and deformations of surface singularities. }
{\bf Inv. Math. 91}  (1988) 299-338.

\bibitem[Kol93]{Kol-Shaf-art}
J. Koll\'ar  {\em Shafarevich maps and plurigenera of algebraic varieties.}
 Invent. Math. 113 (1993), no. 1, 177--215. 

\bibitem[Kol95]{Kollar-Shaf} {\sc  Koll\'ar, J. :} {\em  Shafarevich maps
and automorphic forms. } M. B. Porter Lectures. Princeton University
Press,
Princeton, NJ, (1995) x+201 pp.



\bibitem[Kol13]{Kolhand}
J. Koll\'ar  {\em Moduli of varieties of general type. }
Handbook of moduli. Vol. II, 131--157, Adv. Lect. Math. (ALM), 25, Int. Press, Somerville, MA, (2013).

\bibitem[K-M]{KollarMori}
 Koll\'ar, J., Mori, S. {\it  Birational geometry of algebraic varieties. }
 With the collaboration of C. H. Clemens and A. Corti. 
 Cambridge Tracts in Mathematics, 134. Cambridge University Press, Cambridge, (1998) viii+254 pp.
 
 \bibitem[Kot99]{kotschick}
D. Kotschick, {\em On regularly fibered complex surfaces. } Proceedings
of the Kirbyfest (Berkeley, CA, 1998),  291--298 (electronic), Geom.
Topol. Monogr., 2,
Geom. Topol. Publ., Coventry, (1999).
 
\bibitem[Kur62]{kur1}
M. Kuranishi,
{\em On the locally complete families of complex analytic structures,}
{\bf Ann. Math. (2) 75} (1962), 536-577.

\bibitem[Kur65]{kur2}
M. Kuranishi,
{\em New proof for the existence of locally complete families
of complex structures,}
Proc. Conf. Complex Analysis, Minneapolis 1964 (1965), 142-154.

\bibitem[Kur69]{kur3}
M. Kuranishi,
{\em A note on families of complex structures. }
Global Analysis (Papers in Honor of K. Kodaira), (1969),  pp. 309-313, Princeton Univ. Press and Univ. Tokyo Press, Tokyo.

\bibitem[Kur71]{kur4}
M. Kuranishi,
{\em Deformations of compact complex manifolds.} S\'eminaire de Math\'ematiques Sup\'erieures, No. 39 (\'Et\'e 1969).
 Les Presses de l'Universit\'e de Montr\'eal, Montreal, Que., (1971), 99 pp.
 
 \bibitem[LanBir]{LangeBirkenhake} Lange, H. ; Birkenhake, C. 
{\it 
Complex abelian varieties. }
Grundlehren der
Mathematischen 
Wissenschaften,  302. Springer-Verlag, Berlin, (1992), 
viii+435 
pp.;
  second edition (2004)  xii+635 pp.
  
 \bibitem[Lan01]{Lange}
  Lange, Herbert
{\em Hyperelliptic varieties. }
Tohoku Math. J. (2) 53 (2001), no. 4, 491--510. 

\bibitem[LaRe04]{LR}
Lange, Herbert, Recillas, Sevin
{\em Abelian varieties with group action. }
J. Reine Angew. Math. 575 (2004), 135--155.
 
  \bibitem[Lars73]{Larsen}
  Larsen, Mogens Esrom {\em On the topology of complex projective manifolds}, Invent. Math. 19 (1973), 251--260.
  
  \bibitem[Laz]{laz}
  Lazarsfeld, Robert  {\em Positivity in algebraic geometry. I. Classical setting: line bundles and linear series.  II. Positivity for vector bundles, and multiplier ideals.}
  Ergebnisse der Mathematik und ihrer Grenzgebiete. 3. Folge, 48 and 49. Springer-Verlag, Berlin, (2004) xviii+387 resp. xviii+385 pp.  
  
 \bibitem[LB97]{LeBrunAMS}
   LeBrun, Claude,
   {\em  Twistors for tourists: a pocket guide for algebraic geometers.}
    Algebraic geometry--Santa Cruz 1995, 361--385, Proc. Sympos. Pure Math., 62, Part 2, Amer. Math. Soc., Providence, RI, (1997).
 


 \bibitem[Lef21]{Lefschetz}
 Lefschetz, Solomon {\em On certain numerical invariants of algebraic varieties with application to abelian varieties,}
  Trans. Amer. Math. Soc. 22 (1921), no. 3, 327--406.
  
    \bibitem[Lef24]{Lefschetzbook}
    Lefschetz, S.{\em  L'analysis situs et la g\'eom\'etrie alg\'ebrique. } Gauthier-Villars, Paris, (1924) vi+154 pp..

\bibitem[LW12]{wenfei} Liu, Wenfei, { \it Stable Degenerations of
Surfaces Isogenous to a Product II.}
	Trans. A.M.S. 364 (2012), n. 5, 2411--2427.
	
\bibitem[Liv85]{Liv}
Livingston, C.: {\em Stabilizing surface symmetries}. Mich. Math. J. 32, 249--255 (1985).

\bibitem[LP14]{lp}
L\"onne, Michael, Penegini, Matteo
{\em On asymptotic bounds for the number of irreducible components of the moduli space of
surfaces of general type}, arXiv:1402.6438.

\bibitem [M-W07] {m-w} I. Madsen, M. Weiss, 
{\it The stable moduli space of Riemann surfaces: Mumford's conjecture.} Ann. of Math. (2) 165 (2007), no. 3, 843--941.

\bibitem[MSSV01]{mssv}
Magaard, K., Shaska, T., Shpectorov, S.,  V\"olklein, H. ,
{\em The locus of curves with prescribed automorphism group.} Communications in arithmetic fundamental groups (Kyoto, 1999/2001). S$\bar{u}$rikaisekikenky$\bar{u}$sho K$\bar{o}$ky$\bar{u}$roku No. 1267 (2002), 112--141.

 \bibitem[MKS66]{Magnusks}
 Magnus, Wilhelm; Karrass, Abraham; Solitar, Donald {\em Combinatorial group theory: 
 Presentations of groups in terms of generators and relations.}
  Interscience Publishers [John Wiley - Sons, Inc.], New York-London-Sydney (1966) xii+444 pp. 
  (Dover reprint , (1976)).

\bibitem[Man97]{man3}
M. Manetti, {\em Iterated Double Covers and Connected
   Components of Moduli Spaces},
{\bf  Topology 36, 3} (1997) 745-764.

\bibitem[Man01]{man4}
M. Manetti, {\em On the Moduli Space of diffeomorphic algebraic
surfaces},
{\bf  Inv. Math. 143} (2001), 29-76.

\bibitem[Man04]{Manobs}
M. Manetti, {\em Cohomological constraint on deformations of compact K\"ahler manifolds.}
 Adv. Math. 186 (2004), no. 1, 125--142.
 
 \bibitem[Man09]{man09}
M. Manetti, {\em Differential graded Lie algebras and formal deformation theory.}
 Algebraic geometry-Seattle 2005. Part 2, 785--810, 
Proc. Sympos. Pure Math., 80, Part 2, Amer. Math. Soc., Providence, RI, (2009). 

\bibitem[Mark07]{markovic}
Markovic, Vladimir,
{\em Realization of the mapping class group by homeomorphisms. }
Invent. Math. 168 (2007), no. 3, 523--566.

\bibitem[Math70]{mather}
J. Mather, {\em  Notes on topological stability },
{\bf Harvard University Math. Notes }(1970).

\bibitem[Mats63]{matsumura}
H. Matsumura, {\em
On algebraic groups of birational transformations} 
Atti Accad. Naz. Lincei Rend. Cl. Sci. Fis. Mat. Natur. (8) 34 (1963) 151--155. 

\bibitem[Meer13]{meerssemann}
Laurent Meersseman,
{\em Global moduli space of complex structures, groupoids and analytic stacks},
arXiv 1311.4170v1 .

\bibitem[Migl95]{miglio}
Migliorini, Luca 
{\em A smooth family of minimal surfaces of general type over a curve of genus at most one is trivial.}
 J. Algebraic Geom. 4 (1995), no. 2, 353--361.
     


\bibitem[Mil52]{miller}
Miller, Clair,
{\em  The second homology group of a group; relations among commutators.}
 Proc. Amer. Math. Soc. 3, (1952) 588�595. 

\bibitem[Milne01]{milne}
Milne, J.S. 
{\em Kazhdan's Theorem on Arithmetic Varieties,}   arXiv:math/0106197.

\bibitem[Mil63]{milnorMT}
J. W. Milnor,
{\em  Morse theory. } Based on lecture notes by M. Spivak and R. Wells. Annals of Mathematics Studies, No. 51 Princeton University Press, 
Princeton, N.J. (1963) vi+153 pp. .

\bibitem[Mil68]{milnor}
J. W. Milnor,
{\em Singular points of complex hypersurfaces,}
{\bf Annals of Mathematics Studies 61, Princeton University Press
and the University of
Tokyo Press} (1968),  122 p. .

\bibitem[Mil71]{Milnor}
J. W. Milnor,
{\em Introduction to algebraic K-theory,}
{\bf Annals of Mathematics Studies 72, Princeton University Press} (1971).

 \bibitem[M-S]{milnorCC}
 Milnor, John W.; Stasheff, James D.
{\em Characteristic classes. }
Annals of Mathematics Studies. No.76. Princeton, N.J.: Princeton University Press and University of Tokyo Press. VII, 331 p.  (1974).

 \bibitem[Mil76]{milnorcurv}{ Milnor J.  :} {\it  Curvatures of left invariant metrics on Lie groups.}
Adv. Math. 21 (3) (1976), 293--329.


\bibitem[Moi94]{moi}
B. Moishezon, {\em  The arithmetics of braids and a statement of Chisini},
in `Geometric Topology, Haifa 1992' , {\bf  Contemp. Math. 164,
Amer. Math. Soc., Providence, RI}
(1994), 151-175.

\bibitem[Mok85]{Mok1} { Mok, Ngaiming. :} {\it The holomorphic or anti-holomorphic character of harmonic maps into irreducible compact quotients of polydiscs. } Math.
Ann. 272  , 197--216 (1985).

\bibitem[Mok88]{Mok2} { Mok, Ngaiming. :} {\it Strong rigidity  of  irreducible  quotients of polydiscs of finite volume. } Math.
Ann. 282  , 555-578 (1988).
 
\bibitem[Mok02]{Mok} { Mok, N. :} {\it Characterization of
certain
holomorphic geodesic cycles on quotients of bounded
symmetric domains
in
terms of tangent subspaces. } Compositio Math.
132 (2002), no. 3, 289-309.

\bibitem[Mor12]{moravec}
Moravec, Primoz 
{\em Unramified Brauer groups of finite and infinite groups.}
 Amer. J. Math. 134 (2012), no. 6, 1679--1704.

\bibitem[Mor96]{mor}
   J. W. Morgan,{\em
The Seiberg-Witten equations and applications to the topology of
   smooth four-manifolds.}
{\bf Mathematical Notes  44.
    Princeton Univ. Press}  vi (1996).
    
  \bibitem[M-T]{morgan-tian} 
    Morgan, John; Tian, Gang
{\em Ricci flow and the Poincar\'e conjecture.}
Clay Mathematics Monographs 3. Providence, RI: American Mathematical Society (AMS); Cambridge, MA: Clay Mathematics Institute. xlii, 521~p. (2007).

\bibitem[Most73]{mostow}
G.D. Mostow,
{\em Strong rigidity of locally symmetric spaces,}
{\bf Annals of Mathematics Studies, no. 78,
Princeton University Press, Princeton, N.J.;
University of Tokyo Press, Tokyo} (1973).

\bibitem[Mum62]{mum1}
D. Mumford, {\em The canonical ring of an algebraic surface},
  {\bf   Ann. of Math. (2)  76} (  1962) 612--615,
appendix to
{\em The theorem of Riemann-Roch for high multiples of an effective divisor on
an algebraic surface}, by O. Zariski, ibid. 560-611.

\bibitem[Mum65]{GIT}
D. Mumford, {\em Geometric invariant theory} Ergebnisse der Mathematik und ihrer Grenzgebiete, Neue Folge, Band 34 Springer-Verlag, Berlin-New York (1965) vi+145 pp.

\bibitem[Mum70]{abvar}
D. Mumford, {\em Abelian varieties}, Tata Institute of Fundamental Research Studies in Mathematics, No. 5 
Published for the Tata Institute of Fundamental Research, Bombay; Oxford University Press, London (1970) viii+242 pp.

\bibitem[Mum77a]{EnsMath} Mumford, David {\em Stability of projective varieties.}  Enseignement Math. (2) 23 , no. 1-2, 39--110 (1977).

\bibitem[Mum77b]{Hirzprop} Mumford, D. {\em Hirzebruch's proportionality theorem in the
noncompact case. }{\bf  Invent. Math.  42}  (1977), 239--272.

\bibitem[Mum76]{mumErgebnisse}
Mumford, David
{\em Algebraic geometry. I: Complex projective varieties.} 
Grundlehren der mathematischen Wissenschaften, 221. Berlin-Heidelberg-New York: Springer-Verlag. X, 186 p. (1976).

\bibitem[Mum83]{mumshaf}
D. Mumford,
{\em Towards an enumerative
geometry of the moduli space of curves,}
in `Arithmetic and geometry, Pap. dedic. I. R. Shafarevich,
Vol. II: Geometry', {\bf Prog. Math. 36, Birkh\"auser}  (1983), 271-328.

\bibitem[Mum-Suom72]{Mum-Suom}
D. Mumford, K. Suominen, 
{\em  Introduction to the theory of moduli. } in `Algebraic geometry, Oslo 1970' (Proc. Fifth Nordic
Summer-School in Math.),  pp. 171-222. Wolters-Noordhoff, Groningen, (1972). 

\bibitem[Mum11]{mum-notices}
D. Mumford, {\em Intuition and rigor and Enriques's quest.}
 Notices Amer. Math. Soc. 58 (2011), no. 2, 250--260. 

\bibitem[Nai94]{naie94} Naie, D. {\it Surfaces 
d'Enriques et une 
construction de surfaces de type g\'en\'eral avec 
$p\sb g=0$}.
Math. Z.  215  (1994),  no. 2, 269--280.

\bibitem[Nami80]{Namikawa}
Namikawa, Yukihiko 
{\em Toroidal compactification of Siegel spaces.} Lecture Notes in Mathematics, 812. Springer, Berlin(1980), viii+162 pp..

\bibitem[New-Nir57]{NN}
A. Newlander, L. Nirenberg, 
{\em Complex analytic coordinates in almost complex manifolds.}
Ann. of Math. (2) 65 (1957), 391-404.

\bibitem[Nielsen]{nielsen}
J. Nielsen,  {\em Untersuchungen zur Topologie der geschlossenen zweiseitigen Fl\"achen I, II, III. } Acta Math. 50 (1927), no. 1, 189--358,
53 (1929), no. 1, 1--76, 58 (1932), no. 1, 87�167.

\bibitem[Nielsen43]{nielsen2}
J. Nielsen,  {\em Abbildungsklassen endlicher Ordnung}, Acta Math. 75, (1943) 23--115. 

\bibitem[Nielsen44]{nielsen3}
J. Nielsen,  {\em Surface transformation classes of algebraically finite type,}  Danske Vid. Selsk. Math.-Phys. Medd. 21, (1944). no. 2, 89 pp.


\bibitem[OS01]{Oguiso}
K. Oguiso, J. Sakurai, Jun
{\em Calabi-Yau threefolds of quotient type. } 
Asian J. Math. 5 (2001), no. 1, 43--77. 

\bibitem[Par91]{Pardini}
Pardini, Rita 
{\em Abelian covers of algebraic varieties. }
J. Reine Angew. Math. 417 (1991), 191--213.



\bibitem[Pe-Pol10]{pepo}
M.  Penegini, F. Polizzi {\em On surfaces with $p_g=q=2, K^2=5$ and Albanese map of degree 3},
 Osaka J. Math. 50, n. 3 (2013), 643--686.

\bibitem[Pe-Pol11]{pepo2}
M.  Penegini, F. Polizzi {\em Surfaces with $p_g=q=2, K^2=6$ and Albanese map of degree 2},
 Canad. J. Math. 65, n.1 (2013), 195--221.





\bibitem[Per1]{perelman1}
Grisha Perelman, 
{\em The entropy formula for the Ricci flow and its geometric applications}
arXiv:math/0211159 

\bibitem[Per2]{perelman2}
Grisha Perelman, {\em
Ricci flow with surgery on three-manifolds}
arXiv:math/0303109

\bibitem[Per3]{perelman3}
Grisha Perelman, {\em Finite extinction time for the solutions to the Ricci flow on certain three-manifolds}
arXiv:math/0307245 

\bibitem[PPS87]{PPS}
Peternell, Th.; Le Potier, J.; Schneider, M. 
{\em Vanishing theorems, linear and quadratic normality. } Invent. Math. 87 (1987), no. 3, 573--586.


\bibitem[Pont96]{pontecorvo}
Pontecorvo, Massimiliano, 
{\em On the complex geometry of twistor spaces. } Geometry Seminars, 1994--1995 (Italian) (Bologna), 195--216, Univ. Stud. Bologna, Bologna, (1996). 


\bibitem[Quinn86]{quinn}
F. Quinn, {\em Isotopy of 4-manifolds} J. Differential Geom. 24 (1986), no. 3, 343--372.

\bibitem[Reid87]{youngguide}
M. Reid,
     {\em  Young person's guide to canonical singularities,}
in ' Algebraic geometry, Proc. Summer Res. Inst.,
     Brunswick/Maine 1985', part 1, {\bf Proc. Symp. Pure Math. 46}
(1987), 345-414 .

\bibitem[Ries93]{ries}
Ries, John F. X. , {\em Subvarieties of moduli space determined by finite groups acting on surfaces.}
 Trans. Amer. Math. Soc. 335 (1993), no. 1, 385--406.
 
 
 \bibitem[Rol10]{soenkecomp}
S. Rollenske, {\em 
Compact moduli for certain Kodaira fibrations}, {\bf  Ann. Sc. Norm. Super. Pisa Cl. Sci. IX (4)}, (2010) 851--874.

\bibitem[Roos00]{Roos} {\sc
Roos, G.:} {\it Jordan Triple Systems},
pp.\ 183--282, in J.\
Faraut, S.\ Kaneyuki, A.\ Kor\'{a}nyi, Q.k.~Lu,
G.\ Roos,
Analysis
and geometry on complex homogeneous domains, Progress in
Mathematics, vol.\ \textbf{185}, Birkh\"{a}user, Boston
(2000).

\bibitem[Sam78]{sampson}
Sampson, J.H.
{\em Some properties and applications of harmonic mappings. }
Ann. Sci. \'Ec. Norm. Sup\'er. (4) 11, No. 2, 211--228 (1978). 


\bibitem[Sam86]{sampsonK}
Sampson, J.H.
{\em 
Applications of harmonic maps to K\"ahler geometry. }
Contemp. Math. 49, 125--134 (1986).

\bibitem[Seid08]{seidelcime}
 Seidel, Paul 
 {\em Lectures on four-dimensional Dehn twists.} in: ` Symplectic 4-manifolds and algebraic surfaces', 231--267,
  Lecture Notes in Math., 1938, Springer, Berlin, (2008).

\bibitem[S-T34]{ST}
Seifert, H., Threlfall, W., 
{\em Lehrbuch der Topologie. } 
IV + 353 S. 132 Abb. Leipzig, B. G. Teubner (1934).

\bibitem[Sern06]{sernesi}
E. Sernesi,{\em Deformations of algebraic schemes. }
Grundlehren der Mathematischen Wissenschaften, 334. Springer-Verlag, Berlin, (2006) xii+339 pp.


\bibitem[Ser58]{ser5}
Serre, Jean-Pierre,
{\em Sur la topologie des vari\'et\'es alg\'ebriques en caract\'eristique p. },
 1958 Symposium internacional de topolog\'ia algebraica,  pp. 24--53,  Universidad Nacional Aut\'onoma de M\'exico and UNESCO, Mexico City 

\bibitem[Ser64]{serre}
Serre, Jean-Pierre,
{\em Exemples de vari\'et\'es projectives conjugu\'ees non hom\'eomorphes}.
{\bf C. R. Acad. Sci. Paris 258},  4194--4196, (1964).

\bibitem[Shaf74]{shaf}
 Shafarevich, I. R.
{\em Basic algebraic geometry.} Translated from the Russian by K. A.
Hirsch. Revised printing of
{\bf Grundlehren der mathematischen Wissenschaften, Vol.
213, 1974. Springer Study Edition. Springer-Verlag}, Berlin-New York,
(1977), xv+439 pp.

\bibitem[S-T54]{s-t}
Shephard, G. C., Todd, J. A. 
{\em Finite unitary reflection groups}, Canadian J. Math. 6: 274--304 (1954).

\bibitem[Sie48]{Siegel} {\sc Siegel, C.L.:}
  { \em Analytic
Functions of Several Complex Variables.}
   Notes by P. T. Bateman.
Institute for Advanced Study, Princeton,
N.J., (1950)  ii+200 pp.

\bibitem[Sieg73]{sieg}
C.L. Siegel,
{\em Topics in complex function theory. Vol. III. Abelian functions and
modular functions of several variables. }
Translated from the German by E. Gottschling
and M. Tretkoff. With a preface by Wilhelm Magnus. Reprint of the
1973 original.
{\bf Wiley Classics Library. A Wiley-Interscience Publication. John
Wiley and Sons,
Inc.}, New York (1989), x+244 pp.



\bibitem[Simp92]{simpsonHiggs}
Simpson, Carlos T.,
{\em  Higgs bundles and local systems.} Inst. Hautes \'Etudes Sci. Publ. Math. No. 75, 5--95 (1992).


\bibitem[Simp93]{simpsonTors}
Simpson, Carlos T.,
 {\em Subspaces of moduli spaces of rank one local systems.}  Ann. Sci. \'Ecole Norm. Sup. (4) 26, no. 3, 361--401(1993).
 
\bibitem[Sing72]{singer}
D. Singerman, {\em Finitely maximal Fuchsian groups}, J. London Math. Soc. 6 (1972), 29--38. 


\bibitem[Siu80]{siuannals} Siu, Y. T. {\it
The complex-analyticity of harmonic maps and the strong rigidity of 
compact K\"ahler
manifolds.} Ann. of Math. (2) 112 (1980), no. 1, 73--111.

\bibitem[Siu81]{siu2} Siu, Y. T. {\it Strong rigidity of compact 
quotients of exceptional
bounded symmetric domains.}  Duke Math. J.  48  (1981),  no. 4, 857--871.

\bibitem[Siu82]{siuJDG}
Siu, Yum-Tong
{\em Complex-analyticity of harmonic maps, vanishing and Lefschetz theorems. }
J. Differ. Geom. 17, 55-138 (1982).

\bibitem[Siu87]{siucurves}
 Siu, Yum-Tong
{\em Strong rigidity for K\"ahler manifolds and the construction of bounded holomorphic functions.}
Discrete groups in geometry and analysis, Pap. Hon. G. D. Mostow 60th Birthday, Prog. Math. 67, 124--151 (1987).


\bibitem[Som75]{somm75}
A.J. Sommese, {\em Quaternionic manifolds.} Math. Ann. 212 (1974/75), 191--214.

\bibitem[Span66]{Spanier}
Spanier, Edwin H. , {\em Algebraic topology.} McGraw-Hill Book Co., New York-Toronto, Ont.-London xiv+528 pp, (1966 ).

\bibitem[Stee51]{steenrod}
 Steenrod, Norman, {\em  The Topology of Fibre Bundles. }
 Princeton Mathematical Series, vol. 14. Princeton University Press, Princeton, N. J., viii+224 pp., (1951). 
 
 \bibitem[Taub92]{Taubes}
 Taubes, Clifford Henry, 
{\em The existence of anti-self-dual conformal structures.} 
J. Differential Geom. 36, no. 1, 163--253 , (1992). 


\bibitem[T-Z13]{tianKE}
 Tian Gang, Zhang Zhenlei,
{\em Regularity of the K\"ahler-Ricci flow,} C.R. Math. Acad. Sci. Paris 351, n. 15--16 (2013), 635--638; announcement of
{\em Regularity of K\"ahler-Ricci flows on Fano manifolds}
arXiv:1310.5897. 

\bibitem[Tju70]{tju}
G.N. Tjurina, {\em Resolution of singularities of plane (= flat)
deformations of double rational points },
{\bf Funk.Anal. i Prilozen  4} (1970), 77-83.

\bibitem[To93]{toledo}
Toledo, Domingo,
{\em  Projective varieties with non-residually finite fundamental group.}
 Inst. Hautes \'Etudes Sci. Publ. Math. No. 77 (1993), 103--119.
 
 \bibitem[Tot99]{Totaro}
  Totaro, Burt,
  {\em The Chow ring of a classifying space. } 
  In `Algebraic K-theory' (Seattle, WA, 1997), 249-281, Proc. Symp. Pure Math. 67, A.M.S., Providence R.I. (1999).

\bibitem[Tro96]{tromba}
A. Tromba,
{\em  Dirichlet's energy on Teichm\"uller 's moduli space and the 
Nielsen realization
problem}, {\bf Math. Z. 222} (1996), 451-464.

\bibitem[UY76]{UchidaYoshihara}
Uchida, Koji; Yoshihara, Hisao
{\em Discontinuous groups of affine transformations of $ \CC^3$}. 
Tohoku Math. J. (2) 28,  no. 1, 89--94 (1976).

 \bibitem[Ueno75]{keno}
 Ueno, Kenji 
 {\em Classification theory of algebraic varieties and compact complex spaces.}
  Notes written in collaboration with P. Cherenack. Lecture Notes in Mathematics, Vol. 439. Springer-Verlag, Berlin-New York, (1975). xix+278 pp.


\bibitem[Va06]{murphy} Vakil, R. {\em Murphy's law in algebraic
geometry: badly-behaved deformation spaces.}
{\bf Invent. Math.  164  no. 3} , (2006),  569--590.

\bibitem[Vieh95]{vieh}
E. Viehweg,
{\em Quasi-projective moduli for polarized manifolds,}
{\bf Ergebnisse der Mathematik und ihrer Grenzgebiete, 3. Folge,
30, Springer-Verlag, Berlin.}(1995), viii + 320 pp. .




\bibitem[Vieh10]{lastopus}
E. Viehweg,
{\em Compactifications of smooth families and of
moduli spaces of polarized manifolds}
{\bf Annals of Mathematics, 172},  (2010), 809--910.

\bibitem[ViZu07]{V-Z}
{\sc Viehweg, E. and Zuo, K.:} {\it Arakelov
inequalities and the
uniformization of certain rigid Shimura
varieties.}, J.
Differential
Geom. 77 (2007), no. 2, 291-- 352.

\bibitem[Vois02]{vois}
C. Voisin,
{\em Th\'eorie de Hodge et g\'eom\'etrie alg\'ebrique complexe,}
{\bf Cours Sp\'ecialis\'es , 10,
Soci\'et\'e Math\'ematique de France}, Paris, (2002), viii+595 pp.

\bibitem[Vois04]{voisin1}
C. Voisin,
{\em On the homotopy types of compact K\"ahler and complex projective manifolds,} Invent. Math. 157 (2004), no. 2, 329--343.

\bibitem[Vois06]{voisin2}
C. Voisin,
{\em  On the homotopy types of K\"ahler manifolds and the birational Kodaira problem.} J. Differential Geom. 72 (2006), no. 1, 43--71. 
 
\bibitem[Waj83]{wajnryb}
B. Wajnryb,
{\em  A simple presentation for the mapping class group of an
    orientable surface,}
{\bf Israel J. Math. 45,  no. 2-3} (1983), 157--174.

\bibitem[Waj99]{waj2}
B. Wajnryb,
{\em   An elementary approach to the mapping class group of a surface},
{\bf Geom. Topol. 3}(1999), 405--466.

\bibitem[Wall62]{wall}
C.T.C. Wall, {\em On the orthogonal groups of unimodular quadratic forms},
{\bf   Math. Ann. 147} (1962), 328-338.

 \bibitem[Wal50]{walker}
Walker, R.  J. :
{ \em  Algebraic curves.  }
(Princeton Mathematical Series. 13) Princeton, N. J.: Princeton University Press; London: Oxford University Press X, 201 p. (1950),
reprinted also by Dover and Springer 


\bibitem[Was82]{washington}
Washington, Lawrence C.
{\em  Introduction to cyclotomic fields. } Graduate Texts in Mathematics, 83. Springer-Verlag, New York,  xi+389 pp.,(1982); Second edition. . xiv+487 pp.(1997). 





\bibitem[Wav69]{Wav}
J. Wavrik, {\em  Obstructions to the existence of a space of moduli.} In `Global Analysis' (Papers in Honor of K. Kodaira)
 pp. 403--414 Princeton University and Univ. Tokyo Press, Tokyo (1969).

\bibitem[Weib94]{weibel}
C.A. Weibel, {\em An introduction to homological algebra.}
 Cambridge Studies in Advanced Mathematics, 38. Cambridge University Press, Cambridge, (1994) xiv+450 pp.

\bibitem[Weyl39]{Weyl} 
Weyl,  Hermann  {\em 
 Invariants}, Duke Mathematical Journal 5 (3), (1939),  489--502.
 

\bibitem[Wie]{Wiegold}  Wiegold, J., {\it The {S}chur multiplier: an elementary approach},
{Groups---{S}t. {A}ndrews 1981 ({S}t. {A}ndrews, 1981)},  {London Math. Soc. Lecture Note Ser.},
 {\bf 71},  {137--154}, {Cambridge Univ. Press}, (1982).



\bibitem[Yau77]{yau}
S.T. Yau,
{\em  Calabi's conjecture and some new results in algebraic geometry,}
{\bf Proc. Natl. Acad. Sci. USA 74} (1977), 1798-1799 .

\bibitem[Yau78]{Yau78}
Yau, S. T.:
{\em On the Ricci curvature of a compact K\"ahler manifold and the complex
Monge-Amp\`ere equation. I}.
{\bf Commun. Pure Appl. Math.  31, no. 3} (1978), 339--411.

\bibitem[Yau93]{yau2} {\sc Yau, S.-T:} {\it A splitting theorem and
an algebraic geometric characterization of locally Hermitian symmetric
spaces.}, Comm. Anal. Geom. 1 (1993), no. 3-4, 473 -- 486.

 \bibitem[Zak81]{zak}
 Zak, F. L. {\em Projections of algebraic varieties. }  Mat. Sb. (N.S.) 116(158) (1981), no. 4, 593--602, 608.
 
 \bibitem[Zar35]{Zar} 
 Zariski, Oscar {\em Algebraic surfaces.}With appendices by S. S. Abhyankar, J. Lipman and D. Mumford. 
 Preface to the appendices by Mumford. Ergebnisse der Mathematik und ihrer Grenzgebiete, Band 61, second  edition (1971), Springer-Verlag, 
 reprint of the first edition Band III (1935), xii+270 pp..

 
 \bibitem[Z87]{Zimmer} Zimmermann, B., {\it Surfaces and the second homology of a group},
  Monatsh. Math. {\bf 104} (1987), no. 3, 247--253. 
.

 

\end{thebibliography}
\end{document}